\documentclass[oneside,12pt]{amsart}

\usepackage{amsmath,amssymb,latexsym,eucal,amsthm}

\newtheorem{theorem}{Theorem}[section]
\newtheorem{lemma}[theorem]{Lemma}
\newtheorem{corollary}[theorem]{Corollary}
\newtheorem{proposition}[theorem]{Proposition}

\theoremstyle{definition}

\newtheorem{definition}[theorem]{Definition}
\newtheorem{conjecture}[theorem]{Conjecture}
\newtheorem{remark}[theorem]{Remark}

\newtheorem{notation}[theorem]{Notation}

\theoremstyle{remark}

\newtheorem*{note}{Note}
\newtheorem*{warning}{Warning}
\newtheorem*{question}{Question}

\newenvironment*{subproof}[1][Proof]
{\begin{proof}[#1]}{ \end{proof}}

\newenvironment*{case}[1]
{\textbf{Case #1.  }\itshape }

\newenvironment*{claim}[1][Claim]
{\textbf{#1.  }\itshape }

\newcommand{\setof}[2]{\left\{ \, #1 \, \left| \, #2 \, \right.\right\}}
\newcommand{\lsetof}[2]{\left\{\left. \, #1 \, \right| \, #2 \,  \right\}}
\newcommand{\bigsetof}[2]{\bigl\{ \, #1 \, \bigm | \, #2 \,\bigr\}}

\newcommand{\dotsetof}[2]{\left\{ \, #1 \, : \, #2 \, \right\}}

\newcommand{\sequence}[2]{\left\langle \, #1 \,\left| \, #2 \, \right. \right\rangle}
\newcommand{\lsequence}[2]{\left\langle\left. \, #1 \, \right| \,#2 \,  \right\rangle} 
\newcommand{\bigsequence}[2]{\bigl\langle \,#1 \, \bigm | \, #2 \, \bigr\rangle}

\newcommand{\singleton}[1]{\left\{#1\right\}}
\newcommand{\angles}[1]{\left\langle #1 \right\rangle}

\newcommand{\forcein}[2]{\overset{#2}{\Vert\underset{\!\!\!\!\!#1}%
{\!\!\!\relbar\!\!\!\relbar\!\!\relbar\!\!\relbar\!\!\!\relbar\!\!\relbar\!%
\!\relbar\!\!\!\relbar\!\!\relbar\!\!\relbar\!\!\relbar\!\!\!\relbar\!\!%
\relbar\!\!\relbar}}}

\newcommand{\pre}[2]{{}^{#2}\!{#1}}

\newcommand{\restr}{\!\!\upharpoonright\!}

\newcommand{\intersect}{\cap}
\newcommand{\union}{\cup}

\newcommand{\Union}[1]{\bigcup\limits_{#1}}

\newcommand{\defeq}{=_{\text{def}}}

\newcommand{\lexleq}{\leq_{\text{lex}}}
\newcommand{\lexless}{<_{\text{lex}}}

\newcommand{\R}{\mathbb{R}}
\renewcommand{\P}{\mathbb{P}}
\newcommand{\Q}{\mathbb{Q}}
\newcommand{\Z}{\mathbb{Z}}

\newcommand{\LofR}{L(\R)}
\newcommand{\JofR}[1]{J_{#1}(\R)}
\newcommand{\SofR}[1]{S_{#1}(\R)}
\newcommand{\JalphaR}{\JofR{\alpha}}
\newcommand{\JbetaR}{\JofR{\beta}}
\newcommand{\SalphaR}{\SofR{\alpha}}
\newcommand{\SbetaR}{\SofR{\beta}}
\DeclareMathOperator{\ORD}{OR}

\DeclareMathOperator{\WO}{WO}
\DeclareMathOperator{\OD}{OD}
\DeclareMathOperator{\HC}{HC}
\DeclareMathOperator{\WF}{WF}
\DeclareMathOperator{\HF}{HF}

\newcommand{\ONE}{I}

\newcommand{\TWO}{II}

\DeclareMathOperator{\ZFC}{ZFC}
\DeclareMathOperator{\ZF}{ZF}
\DeclareMathOperator{\AD}{AD}
\DeclareMathOperator{\KP}{KP}
\DeclareMathOperator{\PD}{PD}

\newcommand{\pred}{\text{-pred}}

\DeclareMathOperator{\Det}{Det}
\DeclareMathOperator{\dom}{dom}
\DeclareMathOperator{\ran}{ran}
\DeclareMathOperator{\crit}{crit}
\DeclareMathOperator{\card}{card}
\DeclareMathOperator{\cof}{cof}
\DeclareMathOperator{\rank}{rank}

\DeclareMathOperator{\rud}{rud}
\DeclareMathOperator{\Powerset}{\mathcal{P}}
\DeclareMathOperator{\length}{lh}
\DeclareMathOperator{\lh}{lh}

\DeclareMathOperator{\Ult}{Ult}

\DeclareMathOperator{\Col}{Coll}

\newcommand{\Implies}{\Rightarrow}
\newcommand{\Iff}{\Leftrightarrow}
\newcommand{\AND}{\wedge}
\newcommand{\OR}{\vee}

\newcommand{\surjection}{\xrightarrow{\text{onto}}}
\newcommand{\bijection}{\xrightarrow[\text{onto}]{\text{1-1}}}
\newcommand{\cofmap}{\xrightarrow{\text{cof}}}
\newcommand{\map}{\rightarrow}

\newcommand{\initseg}{\trianglelefteq}
\newcommand{\properseg}{\lhd}
\newcommand{\notinitseg}{\ntrianglelefteq}
\newcommand{\notproperseg}{\ntriangleleft}

\newcommand{\cH}{\mathcal{H}}

\newcommand{\cJ}{\mathcal{J}}

\newcommand{\cL}{\mathcal{L}}
\newcommand{\cM}{\mathcal{M}}
\newcommand{\cN}{\mathcal{N}}

\newcommand{\cP}{\mathcal{P}}
\newcommand{\cQ}{\mathcal{Q}}
\newcommand{\cR}{\mathcal{R}}
\newcommand{\cS}{\mathcal{S}}
\newcommand{\cT}{\mathcal{T}}
\newcommand{\cU}{\mathcal{U}}

\newcommand{\cW}{\mathcal{W}}

\newcommand{\bprime}{b^{\prime}}
\newcommand{\cprime}{c^{\prime}}

\newcommand{\mprime}{m^{\prime}}
\newcommand{\nprime}{n^{\prime}}

\newcommand{\pprime}{p^{\prime}}

\newcommand{\sprime}{s^{\prime}}

\newcommand{\wprime}{w^{\prime}}
\newcommand{\xprime}{x^{\prime}}

\newcommand{\Bprime}{B^{\prime}}

\newcommand{\alphaprime}{\alpha^{\prime}}
\newcommand{\betaprime}{\beta^{\prime}}
\newcommand{\gammaprime}{\gamma^{\prime}}
\newcommand{\deltaprime}{\delta^{\prime}}

\newcommand{\rhoprime}{\rho^{\prime}}
\newcommand{\xiprime}{\xi^{\prime}}
\newcommand{\cMprime}{\cM^{\prime}}

\newcommand{\formulaphi}{\text{\large $\varphi$}\, }

\newcommand{\Sa}[2][\alpha]{\Sigma_{(#1,#2)}}
\newcommand{\Pa}[2][\alpha]{\Pi_{(#1,#2)}}
\newcommand{\Da}[2][\alpha]{\Delta_{(#1,#2)}}
\newcommand{\Aa}[2][\alpha]{A_{(#1,#2)}}
\newcommand{\Ca}[2][\alpha]{C_{(#1,#2)}}
\newcommand{\Qa}[2][\alpha]{Q_{(#1,#2)}}

\newcommand{\San}{\Sa{n}}
\newcommand{\Pan}{\Pa{n}}
\newcommand{\Dan}{\Da{n}}
\newcommand{\Can}{\Ca{n}}
\newcommand{\Qan}{\Qa{n}}
\newcommand{\Aan}{\Aa{n}}

\DeclareMathOperator{\code}{code}

\newcommand{\CodeSet}{\mathfrak{C}}
\newcommand{\Rlanguage}{{\cL}_{\R}}

\newcommand{\BfSigmaOne}{\underset{\sim}{\Sigma}\mbox{${}_1$}}

\newcommand{\BfSigman}{\underset{\sim}{\Sigma}\mbox{${}_n$}}

\newcommand{\BfdeltaTwoOne}{\underset{\sim}{\delta}\mbox{${}^2_1$}}
\newcommand{\BfDeltaTwoOne}{\underset{\sim}{\Delta}\mbox{${}^2_1$}}
\newcommand{\BfSigmaTwoOne}{\underset{\sim}{\Sigma}\mbox{${}^2_1$}}
\newcommand{\BfSan}{\underset{\sim}{\Sigma}\mbox{${}_{(\alpha,n)}$}}

\newcommand{\BfPan}{\underset{\sim}{\Pi}\mbox{${}_{(\alpha,n)}$}}
\newcommand{\BfDan}{\underset{\sim}{\Delta}\mbox{${}_{(\alpha,n)}$}}

\newcommand{\BfDeltaOne}{\underset{\sim}{\Delta}\mbox{${}_1$}}
\newcommand{\BfSigmanplusone}{\underset{\sim}{\Sigma}\mbox{${}_{n+1}$}}

\newcommand{\skipsmall}{\vspace{1em}}
\newcommand{\skipmed}{\vspace{2em}}
\newcommand{\skipbig}{\vspace{3em}}
\newcommand{\skipsmallminus}{\vspace{-1em}}

\begin{document}

\title{Mouse Sets}

\author{Mitch Rudominer}

\address{Department of Mathematics\\
         Florida International University\\
         Miami, FL 33199}

\email{rudomine@fiu.edu}


\begin{abstract}
In this paper we explore a connection between descriptive set
theory and inner model theory. From descriptive set theory, we
will take a countable, definable set of reals, $A$. We will then
show that $A=\R\intersect\cM$, where $\cM$ is a canonical model
from inner model theory. In technical terms, $\cM$ is a ``mouse''.
Consequently, we say that $A$ is a mouse set.
For a concrete example of the type of set $A$ we are working with,
let $\OD^{\omega_1}_n$ be the set of reals which are $\Sigma_n$ 
definable over the model $L_{\omega_1}(\R)$, from an  ordinal
parameter. In this paper we will show  that for all $n\geq 1$,
$\OD^{\omega_1}_n$ is a mouse set.
Our work extends some similar results
due to D.A. Martin, J.R. Steel, and H. Woodin.
Several interesting questions in this area remain open. 
\end{abstract}

\maketitle

\tableofcontents

\setcounter{section}{-1}

%


\section{Introduction}

In this paper we will extend some work by  D. A. Martin, J. R. Steel, and
H. Woodin. We begin by describing this work.

If $\xi$ is a countable ordinal and $x\in\R$  and $n\geq2$, then  we say that
$x\in\Delta^1_n(\xi)$ if for every real $w\in\WO$ such that $|w|=\xi$,
$x\in\Delta^1_n(w)$. Set
$$A_n = \setof{x\in\R}{(\exists\xi<\omega_1)\; x\in\Delta^1_n(\xi)}.$$
In short, $A_n$ is the set of reals which are $\Delta^1_n$
in a countable ordinal.
Under the hypothesis of Projective Determinacy ($\PD$),
the sets $A_n$ are of
much interest to descriptive set theorists.
If $n\geq2$ is even then $A_n = C_n$, the largest
countable $\Sigma^1_n$ set. If $n\geq3$ is odd then
$A_n=Q_n$, the largest countable $\Pi^1_n$ set which is closed downward under
$\Delta^1_n$ degrees. The sets $A_n$ have been studied extensively by
descriptive set theorists. See for example \cite{Ke1}, and \cite{KeMaSo}.

The work of Martin, Steel, and Woodin shows that there is a connection
between the sets $A_n$, and certain inner models of the form
$L[\vec{E}]$, where $\vec{E}$ is a sequence of extenders.
In the paper ``Projectively Wellordered Inner Models'' \cite{MaSt}, Martin,
Steel, and Woodin
prove the following theorem:

\begin{theorem}[Martin, Steel, Woodin]
\label{FirstThm}
Let $n\geq1$ and suppose that there are $n$ Woodin
cardinals with a measurable cardinal above them. Let
$\cM_n$ be the canonical $L[\vec{E}]$ model with $n$ Woodin cardinals. Then
$$\R\intersect\cM_n=A_{n+2}.$$
\end{theorem}

The above theorem is also true with $n=0$. The ``canonical $L[\vec{E}]$ 
model with $0$ Woodin cardinals'' is just the model $L$. If there is
a measurable cardinal, then $A_2=\R\intersect L$. Actually, this result
is true under the weaker hypothesis that $\R\intersect L$ is countable.

Now let $A^*$ be the collection of reals which are ordinal definable in $\LofR$.
To emphasize the parallel with the sets $A_n$ above, note that assuming
that every game in $L(\R)$ is determined ($\AD^{\LofR}$), we have:
$$A^*=\setof{x\in\R}{(\exists\xi<\omega_1)\; x\in(\Delta^2_1)^{\LofR}(\xi)},$$
and $A^*$ is the largest countable $(\Sigma^2_1)^{\LofR}$ set of reals.

In the paper ``Inner Models with Many Woodin Cardinals'' \cite{St2}, Steel
and Woodin prove
the following theorem:

\begin{theorem}[Steel, Woodin]
\label{SecondThm}
 Suppose that there are $\omega$ Woodin
cardinals with a measurable cardinal above them. Let
$\cM_{\omega}$ be the canonical $L[\vec{E}]$ model with $\omega$ Woodin
cardinals. Then
$$\R\intersect\cM_{\omega}=A^*.$$
\end{theorem}

Theorems \ref{FirstThm} and  \ref{SecondThm} are obviously similar to
each other.  The two theorems establish a very deep connection between
descriptive set theory and inner model theory.
Notice that between Theorem \ref{FirstThm} and Theorem \ref{SecondThm} there
is some room.  Our goal in writing this paper was to prove a theorem, 
similar
to the these two theorems, for pointclasses $\Delta$  \emph{between}
the $\Delta^1_n$ pointclasses and the pointclass $(\Delta^2_1)^{\LofR}$.

More specifically, let us say that a set of reals $A$ is a 
\emph{\textbf{mouse set}}
iff $A=\R\intersect\cM$ for some canonical $L[\vec{E}]$ model $\cM$.
(To be more precise, we will take $\cM$ to be a countable, iterable 
\emph{mouse}. A mouse is a model of the form $L_{\alpha}[\vec{E}]$, for
some ordinal $\alpha$. Theorems \ref{FirstThm} and  \ref{SecondThm} above
can be restated in terms of mice.)
Given a pointclass $\Delta$ let us set:
$$A_{\Delta}=\setof{x\in\R}{(\exists\xi<\omega_1)\; x\in\Delta(\xi)}.$$
Our goal in writing this paper was to find some pointclasses $\Delta$ such
that $A_n\subset A_{\Delta}$ for all $n$, and $A_{\Delta}\subset A^*$,
and then to prove that $A_{\Delta}$ is a mouse set.

To accomplish this task, first we need to have some such pointclasses 
$\Delta$. In Section \ref{section:pointclasses} we define some lightface 
pointclasses in $\LofR$
which are similar to the pointclasses of the analytical
hierarchy. For ordinals $\alpha$ such that 
$1\leq\alpha\leq(\BfdeltaTwoOne)^{\LofR}$,
we define a sequence of pointclasses:
$$\bigsequence{\San,\Pan,\Dan}{n\in\omega}.$$
For $\alpha=1$, we will have that
$\Sa[1]{n}=\Sigma^1_n$. So our setting includes the setting
of the analytical hierarchy.
For $\alpha=(\BfdeltaTwoOne)^{\LofR}$, we will have that
$\Sa{0}=(\Sigma^2_1)^{\LofR}$. So our setting also includes the setting
of Theorem \ref{SecondThm} above.

Then we define:
$$\Aan=\setof{x\in\R}{(\exists\xi<\omega_1)\; x\in\Dan(\xi)}.$$
We leave the full technical definition
of $\Aan$ for later. For the sake of this introduction, it is
enough to understand that $\Aan$ is the set of reals which
are ordinal definable in $\LofR$, ``at the $(\alpha,n)^{\text{th}}$ level.''
Assuming ${\AD}^{\LofR}$, these sets are very much like the analytical sets 
$A_n$ defined above.
That is, for even $n$, $\Aan=\Can$, the largest countable $\San$ set.
For odd $n$, $\Aan=\Qan$, the largest countable $\Pan$ set closed 
downward under $\Dan$ degrees. 

Because of this analogy, we conjecture that, assuming there are
$\omega$ Woodin cardinals,
 $\Aan$ is a mouse set,
for all $(\alpha,n)$. To establish this, we need to have some
canonical $L[\vec{E}]$ models $\cM$. In Section \ref{section:bigmice},
 we
study mice $\cM$ which are ``bigger''
than the models $\cM_n$ mentioned in Theorem \ref{FirstThm}, but ``smaller'' 
than the model $\cM_{\omega}$ mentioned in Theorem \ref{SecondThm}.
But how can $\cM$ have ``more than'' $n$ Woodin cardinals, for each $n$,
but ``less than'' $\omega$ Woodin cardinals? The answer of course is
that ``more than'' and ``less than'' here refer to consistency strength.
Our mice $\cM$ will have the property that for each $n$, $\cM$ satisfies
that there is a transitive set model of $\ZFC + \; \exists$ $n$
 Woodin Cardinals.

In this paper we are not able to fully prove our conjecture that
$\Aan$ is a mouse set, for all $(\alpha,n)$. We are able to show
this for some $(\alpha,n)$ though. Our argument has several parts,
and each part works for a different set of $(\alpha,n)$. In the
end we manage to prove that $\Aan$ is a mouse set for for all
$(\alpha,n)$ such that $\alpha\leq\omega_1^{\omega_1}$, and
either $\cof(\alpha)>\omega$, or $\cof(\alpha)\leq\omega$ and
$n=0$. There is nothing special about the ordinal $\omega_1^{\omega_1}$
here. The main property of $\omega_1^{\omega_1}$ which we use is the
fact that  all ordinals 
$\alpha<\omega_1^{\omega_1}$ are simply definable from countable
ordinals. It is easy to see how to extend our results somewhat 
beyond $\omega_1^{\omega_1}$.  
 Eliminating the restriction in our work to ordinals 
$\alpha<\omega_1^{\omega_1}$
does not
appear to be a very difficult problem. On the other hand,
the restriction to $n=0$ in
the case that $\cof(\alpha)\leq\omega$ seems to be more difficult to
overcome.  In fact, at the very first step,
we are not even able to show that the set $\Aa[2]{1}$
is a mouse set.
We believe that the full conjecture, for all $(\alpha,n)$,
 is a very interesting open problem.
In this paper we develop some ideas which should be useful 
towards solving this problem.

To explain the roles of the various sections of this paper, we
will describe our strategy for proving the conjecture that
$\Aan$ is a mouse set.
 Our strategy is a natural generalization of the strategy used
by Martin, Steel, and Woodin in \cite{MaSt}.
The
strategy begins with defining, for each ordinal $\alpha$ and each
$n\in\omega$,  the notion
of a mouse $\cM$ being $(\alpha,n)$-big.  If a mouse is $(\alpha,n)$-big
then it satisfies a certain large cardinal hypothesis which is stronger
than the hypothesis ``$\exists k $ Woodin cardinals'', for all
$k\in\omega$, but weaker than the hypothesis ``$\exists \omega $ Woodin
cardinals.'' If a mouse is not
$(\alpha,n)$-big, we call it $(\alpha,n)$-small. The definition
of $(\alpha,n)$-big is given in Section \ref{section:bigmice}.
Our definition there only works for $\alpha\leq\omega_1^{\omega_1}$.
The next step in our strategy is to prove a theorem which states
that if $\cM$ is a countable, iterable mouse and $\cM$
is  $(\alpha,n)$-big, then $\Aan\subseteq\R\intersect\cM$.
We  prove such a theorem in Section \ref{section:correctness}.

The next step in our strategy is to define, for each ordinal $\alpha$ and 
each $n\in\omega$,  the notion of a mouse $\cM$ being $\Pan$-iterable. 
This notion is an approximation of full iterability, with the added
advantage that it is simply definable: 
The set of (reals coding) countable mice $\cM$
which are $\Pan$-iterable is a $\Pan$ set.
We give the definition
of $\Pan$-iterable in Section \ref{section:iterability}. Our definition
there only works for $\alpha<\kappa^{\R}$, where $\kappa^{\R}$ is the
least $\R$-admissible ordinal.
Next, our strategy is to prove a theorem which states
that if $\cM$ and
$\cN$ are countable mice, and $\cN$ is $\Pan$-iterable, and
$\cM$ is fully iterable and $(\alpha,n+1)$-small, then $\cM$ and
$\cN$ can be compared.  

Unfortunately we are not able to prove this theorem. This accounts
for our inability to prove the above-mentioned conjecture, in the
case $\cof(\alpha)\leq\omega$ and $n>0$. We can prove an approximation
of this theorem though. In Section \ref{section:petitemice}, we define
the notion of a mouse $\cM$ being $(\alpha,n)$-petite. This is
a condition which is slightly more restrictive than $(\alpha,n)$-small.
We are able to prove the above-mentioned comparison theorem with
$(\alpha,n)$-petite in place of $(\alpha,n)$-small. We do this in
Section \ref{section:comparison}, for $\alpha\leq\omega_1^{\omega_1}$.

In the case $\cof(\alpha)>\omega$, $(\alpha,n)$-petite is synonymous
with $(\alpha,n)$-small. So in this case we get a proof of our
conjecture.  The proof goes as follows:
Let $\cM$ be a countable, iterable mouse
such that $\cM$ is  $(\alpha,n)$-big but every proper initial
segment of $\cM$ is $(\alpha,n)$-small. The results of Section
\ref{section:correctness} give us that $\Aan\subseteq\R\intersect\cM$,
and the results of Section  \ref{section:comparison} give us that
$\R\intersect\cM\subseteq\Aan$. Thus $\Aan=\R\intersect\cM$ and so
$\Aan$ is a mouse set. This argument is given in Section 
\ref{section:concl}.

In the case $\cof(\alpha)\leq\omega$, our definition of $(\alpha,n)$-petite 
is not synonymous with $(\alpha,n)$-small. So in this case we do not
get a proof of our conjecture. However we do manage to eke out a proof
of the conjecture in the case $n=0$. That is, we prove that for all
$\alpha\leq\omega_1^{\omega_1}$, $\Aa{0}$ is a mouse set.
Actually, the fact that $\Aa{0}$ is a mouse
set is perhaps our most ``quotable'' result, since $\Aa{0}$ is
a very natural, and easily described set. 
$\Aa{0}$ is the set of reals which are ordinal
definable in some level of $\LofR$ below  the $\alpha^{\text{th}}$
level. For example, with $\alpha=\omega_1$
$$\Aa[\omega_1]{0}=\setof{x\in\R}{(\exists\beta<\omega_1)\; x %
\text{ is ordinal definable in } L_{\beta}(\R)\;}.$$ So in this paper
we prove that, assuming there are $\omega$ Woodin cardinals,
$\Aa[\omega_1]{0}$ is a mouse set.

This paper is based in part on my PhD dissertation. I would like
to thank my thesis advisor, Professor John Steel, for his guidance.

%


\skipbig

\section{Prerequisites}

\label{section:prereq}

In this section we briefly discuss some material which is a
prerequisite for understanding this paper. First we discuss
some material from the paper ``Scales in $\LofR$'' by J. R. 
Steel \cite{St1}. Then we discuss some material from the
papers ``Fine Structure and Iteration Trees'' by W. Mitchell
and J. R. Steel \cite{MiSt}, and 
``Inner Models With Many Woodin Cardinals'' by J. R. Steel \cite{St2}.
 We present
this material here purely as a service to the reader.

\skipmed

\noindent
\textbf{Scales In $\LofR$}

\skipsmall

All of the definitions and results in this sub-section are taken 
word-for-word
from  \cite{St1}.  As we shall heavily use this material we reprint it here
for the reader's convenience.     We work in the theory ZF + DC.

Throughout this paper we will use the symbol
$\R$ to mean $\pre{\omega}{\omega}$
and the word \emph{real} to refer to an element of $\R$.

A {\em pointclass} is a class
of subsets of $\R$ closed under recursive substitution. If $\Gamma$ is a
pointclass then
\begin{align*}
\neg\Gamma &= \setof{\R - A}{A \in \Gamma}\,=\,\text{the dual of $\Gamma$}\\
\exists^{\R}\Gamma &= \setof{\exists^{\R}A}{A \in \Gamma}; \;
\text{where } \exists^{\R}A = \setof{x}{\exists y \,(\angles{x,y} \in A)}\\
\text{and}\qquad
\forall^{\R}\Gamma &= \setof{\forall^{\R}A}{A \in \Gamma}; \;
\text{where } \forall^{\R}A = \setof{x}{\forall y \,(\angles{x,y} \in A)}.
\end{align*}
By $\Det(\Gamma)$ we mean the assertion that all games whose payoff
set is in $\Gamma$ are determined.

Because it is convenient, we shall use the Jensen J-hierarchy for $\LofR$.
Let $\rud(M)$ be the closure of $M \union \singleton{M}$ under the
rudimentary functions. Let $\HF$ be the set of hereditarily
finite sets. Now let
\label{def:JalphaR}
\begin{align*}
\JofR{1} &= \HF\union\,\R\\
\JofR{\alpha + 1} &= \rud(\JalphaR); \quad \text{for $\alpha > 0$}\\
\text{and}\quad
\JofR{\lambda} &= \Union{\alpha < \lambda} \JalphaR;
\quad \text{for $\lambda$ a limit ordinal.}
\end{align*}
Of course $\LofR = \Union{\alpha \in \ORD} \JalphaR$.
For all $\alpha \geq 1$, $\JalphaR$ is transitive.  For $1 \leq
\alpha < \beta$, $\JalphaR \in \JbetaR$.  For $\alpha
> 1$, $\rank{\JalphaR} = \ORD \intersect \JalphaR =
\omega\alpha$.  Also
\begin{multline*}
\JofR{\alpha + 1} \intersect \Powerset(\JalphaR) = \\
\setof{x \subseteq \JalphaR}{x \text{ is definable over $\JalphaR$
from parameters in $\JalphaR$}}.
\end{multline*}

It is useful to refine this hierarchy.  There is a single rudimentary
function $S$ such that for all transitive sets $M$, $S(M)$
is transitive, $M \in S(M)$, and $\rud(M) =
\Union{n \in \omega} S^n(M)$.  Now let
\begin{align*}
\SofR{\omega} &= \JofR{1}\\
\SofR{\alpha + 1} &= S(\SalphaR); \quad \text{for $\alpha > \omega$}\\
\text{and}\quad
\SofR{\lambda} &= \Union{\alpha < \lambda} \SalphaR;
\quad \text{for $\lambda$ a limit ordinal.}
\end{align*}
Then for all $\alpha \geq 1$, $\JalphaR = \SofR{\omega\alpha}$.

If $\alpha \geq 2$ and $X \subseteq \JalphaR$ then
$\Sigma_n(\JalphaR,X)$ is the class of all relations on $\JalphaR$
definable over $(\JalphaR;\: \in)$ by a $\Sigma_n$ formula, from
parameters in $X \union \singleton{\R}$.  Thus we are always
allowed $\R$ itself (but not necessarily its
elements) as a parameter in definitions over $\JalphaR$.
$\Sigma_{\omega}(\JalphaR, X) \defeq \Union{n \in \omega} \Sigma_n
(\JalphaR,X)$.  We write ``$\Sigma_n(\JalphaR)$'' for
``$\Sigma_n(\JalphaR, \emptyset)$,'' and ``$\BfSigman(\JalphaR)$''
for ``$\Sigma_n(\JalphaR,\JalphaR)$.''  Similar conventions apply to
the $\Pi_n$ and $\Delta_n$ notations.

Our first lemma states that for all $\alpha\geq 1$ and all $A\subseteq\R$ if
$A\in\Sigma_1(\JofR{\alpha + 1})$ then $A=\Union{n\in\omega}A_n$
with each $A_n \in \Sigma_{\omega}(\JalphaR)$.
\begin{lemma}
\label{SigmaOneIsUnion}
There is a recursive function which assigns to each $\Sigma_1$ formula
$\formulaphi$ with 2 free variables, a sequence of formulae
$\lsequence{\formulaphi_n}{n\in\omega}$ so that for all $\alpha\geq1$ and
all $x\in\R$, $\JofR{\alpha+1}\models \formulaphi[x,\R]$ iff
$(\exists n\in\omega)\; \JalphaR\models \formulaphi_n[x,\R]$.
\end{lemma}
\begin{proof}
This is essentially proved in Lemma 1.11 of \cite{St1}.  We have\\
$\JofR{\alpha+1}\models \formulaphi[x,\R]$\\
$\Iff \quad (\exists n\in\omega)\;
\SofR{\omega\alpha+n}\models \formulaphi[x,\R]$\\
$\Iff \quad (\exists n\in\omega)\;
S^n(\JalphaR)\models \formulaphi[x,\R]$.\\
Now rudimentary functions are simple and so there are $\Sigma_0$ formulae
$\lsequence{\theta_n}{n\in\omega}$ such that for all transitive sets
$M$ and all $a, b\, \in \, M$\\
$\theta_n(M,a,b) \Iff S^n(M)\models \formulaphi[a,b]$.\\
Finally let $\formulaphi_n$ be formulae such that for all transitive sets
$M$ and all $a, b\, \in \, M$\\
$\theta_n(M,a,b) \Iff M\models \formulaphi_n[a,b]$.
\end{proof}

The following is Lemma 1.1 of \cite{St1}.
\begin{lemma}
\label{SequenceOfLevelsSigma1Definable}
The sequences $\lsequence{\JbetaR}{\beta < \alpha}$ and
$\lsequence{\SbetaR}{\beta < \alpha}$ are uniformly $\Sigma_1(\JalphaR)$
for $\alpha > 1$.
\end{lemma}

If $1 \leq \alpha < \beta$ and $X \subseteq \JalphaR$, then
$\JalphaR \prec^{X}_n \JbetaR$ means that for all $\vec{a} \in
(X \union \singleton{\R})^{<\omega}$, and all $\Sigma_n$
formulae $\formulaphi$, $\JalphaR \models \formulaphi[\vec{a}]$ iff
$\JbetaR \models \formulaphi[\vec{a}]$.  We write
``$\JalphaR \prec_n \JbetaR$'' for
``$\JalphaR \prec^{\JalphaR}_n \JbetaR$.''

Let $\theta$ be the least ordinal not the surjective image of $\R$, and
$\BfdeltaTwoOne$ the least ordinal not the image of $\R$ under
a surjection $f$ such that $\setof{\angles{x,y}}{f(x) \leq f(y)}$
is $\BfDeltaTwoOne$. Then $\theta^{\LofR}$ = the least $\alpha$
such that $\Powerset(\R) \intersect \LofR \subseteq \JalphaR$.  The
following is Lemma 1.2 of \cite{St1}.
\begin{lemma}
\label{CharacterizeDelta-2-1}
Let $\sigma$ be least such that $\JofR{\sigma} \prec^{\R}_1 \LofR$. Then
\begin{itemize}
\item[(a)] $(\BfSigmaTwoOne)^{\LofR} =
\BfSigmaOne(\JofR{\sigma}) \ \intersect \ \Powerset(\R)$
\item[(b)] $(\BfDeltaTwoOne)^{\LofR} =
\JofR{\sigma}  \intersect  \Powerset(\R)$
\item[(c)] $(\BfdeltaTwoOne)^{\LofR} = \sigma$
\end{itemize}
\end{lemma}
\begin{definition}
\em
Let $\alpha, \beta \in \ORD$ with $\alpha \leq \beta$.  The interval
$[\alpha,\beta]$ is a $\Sigma_1$-{\em gap} iff
\begin{itemize}
\item[(i)] $\JalphaR \prec^{\R}_1 \JbetaR$
\item[(ii)] $(\forall \alpha^{\prime} < \alpha) \big[\JofR{\alpha^{\prime}}
{\not\prec}^{\R}_1 \JalphaR\big]$
\item[(iii)] $(\forall \beta^{\prime} > \beta) \big[\JbetaR {\not\prec}^{\R}_1
\JofR{\beta^{\prime}}\big]$
\end{itemize}
\end{definition}

That is, a $\Sigma_1$-gap is a maximal interval of ordinals in which no new
$\Sigma_1$ facts about elements of $\R \union \singleton{\R}$
are verified.  If $[\alpha,\beta]$ is a $\Sigma_1$-gap, we say that
$\alpha$ begins the gap and $\beta$ ends it.
Notice that we allow $\alpha = \beta$.
We shall also allow $[(\BfdeltaTwoOne)^{\LofR},\theta^{\LofR}]$
as a $\Sigma_1$-gap.  Notice that if $[\alpha,\beta]$ is a $\Sigma_1$-gap
and $\alpha < \beta$ then $\JalphaR$ is admissible.
The following is Lemma 2.3
of \cite{St1}.
\begin{lemma}
The $\Sigma_1$-gaps partition $\theta^{\LofR}$.
\end{lemma}

The following is Lemma 1.11 of \cite{St1}.
\begin{lemma}
\label{StartOfGapProjects}
Suppose that $\alpha > 1$ and $\alpha$ begins a $\Sigma_1$-gap. Then
\begin{itemize}
\item[(a)] There is a $\Sigma_1(\JalphaR)$ {\em\bf partial} map
$h\vdots\R\surjection\JalphaR$.
\item[(b)] If $\beta < \alpha$ then there is a {\em\bf total} map
$f:\R\surjection\JbetaR$ such that $f \in \JalphaR$.
\end{itemize}
\end{lemma}

We conclude  with the main theorem from \cite{St1}.  The
following is Theorem 2.1 from that paper.
\begin{lemma}
\label{ScalesInLofR}
If $\alpha > 1$ and $\Det(\JalphaR)$ holds then the pointclass
$\Powerset(\R) \intersect \Sigma_1(\JalphaR)$ has the scale property.
\end{lemma}

\skipmed

\noindent
\textbf{Fine Structure and Iteration Trees}

\skipsmall

In this subsection we briefly review some of the terms we will use
from inner model theory. We would like to give the reader an
intuitive understanding of what a ``mouse'' is.  At the end of
this subsection we
 state the  Condensation Theorem, and then we
prove a corollary to the Condensation Theorem which will be important
to us in the sequel.

As we did for $\LofR$, we
also use the Jensen $J$-hierarchy for models of the form
$L[\vec{E}]$, where $\vec{E}$ is a sequence of extenders.
Recall that one of the properties of a $J$-hierarchy implies
 that for any ordinal $\alpha$,
$\ORD\intersect J^{\vec{E}}_{\alpha}=\omega\alpha$.

A \emph{potential premouse} is a model of the form  
$\cM=(M;\;\in, E^{\cM}, F^{\cM})$,
where $M$ is a transitive set, and $E^{\cM}$, $F^{\cM}$ are amenable
predicates over $M$, and for some sequence 
$\vec{E}=\sequence{E_{\alpha}}{\alpha\in\dom(\vec{E})}$
 we have that:
\begin{itemize}
\item[(i)]$\vec{E}$ is a \emph{good} sequence of  extenders, and
the domain of $\vec{E}$ is some set of ordinals, and
\item[(ii)] for all $\gamma\in\dom(\vec{E})$, $\gamma=\omega\gamma$,
and $E_{\gamma}$ is an extender of length $\gamma$ over the model
$J^{\vec{E}\restriction\gamma}_{\gamma}$, and
\item[(iii)]$M=J^{\vec{E}\restriction\alpha}_{\alpha}$, 
for some $\alpha$, and
\item[(iv)] $E^{\cM}$ codes $\vec{E}\restr\alpha$, and
\item[(v)]  $F^{\cM}$ codes $E_{\alpha}$, if $\alpha\in\dom(\vec{E})$;
otherwise $F^{\cM}=\emptyset$.
\end{itemize}
For the definition of a \emph{good} sequence of  extenders, see
\cite{MiSt}.
If $F^{\cM}=\emptyset$ then we say that $\cM$ is \emph{passive}.
Otherwise we say that $\cM$ is \emph{active}.

If  $\cM = 
(J^{\vec{E}\restriction\alpha}_{\alpha};\;\in, E^{\cM}, F^{\cM})$ 
is a potential
premouse then, for $\beta\leq\alpha$, $\cJ^{\cM}_{\beta}$ is the
premouse 
$\cN = (J^{\vec{E}\restriction\beta}_{\beta};\;\in, E^{\cN}, F^{\cN})$,
where $E^{\cN}$ codes $\vec{E}\restr\beta$ and
$F^{\cN}$ codes $E_{\beta}$, if $\beta\in\dom(\vec{E})$, or
 $F^{\cN}=\emptyset$ otherwise. For $\beta\leq\alpha$,
$\cJ^{\cM}_{\beta}$ is called an \emph{initial segment} of
$\cM$. We write $\cN\unlhd\cM$ to mean that $\cN$ is an initial
segment of $\cM$. $\cN\lhd\cM$ means that $\cN$ is a proper initial
segment of $\cM$.

We must briefly consider the fine structure of definability
over potential premice.
A potential premouse $\cM$ is \emph{$n$-sound}, iff 
$\cM$ is the $\Sigma_n$ hull in $\cM$ of 
$\rho_n(\cM)\union\singleton{p_n(\cM)}$, where $\rho_n(\cM)$ is the
$\Sigma_n$ \emph{projectum} of $\cM$, and $p_n(\cM)$ is the
$n$-th \emph{standard parameter} of $\cM$.  $\cM$ is \emph{sound}
iff for all $n\in\omega$, $\cM$ is $n$-sound.
(I am being a little
sloppy here. See \cite{MiSt} for the real definition of
sound. In particular, when talking about definability over
potential premice, we do not actually use the Levy hierarchy
of $\Sigma_n$ definability, but rather a quasi-Levy hierarchy
called $r\Sigma_n$ definability.)

A \emph{premouse} is a potential premouse, all of whose proper initial
segments are sound. From now on, we shall never need to talk about
potential premice, only premice. Potential premice are only defined
in order to be able to define premice.

\begin{notation}
Let $\cM=(M;\;\in, E^{\cM}, F^{\cM})$ be a premouse.  Then we will use
the symbol $\cM$ to refer both to the \emph{model} $\cM$, and to
the \emph{set} $M$.  From the context it will be apparent to which object
we are referring.
\end{notation}

A premouse $\cM$ is \emph{meek}, iff there is no initial segment
$\cN\unlhd\cM$ such that $\cN$ is active and, letting $\kappa$ be
the critical point of last extender on the $\cN$ sequence, $\cN$
satisfies that there are cofinally many Woodin cardinals below $\kappa$.
See \cite{St2}. Meekness is a property which ensures that a premouse
does not satisfy too large of a large cardinal hypothesis. For
the sake of this paper, it is enough to know that if no initial segment
of $\cM$ has $\omega$ Woodin cardinals, then $\cM$ is meek.

A  premouse $\cM$ is \emph{iterable} iff one can build 
iteration trees on $\cM$ without encountering illfoundedness. 
There
is more than one way to make this notion precise. One particularly
useful way is via the notion of realizability. A premouse $\cM$
is \emph{realizable} iff $\cM$ can be embedded into another premouse
$\bar{\cM}$, where $\bar{\cM}$ has the property that the extenders
on the $\bar{\cM}$ sequence all come from extenders over $V$. 
See \cite{St2} for the precise definition of realizable. In
\cite{St2} it is shown that if a countable, meek premouse 
$\cM$ is realizable, then 
 $\cM$ is iterable enough to do whatever we want with it.
In particular, if $\cM$ and $\cN$ are two countable, meek,
realizable premice, then $\cM$ and $\cN$ can be compared.

A \emph{mouse} is an iterable premouse. Ironically, since there
is not a single, universally useful notion of iterable, the
term mouse is never used. In our setting,
 we could say that a mouse is a  realizable, meek
premouse.

A premouse is an ``$L$-like'' model.  If $\cM$ is a premouse then
$\cM$ is a model of GCH, and in fact, $\cM$ is strongly
acceptable. Thus if $\kappa$ is a cardinal of $\cM$ and
$x\in J^{\cM}_{\kappa}$, then 
$\Powerset(x)\intersect\cM\subset J^{\cM}_{\kappa}$. Another way in which
a premouse $\cM$ is like $L$ is that $\cM$ satisfies the following
Condensation Theorem.

\begin{theorem}[The Condensation Theorem]
Let $\cH$ and $\cM$ be  countable, realizable, meek premice. 
Let $n\geq 0$ and suppose that $\cH$ and $\cM$ are both
$n+1$-sound. Suppose there is a map $\pi:\cH\map\cM$ such that
$\pi$ is $r\Sigma_{n+1}$ elementary, and 
$\crit(\pi)\geq\rho_{n+1}(\cH)$. Then either\\
(a) $\cH$ is a proper initial segment of $\cM$\\
or\\
(b) There is an extender $E$ on the $\cM$ sequence such that
$\length(E)=\rho_{n+1}(\cH)$, and $\cH$ is a proper initial
segment of $\Ult_0(\cM,E)$.
\end{theorem}

The Condensation Theorem, which is due to Mitchell and Steel, is
proved in the union of \cite{MiSt} and \cite{St2}. See in particular
Theorem 8.2 of \cite{MiSt}.

We conclude this section with two corollaries to the Condensation
Theorem. These corollaries concern the $\initseg$-least initial segment
$\cN$ of some given premouse $\cM$ such that $\cN$ satisfies some given
fact about some given parameter.
Our corollaries imply, for instance, that the $\initseg$-least premouse
$\cN$
which satisfies a given \emph{sentence} must project to $\omega$.
We will use these corollaries several times in the sequel.

Before stating our corollaries we will need one lemma.
The following lemma is implicit in section 2 of \cite{MiSt}.
See \cite{MiSt} for a definition of the fine structural terms
used in this lemma. 
\begin{lemma}
\label{WeakEmbedding}
Let $\cM$ be an $n$-sound premouse. Let $\cH=\cH^{\cM}_{n+1}(X)$
where $X\subseteq \cM$ and $u_n(\cM)\in X$. Let $\pi:\cH\map\cM$
be the inverse of the collapse.  Then $\pi$ is a weak $n$-embedding
and $\pi$ is $r\Sigma_{n+1}$ elementary.
\end{lemma}
In the above lemma, the fact that $\pi$ is a weak $n$-embedding
means that:
\begin{itemize}
\item[(a)] $\pi$ is $r\Sigma_{n+1}$ elementary.
\item[(b)] For $1\leq k\leq n-1$, $\pi(\rho_k^{\cH})=\rho_k^{\cM}$.
\item[(c)] $\rho_n^{\cH} =$ the least $\alpha$ such that
$\pi(\alpha)\geq\rho^{\cM}_n$.
\item[(d)] $\cH$ is $n$-sound, and
\item[(e)] $\pi(p_k(\cH))=p_k(\cM)$ for $k\leq n$.
\end{itemize}

Now we state our first corollary to the Condensation Theorem.
\begin{corollary}
\label{PhiMinimalProjects}
Let $\cM$ be a countable, realizable, meek premouse.
Let $\delta$ be a cardinal of
$\cM$, and let $q\in J^{\cM}_{\delta}$. Let $n\geq0$ and let $\formulaphi$ 
be
an $r\Sigma_{n+1}$ formula with 1 free variable.  Suppose there is a
proper initial segment $\cN\lhd\cM$ such that $\delta\in\cN$ and
$\cN\models\formulaphi[q]$.
Let $\cN$ be the $\initseg$-least such initial segment.
Then $\rho_{n+1}(\cN)=\delta$.
\end{corollary}
\begin{proof}
Let $\cH=\cH^{\cN}_{n+1}\big(\delta\union\singleton{u_n(\cN)}\big)$
and let $\pi:\cH\map\cN$ be the inverse of the collapse map. By Lemma
\ref{WeakEmbedding} $\pi$ is a weak $n$-embedding and $\pi$ is
$r\Sigma_{n+1}$ elementary. Thus $\cH$ is $n$-sound and
$\pi(u_n(\cH))=u_n(\cN)$.
So $\cH=\cH^{\cH}_{n+1}\big(\delta\union\singleton{u_n(\cH)}\big)$.
As $\delta$ is a cardinal of $\cM$, $\rho_{n+1}^{\cH} = \delta$,
and so $p_{n+1}(\cH)=\angles{\emptyset,u_n(\cH)}$ and $\cH$ is $n+1$-sound,
and $\rho_{n+1}(\cH)=\delta\leq\crit{\pi}$. By the Condensation Theorem
$\cH\unlhd\cN$.
(Note that alternative (b) of the Condensation Theorem cannot hold
since $\delta=\rho_{n+1}(\cH)$ is a cardinal of $\cM$, and so 
$\cJ^{\cM}_{\delta}$ is passive.)
  Now consider two cases.

\underline{Case 1.}\quad $\delta\in H$.\\ Then as $\cH\models\formulaphi[q]$, we
have $\cH=\cN$ and so we are done.

\underline{Case 2.}\quad $\cH=\cJ^{\cN}_{\delta}$.\\
In this case we must have $u_n(\cN)=\emptyset$.
Let $\cH^{\prime}=
\cH^{\cN}_{n+1}\big(\delta\union\singleton{\delta}\big)$
and let $\pi^{\prime}:\cH^{\prime}\map\cN$ be the inverse of the
collapse map.  By Lemma \ref{WeakEmbedding}, $\pi^{\prime}$ is
a weak $n$-embedding and $\pi^{\prime}$ is $r\Sigma_{n+1}$ elementary.
Thus $\cH^{\prime}$ is $n$-sound and $u_n(\cH^{\prime})=\emptyset$.
As $\cH^{\prime}=
\cH^{\cH^{\prime}}_{n+1}
\big(\delta\union\singleton{\delta}\big)$
and $\delta$ is a cardinal of $\cM$, $\rho^{\cH}_{n+1}=\delta$.

\begin{claim}
$\angles{\delta}$ is $(n+1)$-solid over
$\big(\cH^{\prime},\emptyset\big)$.
\end{claim}
\begin{subproof}[Proof of Claim]
$Th^{\cH^{\prime}}_{n+1}(\delta)=
Th^{\cN}_{n+1}(\delta)=
Th^{\cH}_{n+1}(\delta)$. But
$\cH=\cJ^{\cN}_{\delta}=\cJ^{\cH^{\prime}}_{\delta}$. So
$Th^{\cH^{\prime}}_{n+1}(\delta)
\in\cH^{\prime}$.
\end{subproof}
Thus $p_{n+1}(\cH^{\prime})=\angles{\delta,\emptyset}$ and
$\cH^{\prime}$ is $(n+1)$-sound and $\rho_{n+1}(\cH^{\prime})=
\delta\leq\crit{\pi^{\prime}}$. So we may apply the Condensation Theorem
and conclude that $\cH^{\prime}\unlhd\cN$.
As $\cH^{\prime}\models\formulaphi[q]$, $\cH^{\prime}=\cN$.
\end{proof}

Let $\cM$, $\delta$, $q$, and $\formulaphi$ be as in the previous corollary.
Now let $\cN$ be the $\initseg$-least initial segment of $\cM$ such
that $\cN\models\formulaphi[q]$. By analogy with the behavior of the model
$L$, we would expect that $\ORD^{\cN}<\delta$. Using the previous
Corollary we can see only that $\ORD^{\cN}<(\delta^+)^{\cM}$. It is 
almost certainly true that $\ORD^{\cN}<\delta$, but we are unable to see
how to show this, simply by quoting the \emph{statement}
of the Condensation Theorem. A proof
that $\ORD^{\cN}<\delta$ can probably be given using an argument
similar to the \emph{proof} of the Condensation Theorem, that is, by using a
comparison argument. It would take us too far afield to do such a comparison
argument 
here, and besides we do not need the result in full generality. Below
we give a short argument which only uses the \emph{statement} of
 the Condensation Theorem, in the special case that $\delta=\omega_1^{\cM}$.
What follows is our second corollary to the Condensation Theorem.

\begin{corollary}
\label{PhiMinimalIsCountable}
Let $\cM$ be a countable, realizable, meek premouse.
Suppose $\cM\models$``$\omega_1$ exists,'' and let 
$\delta=\omega_1^{\cM}$.
Let $q\in J^{\cM}_{\delta}$. Let  $\formulaphi$ 
be any  formula with 1 free variable.  Suppose there is a
proper initial segment $\cN\properseg\cM$ such that 
$\cN\models\formulaphi[q]$.
Let $\cN$ be the $\initseg$-least such initial segment.
Then $\ORD^{\cN}<\delta$.
\end{corollary}
\begin{proof}
Suppose that $\ORD^{\cN}\geq\delta$, and we will derive a contradiction.
If $\ORD^{\cN}>\delta$, then by the previous corollary 
$\rho_{\omega}(\cN)=\delta$. If $\ORD^{\cN}=\delta$ then trivially
$\rho_{\omega}(\cN)=\delta$. Let $n\in\omega$ be large enough so that
$\rho_n(\cN)=\delta$, and also so that $\formulaphi$ is an $r\Sigma_n$
formula. 

Fix $\xi<\delta$ such that $q\in J^{\cM}_{\xi}$ and
$\rho_{\omega}(\cJ^{\cM}_{\xi})=\omega$. Let 
$\cH=\cH^{\cN}_{n+5}(\xi+1\union\singleton{\delta, u_{n+5}(\cN)})$, and let
$\pi:\cH\map\cN$ be the inverse of the collapse map. Let 
$\kappa=\crit(\pi)$. A standard argument show that 
$\ran(\pi)\intersect\delta$ is transitive, and so
$\pi(\kappa)=\delta$.
By Lemma \ref{WeakEmbedding} $\pi$ is a weak $(n+4)$-embedding and $\pi$
is $r\Sigma_{n+5}$ elementary.  Since $\pi(q)=q$, we have that
$\cH\models\formulaphi[q]$.
Also, $\cH$ is $n$-sound and
$\rho_n(\cH)=\kappa$. Thus we may apply the Condensation Theorem to
conclude that either\\
(a) $\cH$ is a proper initial segment of $\cM$\\
or\\
(b) There is an extender $E$ on the $\cM$ sequence such that
$\length(E)=\kappa$, and $\cH$ is a proper initial
segment of $\Ult_0(\cM,E)$.\\
In case (a) we have that $\ORD^{\cH}<\delta$, and this contradicts
our assumption that $\cN$ was the $\initseg$-least initial segment of
$\cM$ such that $\cN\models\formulaphi[q]$. So suppose
that case (b) occurs.

Fix an extender $E$ as in (b).  Let
$\cP=\cJ^{\cM}_{\kappa}$. Let $i=i^{\cP}_{E}$.
Thus $i:\cP\map\Ult_0(\cP,E)$.
 Now $\cH$ is a proper initial segment of 
$\Ult_0(\cP,E)$, and so $\Ult_0(\cP,E)$ satisfies the statement:
``There is an initial segment 
$\cQ$ of me such that $\cQ\models\formulaphi[q]$.''
Let $\mu=\crit(E)=\crit(i)$. Since 
$\rho_{\omega}(\cJ^{\cM}_{\xi})=\omega$ and $\mu$ is a cardinal
in $\cP$, we have that $\mu\geq\xi$.
Thus $i(q)=q$. It follows that $\cP$ satisfies the statement:
``There is an initial segment
$\cQ$ of me such that $\cQ\models\formulaphi[q]$.''
Again, this contradicts
our assumption that $\cN$ was the $\initseg$-least initial segment of
$\cM$ such that $\cN\models\formulaphi[q]$.
\end{proof}

%


\skipbig

\section{Some Lightface Pointclasses in $\protect{\LofR}$}

\label{section:pointclasses}

In this section we work in the theory ZF + DC. For each ordinal
$\alpha \geq 2$, we now define an analytical-hierarchy-like sequence of
lightface pointclasses.
\begin{definition}
For $\alpha \geq 2$, let
\begin{align*}
\Sa{0} &= \Powerset(\R) \intersect \Sigma_1(\JalphaR)\\
\text{and for $n \geq 0$ let}\qquad
\Pa{n} &= \neg\Sa{n}\\
\Da{n} &= \Sa{n} \intersect \Pa{n}\\
\Sa{n+1} &= \exists^{\R}\Pa{n}
\end{align*}
\end{definition}

So a set of reals $A$ is in the pointclass $\San$, iff $A$ is
definable over $\JalphaR$ using a formula $\formulaphi$ with a 
certain restricted syntax. We would like to give a name to
this syntax.  We will say that a formula $\formulaphi$ is
$\Sigma_{(\bullet,n)}$, iff $\formulaphi$ is the type of formula
which is used to define a $\San$ set. More specifically,
let $\Rlanguage$ be the language with the constant symbol $\dot{\R}$
in addition to the binary relation symbol $\in$.
A formula in $\Rlanguage$ is interpreted in the model $(\JalphaR ; \in )$ 
by interpreting $\dot{\R}$ as $\R$. A formula
$\formulaphi$ is  $\Sigma_{(\bullet,0)}$ iff $\formulaphi$ is a
$\Sigma_1$ formula of  $\Rlanguage$. 
A formula
$\formulaphi$ is  $\Pi_{(\bullet,0)}$ iff $\formulaphi$ is a
$\Pi_1$ formula of  $\Rlanguage$.  If $n\geq1$, then
$\formulaphi$ is $\Sigma_{(\bullet,n)}$ iff there is a 
$\Pi_{(\bullet,n-1)}$ formula $\psi$ such that
$\formulaphi=(\exists v \in \dot{\R} ) \, \psi$, and
$\formulaphi$ is $\Pi_{(\bullet,n)}$ iff there is a 
$\Sigma_{(\bullet,n-1)}$ formula $\psi$ such that
$\formulaphi=(\forall v \in \dot{\R} ) \, \psi$. Thus,
a set of reals $A$ is in $\San$ just in case
$$A=\setof{x\in\R}{\JalphaR\models\formulaphi[x]}$$
for some $\Sigma_{(\bullet,n)}$ formula $\formulaphi$.

If $x \in \R$ then the relativized pointclasses $\Sa{n}(x)$,
$\Pa{n}(x)$, $\Da{n}(x)$ are defined as usual.  Set $\BfSan$ =
$\Union{x \in \R}\Sa{n}(x)$.  Similarly for $\BfPan$.
Set $\BfDan = \BfSan \intersect \BfPan$.
For completeness we also set $\Sa[1]{n} = \Sigma^1_n$, etc.

Since a $\Sigma_{(\bullet,n)}$ formula is, in particular,
a $\Sigma_{n+1}$ formula, $\San\subseteq\Sigma_{n+1}(\JalphaR)$.
But, in general, $\San\not=\Sigma_{n+1}(\JalphaR)$. The reason
we are interested in the pointclasses $\San$ instead of just
$\Sigma_{n+1}(\JalphaR)$, is that the property of being a
\emph{scaled} pointclass is propagated by the quantifiers
$\exists^{\R}$ and $\forall^{\R}$. It is a little annoying
that our definition yields that 
$\San\subseteq\Sigma_{n+1}(\JalphaR)$ instead of
$\San\subseteq\Sigma_{n}(\JalphaR)$. The reason we made our definition
the way we did is so that, as in the analytical hierarchy, it is
the \emph{even} $\Sigma$ pointclasses which are scaled. This emphasizes
the analogy with the pointclasses of the analytical hierarchy.

We are interested
in the pointclasses $\Sa{n}$ only for certain ordinals $\alpha$.
If $[\alpha,\beta]$
is a $\Sigma_1$-gap and $\alpha < \gamma \leq \beta$, then for
all $x \in \R$ and all $n \in \omega$, $\Sa[\gamma]{n}(x) = \Sa{n}(x)$.
So we are only interested in the pointclass $\San$ in the case that
$\alpha$ begins a $\Sigma_1$-gap.
Furthermore, if $\alpha$ begins a $\Sigma_1$-gap and $\JalphaR$ is
admissible, then $\Sigma_1(\JalphaR)$ is closed under the real quantifiers
$\forall^{\R}$ and $\exists^{\R}$.  So $\forall x \in \R$
\begin{align*}
&\Sa{0}(x) = \Pa{1}(x) = \Sa{2}(x) = \dots\\
\text{and}\qquad
&\Pa{0}(x) = \Sa{1}(x) = \Pa{2}(x) = \dots
\end{align*}
So the collection of pointclasses
$
\lsetof{\Sa{n},\Pa{n}}{n \in \omega}
$
is not particularly interesting.  This leads us to the following.
\begin{definition}
An ordinal $\alpha$ is \emph{projective-like} iff $\alpha$ begins a
$\Sigma_1$-gap and $\JalphaR$ is not admissible.
\end{definition}

\begin{remark}
$\alpha$ projective-like $\Implies$ $[\alpha,\alpha]$ is a
$\Sigma_1$-gap.
\end{remark}

We are only interested in the pointclasses
$
\lsetof{\Sa{n},\Pa{n}}{n \in \omega}
$
in the case that $\alpha$ is projective-like.
In Proposition \ref{NonTrivialPointclasses} below we show that if
$\alpha$ is
projective-like then the collection of pointclasses
$
\lsetof{\Sa{n},\Pa{n}}{n \in \omega}
$
does \emph{not} collapse to two pointclasses as it does in the case that
$\alpha$ is admissible. The importance of being projective-like
is given by the following lemma.
\begin{lemma}
\label{GoodParameterExists}
Let $\alpha \geq 2$ be projective-like. Then $\exists x \in \R$ such that
there is a total function
$f:\R\map\omega\alpha = \ORD\intersect\JalphaR$
with $\ran(f)$ cofinal in $\omega\alpha$ and with
$ f \in \Sigma_1(\JalphaR,\singleton{x}).$
\end{lemma}
\begin{proof}
By definition of non-admissible $\exists \beta < \alpha$ and $\exists
g:\JbetaR\cofmap\omega\alpha$ such that $g \in \BfSigmaOne(\JalphaR)$.
By Lemma \ref{StartOfGapProjects}(b) $\exists h \in \JalphaR
(h:\R\surjection\JbetaR).$ Let $f = g\circ h$.  Then $\ran(f)$ is
cofinal in $\omega\alpha$ and $f \in \BfSigmaOne(\JalphaR)$.  By lemma
\ref{StartOfGapProjects}(a) we may take the parameter in the definition
of $f$ over $\JalphaR$ to be a real.
\end{proof}
\begin{definition}
\label{GoodParameter}
Let $\alpha \geq 2$ be projective-like.  A real $x$ as in Lemma
\ref{GoodParameterExists} above is called a \emph{good parameter}
for $\JalphaR$.
\end{definition}

Let $\alpha \geq 2$ be projective-like.  The collection
of pointclasses $\setof{\Sa{n},\Pa{n}}{n \in \omega}$ has many
of the properties of the analytical hierarchy
$\setof{\Sigma^1_n,\Pi^1_n}{n \in \omega}$.  However there is one very
important property of the analytical hierarchy which may fail for the
collection $\setof{\Sa{n},\Pa{n}}{n \in \omega}$.  This latter collection
\emph{may not form a hierarchy.}  That is

\begin{warning}
For some projective-like $\alpha$ it is not
the case that $\Sa{n} \subseteq \Da{n+1}$.
\end{warning}

Of course it is true that $\Sa{n} \subseteq \Pa{n+1} \subseteq \Sa{n+2}$.
Also we do have the following, which explains why we are interested in
good parameters: Let $\alpha \geq 2$ be projective-like.  Let $x$ be a good
parameter for $\JalphaR$.  Then $(\forall n)\quad \Sa{n}(x) 
\subseteq \Da{n+1}(x)$. That is, relative to a good parameter, our 
pointclasses do form a hierarchy. In fact, we have the following.
\begin{lemma}
\label{GoodParameterImpliesHierarchy}
Let $\alpha \geq 2$ be projective-like, and let $n\geq 0$.  
Let $x$ be a good parameter for $\JalphaR$.  
Then 
$\San(x)=\Powerset(\R) \intersect \Sigma_{n+1}(\JalphaR,\singleton{x})$.
\end{lemma}
\begin{proof}
This is essentially what is shown in Lemma 2.5 of \cite{St1}.
So this proof is due to John Steel.
We argue by induction on $n$. For $n=0$ the lemma is true by
definition of $\Sa{0}$. Now let $n\geq 1$. Even without the parameter
$x$ we already have that
that $\San\subseteq\Sigma_{n+1}(\JalphaR)$. So certainly we have 
$\San(x)\subseteq\Sigma_{n+1}(\JalphaR,\singleton{x})$. 
For the other direction, by induction, 
it suffices to see that
$\Sigma_{n+1}(\JalphaR,\singleton{x}) =
\exists^{\R}\Pi_n(\JalphaR,\singleton{x})$. 
Let $h\vdots\R\surjection\JalphaR$ be a partial $\Sigma_1(\JalphaR)$ map. 
There is such an $h$ since $\alpha$ begins a $\Sigma_1$-gap.
Let $R$ be any $\Pi_n(\JalphaR)$ relation. Then
\begin{gather*}
(\exists v \in \JalphaR)\;R(v)\\
\text{iff}\\
(\exists y\in\R)[y\in\dom{h}\;\AND\;
(\forall v\in\JalphaR)(v=h(y)\implies R(v)\;)\;].
\end{gather*}
If $n\geq2$ then ``$y\in\dom{h}$'' is $\Pi_n(\JalphaR)$ and we are
done. If $n=1$ then we must see that the relation ``$y\in\dom{h}$''
is $\Sa{1}(x)$. We have that the relation
``$y\in\dom{h}$'' is $\Sa{0}$. So it suffices to prove the following:

\begin{claim}
$\Sa{0}\subseteq\Sa{1}(x)$.
\end{claim}
\begin{subproof}[Proof of Claim]
 Let $f:\R\cofmap\omega\alpha$ be a total 
$\Sigma_1(\JalphaR,\singleton{x})$ map.  
Let $R$ be any $\Sa{0}$ set.  So there is a $\Sigma_1$ formula
$\formulaphi$ of $\Rlanguage$ so that 
$R(y) \Iff \JalphaR \models \formulaphi[y]$
\begin{align*}
&\Iff (\exists \beta < \omega\alpha) \SofR{\beta} \models \formulaphi[y]\\
&\Iff(\exists z \in \R) \SofR{f(z)} \models \formulaphi[y]\\
&\Iff(\exists z \in \R)(\forall \beta < \omega\alpha)(\forall S \in
\JalphaR)\big[(\beta = f(z) \AND S = \SofR{\beta}) \implies S \models
\formulaphi[y]\;\big].
\end{align*}
So $R \in \exists^{\R}\Pi_1(\JalphaR,\singleton{x}) = \Sa{1}(x)$.
\end{subproof}
This completes the proof of the lemma.
\end{proof}

\begin{remark}
In particular, if $\alpha\geq 2$ is projective like, then
$\BfSan$ is simply the pointclass $\BfSigmanplusone(\JalphaR)$.
\end{remark}

If $0$ is a good parameter for $\JalphaR$, then for all $n$,
$\San\subseteq \Da{n+1}$,
and so we would expect  the $\San$ pointclasses to behave exactly
like the analytical pointclasses. More generally, if $x$ is a good
parameter for $\JalphaR$, then we would expect the pointclasses
$\San(x)$ to behave exactly like the analytical pointclasses.
But we want to consider the \emph{lightface} pointclasses
$\San$, even when $0$ is not a good parameter for $\JalphaR$.
In this situation,
even though $\Sa{n} \subseteq \Da{n+1}$ may fail, many of the proofs
of various well-known properties of the analytical hierarchy go
through in our setting.  Below we record some of these properties.
We begin with some
trivial observations.
\begin{lemma}
\label{SimplePointclassProperties}
Let $\alpha \geq 2$ be projective-like.
\begin{itemize}
\item[(1)] Every analytical pointset is in $\Da{0}$.
\item[(2)] For all $n \geq 0$, $\Sa{n}$ and $\Pa{n}$ are {\em\bf adequate}.
That is, they are closed under recursive substitution, $\AND$, $\OR$,
$\exists^{\leq}$, and $\forall^{\leq}$. (See \cite{Mo} pg. 158.)
\item[(3)] For all $n \geq 0$, $\Sa{n}$ and $\Pa{n}$ are
$\omega$-parameterized.
\item[(4)] For all $n \geq 0$, $\Sa{n}$ is closed under $\exists^{\R}$
and $\Pa{n}$ is closed under $\forall^{\R}$.
\item[(5)] For all $n \geq 1$, $\Sa{n}$ and $\Pa{n}$ are closed under
both integer quantifiers: $\exists^{\omega}$, $\forall^{\omega}$.
\end{itemize}
\end{lemma}

\begin{remark}
Notice that (5) may fail for $n = 0$ if $\alpha$ is a successor ordinal
or \mbox{$\cof(\alpha) = \omega$}.
\end{remark}

\begin{proof}
These are all easy.  For example, (3) holds for $n = 0$ because the
$\Sigma_1$-satisfaction relation over $\JalphaR$ is $\Sigma_1(\JalphaR)$.
Then (3) holds for $n>0$ by induction.
\end{proof}

Now we can see that if
$\alpha$ is
projective-like then the collection of pointclasses\\
$
\lsetof{\Sa{n},\Pa{n}}{n \in \omega}
$
does not collapse to two pointclasses as it does in the case that
$\alpha$ is admissible.
\begin{proposition}
\label{NonTrivialPointclasses}
Let $\alpha \geq 2$ be projective-like.  Then $(\forall n \in \omega)
\quad\Sa{n+1} \neq \Pa{n}$.
\end{proposition}
\begin{proof}
Suppose for some $n$ that $\Sa{n+1} = \Pa{n}$.  Then for all $x \in \R$,
$\Sa{n+1}(x) = \Pa{n}(x)$.
But then for $x$ a good parameter for $\JalphaR$,
$\Sa{n}(x) \subseteq \Pa{n}(x)$.  But this is impossible as there is a
universal $\Sa{n}$ set.
\end{proof}
\begin{lemma}[Second Periodicity Theorem]
\label{SecondPeriodicity}
Let $\alpha \geq 2$ be projective-like.  Assume $\Det(\JofR{\alpha + 1})$.
Let $n \in \omega$.  Then $\Sa{2n}$ and $\Pa{2n+1}$ have the scale
property.
\end{lemma}
\begin{proof}
This follows from Lemma \ref{ScalesInLofR}, Lemma
\ref{SimplePointclassProperties}, and Theorems 6C.2 and 6C.3 of \cite{Mo}.
\end{proof}

\begin{remark}
In particular, under the hypotheses of Lemma \ref{SecondPeriodicity},
$\Pa{2n+1}$ and $\Sa{2n+2}$ are {\em Spector} pointclasses. (See
section 4C of \cite{Mo}.)
\end{remark}

One very important consequence of the scale property is the
basis property. The next lemma says that our even $\Sigma$ pointclasses
have the basis property.

\begin{lemma}[The Basis Property]
\label{basislemma}
Let $\alpha\geq 2$ be projective-like. Assume $\Det(\JofR{\alpha+1})$.
Let $n\geq 2$ be even. Suppose $A\subset\R$ and $A\in\San$ and
$A$ is not empty. Then there is an $x$ in $A$ such that
$\singleton{x}\in\San$.
\end{lemma}
\begin{proof}
This follows from 4E.7 in \cite{Mo}.
\end{proof}

\begin{definition}
Let $\alpha \geq 2$ be projective-like, and let $n\geq 0$.  Let
$x,y \in \R$.  Then we say that $x \in \Da{n}(y)$ iff
$\singleton{x} \in \Da{n}(y)$, equivalently iff $\singleton{x} \in
\Sa{n}(y)$.  We say that $x \in \Da{n}$ iff $x \in \Da{n}(0)$.
Let  $A\subseteq \R$. Then we say that $A$ is \emph{closed downward under}
$\Dan$ iff whenever $x\in A$ and $y\in\Dan(x)$ we have that
$y\in A$.
\end{definition}

\begin{remark}
\label{rem:ZeroIsDifferent}
If $n \geq 1$ then $x \in \Da{n}(y)$ iff
$\setof{(i,j)}{x(i)=j} \in \Da{n}(y)$. For $n=0$, since $\Sa{0}$ is
not closed under $\forall^{\omega}$ for some $\alpha$,
 this latter condition is not equivalent to
$x \in \Da{0}(y)$.
\end{remark}

Another consequence of the scale property is the bounded quantification
property. This is explained in the next lemma. 

\begin{lemma}[Bounded Quantification]
\label{BoundedQuantification}
Let $\alpha \geq 2$ be projective-like. Let
$n \geq 1$ be odd.  Assume $\Det(\BfPan)$.  Let $Q(x,y,z) \in \Pa{n}$.  Set
$$P(x,y) \Iff \exists z \in \Dan(y)Q(x,y,z).$$
Then $P$ is $\Pa{n}$.
\end{lemma}
\begin{proof}
This follows from 4D.3 of \cite{Mo}.
\end{proof}

Another consequence of the scale property is the following result
due to Kechris. Recall that $\WO$ is the set of reals which code
wellorders of $\omega$, and hence code countable ordinals. 
If $w\in\WO$ then $|w|$ is the countable ordinal coded by $w$.

\begin{lemma}[Boundedness Theorem for Subsets of $\omega_1$]
\label{BoundedSetofOrdinals}
Let $\alpha \geq 2$ be projective-like.  Assume $\Det(\JofR{\alpha + 1})$.
Let $n \geq 1$ be odd.  Let $A \subset \WO$
be $\San$ and suppose that $\sup\dotsetof{|x|}{x \in A} < \omega_1$.
Then there exists a $w\in\WO$ such that $w\in\Dan$ and such that
$|x|<|w|$ for all $x\in A$.
\end{lemma}
\begin{proof}
The proof of Theorem 1.1 from \cite{Ke2} goes through here.
\end{proof}

The study of the $\San$ pointclasses, and in particular the
comparison of these pointclasses with the pointclasses of
the analytical hierarchy, is an interesting endeavor in it
own right. We pursue this study in the papers \cite{Ru1} and
\cite{Ru2}. In this paper however, descriptive set theory is
not our primary focus. So we will discuss  only those topics
from descriptive set theory which will be needed in this paper.
In addition to the several results we mentioned above,
there are two other topics which we will discuss:
the topic of countable sets of ordinal definable reals, and
 the topic of the
\emph{correctness} of a given set of reals. 


\vspace{2em}

\noindent
\textbf{Countable Sets of Ordinal Definable Reals}

\vspace{1em}

\label{section:ordef}

As Remark \ref{rem:ZeroIsDifferent} indicates,
when considering the property:
$$x\in\Dan(y)$$
the case $n=0$ is somewhat
different than the case $n\geq1$. Notice that the next two definitions
are only for $n\geq 1$.

\begin{definition}
Let $x \in \R$, $\xi < \omega_1$. Let $n \geq 1$.  Then we say
that $x \in \Dan(\xi)$ iff $\forall w \in \WO \text{ such that }
|w| =\xi, \quad x \in \Dan(w)$.
\end{definition}

\begin{definition}
Let $\alpha \geq 2$ be projective-like. Let $n \geq 1$.  Then
$$\Aa{n} \defeq
\setof{x \in \R}{(\exists \xi < \omega_1)\; x \in \Dan(\xi)}.$$
\end{definition}

Assuming $\Det(\JofR{\alpha+1})$, the sets $\Aan$ for $n\geq 1$ have
much descriptive set-theoretic interest. For even $n\geq 2$,
$\Aan=\Can$, where $\Can$ is the largest countable $\San$ set.
For odd $n\geq 1$, $\Aan=\Qan$, where $\Qan$ is the largest countable
$\Pan$ set which is closed down under $\Dan$. Even though
$\San\subset\Da{n+1}$ may fail, the sets $\Can$ and $\Qan$ still
share most of the properties of the analogous sets $C_n$ and $Q_n$
in the theory of the analytical hierarchy. For example, for $n\geq 1$
there is
a $\Dan$-good wellorder of $\Aan$.  We will not go into any more
detail here about the descriptive set-theoretic properties of the
sets $\Aan$, for $n\geq 1$. But see \cite{Ru1} and \cite{Ru2} where
these sets are studied in detail.

We wish now to discuss
the basic notion of a real $x$ being 
\emph{ordinal definable} over the model $\JalphaR$. 
For $n\geq 1$, the set of reals which are $\Sigma_n$ ordinal definable
over $\JalphaR$ has a wellordering which is definable over $\JalphaR$.
So assuming $\Det(\JofR{\alpha+1})$, this set is countable. Consequently,
every real which is ordinal definable over $\JalphaR$ is in fact definable
over $\JalphaR$ from a single countable ordinal. This leads us
to the following definition.

\begin{definition}
Let  $n\geq 1$.
$\OD^{\JalphaR}_n$ is the set of reals $x$ such that for some $\Sigma_n$
formula $\formulaphi$ in the language $\Rlanguage$, 
and some countable ordinal $\xi$,
and every real  $w\in\WO$ such that $|w|=\xi$, we have that $x$ is 
the unique real $\xprime$
such that  $\JalphaR\models\formulaphi[\xprime,w]$. Set 
$$\OD^{\JalphaR} = \Union{n}\OD^{\JalphaR}_n.$$
\end{definition}

Let $n\geq 1$, and let us consider the following question.
What is the relationship between $\Aan$ and $\OD^{\JalphaR}_{n+1}$?
Since $\Dan\subseteq\Delta_{n+1}(\JalphaR)$, we expect that
$\Aan\subseteq \OD^{\JalphaR}_{n+1}$. Notice however
that this fact is
not completely immediate, for the following technical reason.
Let $x\in\Aan$. Then there is some $\xi<\omega_1$ such that
$x\in\Dan(w)$ for every $w\in\WO$ such that $|w|=\xi$. But
the $\Dan$ definition of $x$ is allowed to depend on $w$,
so it is not obvious that $x$ is actually definable from $\xi$.
However, assuming $\Det(\JofR{\alpha+1})$, we have that 
$\Aan=\Can$ if $n$ is even, and $\Aan=\Qan$ if $n$ is odd. 
This implies that $\Aan$ is countable,
and $\Aan$ has a $\Dan$-good wellordering. This in turn implies that
there is a single fixed $\Sigma_{(\bullet,n)}$
formula $\formulaphi$, such that for all $x\in\Aan$, 
there is a countable ordinal $\xi$ such that whenever $w\in\WO$ and
$|w|=\xi$, $x$ is the unique real $\xprime$ such that
$\JalphaR\models\formulaphi[\xprime,w]$. Since $\formulaphi$ is,
in particular, a $\Sigma_{n+1}$ formula, we have that
$x\in\OD^{\JalphaR}_{n+1}$. That is, assuming $\Det(\JofR{\alpha+1})$,
we have that for $n\geq1$, $\Aan\subseteq\OD^{\JalphaR}_{n+1}$.
We ask the following two questions.

\begin{question}
\label{aquestion}
Assume $\AD^{\LofR}$. Let $\alpha\geq2$ be projective-like.
Let $n\geq 1$.
Is it true that $\Aan=\OD^{\JalphaR}_{n+1}$?
\end{question}

 A courser question is the following.

\begin{question}
Assume $\AD^{\LofR}$. Let $\alpha\geq2$ be projective-like.
Is it true that $\Union{n} \Aan = \OD^{\JalphaR}$?
\end{question}

That is to say, if $x$ is ordinal definable over $\JalphaR$, does it follow
that $x$ is definable using only quantification over the reals, over 
$\Sigma_1(\JalphaR)$? We  do not know the answer to these questions
in full generality.

We do have some partial results along these lines. 
Notice that if $\beta<\alpha$, then 
$\OD^{\JbetaR}\subseteq\OD^{\JalphaR}$.
Notice also that $\OD^{\JbetaR}=\OD^{\JalphaR}$ is possible. In fact,
assuming $\AD^{\LofR}$, there are only countably many $\alpha$
such that $\OD^{\JalphaR}\not=\Union{\beta<\alpha}\OD^{\JbetaR}$.
Now, given
$A\subseteq\OD^{\JalphaR}$, let us say that
$A$ is \emph{non-trivial} iff there is an $x\in A$ such that
$x\not\in\OD^{\JbetaR}$ for any $\beta<\alpha$.
In \cite{Ru1} and \cite{Ru2} we show that if $\Aan$  is non-trivial,
then $\Aan=\OD^{\JalphaR}_{n+1}$. So the first question can be
rephrased as: if $\OD^{\JalphaR}_{n+1}$ is non-trivial,
must $\Aan$ also be non-trivial. The second, courser question can
be rephrased as: if $\OD^{\JalphaR}$ is non-trivial,
must $\Aan$ also be non-trivial for some $n$. In \cite{Ru1} and
\cite{Ru2} we  show that the answer to the second question is yes,
in the case that $\alpha$ has uncountable cofinality. We are not able
to prove anything more about these questions---using techniques from
descriptive set theory. Somewhat surprisingly, we do get some additional
information using inner model theory.  

In this paper we show that the answer to the second question is yes,
for $\alpha\leq\omega_1^{\omega_1}$. Also, 
in the case that $\cof(\alpha)>\omega$,
we show that the answer to the first question is yes.
This is done as follows. For $\cof(\alpha)>\omega$ we show that there
is a mouse $\cM$ such that $\OD^{\JalphaR}_{n+1}\subseteq\R\intersect\cM$,
and  $\R\intersect\cM\subseteq\Aan$. Thus $\Aan=\OD^{\JalphaR}_{n+1}$.
For $\cof(\alpha)\leq\omega$ we have a gap in our proof. We only get
that there is a mouse $\cM$ such that 
$\OD^{\JalphaR}_{n+1}\subseteq\R\intersect\cM$, and 
$\R\intersect\cM\subseteq\Aa{n+1}$. Thus
$\Aan\subseteq\OD^{\JalphaR}_{n+1}\subseteq\Aa{n+1}$,
and so $\Union{n} \Aan = \OD^{\JalphaR}$.

\skipmed

There is one last point we need to discuss. We would like to extend
the definition of $\Aan$ to the case $n=0$.

\begin{definition}
Let $\alpha\geq 2$. Then $\Aa{0}=\Union{\beta<\alpha} \OD^{\JbetaR}$.
\end{definition}

Assuming determinacy, we have the following alternate characterization
of $\Aa{0}$.

\begin{proposition}
\label{ODone}
Assume $\Det(\Sigma_1(\JalphaR))$. Then
$\Aa{0}=\OD^{\JalphaR}_{1}$.
\end{proposition}
\begin{proof}
Notice that $\Aa{0}$ is a $\Sigma_1(\JalphaR)$ set
and there is a $\Sigma_1(\JalphaR)$ wellorder of $\Aa{0}$.
So assuming $\Det(\Sigma_1(\JalphaR))$, $\Aa{0}$ is countable.
Then, since there is a $\Sigma_1(\JalphaR)$ wellorder of $\Aa{0}$,
we have that $\Aa{0}\subseteq\OD^{\JalphaR}_{1}$. Conversely,
suppose that $x\in\OD^{\JalphaR}_{1}$. Then it is easy to see that
$x$ is ordinal definable in $\JbetaR$ for some $\beta<\alpha$. (If
$\alpha$ is a successor ordinal, then use Lemma 
\ref{SigmaOneIsUnion}.)
\end{proof}

So $\Aa{0}$ is a $\Sa{0}$ set, and assuming determinacy,
it is countable.
In analogy with even $n\geq2$ we might expect that $\Aa{0}$ is the
largest countable $\Sa{0}$ set. It is pretty easy to see that this
is true in the case that $\alpha$ is a limit ordinal. If $\alpha$
is a successor ordinal, we do not know if this is true in general.
In \cite{Ru1} and \cite{Ru2} we study this question and get some
partial results.

Now that we have defined $\Aan$ for all $n\geq 0$, we are ready
to state the main conjecture which motivated
this paper.

\begin{definition}
Let $A\subset \R$. Then we say that $A$ is a \emph{mouse set}
iff $A=\R\intersect\cM$ for some countable, realizable, meek
premouse $\cM$.
\end{definition}

\begin{conjecture}
Let $\alpha\geq2$ be projective-like. Then
$\Aan$ is a mouse set for all $n\geq 0$.
\end{conjecture}

As we mentioned in the introduction, in this paper we will prove the
conjecture in the case $\alpha\leq\omega_1^{\omega_1}$, and either
$\cof(\alpha)>\omega$, or $\cof(\alpha)\leq\omega$ and $n=0$. In
Section \ref{section:bigmice} we begin the process of proving
the conjecture by giving the definition of an
$(\alpha,n)$-big premouse. Then in Section \ref{section:correctness}
we give one half of the proof that $\Aan$ is a mouse set: we show
that if $\cM$ is a countable, iterable premouse, and $\cM$ is
$(\alpha,n)$-big, then $\Aan\subseteq\R\intersect\cM$. 
Actually, we prove the
ostensibly stronger result that 
$\OD^{\JalphaR}_{n+1}\subseteq\R\intersect\cM$. We will then use this
ostensibly stronger result to prove that $\Aan=\OD^{\JalphaR}_{n+1}$.
(See the discussion above.) The proof in Section \ref{section:correctness}
that $\OD^{\JalphaR}_{n+1}\subseteq\R\intersect\cM$ will use that
fact that $\R\intersect\cM$ has a certain \emph{correctness}
property. Below we explain this correctness property. It turns out
that it is related to our notion of good parameters.


\skipmed

\noindent
\textbf{Good Parameters and Correctness}

\skipsmall

\label{section:goodp}

Recall that if $\alpha$ is projective-like, then
a good parameter for $\JalphaR$ is a real $x$ such that there is
a total function $f:\R\cofmap\omega\alpha$ which is
$\Sigma_1(\JalphaR,\singleton{x})$.
In this section we study good parameters and their relationship
to the correctness of a given set of reals.
We work in the theory ZF + DC.

We start with an alternative characterization of good parameters.

\begin{definition}
Let $x\in\R,\alpha\geq2, n\geq 0$. Then we say that
$x$ is an {\em $n$-name for $\alpha$} iff there is a 
formula $\theta$ such that $\alpha$ is the least ordinal so that
$\JalphaR\models\theta[x]$, and $\theta$ is a $\Sigma_{(\bullet, n)}$
formula if $n$ is even, or a $\Pi_{(\bullet, n)}$ formula if
$n$ is odd.
\end{definition}

\begin{proposition}
Let $\alpha\geq 2$ be projective-like. Let $x\in\R$. Then the following
are equivalent:
\begin{itemize}
\item[(1)] $x$ is a good parameter for $\JalphaR$.
\item[(2)] $x$ is a 1-name for $\alpha$
\end{itemize}
\end{proposition}
\begin{proof}

$(1) \Implies (2)$. \quad Let $g:\R\cofmap\omega\alpha$ be a total function
which is $\Sigma_1(\JalphaR,\singleton{x})$. Let $\formulaphi$ be a
$\Sigma_1$ formula so that $g(y) = \beta$ iff
$\JalphaR\models \formulaphi[x,y,\beta]$.  Then $\alpha$ is least
so that 
$\JalphaR\models(\forall y\in\R)(\exists\beta)\formulaphi[x,y,\beta]$.

$(2) \Implies (1)$. \quad  Let $\formulaphi$ be a $\Sigma_1$ formula such 
that
$\alpha$ is least so that $\JalphaR\models(\forall y\in\R)\formulaphi[x,y]$.
Let $g:\R\map\omega\alpha$ be defined by
$$g(y) = (\mu\,\beta)\big[\SbetaR\models%
 \formulaphi[x,(y)_0]\:\big]+(y)_1(0).$$
Then $g\in\Sigma_1(\JalphaR,\singleton{x})$ and $\ran(g)$ is cofinal
in $\omega\alpha$.
\end{proof}

Good parameters play an important role in the
{\em correctness} of sets of reals.  We turn to this topic now.  We
begin with a definition. If $A\subseteq\R$ is any set of reals, and
$\alpha\geq 1$ is any ordinal, then $J_{\alpha}(A)$ is defined similarly
to the way we defined $\JalphaR$ an page \pageref{def:JalphaR}.

\begin{definition}
\label{Correct}
Let $\alpha\geq1, n\geq1, A\subseteq\R, A\not=\emptyset$. Then we
say that $A$ is $\Sigma_n(\JalphaR)$-{\em correct} iff there is
an $\bar{\alpha}\geq1$ such that
\begin{itemize}
\item[(1)] $\R\intersect J_{\bar{\alpha}}(A) = A$\qquad and
\item[(2)] $\exists\;j:(J_{\bar{\alpha}}(A);\,\in,\,A)\map(\JalphaR;\,\in,\,\R)$
such that $j$ is $\Sigma_n$ elementary.
\end{itemize}
\end{definition}

If $\Gamma$ is a pointclass and $A$ a pointset then we say that $A$
is a {\em basis} for $\Gamma$ iff $B\intersect A\not=\emptyset$ for
all $B\in\Gamma$.
\begin{lemma}
Suppose $\alpha$ begins a $\Sigma_1$-gap.  
Let $A$ be a nonempty set of reals.
Then the following are equivalent:
\begin{itemize}
\item[(a)] $A$ is $\Sigma_1(\JalphaR)$-correct.
\item[(b)] Letting
$$H=\bigsetof{a\in\JalphaR}%
{(\exists x\in A)\; \singleton{a}\in\Sigma_1(\JalphaR,\singleton{x})}$$
we have that $H\prec_1\JalphaR$ and $\R\intersect H = A$.
\item[(c)] $A$ is closed under join, and for all $x\in A$, $A$ is
a basis for $\Sigma_1(\JalphaR,\singleton{x})$.
\end{itemize}
\end{lemma}
\begin{proof}
Clearly $(b)\Implies(a)\Implies(c)$.  We will show that $(c)\Implies(b)$.
Let \mbox{$f\vdots\R\surjection\JalphaR$} 
be a partial $\Sigma_1(\JalphaR)$ map.
Let $X = f\text{``}A$.  Then $X \subseteq H$.  But (c) easily implies
that $X\prec_1\JalphaR$.  Thus $X=H$ and so $H\prec_1\JalphaR$.  But
(c) also easily implies that $\R\intersect X = A$.
\end{proof}

\begin{remark}
If $\alpha$ begins a $\Sigma_1$-gap and
$A$ is $\Sigma_1(\JalphaR)$-correct then the $\bar{\alpha}$ and the
$j$ of definition \ref{Correct} are unique. [Proof: Given any such
$\bar{\alpha}$ and $j$ let $X=\ran{j}$. Then $X\prec_1\JalphaR$
so $H\subseteq X$, where $H$ is as in the lemma.  But let
$f\vdots\R\surjection\JalphaR$ be a partial $\Sigma_1(\JalphaR)$ map.
Let $a\in X$. Since
$\JalphaR\models\text{``}(\exists x\in\R)\;f(x) = a\text{''}$,
there is an $x$ in $A$ such that $f(x)=a$.  But then $a\in H$.
So $ X = H$.]  We will use the notation $j_A$ to refer to this
unique j.
\end{remark}
\begin{corollary}
\label{ClosureIsCorrectness}
Suppose $\alpha$ begins a $\Sigma_1$-gap.  Let $n\geq2$.
Let $A$ be a non-empty set
of reals.  Then the following are equivalent:
\begin{itemize}
\item[(a)] $A$ is $\Sigma_n(\JalphaR)$-correct.
\item[(b)] $A$ is $\Sigma_1(\JalphaR)$-correct and $j_A$ is $\Sigma_n$
elementary.
\item[(c)] $A$ is closed under join and $(\forall x\in A)$ $A$ is a
basis for $\Sigma_n(\JalphaR,\singleton{x})$.
\end{itemize}
\end{corollary}
\begin{proof}
Clearly $(b)\Implies(a)\Implies(c)$. Then $(c)\Implies(b)$ follows from the
proof of the previous lemma.
\end{proof}

In the analytical hierarchy, if $n\geq2$ is even and a non-empty set
of reals $A$ is closed under join and closed downward under 
$\Delta^1_n$, then $A$ is $\Sigma^1_n$-correct. 
This is because, for even $n\geq 2$,
 $\Sigma^1_n$ has the basis property.
 In our setting, for even $n\geq 2$,
$\San$ has the basis property. (See Lemma \ref{basislemma}.)
However, the next two lemmas show that
in our setting a good parameter is also required for correctness.
The first lemma handles the case $n=0$, and the second lemma
handles the case $n\geq 2$ even.
\begin{lemma}
\label{OneCorrectness}
Suppose $\alpha$ begins a $\Sigma_1$-gap. Assume $\Det(\JalphaR)$.
Let $A$ be a non-empty set
of reals.  Suppose $A$ is closed under join and closed downward under
$\Da{0}$.
\begin{itemize}
\item[(1)] If $\alpha$ is a limit ordinal then $A$ is
$\Sigma_1(\JalphaR)$-correct.
\item[(2)] If $\alpha=\beta+1$ and there is no $x\in A$ such that
$x$ is a 0-name for $\alpha$,
then $A$ is $\Sigma_1(\JalphaR)$-correct.
\item[(3)] If $\alpha=\beta+1$ and there is an $x\in A$ such that
$x$ is a 0-name for $\alpha$, {\em and}
$\beta$ is projective-like, then the following are equivalent:
\begin{itemize}
\item[(a)] $A$ is $\Sigma_1(\JalphaR)$-correct.
\item[(b)] $A$ is $\Sigma_{\omega}(\JbetaR)$-correct.
\item[(c)] $(\exists z\in A)$ $z$ is a good parameter for $\JbetaR$.
\end{itemize}
\end{itemize}
\end{lemma}
\begin{proof} \hspace*{1in}\\
\noindent
(1) and (2).\qquad Let $x\in A$ and let $\formulaphi$ be a $\Sigma_1$
formula.  Suppose that $\JalphaR\models(\exists y\in\R)\formulaphi[x,y]$.
Let $\gamma$ be the least ordinal such that
$\JofR{\gamma}\models(\exists y\in\R)\formulaphi[x,y]$.  In either case
(1) or case (2) we have that
 $\gamma<\alpha$.  Notice then that $\singleton{\gamma}
\in\Sigma_1(\JalphaR,\singleton{x})$.  Also $\gamma$ is projective-like.
Let $B=\setof{y\in\R}{\JofR{\gamma}\models \formulaphi[x,y]}$.  Then
$B\in\Sa[\gamma]{0}(x)\subseteq\Sa[\gamma]{2}(x)$.
By the Basis Property for $\Sa[\gamma]{2}$,
 $(\exists y\in B)$ with $y\in\Da[\gamma]{2}(x)$. Since
$\singleton{\gamma} \in\Sigma_1(\JalphaR,\singleton{x})$ we have that
$\singleton{y}\in\Sigma_1(\JalphaR,\singleton{x})$. 
So $y\in A$.  This shows that
$(\forall x\in A)$ $A$  is a basis for $\Sigma_1(\JalphaR,\singleton{x})$.
So $A$ is $\Sigma_1(\JalphaR)$-correct.

\noindent
(3).  Let $\alpha,\beta,x$ be as in (3). Notice that
$\singleton{\beta}\in\Sigma_1(\JalphaR,\singleton{x})$.

\noindent
$(a)\Implies(b)$.\quad
Since $\singleton{\beta}\in\Sigma_1(\JalphaR,\singleton{x})$,
we have that
$\Sigma_{\omega}(\JbetaR)\subseteq\Sigma_1(\JalphaR,\singleton{x})$. So
if $A$ is $\Sigma_1(\JalphaR)$-correct, then
$(\forall w\in A)$ $A$ is a basis for $\Sigma_1(\JalphaR,\singleton{x,w})$,
so $(\forall w\in A)$
$A$ is a basis for $\Sigma_{\omega}(\JbetaR,\singleton{w})$.

\noindent
$(b)\Implies(c)$.\quad This follows from the fact that
$$\JbetaR\models(\exists z)(z \text{ is a good parameter for }\JbetaR).$$

\noindent
$(c)\Implies(a)$.\quad Let $z$ be as in (c).  Let $w\in A$. Let
$B$ be a non-empty set of reals with $B\in\Sigma_1(\JalphaR,\singleton{w})$.
Then there is a $\Bprime\subseteq B$, $\Bprime\not=\emptyset$,
$\Bprime\in\Sigma_{\omega}(\JbetaR,\singleton{w})$.  Let $n\in\omega$
with $\Bprime\in\Sigma_{n}(\JbetaR,\singleton{w})$, and $n$ odd.
Then $\Bprime\in\Sa[\beta]{n-1}(\angles{z,w})$.
As $n-1$ is even, by the Basis Property for $\Sa[\beta]{n-1}$,
$(\exists y\in\Bprime)$
$y\in\Da[\beta]{n-1}(\angles{z,w})$.  As
$\singleton{\beta}\in\Sigma_1(\JalphaR,\singleton{x})$, we have that
$y\in\Sigma_1(\JalphaR,\singleton{x,z,w})$. So $y\in A$.
This shows that
$(\forall w\in a)$ $A$  is a basis for $\Sigma_1(\JalphaR,\singleton{w})$.
So $A$ is $\Sigma_1(\JalphaR)$-correct.
\end{proof}

\begin{remark}
Notice that the one case not covered by Lemma \ref{OneCorrectness}
is the case in which
$\alpha=\beta+1$, and there is an $x\in A$ such that $x$ is a 0-name
for $\alpha$,  but $\beta$ is not
projective-like.  For example, let $\beta=\kappa^{\R}=$ the least ordinal
$\kappa$ such that $\JofR{\kappa}$ is admissible, and let
$\alpha=\beta+1$. Then 0 is a 0-name for $\alpha$ and $\beta$ is
not projective-like. In this case there is a set of reals $A$ such that
$A$ is closed under join and closed downward under $\Da{0}$, and $A$
is not $\Sigma_1(\JalphaR)$-correct. Since we are not interested in
this case in this paper, we do not pursue this here.
\end{remark}
\begin{lemma}
\label{nCorrectness}
Let $\alpha\geq 2$ be projective-like. Assume $\Det(\JofR{\alpha+1})$.
Let $n\geq2$ be even.  Let $A$ be a non-empty set of
reals, closed under join and closed downward under $\Dan$. Then the
following are equivalent.
\begin{itemize}
\item[(a)] $A$ is $\Sigma_{n+1}(\JalphaR)$-correct.
\item[(b)] $(\exists x\in A)$ $x$ is a good parameter for $\JalphaR$.
\end{itemize}
\end{lemma}
\begin{proof} \hspace*{1in}\\
\noindent
$(a)\Implies(b)$.\qquad This follows from the fact that the set of good
parameters for $\JalphaR$ is $\Pi_2(\JalphaR)$.

\noindent
$(b)\Implies(a)$.\qquad Let $x$ be as in (b). Now let $z\in A$
and let $B$ be a non-empty set of reals with
$B\in\Sigma_{n+1}(\JalphaR,\singleton{z})$.  Then
$B\in\San(\angles{x,z})$. By the Basis Property for
$\San$,
 $(\exists y\in B)\;y\in\Dan( \angles{x,z})$.
So $B\intersect A\not=\emptyset.$
\end{proof}

Suppose that $A\subset\R$ and $A$ is $\Sigma_n(\JalphaR)$-correct.
Then there is some ordinal $\bar{\alpha}\geq1$ such that
\begin{itemize}
\item[(1)] $\R\intersect J_{\bar{\alpha}}(A) = A$\qquad and
\item[(2)] $\exists\;j:(J_{\bar{\alpha}}(A);\,\in,\,A)\map(\JalphaR;\,\in,\,\R)$
such that $j$ is $\Sigma_n$ elementary.
\end{itemize}
But how can we find the ordinal $\bar{\alpha}$ given the
set $A$? It turns out that this question will be important to us
when we apply our results about correctness in Section 
\ref{section:correctness}. In this paper we will only be dealing
with ordinals $\alpha<\omega_1^{\omega_1}$.  In this special case, the
problem of finding $\bar{\alpha}$  is simpler. This is because
every ordinal $\alpha<\omega_1^{\omega_1}$ is coded by a finite
set of countable ordinals, and every countable ordinal is coded by
a real. We turn to a discussion of this topic now.


\skipmed

\noindent
\textbf{A Coding Scheme}

\skipsmall

\label{section:coding}

We will need a scheme for coding an ordinal $\alpha<\omega_1^{\omega_1}$
using countable ordinals.  Such an $\alpha$ is easily 
coded by a finite function on the countable ordinals.

\begin{definition}
\label{def:CodeSet}
Let $\CodeSet$ be the set of all finite partial functions from $\omega_1$
to $\omega_1-\singleton{0}$. $\CodeSet$ will be our set of codes for
ordinals $\alpha<\omega_1^{\omega_1}$. If $s\in\CodeSet$ let $\bar{s}$
denote the total function $\bar{s}:\omega_1\map\omega_1$ extending
$s$ such that $\bar{s}(\xi)=0$ if $\xi\not\in\dom(s)$. For $s,t\in\CodeSet$
say that $s<t$ if $\bar{s}(\xi)<\bar{t}(\xi)$, where $\xi$ is  greatest
so that $\bar{s}(\xi)\not=\bar{t}(\xi)$. It is easy to see that
$(\CodeSet,<)$ is a wellorder of order type $\omega_1^{\omega_1}$.
For $s\in\CodeSet$ let $|s|$ = the rank of $s$ in $(\CodeSet,<)$.
For $\alpha<\omega_1^{\omega_1}$ let $\code(\alpha)$ = the unique
$s\in\CodeSet$ such that $|s|=\alpha$.  One more piece of notation:
For $s\in\CodeSet$ let $\max(s)$ be the largest countable ordinal
in $\dom(s)\union\ran(s)$.
\end{definition}

Below we record some trivial facts about our coding scheme which we
will need to use later.

\begin{lemma}\hspace*{1in}\\
\label{CodeFacts}
\skipsmallminus
\begin{itemize}
\item[(a)] Let $P$ be the relation on $\CodeSet$ defined by
$P(s,t)\Iff |s|<|t|$. Then $P$ is $\Sigma_0$ definable.
\item[(b)] Let $R$ be the relation on $\CodeSet$ defined by
$R(s)\Iff \cof(|s|)>\omega$. Then $R$ is $\Sigma_0$ definable.
\item[(c)] Let $\alpha<\omega_1^{\omega_1}$. For all 
$\xi<\max(\code(\alpha))$ there is an $\alpha_0<\alpha$ such
that for all $\alphaprime$ with $\alpha_0\leq\alphaprime<\alpha$,
$\xi\leq\max(\code(\alphaprime))$.
\item[(d)] $\max(\code(\alpha+1))\leq\max(\code(\alpha))+1$.
\item[(e)] Let $\alpha<\omega_1^{\omega_1}$ be a limit ordinal
of cofinality $\omega$.  Then there are cofinally many $\alphaprime<\alpha$
such that $\max(\code(\alphaprime))\leq\max(\code(\alpha))$.
\end{itemize}
\end{lemma}
\begin{proof}\hspace*{1in}\\

\skipsmallminus

(a) $P(s,t)\Iff|s|<|t|\Iff s<t$ $\Iff$ 
$\xi\not\in\dom(s) \OR s(\xi)<t(\xi)$,
where $\xi$ is the largest ordinal in $\dom(s)\union\dom(t)$ such that
either $\xi\not\in\dom(s)\intersect\dom(t)$ or
$s(\xi)\not=t(\xi)$.

(b) $R(s)$ $\Iff$  $\xi$ is a successor ordinal $\AND$ $s(\xi)$ is a 
successor ordinal, where $\xi$ is the least ordinal in $\dom(s)$.

(c) Let $s=\code(\alpha)$. Fix $\xi<\max(s)$.  
We may assume that $\xi>0$. Let $\eta$ be
least such that $\eta\in\dom(s)$. 

\begin{case}{1}
$\xi<\eta$.
\end{case}
Let $s_0=s-\singleton{\bigl(\eta,s(\eta)\bigr)}\union
\singleton{(\xi,1)}$
Let $\alpha_0=|s_0|$. This works.

\begin{case}{2}
$\eta\leq\xi<s(\eta)$.
\end{case}
Let $s_0=s-\singleton{\bigl(\eta,s(\eta)\bigr)}\union
\singleton{(\eta,\xi)}$
Let $\alpha_0=|s_0|$. This works.

\begin{case}{3}
$\eta\leq\xi, s(\eta)\leq\xi$.
\end{case} 
Let $s_0=s-\singleton{\bigl(\eta,s(\eta)\bigr)}$. Let 
$\alpha_0=|s_0|$. This works.

(d) This is a particular case of (c).

(e) Let $s\in\CodeSet$ be a code. Then $\alpha = |s|$ is a limit ordinal of
cofinality $\omega$ iff $\xi$ is a limit ordinal or $s(\xi)$ is a
limit ordinal, where $\xi$ is the least ordinal in $\dom(s)$.
If $s(\xi)$ is a limit ordinal then
$\dotsetof{|s-\singleton{\bigl(\xi,s(\xi)\bigr)}\union%
\singleton{(\xi,\eta)}|}{\eta<s(\xi)}$ is a cofinal set of ordinals
$\alphaprime<\alpha$ such that
$\max(\code(\alphaprime))\leq\max(\code(\alpha))$.
If $\xi$ is a limit ordinal and $s(\xi)=\eta+1$ then
$\dotsetof{|s-\singleton{\bigl(\xi,s(\xi)\bigr)}\union%
\singleton{\bigl(\xi,\eta\bigr)}\union
\singleton{(\gamma,1)}|}{\gamma<\xi}$ works.
\end{proof}

We conclude this section with a lemma about correctness, in the
special case that $\alpha<\omega_1^{\omega_1}$. This is the
only result about correctness which we will need in the sequel.

 If $M$ is a transitive model of $\ZF$ and $s\in\CodeSet$ and
$\max(s)<\omega_1^{M}$, then
we may interpret the term ``$|s|$'' in $M$. Clearly
$|s|^{M}$ is just the rank of $s$ in the wellorder
$(\CodeSet^{M},<)$, where $\CodeSet^{M}$ is
$\setof{s\in\CodeSet}{\max(s)<\omega_1^{M}}$.  We will use this
notion in the next lemma.

\begin{lemma}
\label{MaxCodeIsGoodParameter}
Let $\alpha<\omega_1^{\omega_1}$, and let $n\geq0$ be even.
Assume $\Det(\JofR{\alpha+1})$. Let $A\subseteq\R$ be such that
$A$ is closed downward under $\Dan$, and $A=\R\intersect L(A)$.
Let $s=\code(\alpha)$ and suppose that $\max(s)<\omega_1^{L(A)}$. 
Let $\bar{\alpha}=|s|^{L(A)}$. Then there exists a $\Sigma_{n+1}$ elementary
embedding
$$j:J_{\bar{\alpha}}(A)\map\JalphaR$$
\end{lemma}
\begin{proof}
We begin by showing that the conclusion of the lemma is true for \emph{some}
ordinal $\bar{\alpha}$. Then we will show that it is true for our particular
$\bar{\alpha}$.  The fact that the conclusion of the lemma is true for
some $\bar{\alpha}$ is just the fact that $A$ is 
$\Sigma_{n+1}(\JalphaR)$-correct. So we begin by showing this. Since
we are assuming that $n$ is even and $A$ is closed under join and
closed downward under
$\Dan$, Lemmas \ref{OneCorrectness} and \ref{nCorrectness} above
 indicate that we must show that there is some good parameter  in $A$.
 We will show that there is one---namely any real
coding $s$.

Recall that $\WO$ is the  set of reals coding countable
ordinals. $\WO$ is a $\Pi^1_1$ set. If $w\in\WO$ then $|w|$ is the
countable ordinal coded by $w$. (Do not be confused by the fact
that if $s\in\CodeSet$ then we are also using the notation $|s|$ to
refer to the ordinal $\alpha<\omega_1^{\omega_1}$ such that
$\code(\alpha)=s$.)
Let $\CodeSet^*$ be
the set of finite partial functions $h$ from $\WO$ to $\WO$ so that for
no $w$ in the range of $h$ is it true that $|w|=0$.
So $\CodeSet^*$ is essentially a set of reals, and  $\CodeSet^*$
is a collection of codes for the elements of $\CodeSet$.
Let $\pi:\CodeSet^*\map\CodeSet$ be the obvious decoding function.
For $x,y\in\CodeSet^*$, say that
$x<^*y \Iff \pi(x) < \pi(y)$. Thus $(\CodeSet^*,<^*)$ is a pre-wellorder
of ordertype $\omega_1^{\omega_1}$.
Notice that $(\CodeSet^*,<^*)$ is lightface definable over
$\JofR{1}$.
Let 
$\operatorname{Rnk}:\CodeSet^*\map\omega_1^{\omega_1}$  be the rank 
function.
So $\operatorname{Rnk}(x)=|\pi(x)|$.
For $\xi<\omega_1^{\omega_1}$ let us use the abbreviation 
$\operatorname{Rnk}\restr\xi$ to 
mean $\operatorname{Rnk}\restr\setof{x}{\operatorname{Rnk}(x)<\xi}$.

\begin{claim}[Claim 1]
Let $\beta\geq 2$.  Then
\begin{itemize}
\item[(a)] $\operatorname{Rnk}\restr\ORD^{\JbetaR}$ is $\Sigma_1(\JbetaR)$.
\item[(b)] For all $\xi\in\JbetaR$, $\operatorname{Rnk}\restr\xi \in \JbetaR$.
\end{itemize}
\end{claim}
\begin{subproof}[Proof of Claim 1] 
This is a standard induction on $\beta$.
\end{subproof}

Since $\max(s)<\omega_1^{L(A)}$, there is an $x_0\in L(A)$ such
that $x_0\in\CodeSet^*$ and $\pi(x_0)=s$. Thus
$\operatorname{Rnk}(x_0)=\alpha$.
By  a simple recursive coding,
we may think of $x_0$ as a real. Thus $x_0\in A$.
Fix such a real $x_0$.

\begin{claim}[Claim 2]
$x_0$ is a good parameter for $\JalphaR$.
\end{claim}
\begin{subproof}[Proof of Claim 2]
$\alpha$ is least such that $\JalphaR$ satisfies:
$$(\forall y\in \CodeSet^*)\bigl[y<^*x_0 \Implies 
(\exists \beta,\; \exists X)
(\beta=\operatorname{Rnk}(y) \AND X=\JbetaR)\,\bigr].$$
Thus $x_0$ is a 1-name for $\alpha$ and so $x_0$ is a good parameter for
$\JalphaR$.
\end{subproof}

We want to show that $A$ is $\Sigma_{n+1}(\JalphaR)$-correct.  
First suppose that $n\geq 2$. Then we have that $A$ is closed downward
under $\Dan$, and $A$ is closed under join,
 and there is a good parameter for $\JalphaR$ in $A$.
It follows from lemma \ref{nCorrectness}
that $A$ is $\Sigma_{n+1}(\JalphaR)$-correct.  
Next suppose
that $n=0$.  We have that
$A$ is closed downward under $\Da{0}$ and $A$ is closed under join.
If $\alpha$ is a limit ordinal, then  it follows from part (a) of
Lemma \ref{OneCorrectness} that $A$ is $\Sigma_1(\JalphaR)$-correct.
Suppose then that $\alpha$ is a successor
ordinal. Since $\max(\code(\alpha))< \omega_1^{L(A)}$, it is also
true that $\max(\code(\alpha-1))< \omega_1^{L(A)}$. Thus there is
an $x\in L(A)$ such
that $x\in\CodeSet^*$ and $\operatorname{Rnk}(x)=\alpha-1$. As in the
claim above, $x$ is a good parameter  for $\JofR{\alpha-1}$.
It follows from parts (b) and (c) of Lemma \ref{OneCorrectness}
that $A$ is $\Sigma_1(\JalphaR)$-correct. 

So we now have that
in all cases $A$ is $\Sigma_{n+1}(\JalphaR)$-correct.
Thus there
is a $\Sigma_{n+1}$ embedding:
$$j:J_{\gamma}(A)\map\JalphaR.$$
We only need to see that $\gamma=\bar{\alpha}$.
Recall that  $s=\code(\alpha)$, and  $\bar{\alpha}=|s|^{L(A)}$,
and $x_0\in A$, and $\pi(x_0)=s$.

We split into three cases.

\begin{case}{1}
$\alpha<\omega_1$
\end{case}

Then $s=\singleton{(0,\alpha)}$. Thus $\bar{\alpha}=\alpha$.
Also, $x_0$ is essential an element of $\WO$, with $|x_0|=\alpha$.
Notice that if $w\in\WO\intersect A$, and $|w|=\xi\in J_{\gamma}(A)$,
then $|w|^{J_{\gamma}(A)}=\xi$. Furthermore, the 
relation $\xi=|w|$ is $\Sigma_1(J_{\gamma}(A))$.
Since $x_0\in A$,  
and $j$ is at least $\Sigma_1$ elementary,
it is easy to see that $\gamma=\alpha$.

\begin{case}{2}
$\alpha=\omega_1$
\end{case}

Then $s=\singleton{(1,1)}$ and so $\bar{\alpha}=\omega_1^{L(A)}$.
Since $j$ is at least $\Sigma_1$ elementary, it is easy to see
that $\gamma=(\omega_1)^{J_{\gamma}(A)}=\omega_1^{L(A)}$.

\begin{case}{3}
$\alpha>\omega_1$
\end{case}

Then $\gamma>\omega_1^{J_{\gamma}(A)}=\omega_1^{L(A)}$.
As $\max(s)<\omega_1^{L(A)}$, $s\in J_{\gamma}(A)$.
Notice that  $j(s)=s$. First suppose towards a contradiction that
$\gamma>\bar{\alpha}$. In particular $\bar{\alpha}\in J_{\gamma}(A)$.
Then $\bar{\alpha}=|s|^{J_{\gamma}(A)}$,
and so $j(\bar{\alpha})=\alpha$. But this
is a contradiction as $J_{\gamma}(A)$ satisfies that
$\JofR{\bar{\alpha}}$ exists. So $\gamma\leq\bar{\alpha}$.
Now suppose towards a contradiction that $\gamma<\bar{\alpha}$.
Then there is a $t\in\CodeSet$ with $\max(t)<\omega_1^{J_{\gamma}(A)}$ 
such that $|t|^{J_{\gamma}(A)}=\gamma$. 
Notice that $j(t)=t$. Let $\delta=|t|$ in $V$.
Thus $j(\gamma)=\delta$.
But this contradicts the fact that $\JalphaR$ satisfies that
$\JofR{\delta}$ exists.
\end{proof}

%


\skipbig

\section{Big Premice}

\label{section:bigmice}

In this section we finally get to inner model theory. This section
marks the beginning of the main portion of this paper. 

In the previous section we defined the sets $\Aan$ and we stated
the conjecture that $\Aan$ is a mouse set. We now begin the process
of proving this conjecture. In this section we will define the notion
of a premouse $\cM$ being $(\alpha,n)$-big.
 In the next section we will show that if $\cM$ is a countable,
iterable premouse, and $\cM$ is $(\alpha,n)$-big, then
$\Aan\subseteq\R\intersect\cM$. This will give us one half of the
proof that $\Aan$ is a mouse set. 

Our approach here is to generalize some work of
Martin and Steel. We begin by describing this work.
The following is (equivalent to) a definition from \cite{MaSt}. 
We repeat it here for the readers convenience.

\begin{definition}
\label{alphaisone}
Let $\cM$ be a premouse and $\beta\in\ORD^{\cM}$. Let $n\geq0$.  Then
$\cM$ is \emph{$n$-big above $\beta$} iff there is an \emph{active}
initial segment
$\cN\unlhd\cM$ such that, letting $\kappa$ be the critical point of the
last extender from the $\cN$-sequence, we have that
$\beta<\kappa$ and there are
$n$ ordinals $\delta_1, \dots ,\delta_{n}$,
with $\beta<\delta_1 < \dots < \delta_n < \kappa$, 
such that each $\delta_i$ is a Woodin cardinal of  $\cN$.

$\cM$ is \emph{$n$-small} above $\beta$ iff $\cM$ is not $n$-big
above $\beta$. We will say that $\cM$ is $n$-big iff $\cM$ is
$n$-big above $0$, equivalently iff $\cM$ is $n$-big above
$\omega$.  Similarly with $n$-small.
\end{definition}

In order to parallel our definition of $\Aan$ for $\alpha\geq 2$,
we make the following definition with $\alpha=1$.
For $n\geq 0$ we set
$$\Aa[1]{n}=\setof{x\in\R}{(\exists\xi<\omega_1)\,x\in\Delta^1_{n+2}(\xi)}.$$
So $\Aa[1]{n}$ is the same set that we were calling $A_{n+2}$ in 
the introduction. We have changed the indexing here in order to
better express the analogy between $\alpha=1$ and $\alpha\geq 2$.
The following theorem and its proof are implicit in  \cite{MaSt}.

\begin{theorem}[Martin, Steel, Woodin]
\label{SteelThm1}
Let $\cM$ be a countable, iterable premouse and 
$\beta\in\ORD^{\cM}$. Let $n\geq0$, and suppose that $\cM$ is
$n$-big above $\beta$. Let $\P\in\cM$ be a partial order, with
$(\card(\P)^+)^{\cM}\leq\beta$.
 Suppose that $G$ is $\cM$-generic over 
$\P$. Then for all $x\in\R\intersect\cM[G]$,
$\Aa[1]{n}(x)\subseteq\R\intersect\cM[G]$. 
In particular, if $\cM$ is $n$-big
then $\Aa[1]{n}\subseteq\R\intersect\cM$.
\end{theorem}

Our goal in this section is to give a definition which extends Definition
\ref{alphaisone} above.  In order to emphasize this, let us give
a synonym for the term $n$-big.  For $n\geq 0$ we will say that a 
premouse is $(1,n)$-big above $\beta$ iff it is $n$-big above $\beta$.
Now let $\cM$ be a premouse and $\beta\in\ORD^{\cM}$. Let $\alpha\geq 2$, 
and $n\geq 0$. We wish  to define the notion: $\cM$ is 
\emph{$(\alpha,n)$-big above $\beta$}.
To motivate the definition, we list five criteria that we would like for
the definition to satisfy.  The definition of $n$-big (Definition
\ref{alphaisone} above) satisfies the analogous criteria.
\begin{itemize}
\item[(1)] If a countable, iterable premouse $\cM$ is $(\alpha,n)$-big above
$\beta$, and if $\P\in\cM$ is a partial order with
$(\card(\P)^+)^{\cM}\leq\beta$, and if
$G$ is $\cM$-generic over $\P$, then for all
$x\in\R\intersect\cM[G]$,  $\Aa{n}(x)$ should be contained
in $\cM[G]$.
\item[(2)] If $\cN$ is $(\alpha,n)$-big above $\beta$ and $\cN\unlhd\cM$
then $\cM$ should be $(\alpha,n)$-big above $\beta$.
\item[(3)] If $\cM$ is $(\alpha,n)$-big above $\beta$ and $\betaprime<\beta$
then $\cM$ should be $(\alpha,n)$-big above $\betaprime$.
\item[(4)] If $\cM$ is $(\alpha,n)$-big above $\beta$ and 
$(\alphaprime,\nprime)\lexless(\alpha,n)$ then
$\cM$ should be $(\alphaprime,\nprime)$-big above $\beta$.
\item[(5)] The property of being $(\alpha,n)$-big should be preserved by 
non-dropping iterations.  That is, if $\cM$ is $(\alpha,n)$-big 
above $\beta$, and if we iterate $\cM$ without dropping to yield
$i:\cM\map\cMprime$, then $\cMprime$ should be $(\alpha,n)$-big
above $i(\beta)$.
\end{itemize}
The fifth criterion is motivated by the proof of Theorem \ref{SteelThm1}
above, which we hope to extend.  This criterion turns out to be the
most problematic for us.  In order to achieve criterion (5), we will aim
for the following stronger criterion:
\begin{itemize}
\item[$(5^{\prime})$] If $\cM$ is $(\alpha,n)$-big above $\beta$ then $\cM$
should have a code for $\alpha$, and the statement ``I am $(\alpha,n)$-big
above $\beta$'' should be definable in $\cM$.
\end{itemize}
Right away we see that there is a conflict between criterion $(5^{\prime})$
and criterion $(4)$. For example, suppose $\cM$ is a countable premouse
and $\cM$ is $(\omega_1,0)$-big above $\beta$.  By (4), $\cM$ should be 
$(\alpha,0)$-big above $\beta$ for every countable ordinal $\alpha$.
Then by $(5^{\prime})$, $\cM$ should have a code for every countable
ordinal $\alpha$.  But this is impossible as $\cM$ is countable.

Our solution to this conflict is the following.  We will define an
auxiliary notion: $\cM$ is \emph{explicitly $(\alpha,n)$-big above $\beta$}.
For this auxiliary notion we will drop criterion (4).  So a premouse
can be explicitly $(\omega_1,0)$-big without being explicitly 
$(\alpha,0)$-big for every countable $\alpha$. This auxiliary notion will
satisfy criteria (1), (2), (3), and $(5^{\prime})$.
In place of criterion (4)
we will have the following:
\begin{itemize}
\item[$(4^{\prime})$] If $\cM$ is explicitly $(\alpha,n)$-big above 
$\beta$ and $(\alphaprime,\nprime)\lexless(\alpha,n)$ then there is
\emph{an iterate}
of $\cM$ which is explicitly $(\alphaprime,\nprime)$-big above $\beta$.
\end{itemize}
Then, to recapture criterion (4), we will define as our primary notion: 
$\cM$ is $(\alpha,n)$-big above $\beta$ iff for some
$(\alphaprime,\nprime)$ with $(\alpha,n)\lexleq(\alphaprime,\nprime)$,
$\cM$ is explicitly $(\alphaprime,\nprime)$-big above $\beta$. Then this
primary notion will satisfy criteria (1) through (5). 
(but not $(5^{\prime})$)

There is one major drawback with the programme we have outlined above: we
are only able to see how to carry it out for small ordinals
$\alpha$.  In this paper, for the sake of concreteness, we will carry
out the programme for $\alpha<\omega_1^{\omega_1}$. The ordinal 
$\omega_1^{\omega_1}$ has no special significance here.
 It is easy to see how the programme can be extended somewhat beyond
$\omega_1^{\omega_1}$.  We are not able to go very far beyond though.
The difficulty
seems to be with criterion $(5^{\prime})$ above.  We are unable to see
how to extend the definition of ``explicitly $(\alpha,n)$-big''
to arbitrarily
large ordinals $\alpha$, so as to satisfy criterion $(5^{\prime})$.
We do not know if this is a limitation of the approach, or merely of
the author.  We suspect that it is a little of both.

In order to satisfy criterion $(5^{\prime})$ (described above),
we will need a scheme for coding ordinals which yields a code for
the uncountable ordinal $\alpha$ that may be found inside of a
countable premouse. Recall that in the previous section we
 defined the coding scheme
$(\CodeSet, <)$ for ordinals $\alpha<\omega_1^{\omega_1}$.

We are now ready to begin our definition. For $\cM$ a countable 
premouse and $\beta\in\ORD^{\cM}$, and $2\leq\alpha<\omega_1^{\omega_1}$, 
and $n\geq 0$, we wish  to describe the notion: $\cM$ is 
\emph{explicitly $(\alpha,n)$-big above $\beta$}.
Our definition is by induction on $(\alpha,n)$.  As the basis of
the induction, for $\alpha=1$ and $n\geq0$, we will say that $\cM$
is \emph{explicitly $(1,n)$-big above $\beta$} iff $\cM$ is $n$-big
above $\beta$. (See Definition \ref{alphaisone} above.) Now we turn
to the inductive step.
There are four different cases of the induction, 
and for ease of readability we split these up into four definitions below.

We begin with the case $n=0$. In this context, recall that we have defined
$$\Aa{0} =\setof{x\in\R}%
{(\exists\alphaprime<\alpha)\, x\in\OD^{\JofR{\alphaprime}}}.$$
Naively, in order to inductively arrange criterion (1) above,
we would like to say that 
$\cM$ is explicitly $(\alpha,0)$-big above $\beta$ iff for all 
$\alphaprime<\alpha$ and all $n$, there is some initial segment
$\cN_{(\alphaprime,n)}\unlhd\cM$ such that $\cN_{(\alphaprime,n)}$ 
is explicitly $(\alphaprime,n)$-big above
$\beta$. (In fact this is what we will say for countable ordinals
$\alpha$. For example in the case $\alpha=2$ we will have that 
$\cM$ is explicitly $(2,0)$-big iff for all $n$ there is an initial 
segment of $\cM$
which is $n$-big.) But, as was discussed
above, we can not do this for uncountable ordinals $\alpha$. 
So our actual definition will be an approximation of this. 
There are two cases, depending
on whether or not $\cof(\alpha)>\omega$.  The case 
$\cof(\alpha)\leq\omega$ is easy.

\begin{definition}
\label{cofomeganiszero}
Let $\cM$ be a premouse and $\beta\in\ORD^{\cM}$. Let 
$2\leq\alpha<\omega_1^{\omega_1}$, and
suppose that $\alpha$ is a successor ordinal, or  a limit 
ordinal of cofinality $\omega$. Then $\cM$ is 
\emph{explicitly $(\alpha,0)$-big above $\beta$} iff for cofinally many 
$\alphaprime<\alpha$, and all $n$, there is a proper initial segment
$\cN_{(\alphaprime,n)}\lhd\cM$ such that
$\cN_{(\alphaprime,n)}$ is explicitly $(\alphaprime,n)$-big above $\beta$.
\end{definition}

We next turn to the case $n=0$ and $\cof(\alpha)>\omega$.  
The following definition is motivated by criteria $(4^{\prime})$
and $(5^{\prime})$ above.

\begin{definition}
\label{cofnotomeganiszero}
Let $\cM$ be a premouse and $\beta\in\ORD^{\cM}$. Suppose that 
$\alpha<\omega_1^{\omega_1}$ 
is a limit ordinal of cofinality $\omega_1$. Then $\cM$ is 
\emph{explicitly $(\alpha,0)$-big above $\beta$} iff there is an 
\emph{active} initial segment
$\cN\unlhd\cM$ such that, letting $\kappa$ be 
the critical point of the last extender from the $\cN$-sequence,
we have that
$\beta<\kappa$ and $\code(\alpha)\in J^{\cN}_{\kappa}$, and for all 
$\alphaprime<\alpha$, if $\code(\alphaprime)\in J^{\cN}_{\kappa}$, 
then for all $n$, there is a proper initial segment 
$\cN_{(\alphaprime,n)}\lhd J^{\cN}_{\kappa}$ such that
$\cN_{(\alphaprime,n)}$ is  explicitly $(\alphaprime,n)$-big 
above $\beta$.
\end{definition}

\begin{note}
Let $\cN$ be as in the above definition. 
Then we will say that $\cN$ 
witnesses that $\cM$ is explicitly $(\alpha,0)$-big above $\beta$.
\end{note}

We next turn our attention to the case $n>0$.  This is the successor case
of our inductive definition.  The procedure for passing from explicitly
$(\alpha,n-1)$-big to explicitly $(\alpha,n)$-big is the same in our
context as it is in the case $\alpha=1$. That is, we add a Woodin 
cardinal.

\begin{note}
The case $\cof(\alpha)\leq\omega$, and $n=1$ is more complicated and
is handled separately. We start with the easy case.
\end{note}

\begin{definition}
\label{nplusone}
Let $\cM$ be a premouse and $\beta\in\ORD^{\cM}$. Let 
$2\leq\alpha<\omega_1^{\omega_1}$. 
If $\cof(\alpha)>\omega$ then suppose that
$n\geq 1$. 
If $\cof(\alpha)\leq\omega$ then suppose that $n\geq 2$.
Then $\cM$
is \emph{explicitly $(\alpha,n)$-big above $\beta$} 
iff there is an initial segment
$\cN\unlhd\cM$, and an ordinal $\delta\in\cN$
with $\beta<\delta$, such that $\delta$ is a Woodin cardinal of $\cN$, and
$\cN$ is explicitly $(\alpha,n-1)$-big above $\delta$.
\end{definition}

To emphasize the similarity with Definition \ref{alphaisone} we
state the following easy lemma.

\begin{lemma}
\label{nWoodins}
Let $\cM$ be a premouse and $\beta\in\ORD^{\cM}$. Let 
$2\leq\alpha<\omega_1^{\omega_1}$. 

First suppose that $\cof(\alpha)>\omega$ and $n\geq 1$. 
Then $\cM$ 
is \emph{explicitly $(\alpha,n)$-big above $\beta$} 
iff 
there is an initial segment
$\cN\unlhd\cM$, and $n$ ordinals 
$\delta_1 < \dots <\delta_{n}\in\cN$,
with $\beta<\delta_1$, such that each $\delta_i$ is a Woodin cardinal of 
$\cN$, and $\cN$ is explicitly $(\alpha,0)$-big above $\delta_{n}$.

Next suppose that $\cof(\alpha)\leq\omega$ and $n\geq 2$.
Then $\cM$ 
is \emph{explicitly $(\alpha,n)$-big above $\beta$} 
iff 
there is an initial segment
$\cN\unlhd\cM$, and  $n-1$ ordinals 
$\delta_1 < \dots <\delta_{n-1}\in\cN$,
with $\beta<\delta_1$, such that each $\delta_i$ is a Woodin cardinal of 
$\cN$, and $\cN$ is explicitly $(\alpha,1)$-big above $\delta_{n-1}$.
\end{lemma}

To complete our definition then, we must define the notion of
explicitly $(\alpha,1)$-big in the case that $\alpha$ is a successor 
ordinal,
or a limit ordinal of cofinality $\omega$. It is instructive to consider
the definition that we would have made if we had not excluded this
case from Definition \ref{nplusone} above.  Considering Definition
\ref{cofomeganiszero} above, we see that we would have defined
$\cM$ is explicitly $(\alpha,1)$-big above $\beta$ iff there is an 
initial segment
$\cN\unlhd\cM$, and an ordinal $\delta\in\cN$
with $\beta<\delta$, such that $\delta$ is a Woodin cardinal of $\cN$, and
for cofinally many $\alphaprime<\alpha$, and all $n$, there is
a proper initial segment $\cN_{(\alphaprime,n)}\lhd\cN$ such that
$\cN_{(\alphaprime,n)}$ is explicitly $(\alphaprime,n)$-big 
above $\delta$. This is, at
first sight, a reasonable definition.
But it turns out that a premouse with this property is ``bigger'' 
than necessary.  (We are trying to give the weakest possible definition
which still meets our criteria.) So our
actual definition will be a weakening of this property.  The idea is
that we do not need a single fixed $\delta$ for all of the pairs
$(\alphaprime,n)$.
We can have a different $\delta_{(\alphaprime,n)}$ 
for each pair $(\alphaprime,n)$, and  these
$\delta_{(\alphaprime,n)}$ do not have to be Woodin in all of 
$\cN$, but only in the initial
segment $\cN_{(\alphaprime,n)}\lhd\cN$.  
In exchange for this weakening we must demand that each 
$\delta_{(\alphaprime,n)}$ is a cardinal of $\cN$, and that $\cN$ is a 
model of $\KP$. This brings us to the next two definitions. 

\begin{warning}
Readers of \cite{Ru1} should note that the indexing here is different
than in that paper.  For $\alpha$ a countable ordinal, what we are
here calling ``explicitly $(\alpha,1)$-big'' was called
``$(\alpha,0)$-big'' in \cite{Ru1}. What we are here calling
``explicitly $(\alpha,0)$-big'' did not have a name in \cite{Ru1}.
In \cite{Ru1} we did not deal with
uncountable ordinals $\alpha$, and when one does the indexing used
here is better.
\end{warning}

\begin{definition}
Let $\cM$ be a premouse, $\alpha<\omega_1^{\omega_1}$,
and $n\in\omega$. Let $\delta\in\ORD^{\cM}$. 
Then we will say that $\delta$ is \emph{locally $(\alpha,n)$-Woodin}
in $\cM$ iff there is a proper initial segment 
$\cJ^{\cM}_{\gamma}\lhd\cM$ such that $\delta\in\cJ^{\cM}_{\gamma}$
and
\begin{itemize}
\item[(a)] $\delta$ is a Woodin cardinal of $\cJ^{\cM}_{\gamma}$, and
\item[(b)] $\cJ^{\cM}_{\gamma}$ is explicitly $(\alpha,n)$-big above
           $\delta$.
\end{itemize}
\end{definition}

\begin{note}
In the above situation we will say that $\cJ^{\cM}_{\gamma}$ witnesses
that $\delta$ is locally $(\alpha,n)$-Woodin in $\cM$.
\end{note}

\begin{definition}
\label{cofomeganisone}
Let $\cM$ be a premouse and $\beta\in\ORD^{\cM}$. Suppose that 
$2\leq\alpha<\omega_1^{\omega_1}$ and that $\alpha$ is a successor
ordinal, or a limit ordinal of cofinality $\omega$. Then $\cM$ is
\emph{explicitly $(\alpha,1)$-big above $\beta$} iff 
there  is an initial segment $\cN\unlhd\cM$ with $\beta\in\cN$, and an 
ordinal $\rho\in\cN$ with
$\beta<\rho$, such that $\rho$ is a limit cardinal of $\cN$ and:
\begin{itemize}
\item[(i)] $\cN$ is an admissible set, and
\item[(ii)] $\code(\alpha)\in J^{\cN}_{\rho}$, and
\item[(iii)] for cofinally many $\alphaprime<\alpha$, and all $n\in\omega$,
and cofinally many $\delta<\rho$, we have that:
\begin{itemize}
\item[(a)] $\delta$ is a cardinal of $\cN$, and
\item[(b)] $\delta$ is locally $(\alphaprime,n)$ Woodin 
            in $\cJ^{\cN}_{\rho}$.
\end{itemize}
\end{itemize}
\end{definition}

\begin{note}
Let $\cN$ and $\rho$ be as in the above definition. 
Then we will say that $(\cN, \rho)$ 
witnesses that $\cM$ is explicitly $(\alpha,1)$-big above $\beta$.
\end{note}

\begin{figure}[hbt]
\begin{center}
\input{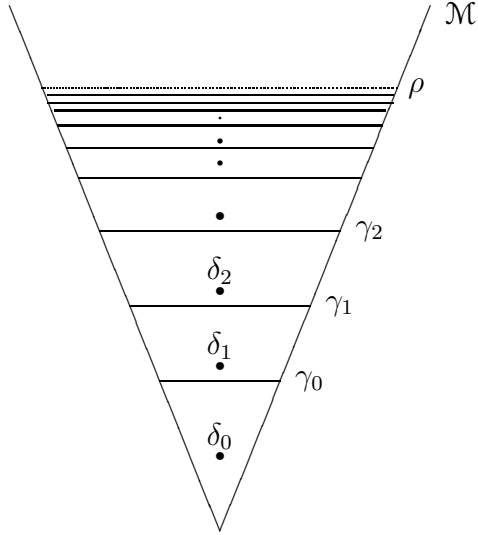}
\caption{$\alpha\geq2$ is a successor and $\cM$ is explicitly 
$(\alpha,1)$-big}
\label{fig:bigpremouse}
\end{center}
\end{figure}

To understand the above definition, we invite the reader to contemplate
what an explicitly
$(2,1)$-big premouse looks like.  In Remark \ref{rem:twoonebig} 
and Figure \ref{fig:bigpremouse} below
we explore the definition of explicitly $(\alpha,1)$-big in the
case that $\alpha$ is a successor ordinal.

\begin{remark}
\label{rem:twoonebig}
Suppose $\alpha=\alphaprime+1$ with $\alphaprime\geq1$.
Suppose that $\cM$ is a countable, iterable premouse and 
$\cM$ is explicitly $(\alpha,1)$-big (above $0$) but no proper initial
segment of $\cM$ is. We will examine what $\cM$ looks like.
(See Figure \ref{fig:bigpremouse} above.) For a particularly concrete
example, consider the case $\alpha=2$.
Because no
proper initial segment of $\cM$ is explicitly $(\alpha,1)$-big,
there is an ordinal $\rho\in\cM$ such that $(\cM,\rho)$ witnesses
that $\cM$ is explicitly $(\alpha,1)$-big. Again, because no
proper initial segment of $\cM$ is explicitly $(\alpha,1)$-big,
we have the following: 
$\rho$ is the largest cardinal of $\cM$,
$\cM$ is a passive premouse, $\cM=J^{\vec{E}}_{\gamma}$, with
$\vec{E}\subset J^{\cM}_{\rho}$, and $\gamma$ is the first admissible ordinal
over $\cJ^{\cM}_{\rho}$. There is a sequence of ordinals 
$\delta_0<\gamma_0<\delta_1<\gamma_1<\delta_2<\dots$ such that
$\max(\code(\alpha))<\delta_0$, and each
$\delta_n$ is a cardinal of $\cM$ 
and a Woodin cardinal of $\cJ^{\cM}_{\gamma_n}$, and
$\cJ^{\cM}_{\gamma_n}$ is explicitly $(\alphaprime,n)$-big above $\delta_n$.
We have that $\sup\singleton{\delta_n}=\sup\singleton{\gamma_n}=\rho$.
(For otherwise let $\rhoprime=\sup\singleton{\delta_n}$. Then, letting
$\cN$ be the first admissible structure over $\cJ^{\cM}_{\rhoprime}$,
we have that $\cN\lhd \cJ^{\cM}_{\rho}$ and
$(\cN,\rhoprime)$ witnesses that $\cN$ is explicitly $(\alpha,0)$-big
above $\beta$.) 

Later we will show that the sequence 
$\lsequence{\delta_n,\gamma_n}{n\in\omega}$ is in $\cM$. This turns out
to be very important. We call this
sequence an \emph{$\alpha$-ladder}.  See Definition \ref{def:ladder}.

Because $\delta_n$ is a (strong limit) cardinal of $\cM$ and $\cM$ is
strongly acceptable, $V^{\cM}_{\delta_n}\subset\cJ^{\cM}_{\gamma_n}$, 
and so $\delta_n$ is actually Woodin in $\cM$ with
respect to all functions $f:\delta_n\map\delta_n$ which are in
$\cJ^{\cM}_{\gamma_n}$. But which functions are in
$\cJ^{\cM}_{\gamma_n}$? The fact that $\cJ^{\cM}_{\gamma_n}$ 
is $(\alphaprime,n)$-big above $\delta_n$ (and that $\cM$ is iterable) 
implies that  $\cJ^{\cM}_{\gamma_n}$ contains all functions 
$f:\delta_n\map\delta_n$ which are (in some appropriate sense) 
$\Delta_{(\alphaprime,n)}$ 
definable with $\cJ^{\cM}_{\delta_n}$ as a parameter. This follows from
the results in the next section.
\end{remark}

This completes the inductive definition of the notion:
$\cM$ is explicitly $(\alpha,n)$-big above $\beta$, for 
$\alpha<\omega_1^{\omega_1}$.
As we remarked earlier, there is no special significance to the
ordinal $\omega_1^{\omega_1}$, and it is easy to see how to extend
our definition somewhat beyond this ordinal.  In fact it will 
actually be
useful for us to extend the definition by one more step. We begin by
extending our coding scheme by one more step.  For uniformity, let us
give a code to the ordinal $\omega_1^{\omega_1}$. We will say
$\code(\omega_1^{\omega_1})=\emptyset$. Compare the
following with Definition \ref{cofnotomeganiszero} above.

\begin{definition}
Let $\cM$ be a premouse and $\beta\in\ORD^{\cM}$. Then $\cM$ is
\emph{explicitly $(\omega_1^{\omega_1},0)$-big above $\beta$} iff
there is an \emph{active} initial segment $\cN\unlhd\cM$ such that, letting
$\kappa$ be the critical point of the last extender from the $\cN$-sequence,
we have that $\beta<\kappa$, and for all $\alpha<\omega_1^{\omega_1}$,
if $\code(\alpha)\in J^{\cN}_{\kappa}$, then for all $n$, there is a
proper initial segment $\cN_{(\alpha,n)}\lhd \cJ^{\cN}_{\kappa}$ such that
$\cN_{(\alpha,n)}$ is explicitly $(\alpha,n)$-big above $\beta$.
\end{definition}

\begin{definition}
Suppose $(1,0)\lexleq(\alpha,n)\lexleq(\omega_1^{\omega_1},0)$.
We will say that $\cM$ is \emph{$(\alpha,n)$-big above $\beta$} iff
$\cM$ is explicitly $(\alphaprime, \nprime)$-big above $\beta$,
for some $(\alphaprime, \nprime)$  with 
$(\alpha,n)\lexleq(\alphaprime, \nprime)$.
We will say that $\cM$ is \emph{$(\alpha,n)$-small} above $\beta$
iff $\cM$ is not $(\alpha,n)$-big above $\beta$.
We will say that $\cM$ is (explicitly) $(\alpha,n)$-big iff $\cM$ is 
(explicitly) $(\alpha,n)$-big above 0, equivalently above $\omega$.
Similarly with $(\alpha,n)$-small.
\end{definition}

In the rest of this section we will study the notions which we have
just defined. In particular we will show that the notion of $(\alpha,n)$-big
satisfies criteria (2) through (5) mentioned above, and that the notion
of explicitly $(\alpha,n)$-big satisfies criteria (2), (3), $(4^{\prime})$,
and $(5^{\prime})$.
The proof 
that criteria (1) is satisfied will be given in the next section.
We begin with the following  lemma which follows immediately from the
definitions.

\begin{lemma}

\textsc{criterion 2.} If $\cN$ is $(\alpha,n)$-big above $\beta$ and $\cN\unlhd\cM$
then $\cM$ is $(\alpha,n)$-big above $\beta$. (Similarly with
explicitly $(\alpha,n)$-big.)

\nopagebreak

\textsc{criterion 3.} If $\cM$ is $(\alpha,n)$-big above $\beta$ and $\betaprime<\beta$
then $\cM$ is $(\alpha,n)$-big above $\betaprime$. (Similarly with
explicitly $(\alpha,n)$-big.)

\nopagebreak

\textsc{criterion 4.} If $\cM$ is $(\alpha,n)$-big above $\beta$ and 
$(\alphaprime,\nprime)\lexless(\alpha,n)$ then
$\cM$ is $(\alphaprime,\nprime)$-big above $\beta$. (Not true with
explicitly $(\alpha,n)$-big.)
\end{lemma}

Although the notion of explicitly $(\alpha,n)$-big does not satisfy
criterion (4), the following two lemmas show that it does ``locally."

\begin{lemma}
\label{biggoesdownbyone}
If $\cM$ is explicitly $(\alpha,n+1)$-big above $\beta$ then
$\cM$ is explicitly $(\alpha,n)$-big above $\beta$.
\end{lemma}
\begin{proof}
If $\cof(\alpha)>\omega$ or $n\geq1$ this follows immediately from
Definition \ref{nplusone}. If $\cof(\alpha)\leq\omega$
and $n=0$,  then
it follows immediately from  Definitions \ref{cofomeganiszero} and
\ref{cofomeganisone}. 
\end{proof}

\begin{lemma}
\label{biggoesdownbycountable}
Suppose $2\leq\alpha<\omega_1^{\omega_1}$. 
Let $\xi$ be a countable ordinal greater than $0$.
If $\cM$ is explicitly $(\alpha+\xi,0)$-big above $\beta$ then
$\cM$ is explicitly $(\alpha,n)$-big above $\beta$
for all $n$.
\end{lemma}
\begin{proof}
By induction on $\xi$. First suppose $\xi=1$, and
$\cM$ is explicitly $(\alpha+1,0)$-big above $\beta$. By Definition
\ref{cofomeganiszero} $\cM$ is explicitly $(\alpha,n)$-big above $\beta$
for all $n$. The inductive step is similar.
\end{proof}

In addition to the above two lemmas, it is also true that
if $\cM$ is explicitly $(\alpha,n)$-big above $\beta$ and
$(\alphaprime,\nprime)\lexless(\alpha,n)$ and $\max(\code(\alphaprime))$
is sufficiently small, then $\cM$ is explicitly $(\alphaprime,\nprime)$-big
above $\beta$.  This idea is expressed in
Lemmas \ref{biggoesdownbelowbeta}, and \ref{biggoesdownbelows}  below. 

\begin{lemma}
\label{biggoesdownbelowbeta}
Let $2\leq\alpha\leq\omega_1^{\omega_1}$, 
and suppose that $\cM$ is explicitly $(\alpha, 0)$-big above $\beta$. 
Let $\alphaprime<\alpha$ and suppose that $\max(\code(\alphaprime))<\beta$.
Then for all $n$, $\cM$ is explicitly $(\alphaprime,n)$-big above $\beta$.
\end{lemma}

\begin{proof}
First suppose that $\cof(\alpha)>\omega$.
Then there is an active initial segment
$\cN\unlhd\cM$ such that, letting $\kappa$ be
the critical point of the last extender from the $\cN$-sequence,
we have that $\beta<\kappa$, and for all 
$\hat{\alpha}<\alpha$, if $\code(\hat{\alpha})\in J^{\cN}_{\kappa}$
then for all $n$, $\cM$ is explicitly $(\hat{\alpha},n)$-big above $\beta$.
Since $\max(\code(\alphaprime))<\beta<\kappa$ we have that
$\code(\alphaprime)\in J^{\cN}_{\kappa}$, and so for all $n$,
$\cM$ is explicitly $(\alphaprime,n)$-big above $\beta$.

Next suppose that $\cof(\alpha)\leq\omega$.  We argue by induction
on $\alpha$. By Lemma \ref{biggoesdownbycountable} we may assume
that $\alphaprime+1<\alpha$. Then, by definition,
there is an $\hat{\alpha}$ with
$\alphaprime<\hat{\alpha}<\alpha$ such that $\cM$ is explicitly
$(\hat{\alpha},0)$-big above $\beta$. So, by the lemma we are
proving, for $\hat{\alpha}$,
$\cM$ is $(\alphaprime,n)$-big above $\beta$ for all $n$.
\end{proof}

\begin{lemma}
\label{biggoesdownbelows}
Let $2\leq\alpha<\omega_1^{\omega_1}$, and suppose that
that $\cM$ is explicitly $(\alpha, 0)$-big above $\beta$. 
Let $\alphaprime<\alpha$ and suppose that $\max(\code(\alphaprime))<
\max(\code(\alpha))$.
Then for all $n$, $\cM$ is explicitly $(\alphaprime,n)$-big above $\beta$.
\end{lemma}

\begin{proof}
Suppose first that $\cof(\alpha)=\omega_1$.
Then there is an active initial segment
$\cN\unlhd\cM$ such that, letting $\kappa$ be
the critical point of the last extender from the $\cN$-sequence,
we have that $\beta<\kappa$, 
and $\code(\alpha)\in J^{\cN}_{\kappa}$,
and for all $\hat{\alpha}<\alpha$, if 
$\code(\hat{\alpha})\in J^{\cN}_{\kappa}$
then for all $n$, $\cM$ is explicitly $(\hat{\alpha},n)$-big above $\beta$.
Since $\max(\code(\alphaprime))<\max(\code(\alpha))$ we have that
$\code(\alphaprime)\in J^{\cN}_{\kappa}$ and so for all $n$,
$\cM$ is explicitly $(\alphaprime,n)$-big above $\beta$.

Next suppose that $\alpha=\alpha_0+1$ is a successor ordinal,
where $\cof(\alpha_0)=\omega_1$. 
By Lemma \ref{biggoesdownbycountable}, $\cM$ is explicitly
$(\alpha_0,0)$-big above $\beta$. The argument is similar
to the previous paragraph, using $\alpha_0$ in place of $\alpha$, 
and using  the fact that
by part (d) of Lemma \ref{CodeFacts},
$\max(\code(\alpha))\leq\max(\code(\alpha_0))+1$. 

Now suppose that $\alpha=\alpha_0+1$ is a successor ordinal,
where $\cof(\alpha_0)\leq\omega$.
By Lemma  \ref{biggoesdownbycountable}, $\cM$ is explicitly 
$(\alpha_0,1)$-big above $\beta$, and we may assume that
$\alphaprime+1<\alpha_0$.
Let $(\cN,\rho)$ witness that
$\cM$ is explicitly $(\alpha_0,1)$-big above $\beta$. 
By definition, $\max(\code(\alpha_0))<\rho$.
As $\rho$ is a limit cardinal of $\cN$,
$\max(\code(\alpha_0))+1<\rho$.
By part (d) of Lemma \ref{CodeFacts},
$\max(\code(\alpha))\leq\max(\code(\alpha_0))+1$.
So $\max(\code(\alphaprime))<\max(\code(\alpha))<\rho$.
Examining Definition \ref{cofomeganisone}  we see that there is an
ordinal $\delta<\rho$ and an ordinal $\hat{\alpha}<\alpha_0$
such that $\max(\code(\alphaprime))<\delta$ and 
$\alphaprime<\hat{\alpha}$,
and there is an ordinal $\gamma$ with $\delta<\gamma<\rho$ such that
$\cJ^{\cN}_{\gamma}$ is explicitly
$(\hat{\alpha},0)$-big above $\delta$.
By Lemma \ref{biggoesdownbelowbeta}, for $\hat{\alpha}$,
$J^{\cN}_{\gamma}$ is explicitly 
$(\alphaprime,n)$-big above $\delta$ for all $n$.

Finally suppose that $\alpha$ is a limit ordinal of cofinality $\omega$.
Let $\xi=\max(\code(\alphaprime))$. So $\xi<\max(\code(\alpha))$.
By part (c) of Lemma \ref{CodeFacts}, there is an $\alpha_0<\alpha$
such that for all $\hat{\alpha}$ with $\alpha_0<\hat{\alpha}<\alpha$ 
we have that
$\xi<\max(\code(\hat{\alpha}))$. Let $\hat{\alpha}$ be least such
that $\alphaprime<\hat{\alpha}$ and $\alpha_0<\hat{\alpha}$ and
$\cM$ is explicitly $(\hat{\alpha},0)$-big
above $\beta$. Then $\hat{\alpha}$ is
not a limit ordinal of cofinality $\omega$. So by the lemma we are
proving, for $\hat{\alpha}$,
$\cM$ is explicitly $(\alphaprime,n)$-big above $\beta$ for all $n$.
\end{proof}

Suppose that $\cM$ is explicitly $(\alpha,n)$-big above $\beta$. 
For which $\betaprime>\beta$ is it also true that $\cM$ is explicitly
$(\alpha,n)$-big above $\betaprime$? For example, if
 $n=k+1$ and there is an ordinal $\delta>\beta$ such that $\delta$ is
Woodin in $\cM$ and $\cM$ is explicitly $(\alpha,k)$-big above
$\delta$, then $\cM$ is explicitly $(\alpha,n)$-big above $\betaprime$
for all $\betaprime<\delta$. Also,
 if $\cof(\alpha)\leq\omega$ and $(\cM,\rho)$ witnesses that
$\cM$ is explicitly $(\alpha,1)$-big above $\beta$, then in fact
$(\cM,\rho)$ witnesses that $\cM$ is explicitly $(\alpha,1)$-big above
$\betaprime$, for all $\betaprime<\rho$. These considerations will crop
up in  the next section. There is one case which is just
convoluted enough to require a short proof. We describe this in the
following proposition.

\begin{proposition}
\label{BigUpToKappa}
Let $\alpha\leq\omega_1^{\omega_1}$ and suppose $\cof(\alpha)=\omega_1$.
Suppose that $\cM$ witnesses that $\cM$ is explicitly $(\alpha,0)$-big
above some ordinal $\beta$. By definition $\cM$ is active. Let $\kappa$
be the critical point of the last extender from the $\cM$ sequence.
 Then $\cM$ witnesses that $\cM$ is
explicitly $(\alpha,0)$-big above $\betaprime$, for all $\betaprime<\kappa$.
\end{proposition}
\begin{proof}
By definition, $\beta<\kappa$, and $\max(\code(\alpha))<\kappa$.
Fix $\betaprime$ with $\beta<\betaprime<\kappa$. 
Let $\alphaprime<\alpha$ with $\max(\code(\alphaprime))<\kappa$. 
Let $n\in\omega$. We must show that there is a proper initial
segment $\cN\lhd J^{\cM}_{\kappa}$ such that
$\cN$ is explicitly $(\alphaprime,n)$-big
above $\betaprime$. By increasing $\betaprime$ we may assume that
$\max(\code(\alphaprime))<\betaprime$.

Let $\bar{\alpha}=\alphaprime+\betaprime+1$.
Notice that $\betaprime<\max(\code(\bar{\alpha}))\leq
\max(\code(\alphaprime))+\betaprime+1$.
As $\cof(\alpha)=\omega_1$, $\bar{\alpha}<\alpha$. As $\kappa$ is
inaccessible in $\cM$, $\max(\code(\bar{\alpha}))<\kappa$. 
Thus there is a proper initial segment $\cN\lhd J^{\cM}_{\kappa}$ such that
$\cN$ is explicitly $(\bar{\alpha},1)$-big above $\beta$. Fix the
$\unlhd$-least such $\cN$. Since $\bar{\alpha}$ is a successor ordinal,
there is some $\rho$ in $\cN$ such that $(\cN,\rho)$ witnesses that
$\cN$ is explicitly $(\bar{\alpha},1)$-big above $\beta$. 
By definition $\max(\code(\bar{\alpha}))<\rho$. 
 So $\betaprime<\rho$.
Thus $\cN$ is explicitly $(\bar{\alpha},1)$-big above $\betaprime$.

Now $\alphaprime<\bar{\alpha}$.
Also, we are assuming that $\max(\code(\alphaprime))<\betaprime$.
So by Lemma \ref{biggoesdownbelowbeta}, $\cN$ is explicitly 
$(\alphaprime,n)$-big above $\betaprime$.
\end{proof}

We now turn to criterion $(5^{\prime})$. That is, we wish to show that
the statement ``I am explicitly $(\alpha,n)$-big above $\beta$'' is
definable.  To begin with we need the following.

\begin{lemma}
\label{bigcontainscode}
If $\cM$ is explicitly $(\alpha,n)$-big above $\beta$, then 
$\code(\alpha)\in\cM$. \textup{(}Unless $\alpha$ is a limit ordinal
of cofinality $\omega$ and $n=0$.  In this case we can only say
that $\max(\code(\alpha))\leq\ORD\intersect \cM$.\textup{)}
\end{lemma}

\begin{proof}
First suppose $\cof(\alpha)>\omega$. 
If $n=0$ then the lemma is true by definition,
and if $n>0$ then we apply lemma \ref{biggoesdownbyone} above.

Now suppose $\cof(\alpha)\leq\omega$.
If $n=1$ then the lemma is true by definition,
and if $n>1$ then we apply lemma \ref{biggoesdownbyone} above.

Next suppose $\alpha=\gamma+1$ is a successor ordinal and
$n=0$, and $\cM$ is explicitly $(\alpha,0)$-big above $\beta$.
Then
$\cM$ is explicitly $(\gamma,n)$-big above $\beta$ for
all $n$, so $\code(\gamma)\in\cM$. By part (d) of Lemma
\ref{CodeFacts}, $\code(\alpha)\leq\code(\gamma)+1$.
So $\code(\alpha)\in\cM$.

Finally, suppose that $\alpha$ is a limit ordinal of cofinality $\omega$
and $n=0$, and $\cM$ is explicitly $(\alpha,0)$-big above $\beta$.
We argue by induction on $\alpha$. 
For cofinally many $\alphaprime<\alpha$,  there is a
proper initial segment $\cN_{\alphaprime}\lhd\cM$ such that
$\cN_{\alphaprime}$ is explicitly $(\alphaprime,0)$-big above $\beta$.
By induction, for cofinally many $\alphaprime<\alpha$,
$\code(\alphaprime)\in\cM$. By part (c) of Lemma \ref{CodeFacts},
$\max(\code(\alpha))\leq\ORD\intersect\cM$.
\end{proof}

Suppose that $\cM$ is a countable premouse and $\beta\in\ORD^{\cM}$,
and $\alpha<\omega_1^{\omega_1}$ and $s=\code(\alpha)$ and $s\in\cM$. 
We want to show that the statement
``I am explicitly $(|s|,n)$-big above $\beta$'' is definable in $\cM$.
 Notice that as $\cM$ is countable,
  the term ``$|s|$'' cannot mean the same
thing in $\cM$ as it does in $V$. So we must be careful about what
we mean. To be precise, we want to show that
there is a formula $\psi$ in the language of premice such
that for all premice $\cM$, and all codes $s\in\CodeSet\intersect\cM$,
and all $\beta$ in $\ORD^{\cM}$, and all $n\in\omega$, $\cM$ is explicitly
$(|s|,n)$-big above $\beta$ iff $\cM\models\psi[s,n,\beta]$.

It suffices to see that a premouse $\cM$ can tell when
its proper initial segments are big. This is true because the definition
involves a ``local $\Sigma_1$ induction.'' We spell out all of the details
below.

\begin{definition}
Let $\cM=\cJ^{\cM}_{\gamma}$ be a premouse. Then $H^{\cM}$ is the function
with domain equal to the set of all $4$-tuples $(\xi,\beta,s,n)$ 
such that $\xi\leq\gamma$, and $\beta\in\cJ^{\cM}_{\xi}$,
and $s\in\CodeSet$ is a code, with 
$\max(s)\leq\ORD\intersect\cJ^{\cM}_{\xi}$, and $n\in\omega$, and
$$H^{\cM}:\dom(H^{\cM})\map\singleton{0,1}$$
is defined by
$H^{\cM}(\xi,\beta,s,n)=1$ iff 
$\cJ^{\cM}_{\xi}$ is explicitly $(|s|,n)$-big above $\beta$. For 
$\rho\leq\gamma$ we will write $H^{\cM}\restr\rho$ to mean
$$\Union{\xi<\rho}H^{\cJ^{\cM}_{\xi}}.$$
\end{definition}

\begin{lemma}
Let $\cM=\cJ^{\cM}_{\gamma}$ be a premouse.
\begin{itemize}
\item[(a)] For all $\rho\leq\gamma$,
$H^{\cM}\restr\rho$ is $\Sigma_1$ definable over
$\cJ^{\cM}_{\rho}$, uniformly for all $\cM$ and all $\rho$.
\item[(b)] For all $\rho<\gamma$,
$H^{\cM}\restr\rho\in\cJ^{\cM}_{\rho+1}$.
\item[(c)] For all $\rho<\gamma$,
$H^{\cM}\restr(\rho+1)=H^{\cJ^{\cM}_{\rho}}\in\cJ^{\cM}_{\rho+1}$
\end{itemize}
\end{lemma}
\begin{proof}
First observe that $H^{\cN}$ is 
$\Sigma_0$ definable from the parameter $\cN$, uniformly for all
premice $\cN$.
That is, there is a $\Sigma_0$ formula
$\theta$ such that for all $\cN$,
$h=H^{\cN}$ iff $\theta[h,\cN]$ holds.
$\theta(v_1,v_2)$ just says that $v_1$ is a function which
obeys the appropriate recursive definition, with $v_2$ as a parameter. 
The fact that $\theta$ can
be taken to be $\Sigma_0$ uses certain facts about our coding scheme:
(See Lemma \ref{CodeFacts}.) 
\begin{itemize}
\item[1.] The relation $|s|<|t|$ is $\Sigma_0$ definable.
\item[2.] The relation $\cof(|s|)\leq\omega$ is $\Sigma_0$ definable. 
\item[3.] If $s\in\CodeSet$ and $|s|$ is a limit ordinal of cofinality 
$\omega$, and if we set $\CodeSet\restr s= \bigl\{\sprime\in\CodeSet :
|\sprime|<|s|\  \AND \ \max(\sprime)\leq\max(s)\bigr\}$, then
we have that $\dotsetof{|\sprime|}{\sprime\in\CodeSet\restr s}$ is
cofinal in $|s|$.
\end{itemize}

Now let $\formulaphi(\xi,\beta,s,n,b)$ be the $\Sigma_1$ formula
in the language of premice which, when interpreted in the premouse $\cM$, 
asserts:
\begin{quote}
$b\in\singleton{0,1}$ and
there exists a function $h$ such that $\theta(h,\cJ^{\cM}_{\xi})$
and $h(\xi,\beta,s,n)=b$. 
\end{quote} 

We now prove (a), (b), and (c) by simultaneous induction on $\rho$.

Given (c) for $\rho^{\prime}<\rho$, it is easy to see that 
$H^{\cM}\restr\rho(\xi,\beta,s,n)=b$ iff
$\cJ^{\cM}_{\rho}\models\formulaphi[\xi,\beta,s,n,b]$. This proves (a).
 
(b) follows immediately from (a).

We turn to (c).
Let $h=H^{\cM}\restr(\rho+1)-H^{\cM}\restr\rho$. By (b) it suffices
to show that $h\in J^{\cM}_{\rho+1}$. Let 
$$X=\setof{(\beta,s,n)\in\cJ^{\cM}_{\rho}}%
{\cJ^{\cM}_{\rho}\text{ is explicitly }(|s|,n)\text{-big above }\beta}.$$
It suffices to show that $X$ is definable over $\cJ^{\cM}_{\rho}$.
Let$$Y=X\intersect
\setof{(\beta,s,n)}{(\forall\xi<\rho)\,%
\cJ^{\cM}_{\xi}\text{ is \emph{not} explicitly }%
(|s|,n)\text{-big above }\beta}.$$
By (a) it suffices to show that $Y$ is definable over $\cJ^{\cM}_{\rho}$.
But to determine whether or not $(\beta,s,n)\in Y$ only requires looking
at some simple first order properties of $\cJ^{\cM}_{\rho}$ (like whether
or not it is active, whether or not it is admissible, and whether or not
it has $n$ Woodin cardinals) and then looking at whether or not certain
proper initial segments of $\cJ^{\cM}_{\rho}$ are big.  By (a) this
last task can be done definably over $\cJ^{\cM}_{\rho}$. So
$Y$ is definable over $\cJ^{\cM}_{\rho}$.
\end{proof}

\begin{corollary}
\label{BigIsDefinable}
\textsc{criterion} $(5^{\prime})$.
There is a formula $\psi$ in the language of premice such
that for all premice $\cM$, and all codes $s\in\CodeSet\intersect\cM$,
and all $\beta$ in $\ORD^{\cM}$, and all $n\in\omega$, $\cM$ is explicitly
$(|s|,n)$-big above $\beta$ iff $\cM\models\psi[s,n,\beta]$.
\end{corollary}
\begin{proof}
This follows from the proof of (c) in the previous lemma.
\end{proof}

\begin{note}
If $\psi$ is as above, then we will express the fact that
$\cM\models\psi[s,n,\beta]$ by saying that $\cM\models$``I am
explicitly $(\alpha,n)$-big above $\beta$.''
\end{note}

We will use the above corollary heavily.  
Our first application of 
the corollary is to prove a series of lemmas which say that the property of
being explicitly $(\alpha,n)$-big can be reflected downward. We will then
go on to use these reflection lemmas to further study the structure
of $(\alpha,n)$-big mice. Among other things, we will learn that the
mice we are interested in studying occur below the first mouse 
with $\omega$ Woodin cardinals. (See Proposition \ref{OmegaWoodinsIsBig})

\begin{lemma}
\label{FirstBigIsBelowLambda}
Let $\cM$ be a countable, realizable, meek premouse and suppose that
$\lambda$ is a limit cardinal of $\cM$. 
Let $\alpha<\omega_1^{\omega_1}$ with 
$\code(\alpha)\in J^{\cM}_{\lambda}$. Let $\beta<\lambda$,
and let $n\in\omega$.
Suppose that there is a proper initial segment $\cN\lhd\cM$ such that 
$\cN$ is explicitly $(\alpha,n)$-big
above $\beta$. Then,  there is a proper initial segment 
$\cN\lhd J^{\cM}_{\lambda}$ such that 
$\cN$ is explicitly $(\alpha,n)$-big
above $\beta$.
\end{lemma}
\begin{proof}
Let $\kappa<\lambda$ be a cardinal of $\cM$ such that
$\code(\alpha)\in J^{\cM}_{\kappa}$ and $\beta<\kappa$.
Suppose that $\cN\lhd\cM$ and $\kappa\in\cN$ and
$\cN$ is explicitly $(\alpha,n)$-big above $\beta$, and $\cN$ is
$\unlhd$-least with these properties. By Corollary \ref{BigIsDefinable},
 there is a formula
$\theta$ such that $\cN$ is the $\unlhd$-least initial segment of
$\cM$ such that $\kappa\in\cN$ and $\cN\models\theta[\code(\alpha),\beta]$.
Then, by our first corollary to the Condensation Theorem,
Corollary \ref{PhiMinimalProjects} on page \pageref{PhiMinimalProjects},
  $\rho_{\omega}(\cN)=\kappa$.
So $\cN\lhd J^{\cM}_{\lambda}$.
\end{proof}

The next lemma says that the least cardinal which is locally
$(\alpha,n)$-Woodin is not regular.

\begin{lemma}
\label{LeastLocalNotRegular}
Let $\cM$ be a countable, realizable, meek premouse and suppose that
$\delta$ is an inaccessible cardinal of $\cM$. 
Let $\alpha<\omega_1^{\omega_1}$ with 
$\code(\alpha)\in J^{\cM}_{\delta}$, and let $n\in\omega$.
Suppose that $\delta$ is locally $(\alpha,n)$-Woodin in $\cM$.
Then, there are cofinally many $\bar{\delta}<\delta$ such that 
$\bar{\delta}$ is a cardinal of $\cM$ and $\bar{\delta}$ is locally 
$(\alpha,n)$-Woodin in $\cJ^{\cM}_{\delta}$.
\end{lemma}
\begin{proof}
Let $\cN\lhd\cM$ be $\unlhd$-least such that
$\cN$ witnesses that $\delta$ is locally $(\alpha,n)$-Woodin in $\cM$.
So $\delta$ is Woodin in $\cN$ and $\cN$ is the $\unlhd$-least initial
segment of $\cM$ which is explicitly $(\alpha,n)$-big above $\delta$. 

\begin{claim}
$\rho_{\omega}(\cN)=\delta$.
\end{claim}
\begin{subproof}[Proof of Claim]
We want to apply our first corollary to the Condensation Theorem,
Corollary \ref{PhiMinimalProjects}, and this corollary has
a rather carefully-worded hypothesis. We need to see that
there is some formula $\theta$ and some parameter $q\in J^{\cM}_{\delta}$,
so that $\cN$ is the $\unlhd$-least initial segment of $\cM$ such that 
$\delta\in\cN$ and $\cN\models\theta[q]$.
We can take $q=\code(\alpha)$ and we can take
$\theta[\code(\alpha)]$ to say:
``There is some ordinal $\deltaprime$ such that I am 
explicitly $(\alpha,n)$-big above $\deltaprime$ and no proper initial 
segment of me is.'' Clearly $\cN$ satisfies this with $\deltaprime=\delta$.
Suppose towards a contradiction that there is some $\cN^{\prime}$ so that
$\delta\in\cN^{\prime}\lhd\cN$, and
$\cN^{\prime}\models\theta[\code(\alpha)]$. Fix $\deltaprime$ witnessing
this. Then $\deltaprime<\delta$. Let $\kappa$ be a cardinal of $\cM$
such that $\kappa<\delta$ and the maximum of
$\max(\code(\alpha))$ and $\deltaprime$ is below $\kappa$. 
Then
$\cN^{\prime}$ is the $\unlhd$-least initial segment of $\cM$ such
that $\kappa\in\cN^{\prime}$ and $\cN^{\prime}$ satisfies ``I am
explicitly $(\alpha,n)$-big above $\deltaprime$.'' By Corollary
\ref{PhiMinimalProjects}, $\rho_{\omega}(\cN^{\prime})=\kappa$.
But this is impossible since $\delta$ is a cardinal of $\cM$.
This shows that $\cN$ is the $\unlhd$-least initial segment of $\cM$ such 
that $\delta\in\cN$ and $\cN\models\theta[\code(\alpha)]$.
By Corollary \ref{PhiMinimalProjects}, $\rho_{\omega}(\cN)=\delta$.
\end{subproof}

Now fix $\beta<\delta$. We will find $\bar{\delta}$ with $\beta<\bar{\delta}
<\delta$ such that $\bar{\delta}$ is locally $(\alpha,n)$-Woodin 
in $\cJ^{\cM}_{\delta}$.
Working in $\cM$ where $\delta$ is inaccessible, let $X\prec\cN$ be
a fully elementary submodel of $\cN$ such that $|X|<\delta$ and
$X\intersect\delta$ is a cardinal greater than the maximum of
$\max(\code(\alpha))$ and $\beta$. Let $\pi:\bar{\cN}\map X$ be the
inverse of the collapse map. Let $\bar{\delta}=\crit(\pi)$.
Then $\pi(\bar{\delta})=\delta$ and $\bar{\delta}$ is a cardinal of $\cM$.
So $\rho_{\omega}(\bar{\cN})=\bar{\delta}$.
So by the Condensation Theorem, $\bar{\cN}\lhd J^{\cM}_{\delta}$.
(Case (b) of the Condensation Theorem cannot occur because $\bar{\delta}$
is a cardinal of $\cM$, and so $\cJ^{\cM}_{\bar{\delta}}$ is passive.)
We have: $\beta<\bar{\delta}$ and $\bar{\delta}$ is a cardinal of $\cM$ 
and $\bar{\delta}$ is a Woodin cardinal of $\bar{\cN}$ and 
$\bar{\cN}$ is explicitly $(\alpha,n)$-big
above $\bar{\delta}$. So $\bar{\cN}$ witnesses that  $\bar{\delta}$ is 
locally $(\alpha,n)$-Woodin in $\cJ^{\cM}_{\delta}$.
\end{proof}

Suppose $\cof(\alpha)\leq\omega$.
The next lemma say that if $\cM$ has a single fixed cardinal $\delta$,
such that for cofinally many $\alphaprime<\alpha$ and all $n$, 
$\delta$  is locally $(\alphaprime,n)$-Woodin in $\cM$, then there is
a proper initial segment of $\cM$ which is explicitly $(\alpha,1)$-big.
Compare this with the discussion following Lemma \ref{nWoodins} on
page \pageref{nWoodins}.

\begin{lemma}
\label{FixedDeltaIsTooBig}
Let $\cM$ be a countable, realizable, meek premouse and suppose that
$\delta$ is an inaccessible cardinal of $\cM$. 
Let $\alpha<\omega_1^{\omega_1}$ with 
$\code(\alpha)\in J^{\cM}_{\delta}$. Suppose that $\cof(\alpha)\leq\omega$
and that for cofinally many $\alphaprime<\alpha$ and all $n$,
$\delta$ is locally $(\alpha,n)$-Woodin in $\cM$. Then, for all 
$\beta<\delta$, $J^{\cM}_{\delta}$ is explicitly $(\alpha,1)$-big above 
$\beta$.
\end{lemma}
\begin{proof}
Fix $\beta<\delta$.  Let $s=\code(\alpha)$. 
So $s\in J^{\cM}_{\delta}$ by assumption.
Recall that $\CodeSet$ is our set of codes for ordinals less than
$\omega_1^{\omega_1}$.
Set $\CodeSet\restr s=\bigl\{\sprime\in\CodeSet :
|\sprime|<|s| \AND \max(\code(\sprime))\leq\max(\code(s)) \bigr\}$.

\begin{claim}
For all $\sprime\in\CodeSet\restr s$, and all $n\in\omega$, and
all $\betaprime<\delta$, there is an ordinal $\deltaprime$
such that $\betaprime<\deltaprime<\delta$ and $\deltaprime$ is
a cardinal of $\cM$ and $\deltaprime$ is locally $(\alphaprime,n)$-Woodin
in $J^{\cM}_{\delta}$.
\end{claim}
\begin{subproof}
Fix $\sprime$, $n$, and $\betaprime$.
By hypothesis there is an ordinal $\hat{\alpha}$ such that
$|\sprime|<\hat{\alpha}<\alpha$ and there is a proper initial segment
$\cN\lhd\cM$ such that $\delta\in\cN$, and $\delta$
is a Woodin cardinal of $\cN$, and $\cN$ is explicitly 
$(\hat{\alpha},n)$-big above $\delta$. By Lemma \ref{biggoesdownbelowbeta},
$\cN$ is explicitly $(|\sprime|,n)$-big above $\delta$. Then the conclusion
of the claim follows from the previous lemma.
\end{subproof}

Since $\CodeSet\restr s\in J^{\cM}_{\delta}$ and 
$J^{\cM}_{\delta}\models\ZFC$, there is a $\rho<\delta$ such that
$\beta<\rho$ and $\max(s)<\rho$ and the above claim is true with
$\rho$ in place of $\delta$.
By part (e) of Lemma \ref{CodeFacts}, 
$\dotsetof{|\sprime|}{\sprime\in\CodeSet\restr s}$ is cofinal in $\alpha$.
Thus $(J^{\cM}_{\delta},\rho)$ witnesses that $J^{\cM}_{\delta}$ is
explicitly $(\alpha,1)$-big above $\beta$.
\end{proof}

Because the definition of $(\alpha,n)$-big involves a fairly complicated
induction, it may not be clear that there are \emph{any} premice which are
$(\alpha,n)$-big. We will now show that if $\cM$ is a countable, iterable
premouse and $\cM$ has $\omega$ Woodin cardinals, then $\cM$ is
$(\omega_1^{\omega_1},0)$-big.
The contrapositive of this implies, in particular, that if $\cM$ is 
$(\omega_1^{\omega_1},0)$-small, then $\cM$ does not have $\omega$ Woodin
cardinals, and in particular $\cM$ is meek.

\begin{proposition}
\label{OmegaWoodinsIsBig}
Let $\cM$ be a countable, realizable, premouse.  
For $i\in\omega$, suppose
that $\delta_i$ is a Woodin cardinal of $\cM$, and $\delta_i<\delta_{i+1}$.
Then $\cM$ is $(\omega_1^{\omega_1},0)$-big.
\end{proposition}
\begin{proof}
We may assume that no proper initial segment of $\cM$ satisfies our
hypotheses on $\cM$. Thus the $\delta_i$ are cofinal in the ordinals
of $\cM$, and $\cM$ is, in particular, meek. 

\begin{claim}
For all $\alpha<\omega_1^{\omega_1}$ such that $\code(\alpha)\in\cM$,
if $\code(\alpha)\in J^{\cM}_{\delta_i}$, then for all $\beta<\delta_i$,
$J^{\cM}_{\delta_i}$ is explicitly $(\alpha,n)$-big above
$\beta$ for all $n$.
\end{claim}
\begin{subproof}[Proof of Claim]
By Lemmas \ref{biggoesdownbyone} and \ref{biggoesdownbycountable},
it suffices to prove the claim with ``for $n=1$'' in place
of ``for all $n$''. We do this by induction on $\alpha$.
So let $\alpha$ be the least counterexample. 
Fix $i$ such that $\code(\alpha)\in J^{\cM}_{\delta_i}$, and fix
$\beta<\delta_i$.

First suppose that $\cof(\alpha)\leq\omega$. To conclude that
$J^{\cM}_{\delta_i}$ is explicitly $(\alpha,1)$-big above $\beta$,
we will apply the previous lemma with $\delta=\delta_i$. 
So we must see that the hypotheses
of this lemma are satisfied. By Lemma \ref{CodeFacts},
there are cofinally many $\alphaprime<\alpha$ such that
$\code(\alphaprime)\in J^{\cM}_{\delta_i}$. Fix such an $\alphaprime$
and fix $n\geq1$. By induction, $J^{\cM}_{\delta_{i+n}}$ is
explicitly $(\alphaprime,1)$-big above $\delta_{i+n-1}$. So
$J^{\cM}_{\delta_{i+n}}$ is
explicitly $(\alphaprime,n)$-big above $\delta_{i}$. Thus the hypotheses
of the previous lemma are satisfied, and our claim is proved.

Now suppose that $\cof(\alpha)=\omega_1$. Let $\kappa$ be a
measurable cardinal of $\cM$ such that 
$\delta_i <\kappa<\delta_{i+1}$. 
Then for all $\alphaprime<\alpha$,
if $\code(\alphaprime)\in J^{\cM}_{\kappa}$, then by induction,
 $J^{\cM}_{\delta_{i+1}}$ is explicitly $(\alphaprime+1,1)$-big
above $\delta_i$. By Lemma \ref{biggoesdownbyone}, for all $n$,
$J^{\cM}_{\delta_{i+1}}$ is explicitly $(\alphaprime,n)$-big
above $\delta_i$. By Lemma \ref{FirstBigIsBelowLambda}, for all
$\alphaprime<\alpha$, if $\code(\alphaprime)\in J^{\cM}_{\kappa}$,
then for all $n$ there is
a proper initial segment $\cN\lhd J^{\cM}_{\kappa}$ such that
$\cN$ is explicitly $(\alphaprime,n)$-big above $\delta_i$.
Thus $J^{\cM}_{\delta_{i+1}}$ is explicitly $(\alpha,0)$-big
above $\delta_i$. Thus $J^{\cM}_{\delta_{i+1}}$ is explicitly 
$(\alpha,1)$-big above $\beta$. By Lemma \ref{FirstBigIsBelowLambda},
there is a proper initial segment $\cN\lhd J^{\cM}_{\delta_i}$ such
that $\cN$ is explicitly $(\alpha,1)$-big above $\beta$.
\end{subproof}

Now, let $\kappa$ be the least measurable cardinal
of $\cM$. Let $\alpha<\omega_1^{\omega_1}$ be such that
$\code(\alpha)\in J^{\cM}_{\kappa}$.  Then for all $n$,
$J^{\cM}_{\delta_0}$ is explicitly $(\alpha,n)$-big (above 0). By Lemma
\ref{FirstBigIsBelowLambda}, for all $n$ there is a proper initial segment
$\cN\lhd J^{\cM}_{\kappa}$ such that $\cN$ is explicitly
$(\alpha,n)$-big. Thus $\cM$ is explicitly $(\omega_1^{\omega_1},0)$-big.
\end{proof}
The next lemma says that the least mouse which is $(\alpha,n)$-big
is not $(\alpha,n+1)$-big, at least if $\max(\code(\alpha))$ is
sufficiently small. This is a somewhat technical lemma which
will be used several times in the sequel.

\begin{lemma}
\label{ProperSegmentsAreBig}
Let $\cM$ be a countable, realizable, meek premouse, and 
$\beta\in\ORD^{\cM}$.
Let $\alpha<\omega_1^{\omega_1}$, let $n\in\omega$,
 and suppose that $\cM$ is explicitly
$(\alpha,n)$-big above $\beta$.  Suppose that
$(\alphaprime,\nprime)\lexless(\alpha,n)$ and 
$\max(\code(\alphaprime))<\beta$. Then there is a
\emph{proper} initial segment $\cN\lhd\cM$ such that
$\cN$ is explicitly $(\alphaprime,\nprime)$-big above $\beta$.
\end{lemma}
\begin{proof}
If $\cof(\alpha)>\omega$ and $n=0$ this follows from the proof
of Lemma \ref{biggoesdownbelowbeta}. If $\cof(\alpha)\leq\omega$
and $n=0$ the statement follows from the definition of
explicitly $(\alpha,0)$-big, and Lemma \ref{biggoesdownbelowbeta}.

So now let $n\geq 0$ and
suppose that $\cM$ is explicitly $(\alpha,n+1)$-big above
$\beta$, and $\max(\code(\alpha))<\beta$,
 and we will show that there is a proper initial segment $\cN\lhd\cM$
such that $\cN$ is explicitly $(\alpha,n)$-big above $\beta$.
Clearly we may assume that no proper initial segment of $\cM$ is
explicitly $(\alpha,n+1)$-big above $\beta$.

If $\cof(\alpha)\leq\omega$ and $n=0$ this follows from the definitions.
So suppose otherwise. Then there is an ordinal $\delta\in\cM$ with
$\beta<\delta$ such that $\delta$ is a Woodin cardinal of $\cM$, and
$\cM$ is explicitly $(\alpha,n)$-big above
$\delta$. We will show that there is a proper initial segment
$\cN\lhd J^{\cM}_{\delta}$
such that $\cN$ is explicitly $(\alpha,n)$-big above $\beta$.

First suppose that $\cof(\alpha)\leq\omega$. Then by case hypothesis we
have that $n\geq 1$. 
There are $n$ ordinals, $\delta_1,\dots,\delta_n$, with
$\delta=\delta_1<\delta_2<\dots<\delta_n$ such that each $\delta_i$
 is a Woodin cardinal of $\cM$, and $\cM$ is explicitly $(\alpha,1)$-big
above $\delta_n$.  In particular, $\cM$ is explicitly $(\alpha,0)$-big
above $\delta_n$.Thus,
for cofinally many $\alphaprime<\alpha$ and all $\nprime$, $\delta_n$
is locally $(\alphaprime,\nprime)$-Woodin in $\cM$. By Lemma
\ref{FixedDeltaIsTooBig}, 
$J^{\cM}_{\delta_n}$ is explicitly $(\alpha,1)$-big above
$\betaprime$, for all $\betaprime<\delta_n$. Thus $J^{\cM}_{\delta_n}$
is explicitly $(\alpha,n)$-big above $\beta$. By Lemma 
\ref{FirstBigIsBelowLambda},
there is a proper initial segment
$\cN\lhd J^{\cM}_{\delta_1}$
such that $\cN$ is explicitly $(\alpha,n)$-big above $\beta$.

Next suppose that $\cof(\alpha)>\omega$. Then
there are $n+1$ ordinals, $\delta_1,\dots,\delta_{n+1}$, with
$\delta=\delta_1<\delta_2<\dots<\delta_{n+1}$ such that each $\delta_i$
 is a Woodin cardinal of $\cM$, and $\cM$ is explicitly $(\alpha,0)$-big
above $\delta_{n+1}$. By definition, 
$\cM$ is an active premouse, and letting
$\kappa$ be the critical point of the last extender from the $\cM$ sequence
we have that $\delta_{n+1}<\kappa$. 

To simplify the notation, let us assume that $n>0$.
Let $\mu$ be a measurable cardinal of
$\cM$ such that $\delta_{n}<\mu<\delta_{n+1}$.
 Let $\cP\lhd J^{\cM}_{\delta_{n+1}}$ be such that
$\cP$ is active, and $\mu$ is the critical point of the last extender
on the $\cP$-sequence. We will show that $\cP$ is explicitly
$(\alpha,0)$-big above $\delta_n$. Let $\alphaprime<\alpha$ with
$\max(\code(\alphaprime))<\mu$, and let $\nprime\in\omega$.
 We must show that there is an initial
segment $\cQ\lhd J^{\cM}_{\mu}$ such that $\cQ$ is explicitly
$(\alphaprime,\nprime)$-big above $\delta_{n}$. But $\cM$ is explicitly
$(\alpha,0)$-big above $\delta_{n+1}$, so there is a proper initial
segment of  $J^{\cM}_{\kappa}$ which is explicitly
$(\alphaprime,\nprime)$-big above $\delta_{n+1}$, and so in particular
above $\delta_n$. So by Lemma \ref{FirstBigIsBelowLambda},
there is an initial
segment $\cQ\lhd J^{\cM}_{\mu}$ such that $\cQ$ is explicitly
$(\alphaprime,\nprime)$-big above $\delta_{n}$.

Thus $\cP$ is explicitly $(\alpha,0)$-big above $\delta_n$, and so
$\cP$ is explicitly $(\alpha,n)$-big above $\beta$. By Lemma 
\ref{FirstBigIsBelowLambda},
there is a proper initial segment
$\cN\lhd J^{\cM}_{\delta_1}$
such that $\cN$ is explicitly $(\alpha,n)$-big above $\beta$.
\end{proof}

The next lemma says that if $\alpha$ is greatest so that $\cM$ is
$(\alpha,0)$-big, then $\max(\code(\alpha))$ is below the least measurable
cardinal of $\cM$. This is another technical lemma which will be
used several times in the sequel.

\begin{lemma}
\label{BignessIsBelowKappa}
Let $\cM$ be a countable, realizable premouse, and $\beta\in\ORD^{\cM}$.
Let $\alpha<\omega_1^{\omega_1}$ and suppose that $\cM$ is explicitly
$(\alpha,0)$-big above $\beta$. Let $\kappa\in\ORD^{\cM}$ and
suppose that $\beta<\kappa$, and $\kappa$ is a cardinal of $\cM$,
and $\kappa$ is the critical point of an extender from the $\cM$-sequence.
Then either $\max(\code(\alpha))<\kappa$, or $\cM$ is explicitly
$(\alpha^*,0)$-big above $\beta$, for some $\alpha^*>\alpha$ such
that $\cof(\alpha^*)=\omega_1$.
\end{lemma}
\begin{proof}
Suppose $\max(\code(\alpha))\geq\kappa$. Let $\alpha^*\leq\omega_1^{\omega_1}$
be least such that $\alpha^*>\alpha$ and $\max(\code(\alpha^*))<\kappa$.
(There is such an $\alpha^*$ since we have set $\code(\omega_1^{\omega_1})=\emptyset$,
and $\max(\emptyset)=0$.)
Then, by Lemma \ref{CodeFacts},
$\alpha^*$ is a limit ordinal of uncountable cofinality.
We will show that $\cM$ is explicitly $(\alpha^*,0)$-big above
$\beta$.
For all $\alphaprime<\alpha^*$, if
$\max(\code(\alphaprime))<\kappa$, then $\alphaprime<\alpha$,
and $\max(\code(\alphaprime))<\max(\code(\alpha))$. By Lemma
\ref{biggoesdownbelows}, for all $\alphaprime<\alpha^*$, if
$\max(\code(\alphaprime))<\kappa$, then $\cM$ is explicitly
$(\alphaprime,n)$-big above $\beta$ for all $n$. Since we may assume that
$\cM$ does
not have $\omega$ Woodin cardinals, we must have that 
for all $\alphaprime<\alpha^*$, if $\max(\code(\alphaprime))<\kappa$, 
then for all $n$, there is a proper initial segment of $\cM$ 
that is explicitly $(\alphaprime,n)$-big above $\beta$. By Lemma
\ref{FirstBigIsBelowLambda}, for all $\alphaprime<\alpha^*$, if 
$\max(\code(\alphaprime))<\kappa$,  then for all $n$, there is a proper 
initial segment of $\cJ^{\cM}_{\kappa}$  that is explicitly 
$(\alphaprime,n)$-big above  $\beta$. Thus letting $\cN\unlhd\cM$ be such
that $\cN$ is active and $\kappa$ is the critical point of the last extender
from the $\cN$ sequence, we have that $\cN$ witnesses that $\cM$ is
explicitly $(\alpha^*,0)$-big above $\beta$.
\end{proof}

In Corollary \ref{BigIsDefinable} above, we showed that the statement ``I
am explicitly $(\alpha,n)$-big above $\beta$'' is definable.
 We then examined some consequences
of this. We now turn to another consequence: If $\cof(\alpha)\leq\omega$
and $\cM$ is explicitly $(\alpha,1)$-big, then there is a sequence of 
ordinals witnessing this. We call this sequence an $\alpha$-ladder. (See
Figure \ref{fig:bigpremouse} on page \pageref{fig:bigpremouse}.) 
It turns out to be important that this
$\alpha$-ladder is an element of $\cM$. We turn to this now.
The technical details concerning the organization of the $\alpha$-ladder
in the following definition
are motivated by our use of the $\alpha$-ladder in the next section.

\begin{definition}
\label{def:ladder}
Let $\cM$ be a premouse and $\beta\in\ORD^{\cM}$. Suppose that 
$2\leq\alpha<\omega_1^{\omega_1}$ and that $\cof(\alpha)\leq\omega$.
Let $f$ be a function. Then
$f$ is an \emph{$\alpha$-ladder above $\beta$ relative to $\cM$}
iff the following conditions hold.
There are two cases. First suppose that $\alpha=\alphaprime+1$ is
a successor ordinal. (See Figure \ref{fig:bigpremouse} on page 
\pageref{fig:bigpremouse}.) Then:
\begin{itemize}
\item[(i)] $\dom(f)=\omega$, and
\item[(ii)] for all $n$,  $f(n)=(\delta_n, \gamma_n)$, 
with $\beta<\delta_n<\gamma_n<\delta_{n+1}<\ORD^{\cM}$, and
\item[(iii)]
$\delta_n$ is a cardinal of $\cM$, and
$J^{\cM}_{\gamma_n}\models$ ``$\delta_n$ is a Woodin cardinal,''
and $\gamma_n$ is least such that
$J^{\cM}_{\gamma_n}$ is explicitly $(\alphaprime,2n)$-big 
above $\delta_n$, and
\item[(iv)] $\max(\code(\alpha))<\delta_0$.
\end{itemize}
Next suppose that $\alpha$ is a limit ordinal of cofinality $\omega$.
Then:
\begin{itemize}
\item[(i)] $\dom(f)$ is a countable limit ordinal 
$\mu\leq\max(\code(\alpha))$, and 
\item[(ii)] for all $\xi\in\dom(f)$ we have that:
$f(\xi)=(\delta_{\xi}, \gamma_{\xi}, s_{\xi})$ 
with $\beta<\delta_{\xi}<\gamma_{\xi}<\ORD^{\cM}$, and
$\delta_{\xi}>\sup\setof{\gamma_{\eta}}{\eta<\xi}$, and
$s_{\xi}\in\CodeSet$ is a code for an ordinal less than $\omega_1^{\omega}$,
 and, letting $\alpha_{\xi}=|s_{\xi}|$,
we have that $\alpha_{\xi}<\alpha$ and
\item[(iii)]
$\delta_{\xi}$ is a cardinal of $\cM$, and
$J^{\cM}_{\gamma_{\xi}}\models$ ``$\delta_{\xi}$ is a Woodin cardinal,''
and $\gamma_{\xi}$ is least such that
$J^{\cM}_{\gamma_{\xi}}$ is explicitly $(\alpha_{\xi},2)$-big 
above $\delta_{\xi}$, and
\item[(iv)] for all $\xi$, $\max(s_{\xi})\leq\max(\code(\alpha))<\delta_0$, and
\item[(v)] $\setof{\alpha_{\xi}}{\xi\in\dom(f)}$ is cofinal in $\alpha$.
\end{itemize}
\end{definition}

\begin{definition}
Suppose $f$ is an $\alpha$-ladder above $\beta$ relative to $\cM$.
Then, using the notation from the previous definition, we define
$$\sup(f)=\sup\setof{\gamma_{\xi}}{\xi\in\dom(f)}=
\sup\setof{\delta_{\xi}}{\xi\in\dom(f)}.$$
\end{definition}

\begin{lemma}
\label{LadderIsInMouse}
Let $\cM$ be a countable, realizable, meek
premouse, and $\beta\in\ORD^{\cM}$. Suppose that
$2\leq\alpha<\omega_1^{\omega_1}$ and that $\alpha$ is a successor
ordinal, or a limit ordinal of cofinality $\omega$. Suppose that
$\cM$ is explicitly $(\alpha,1)$-big above $\beta$, and that $\cM$ itself
witnesses this.  Let $\rho$ be the
least  ordinal such that $(\cM,\rho)$ witnesses this. 
Then there is a function $f$ such that:
\begin{itemize}
\item[(i)] $f$ is an $\alpha$-ladder above $\beta$ relative to $\cM$, and
\item[(ii)] $f\in\cM$, and
\item[(iii)] $\rho=\sup(f)$.
\end{itemize}
\end{lemma}
\begin{proof}
We will assume that $\alpha$ is a limit 
ordinal  of cofinality $\omega$. If $\alpha$ is a successor ordinal then 
the proof is similar but simpler. 

Recall that $\rho$ is a limit cardinal
of $\cM$, and
notice that there are cofinally many $\gamma<\rho$ such that
$J^{\cM}_{\gamma}\models\ZFC$.
Let $s=\code(\alpha)$. Then $s\in J^{\cM}_{\rho}$ by definition.
Recall that $\CodeSet$ is our set of codes for ordinals below
$\omega_1^{\omega_1}$.
Set $\CodeSet\restr s=\bigl\{\sprime\in\CodeSet :
|\sprime|< |s| \AND \max(\sprime)\leq\max(s) \bigr\}$.
Fix some $\gamma<\rho$ such that $\CodeSet\restr s \in J^{\cM}_{\gamma}$
and $J^{\cM}_{\gamma}\models\ZFC$. We work in $J^{\cM}_{\gamma}$ for a
moment. Let $\mu$ be the cofinality, in $J^{\cM}_{\gamma}$, 
of the wellorder $(\CodeSet\restr s, <)$. It is easy to see that
$\mu\leq\max(s)$. Fix a sequence $\lsequence{s_{\xi}}{\xi<\mu}
\in J^{\cM}_{\gamma}$ which is increasing and cofinal in $\CodeSet\restr s$.

Now we work in $J^{\cM}_{\rho}$. For $\xi<\mu$ let $\delta_{\xi}$ be the
least cardinal of $\cM$ greater than $\max(s)$
which is locally $(|s_{\xi}|,2)$-Woodin in $J^{\cM}_{\rho}$. 
By assumption on $\cM$, there is such a $\delta_{\xi}<\rho$.
Let 
$\gamma_{\xi}$ be the least ordinal greater than $\delta_{\xi}$ such that 
$\cJ^{\cM}_{\gamma_{\xi}}$ witnesses that $\delta_{\xi}$ is locally 
$(|s_{\xi}|,2)$-Woodin.  So $\delta_{\xi}$ is Woodin in 
$\cJ^{\cM}_{\gamma_{\xi}}$ and $\cJ^{\cM}_{\gamma_{\xi}}$ is explicitly 
$(|s_{\xi}|,2)$-big above $\delta_{\xi}$. By assumption on $\cM$, 
$\delta_{\xi}<\gamma_{\xi}<\rho$.
By Corollary \ref{BigIsDefinable},
the function $f=\lsequence{\delta_{\xi},\gamma_{\xi},s_{\xi}}{\xi<\mu}$
is definable over $J^{\cM}_{\rho}$. So $f\in\cM$.  We will show that
$f$ is an $\alpha$-ladder above $\beta$ relative to $\cM$, and that
$\sup(f)=\rho$. We start by showing that the sequence of
$\angles{\delta_{\xi}}$ is strictly increasing.

\begin{claim}[Claim 1]
Let $\eta<\xi<\mu$. Then $\delta_{\eta}<\gamma_{\eta}<\delta_{\xi}$.
\end{claim}
\begin{subproof}[Proof of Claim 1]
Since $\max(s)<\delta_{\xi}$, it follows from Lemma
\ref{ProperSegmentsAreBig} that for all $\sprime\in\CodeSet\restr s$,
if $|\sprime|<|s_{\xi}|$ then
 there is a \emph{proper} initial segment of 
$\cJ^{\cM}_{\gamma_{\xi}}$ which is explicitly $(|\sprime|,2)$-big above 
$\delta_{\xi}$.  Since $\delta_{\xi}$ is inaccessible in
$\cJ^{\cM}_{\gamma_{\xi}}$, Lemma \ref{LeastLocalNotRegular} for 
$\cJ^{\cM}_{\gamma_{\xi}}$ and $\delta_{\xi}$
implies that
$\delta_{\eta}<\gamma_{\eta}<\delta_{\xi}$.
\end{subproof}

It then follows from the leastness of $\rho$ that $\sup(f)=\rho$. 
To complete our proof we must show that the sequence of $\delta_{\xi}$
does not contain any of its limit point.

\begin{claim}[Claim 2]
Let $\xi<\mu$ and suppose that $\xi$ is a limit ordinal. Then
$\delta_{\xi}>\sup\setof{\delta_{\eta}}{\eta<\xi}$.
\end{claim}
\begin{subproof}[Proof of Claim 2]
We have that
$\xi<\mu\leq\max(s)<\delta_{\xi}$, and $\delta_{\xi}$ is
regular in $\cJ^{\cM}_{\gamma_{\xi}}$, and 
$f\restr\xi \in \cJ^{\cM}_{\gamma_{\xi}}$.
\end{subproof}

This completes the proof of the lemma.
\end{proof}

We now turn to yet another consequence of the fact that 
the statement ``I am explicitly $(\alpha,n)$-big'' is definable.
We have mentioned several times that the property of being
explicitly $(\alpha,n)$-big does not satisfy our criterion (4).
That is to say, if $\cM$ is explicitly $(\alpha,n)$-big it does not follow
that $\cM$ is explicitly $(\alphaprime,\nprime)$-big for all
$(\alphaprime,\nprime)\lexless(\alpha,n)$. We proved four lemmas above
which gave approximations to this condition though. (See
Lemmas \ref{biggoesdownbyone} through \ref{biggoesdownbelows}.)
Another approximation is that our criterion $(4^{\prime})$ is
satisfied. That is the content of the next  lemma.

\begin{lemma}
\label{CriterionFourPrime}
\textsc{criterion} $(4^{\prime})$.
Suppose that $\cM$ is a countable, realizable premouse, 
and $(1,0)\lexless(\alpha_0,n_0)\lexleq(\omega_1^{\omega_1},0)$,
and $\beta\in\cM$, and  $\cM$ is explicitly
$(\alpha_0,n_0)$-big above $\beta$. Then for all  
$(\alphaprime,\nprime)\lexless(\alpha_0,n_0)$, 
there is a countable premouse $\cN^{\prime}$ such that
\begin{itemize}
\item[(a)] $\cN^{\prime}$ is a linear iterate of some initial
segment $\cN\unlhd\cM$,
\item[(b)] all critical points used in the iteration  to
$\cN^{\prime}$ are above $\beta$, 
\item[(c)] $\cN^{\prime}$ is realizable, and
\item[(d)]$\cN^{\prime}$ is explicitly $(\alphaprime,\nprime)$-big above 
$\beta$.
\end{itemize}
\end{lemma}
\begin{proof}
Let $(\alpha,n)$ be the lexicographically least pair such that 
$(\alphaprime,\nprime)\lexleq(\alpha,n)$
and $\cM$ is explicitly $(\alpha,n)$-big above $\beta$. 
If $(\alphaprime,\nprime)=(\alpha,n)$ then we may set
$\cN^{\prime}=\cM$ and we are done.
So suppose that $(\alphaprime,\nprime)\lexless(\alpha,n)$.
Then by Lemma \ref{biggoesdownbyone} $n=0$, and so by definition 
$\alpha$ is not a limit ordinal of
cofinality $\omega$, and by Lemma \ref{biggoesdownbycountable}
$\alpha$ is not a successor ordinal. So $\alpha$ is a limit
ordinal of cofinality $\omega_1$.  Let $\cN\unlhd\cM$ witness that
$\cM$ is explicitly $(\alpha,0)$-big above $\beta$. By definition,
$\cN$ is active. Let $\kappa$ be the
critical point of the last extender from the $\cN$-sequence.
By definition, $\beta<\kappa$.
Let $\cN^{\prime}$ be the result of
iterating $\cN$ linearly, via its last extender,
$\alphaprime+1$-times. We will show that $\cN^{\prime}$ is explicitly
$(\alphaprime,\nprime)$-big above $\beta$.
Let $i:\cN\map\cN^{\prime}$ be the iteration map.  So
$\kappa=\crit(i)$. If $\alpha<\omega_1^{\omega_1}$ then let
$s=\code(\alpha)$. 
By definition,  $s\in J^{\cN}_{\kappa}$.
Let $\formulaphi$ be a $\Sigma_1$ formula such that
for any premouse $\cP$, and any $\xi,\bar{\beta},\bar{s},\bar{n}$,
$\cP\models\formulaphi[\xi,\bar{\beta},\bar{s},\bar{n}]$
iff $\cJ^{\cP}_{\xi}$ is explicitly $(|\bar{s}|,\bar{n})$-big above 
$\bar{\beta}$.

First suppose that $\alpha<\omega_1^{\omega_1}$. Then
$$\cN\models(\forall\bar{s}\in J^{\cN}_{\kappa},\forall \bar{n})
[|\bar{s}|<|s|\implies
(\exists\xi<\kappa)\formulaphi(\xi,\beta,\bar{s},\bar{n})].$$
So
$$\cN^{\prime}\models
(\forall\bar{s}\in J^{\cN^{\prime}}_{i(\kappa)},\forall \bar{n})
[|\bar{s}|<|s|\implies
(\exists\xi<i(\kappa))\formulaphi(\xi,\beta,\bar{s},\bar{n})].$$
Let $\sprime=\code(\alphaprime)$. 
Then $\sprime\in J^{\cN^{\prime}}_{i(\kappa)}$.
So $\cN^{\prime}$ is explicitly
$(\alphaprime,\nprime)$-big above $\beta$.

If $\alpha=\omega_1^{\omega_1}$ then
$$\cN\models(\forall\bar{s}\in J^{\cN}_{\kappa},\forall \bar{n})
[(\exists\xi<\kappa)\formulaphi(\xi,\beta,\bar{s},\bar{n})]$$
and the argument is similar.
\end{proof}

Finally, we turn to criterion 5.  That is, we will show that the property of
 being $(\alpha,n)$-big is preserved by non-dropping iterations. 
As in the proof of the previous lemma, it is easy to see that the
property of having a proper initial segment which is explicitly
$(\alpha,n)$-big is preserved by non-dropping iterations. So we
want to see that if $\cM$ itself witnesses that $\cM$ is explicitly
$(\alpha,n)$-big, then the same is true of a non-dropping iterate of $\cM$.
We will focus on the one difficult case, which is the case $\cof(\alpha)\leq\omega$ and $n=1$.
This turns out to be somewhat tricky 
because it involves the preservation of admissibility under $\Sigma_0$
ultrapowers. We do not know whether admissibility is preserved under 
$\Sigma_0$ ultrapowers in general, but the next lemma says that it is preserved given a
few extra assumptions.

The following lemma is due jointly to the author and his PhD thesis advisor,
John Steel.
\begin{lemma}
\label{UltrapowersPreserveAdmissibility}
Let $M=J_{\xi}^A$, where $A$ is some set in $M$. Suppose that $M\models\KP$.
Suppose that there is an ordinal $\rho\in M$ such that
\begin{itemize}
\item[(a)] $\rho$ is the largest cardinal of $M$.
\item[(b)] For all $x\in J^{A}_{\rho}$,  
$\cP(x)\intersect\cM\subset J^{A}_{\rho}$.
\end{itemize}
Let $G$ be a $(\kappa,\lambda)$-extender over $M$, with $\kappa<\rho$.
Suppose that for all $a\in[\lambda]^{<\omega}$, $G_a$ is $\BfSigmaOne(M)$.
Suppose that $\Ult_0(M,G)$ is wellfounded, and let $i=i_G$ so that
$\Ult_0(M,G)=J_{\eta}^{i(A)}$ for some $\eta$. Finally suppose that
$i$ is continuous at $\rho$. Then $\Ult_0(M,G)\models\KP$.
\end{lemma}
\begin{proof}
Write $\Ult=\Ult_0(M,G)$. Recall that $i$ is a cofinal $\Sigma_1$
embedding.
Let $\formulaphi$ be a $\Sigma_0$ formula of the
Language of Set Theory, and let $p\in\Ult$. Suppose that
$(\forall x\in J^{i(A)}_{i(\rho)})(\exists y \in \Ult)$ such that
$\Ult\models\formulaphi[x,y,p]$.  It suffices to show that there is an 
ordinal $\gamma\in\Ult$ such that
$(\forall x\in J^{i(A)}_{i(\rho)})(\exists y \in J^{i(A)}_{\gamma} )$ such that
$\Ult\models\formulaphi[x,y,p]$. (To see that this is sufficient, use the 
fact that, since $i(\rho)$ is the largest cardinal of $\Ult$, for all
domains $D\in\Ult$, there is a function $h\in\Ult$ such that
$h:J^{i(A)}_{i(\rho)}\surjection D$.)

Let $a\in[\lambda]^{<\omega}$ and let $h\in M$ be a function with
$\dom{h}=[\kappa]^{|a|}$ so that $p=[a,h]^M_G$.  First we will show that
there is an ordinal $\gamma\in\Ult$ which works for all $x$ which are
generated by $a$. Let
$$B=\setof{f\in J^{A}_{\rho}}%
{f\text{ is a function with }\dom{f}=[\kappa]^{|a|}}.$$
Then $B\in M$, and, using $\Sigma_0$ Los's Theorem and the fact that
the range of $i$ is cofinal in $\ORD^{\Ult}$ we have
$$(\forall f\in B)(\exists\gamma\in M)\big[\setof{z\in[\kappa]^{|a|}}%
{(\exists y\in J^A_{\gamma})M\models\formulaphi[f(z),y,h(z)]}\in G_a\big].$$

Now, using the fact the $G_a$ is $\BfSigmaOne(M)$ and the fact that $M\models\KP$,
we get that there is a fixed ordinal $\gamma\in M$ such that
$$(\forall f\in B)\big[\setof{z\in[\kappa]^{|a|}}%
{(\exists y\in J^A_{\gamma})M\models\formulaphi[f(z),y,h(z)]}\in G_a\big].$$

Fix such an ordinal $\gamma$.
Now we will show that this same ordinal works for elements of $\Ult$ which
are generated by coordinates other than $a$.

\begin{claim}
Let $b\in[\kappa]^{<\omega}$.
Let $f$ be any function in $J^A_{\rho}$ with $\dom{f}=[\kappa]^{|a\union b|}$.
Then,
$$\setof{z\in[\kappa]^{|a\union b|}}%
{(\exists y\in J^A_{\gamma}) M\models\formulaphi[f(z),y,h(z^a_{a\union b})}\in G_{a\union b}.$$
\end{claim}
\begin{subproof}[Proof of Claim]
Suppose otherwise. Let
$$X=\setof{z\in[\kappa]^{|a \union b|}}%
{(\forall y\in J^A_{\gamma}) M\models\neg\formulaphi[f(z),y,,h(z^a_{a\union b})}.$$
So $X\in G_{a\union b}$. Let $Y=\setof{z^a_{a\union b}}{z\in X}$. Then $Y\in G_a$.
Let $\pi:Y\map X$ be defined by
$\pi(w)=$ the lexicographically least $z\in X$ such that $z^a_{a\union b}=w$.
Now let $\bar{f}:[\kappa]^{|a|}\map J^A_{\rho}$ be defined by
$$\bar{f}(w)=
\begin{cases}
f(\pi(w))& \text{ if } w\in Y;\\
\emptyset& \text{ otherwise.}
\end{cases}
$$
Then $\bar{f}\in J^A_{\rho}$ so $\bar{f}\in B$. Thus
$$\setof{w\in[\kappa]^{|a|}}%
{(\exists y\in J^A_{\gamma})M\models\formulaphi[\bar{f}(w),y,h(w)]}\in G_a.$$
Fix $w\in Y$ such that
$(\exists y\in J^A_{\gamma})M\models\formulaphi[\bar{f}(w),y,h(w)]$.
Let $z=\pi(w)$. Then $z\in X$. But this contradicts the fact that
$$(\exists y\in J^A_{\gamma})M\models\formulaphi[f(z),y,h(z^a_{a\union b})].$$
\end{subproof}

To finish the proof, let $x\in J^{i(A)}_{i(\rho)}$.  We will show that
$(\exists y\in J^{i(A)}_{i(\gamma)})\Ult\models\formulaphi[x,y,p].$

Using our assumption that $i$ is continuous at $\rho$, we have that
$x\in J^{i(A)}_{i(\beta)}$ for some $\beta<\rho$.  Thus, there is some
function $f\in M$ such that $\ran(f)\subset J^{A}_{\beta}$,
and some $b\in[\lambda]^{<\omega}$ so that
$x=[a\union b, f]^M_{G}$. Using our assumption that for all 
$w\in J^{A}_{\rho}$,  $\cP(w)\intersect\cM\subset J^{A}_{\rho}$,
we have that $f\in J^{A}_{\rho}$.
Then, using the Claim and $\Sigma_0$ Los's Theorem,
we are done.
\end{proof}

The next lemma states that admissibility is preserved under direct
limits.
\begin{lemma}
\label{DirectLimitsPreserveAdmissibilty}
Let $\lambda$ be a limit ordinal. Let
$\lsequence{M_{\xi}}{\xi<\lambda}$ be a sequence of admissible sets.
For $\xi<\eta<\lambda$, suppose $i_{\xi,\eta}$ is a $\Sigma_1$ embedding
of $M_{\xi}$ into $M_{\eta}$. Suppose that for $\xi<\eta<\nu<\lambda$,
$i_{\xi,\nu}=i_{\eta,\nu}\circ i_{\xi,\eta}$. Suppose that the direct
limit of this system is wellfounded, and let $M_{\lambda}$ be the
(transitive collapse of) the direct limit.
Then $M_{\lambda}$ is admissible.
\end{lemma}
\begin{proof}
For $\xi<\lambda$, let $i_{\xi,\lambda}$ be the direct limit map.
Let $a,p\in M_{\lambda}$ and let $\formulaphi$ be a $\Sigma_1$ formula.
Suppose that $M\models(\forall x\in a)(\exists y)\formulaphi[x,y,p]$.
Let $\xi<\lambda$ be large enough that both $a$ and $p$ are in
the range of $i_{\xi,\lambda}$. Let $\bar{a}$ be such that
$i_{\xi,\lambda}(\bar{a})=a$ and let $\bar{p}$ be such that
$i_{\xi,\lambda}(\bar{p})=p$.
Then $M_{\xi}\models
(\forall x\in\bar{a})(\exists y)\formulaphi[x,y,\bar{p}]$. Since
$M_{\xi}$ is admissible, there is a $\bar{b}\in M_{\xi}$ such that
$M_{\xi}\models(\forall x\in \bar{a})(\exists y\in \bar{b})\formulaphi[x,y,\bar{p}]$.
Let $b=i_{\xi,\lambda}(\bar{b})$.
Then $M_{\lambda}\models(\forall x\in a)(\exists y\in b)\formulaphi[x,y,p]$.
Thus $M_{\lambda}$ is admissible.
\end{proof}

Recall that our current goal is to prove that if $\cM$ is explicitly
$(\alpha,n)$-big above $\delta$, and $i:\cM\map\cM^{\prime}$ is a
non-dropping iteration map, then $\cM^{\prime}$ is explicitly
$(\alpha,n)$-big above $i(\delta)$.  This is pretty easy to see
in all but one case. The difficult case is the case in which
$\cof(\alpha)\leq\omega$, and $n=1$, and no proper initial segment
of $\cM$ is explicitly $(\alpha,1)$-big above $\delta$. In this case
we have that $\cM$ is admissible, and we must see that $\cM^{\prime}$
is admissible. The previous two lemmas will give this to us, if we add
one more hypothesis. The following lemma
gives the details.

\begin{lemma}
\label{IteratingPreservesBigness}
Let $\cM$ be a countable, realizable,
pre-mouse, and $\delta\in\ORD^{\cM}$.
Let $\alpha<\omega_1^{\omega_1}$ and suppose that $\cof(\alpha)\leq\omega$.
Suppose that $\cM$ is explicitly $(\alpha,1)$-big above $\delta$ and
no proper initial segment of $\cM$ is explicitly 
$(\alpha,1)$-big above $\delta$. Furthermore, suppose that there is some
ordinal $\beta<\delta$ such that $\max(\code(\alpha))<\beta$.
Let $\cT$ be a $0$-maximal iteration tree on $\cM$, and suppose
that $\cT$ is above $\beta$. Then $\forall\xi<\length(\cT)$, if $i_{0,\xi}$
exists, then $\cM^{\cT}_{\xi}$ is explicitly $(\alpha,1)$-big above 
$i_{0,\xi}(\delta)$.
\end{lemma}
\begin{proof}
By induction on $\xi$. For $\xi<\lh{\cT}$, we will write 
$\cM_{\xi}$ instead of $\cM^{\cT}_{\xi}$.  
Suppose that $i_{0,\xi}$ exists.

By hypothesis, there is a an ordinal $\rho\in\cM$ such that
$\delta<\rho$ and $(\cM,\rho)$ witnesses that $\cM$ is explicitly
$(\alpha,1)$-big above $\delta$. Recall that this means that
$\rho$ is a limit cardinal of $\cM$ and that:
\begin{itemize}
\item[(i)] $\cM$ is an admissible set, and
\item[(ii)] $\code(\alpha)\in J^{\cM}_{\rho}$, and
\item[(iii)] for cofinally many $\alphaprime<\alpha$, and all 
$\nprime\in\omega$, and cofinally many $\deltaprime<\rho$ we have that:
\begin{itemize}
\item[(a)] $\deltaprime$ is a cardinal of $\cM$, and
\item[(b)] $\deltaprime$ is locally $(\alphaprime,\nprime)$-Woodin 
            in $\cJ^{\cM}_{\rho}$.
\end{itemize}
\end{itemize}

Let $\kappa_0 = \crit(i_{0,\xi})$.  Let $s=\code(\alpha)$.
So we are assuming that $\max(s)<\beta<\kappa_0$.

We will show that $(\cM_{\xi},i_{0,\xi}(\rho))$ witnesses that
$\cM_{\xi}$ is explicitly $(\alpha,1)$-big above $i_{0,\xi}(\delta)$.
It is easy to see that:
\begin{itemize}
\item[(ii)] $s\in J^{\cM_{\xi}}_{i_{0,\xi}(\rho)}$, and
\item[(iii)] for cofinally many $\alphaprime<\alpha$, and all 
$\nprime\in\omega$, and cofinally many  $\deltaprime<i_{0,\xi}(\rho)$ 
we have that:
\begin{itemize}
\item[(a)] $\deltaprime$ is a cardinal of $\cM_{\xi}$, and
\item[(b)] $\deltaprime$ is locally $(\alphaprime,\nprime)$ Woodin 
            in $J^{\cM_{\xi}}_{i_{0,\xi}(\rho)}$.
\end{itemize}
\end{itemize}

So we have only to see that $\cM_{\xi}$ is an admissible set.

By induction, for all $\nu\in[0,\xi)_{T}$, 
$\cM_{\nu}$ is an admissible set.

If $\xi$ is a limit ordinal then $\cM_{\xi}$ is an admissible set by Lemma
\ref{DirectLimitsPreserveAdmissibilty}.  So suppose that $\xi=\eta+1$. Let
$\nu=T\pred(\xi)$. By induction $\cM_{\nu}$ is admissible.
Let $G=E^{\cT}_{\eta}$, and let $i=i_{\nu,\xi}$, and let
$\kappa=\crit(i)$.
We can conclude from Lemma 
\ref{UltrapowersPreserveAdmissibility}
that $\cM_{\xi}$ is admissible, if we can see that $\cM_{\nu}$ and $G$
and $i$ and $\kappa$
satisfy the hypotheses of that lemma. Let 
$\rhoprime=i_{(0,\nu)}(\rho)$. Then $\rhoprime$ is a limit cardinal
of $\cM_{\nu}$, and since no proper initial segment of $\cM_{\nu}$
is explicitly $(\alpha,1)$-big above
$i_{(0,\nu)}(\delta)$, we have that $\cM_{\nu}$ is the first admissible 
structure over $\cJ^{\cM_{\nu}}_{\rhoprime}$. Thus $\rhoprime$ is the 
largest cardinal of $\cM_{\nu}$ and there are no extenders on the
$\cM_{\nu}$ sequence with index above $\rhoprime$. This implies that
$\kappa<\rhoprime$.
Since premice are strongly acceptable we have
that for all $x\in J^{\cM_{\nu}}_{\rhoprime}$,  
$\cP^{\cM_{\nu}}(x)\in J^{\cM_{\nu}}_{\rhoprime}$.
Since $\cT$ is a $0$-maximal tree, $G$ is close to $\cM_{\nu}$, and
so for all $a\in[\length(G)]^{<\omega}$, 
$G_a$ is $\BfSigmaOne(\cM_{\nu})$. So to complete our proof it suffices
to see that  $i$ is continuous at $\rhoprime$. This follows from the
next claim.

\begin{claim}[Claim 3]
$\cof^{\cM_{\nu}}(\rhoprime)<\kappa$.
\end{claim}
\begin{subproof}[Proof of Claim 3]
As there is no proper initial segment of $\cM$ which is explicitly
$(\alpha,1)$-big above $\delta$, $\rho$ is the least ordinal such
that $(\cM,\rho)$ witnesses that $\cM$ is explicitly $(\alpha,1)$-big
above $\delta$.
By Lemma \ref{LadderIsInMouse}, there is a function $f$ such that:
\begin{itemize}
\item[(i)] $f$ is an $\alpha$-ladder above $\delta$
relative to $\cM$, and
\item[(ii)] $f\in\cM$, and
\item[(iii)] $\rho=\sup(f)$.
\end{itemize}
By definition, $\dom(f)\leq\max(\code(\alpha))$.
Thus $\cof^{\cM}(\rho)\leq\max(s)<\kappa_0$. So
$\cof^{\cM_{\nu}}(\rhoprime)<\kappa_0<\kappa$.
\end{subproof}

This completes the proof of the lemma.
\end{proof}

%


\skipbig

\section{A Correctness Theorem}

\label{section:correctness}

In this section we prove that if a countable, iterable premouse $\cM$ is
is $(\alpha,n)$-big, then $\Aan\subseteq\R\intersect\cM$. This gives
us one half of the proof that $\Aan$ is a mouse set.  We work in the theory ZFC.

We shall use heavily a result due to H. Woodin concerning genericity over
$L[\vec{E}]$ models. This result is often referred  to as ``iterating
to make a real generic.'' We state the result here in the form in which
we shall need it. The following is also Theorem 4.3 from \cite{St2}.
\begin{lemma}[Woodin]
\label{MakeRealGeneric}
Let $\cM$ be a countable, realizable,  premouse.  Suppose
$\kappa<\delta<\ORD^{\cM}$. Suppose
$\cM\models\text{``}\delta$ is a Woodin cardinal.''
Let $\P\subseteq J^{\cM}_{\kappa}$ be a partial
order in $\cM$ and let $G$ be $\cM$-generic over $\P$. Then there
is a partial order $\Q\subseteq J^{\cM}_{\delta}$, with
$\Q\in \cM$, such that for any real $w$, there is a $0$-maximal
iteration tree $\cT$ on $\cM$ of countable length $\theta+1$ such that
\begin{itemize}
\item[(a)] $\cM^{\cT}_{\theta}$ is realizable, and
\item[(b)] $D^{\cT}=\emptyset$ so that $i^{\cT}_{0,\theta}$ is defined,
and
\item[(c)] $\length(E^{\cT}_{\xi})<i^{\cT}_{0,\xi}(\delta)$ for all
$\xi<\theta$, and
\item[(d)] $\crit(E^{\cT}_{\xi}) > \kappa$ for all $\xi<\theta$
(so $G$ is $\cM^{\cT}_{\theta}$-generic over $\P$), and
\item[(e)] $w$ is $\cM^{\cT}_{\theta}[G]$-generic over
$i^{\cT}_{0,\theta}(\Q)$.
\end{itemize}
\end{lemma}

For an idea of how to prove this lemma, see the exercises at the
end of chapters 6 and 7 in \cite{Ma2}.
We shall also need a slightly different form of the above lemma.
\begin{lemma}
\label{MakeRealGenericTwo}
Let $\cM$ be a countable, realizable premouse. Suppose
$\kappa<\delta<\gamma<\ORD^{\cM}$,
and suppose that
$\delta$ is a cardinal of $\cM$, and $J^{\cM}_{\gamma}\models\text{``}
\delta$
is a Woodin cardinal.'' Let $\P\subseteq J^{\cM}_{\kappa}$ be a partial
order in $\cM$ and let $G$ be $\cM$-generic over $\P$. Then there
is a partial order $\Q\subseteq J^{\cM}_{\delta}$, with
$\Q\in J^{\cM}_{\gamma}$, such that for any real $w$, there is a
$0$-maximal
iteration tree $\cT$ on $\cM$ of countable length $\theta+1$ such that
\begin{itemize}
\item[(a)] $\cM^{\cT}_{\theta}$ is realizable, and
\item[(b)] $D^{\cT}=\emptyset$ so that $i^{\cT}_{0,\theta}$ is defined,
and
\item[(c)] $\length(E^{\cT}_{\xi})<i^{\cT}_{0,\xi}(\delta)$ for all
$\xi<\theta$, and
\item[(d)] $\crit(E^{\cT}_{\xi}) > \kappa$ for all $\xi<\theta$
(so $G$ is $\cM^{\cT}_{\theta}$-generic over $\P$), and
\item[(e)] $w$ is $i^{\cT}_{0,\theta}(J^{\cM}_{\gamma})[G]$-generic over
$i^{\cT}_{0,\theta}(\Q)$.
\end{itemize}
\end{lemma}

\begin{remark}
Lemma \ref{MakeRealGenericTwo} does not seem to follow
from Lemma \ref{MakeRealGeneric} applied to the
premouse $J^{\cM}_{\gamma}$.  Since $\delta$ is not necessarily a
regular cardinal of $\cM$, if $w$ is a real,
and if $\cT^{\prime}$ is the iteration tree
on $J^{\cM}_{\gamma}$ given by Lemma \ref{MakeRealGeneric}, and if
$\cT$ is the tree on $\cM$ induced by $\cT^{\prime}$,
and if $\length(\cT)=\theta+1$, then it is possible that
$i^{\cT^{\prime}}_{0,\theta}(\delta)<i^{\cT}_{0,\theta}(\delta)$.
This means that $w$ may not be
$i^{\cT}_{0,\theta}(J^{\cM}_{\gamma})[G]$-generic.

However Lemma \ref{MakeRealGenericTwo} does follow from the proof
of Lemma \ref{MakeRealGeneric}. The idea is that in the situation described
in the above paragraph, we would simply continue to iterate.
We sketch this idea below.
\end{remark}

\begin{proof}[Proof of Lemma  \ref{MakeRealGenericTwo}]
We assume that the reader is familiar with the proof of
Lemma \ref{MakeRealGeneric}.
Let $\Q$ be the partial order given by Lemma \ref{MakeRealGeneric}
applied to the premouse $J^{\cM}_{\gamma}$. Fix any real $w$.
The proof of
Lemma \ref{MakeRealGeneric} proceeds by building an iteration tree
$\cT^{\prime}$ on $J^{\cM}_{\gamma}$ such that for all
$\xi<\length{\cT^{\prime}}$, $E^{\cT^{\prime}}_{\xi}$ is the ``least
witness of the non-genericity of $w$.'' (That is $E^{\cT^{\prime}}_{\xi}$
is the shortest extender $E$ from the $\cM^{\cT^{\prime}}_{\xi}$
sequence such that there is a set $X\subseteq\crit{E}$ such that
$w$ is not consistent with $X$ but $w$ is consistent with $i_{E}(X)$.)
To prove Lemma \ref{MakeRealGenericTwo}, we can build a tree $\cT$ on
$\cM$ such that for all
$\xi<\length{\cT}$, $E^{\cT}_{\xi}$ is the least
witness that $w$ is not generic over
$i^{\cT}_{0,\xi}(J^{\cM}_\gamma)[G]$. (Note that since
$\delta$ is a cardinal of $\cM$, any extender in $J^{\cM}_\delta$
which is total on $J^{\cM}_\gamma$ is also total on $\cM$.) The proof
that this process must terminate goes through in our setting even
though $\delta$ is not Woodin in $\cM$.
\end{proof}

We are now ready to state the main theorem of this section. This theorem
says that our Criterion 1
from the previous section is satisfied. First we give a definition
which will help us to state this criterion.

\begin{definition}
\label{ClosedUnder}
Let $\alpha\geq 2$ be projective-like. Let $n\in\omega$.  
Let $X\subset\R$.  Then we
say that $X$ is \emph{closed under $\Aan$} iff for all $y\in X$,
$\Aan(y)\subseteq X$.
\end{definition}

Using this definition we can state Criterion 1 from the previous
section as follows:
If a countable, iterable premouse $\cM$ is $(\alpha,n)$-big above
$\beta$, and if $\P\in\cM$ is a partial order with
$(\card(\P)^+)^{\cM}\leq\beta$, and if
$G$ is $\cM$-generic over $\P$, then $\R\intersect\cM[G]$
should be closed under $\Aan$. The motivation for this  is
that the  criterion is met in the case $\alpha=1$. This was proved
by Steel and Woodin in \cite{MaSt}. Our proof below borrows many ideas
from the proof in \cite{MaSt}.

In theorem \ref{BigMiceAreClosed} below, we will actually prove something 
which is ostensibly stronger than the statement of Criterion 1 in the
previous paragraph. Recall that on page \pageref{section:ordef} we discussed
Countable Sets of Ordinal Definable Reals, in $\LofR$. In that discussion
we
defined the sets $\OD^{\JalphaR}_n$ for $n\geq 1$,
and we pointed out that for $n\geq 0$,
$\Aan\subseteq\OD^{\JalphaR}_{n+1}$. We also asked the question: Is it
true that $\Aan=\OD^{\JalphaR}_{n+1}$? As we mentioned before, we
are not able to answer this question using techniques from descriptive
set theory, but we are able to shed some light on the question
using inner model theory. The way we will do this is by replacing
$\Aan$ with $\OD^{\JalphaR}_{n+1}$ in Criterion 1, from the previous
paragraph. We start with the
 obvious analogue of Definition \ref{ClosedUnder} above.

\begin{definition}
Let $\alpha\geq 2$ be projective-like. Let $n\in\omega$.  
Let $X\subset\R$.  Then we
say that $X$ is \emph{closed under} 
$\OD^{\JalphaR}_{n+1}$
iff for all $y\in X$, $\OD^{\JalphaR}_{n+1}(y)\subseteq X$.
\end{definition}

 Now we sate the main
theorem of this section.

\begin{theorem}
\label{BigMiceAreClosed}
Assume $\Det(\BfSigmaOne(\JofR{\omega_1^{\omega_1}}))$.
Let $(\alpha,n)$ be such that 
$(1,0)\lexleq(\alpha,n)\lexleq(\omega_1^{\omega_1},0)$.
Let $\cM$ be a countable, realizable premouse. Let
$\P$ be a poset in $\cM$ and let $\kappa\geq\omega$ be such that 
$\P\subseteq J^{\cM}_{\kappa}$.
Suppose $\cM\models\kappa^+$ exists, and let $\beta=(\kappa^+)^{\cM}$.
Suppose that $\cM$ is $(\alpha,n)$-big above $\beta$.
Let $G$ be $\cM$-generic over $\P$.
 Then $\R\intersect\cM[G]$ is closed under $\OD^{\JalphaR}_{n+1}$.
\end{theorem}

\begin{note}
The statement of the  theorem includes a determinacy
hypothesis. We have not assumed the optimal determinacy hypothesis,
but rather an easily statable one.
It is possible that the optimal determinacy hypothesis actually
follows from the other assumptions of the theorem, and therefore could
be omited. We do not know if this is true.
\end{note}

Before proving Theorem \ref{BigMiceAreClosed}, we will give a corollary
to the theorem. The following is the only result from this section which
we will actually use in the rest of the paper.

\begin{corollary}
\label{BigMouseContainsAan}
Assume $\Det(\BfSigmaOne(\JofR{\omega_1^{\omega_1}}))$.
Let $(\alpha,n)$ be such that 
$(1,0)\lexleq(\alpha,n)\lexleq(\omega_1^{\omega_1},0)$.
Let $\cM$ be a countable, realizable premouse.
Suppose that $\cM$ is $(\alpha,n)$-big (above $0$.) Then 
$\OD^{\JalphaR}_{n+1}\subseteq\R\intersect\cM$.
\end{corollary}
\begin{proof}
Let $(\alphaprime,\nprime)$ be the lexicographically least pair such
that $(\alpha,n)\lexleq(\alphaprime,\nprime)$, and $\cM$ is 
\emph{explicitly} $(\alphaprime,\nprime)$-big. Let $\cN\initseg\cM$
be the $\initseg$-least initial segment of $\cM$ which is explicitly
$(\alphaprime,\nprime)$-big. The case $\nprime=0$ and
$\cof(\alphaprime)\leq\omega$ is a bit awkward, so let us consider the
other cases first. In all other cases, $\cN\models$``$\omega_1$ exists'',
and $\cN$ is explicitly $(\alphaprime,\nprime)$-big above $\omega_1^{\cN}$.
Then we may apply Theorem \ref{BigMiceAreClosed} to $\cN$, with
$\P=\emptyset$ and $\beta=\omega_1^{\cN}$, to conclude that
$\R\intersect\cN$ is closed under $\OD^{\JofR{\alphaprime}}_{\nprime+1}$.
In particular, $\OD^{\JalphaR}_{n+1}\subseteq\R\intersect\cM$.

Now suppose that $\cof(\alphaprime)\leq\omega$ and $\nprime=0$. In this
case it may not be true the  $\cN\models$``$\omega_1$ exists.'' 
Nevertheless, we
 will show that $\OD^{\JofR{\alphaprime}}_{1}\subseteq\R\intersect\cN$.
By proposition \ref{ODone} we have that
$\OD^{\JofR{\alphaprime}}_{1}$ = $A_{(\alphaprime,0)}$ =
$$\setof{x\in\R}{(\exists\gamma<\alphaprime)\, x\in\OD^{\JofR{\gamma}}}.$$
Fix $\gamma<\alphaprime$. We will show that 
$\OD^{\JofR{\gamma}}\subseteq\R\intersect\cN$. Fix $n_0\geq1$.
We will show that $\OD^{\JofR{\gamma}}_{n_0+1}\subseteq\R\intersect\cN$.
By definition of explicitly $(\alphaprime,0)$-big, there is an 
$\alpha_0$ such that $\gamma\leq\alpha_0<\alphaprime$, and $\cN$ is
explicitly $(\alpha_0,n_0)$-big. Let $\cP\initseg\cN$ be the
$\initseg$-least initial segment of $\cN$ such that $\cP$ is
explicitly $(\alpha_0,n_0)$-big. Since $n_0\geq 1$, we have as above that
 $\cP\models$``$\omega_1$ exists'',
and $\cP$ is explicitly $(\alpha_0,n_0)$-big above $\omega_1^{\cP}$.
Then we may apply Theorem \ref{BigMiceAreClosed} to $\cP$, with
$\P=\emptyset$ and $\beta=\omega_1^{\cP}$, to conclude that
$\R\intersect\cP$ is closed under $\OD^{\JofR{\alpha_0}}_{n_0+1}$.
In particular $\OD^{\JofR{\gamma}}_{n_0+1}\subseteq\R\intersect\cN$.
\end{proof}

We now turn to the proof of Theorem \ref{BigMiceAreClosed}.
We want to prove the theorem  by induction on
$(\alpha,n)$. It turns out that the induction is easier to do
if we add a few extra assumptions. What we will do is the following.
First we will state a lemma which is almost identical to the
theorem, but with a few extra assumptions. Then we will derive
the theorem from the lemma. Then we will prove the lemma by induction
on $(\alpha,n)$. 

The following lemma is almost identical to Theorem \ref{BigMiceAreClosed},
with a few exceptions. Notice, in particular, that below
we say ``explicitly $(\alpha,n)$-big'' and in the theorem above
we say ``$(\alpha,n)$-big.'' Also, in the lemma below we put a requirement
on $\max(\code(\alpha))$.

\begin{lemma}
\label{TechnicalVersion}
Let $(\alpha,n)$ be such that 
$(1,0)\lexleq(\alpha,n)\lexleq(\omega_1^{\omega_1},0)$.
Assume $\Det(\JofR{\alpha+1})$.
Let $\cM$ be a countable, realizable premouse. Let
$\P$ be a poset in $\cM$ and let $\kappa$ be such that 
$\P\subseteq J^{\cM}_{\kappa}$.
Suppose $\cM\models\kappa^+$ exists, and let $\beta\geq(\kappa^+)^{\cM}$.
Suppose that $\cM$ is explicitly $(\alpha,n)$-big above $\beta$,
and that, if $\alpha<\omega_1^{\omega_1}$,
$\max(\code(\alpha))<\beta$.
Let $G$ be $\cM$-generic over $\P$.
Then $\R\intersect\cM[G]$ is closed under $\OD^{\JalphaR}_{n+1}$.
\end{lemma}

First we will derive the theorem from the lemma.

\begin{proof}[Proof of Theorem \ref{BigMiceAreClosed}]
Fix $(\alpha,n)$, $\cM$, $\P$, $\kappa$, $\beta$, and $G$ as in the
statement of the theorem. Let $(\alphaprime,\nprime)$ be the
lexicographically greatest pair such that
$(\alpha,n)\lexleq(\alphaprime,\nprime)$ and $\cM$ is explicitly
$(\alphaprime,\nprime)$-big above $\beta$. It is easy to see that
there is such a pair. (Let $X$ be the set of pairs $(\alphaprime,\nprime)$
such that that $\cM$ is explicitly $(\alphaprime,\nprime)$-big above 
$\beta$. Notice that $X$ is countable. Suppose towards a contradiction
that $(X,\lexleq)$ has limit order type. Let $(\alpha_0,0)$ be the
supremum of $X$. Then by lemma \ref{biggoesdownbyone},
for cofinally many $\alphaprime<\alpha_0$, and
all $\nprime$, $\cM$ is explicitly $(\alphaprime,\nprime)$-big 
above $\beta$. If $\cM$ had $\omega$ Woodin cardinals above $\beta$,
then $\cM$ would be explicitly $(\omega_1^{\omega_1},0)$-big above
$\beta$, and so $(\omega_1^{\omega_1},0)$ would be the greatest pair
in $X$. So suppose that $\cM$ does not have $\omega$ Woodin cardinals.
Then for cofinally many $\alphaprime<\alpha_0$, and
all $\nprime$, some \emph{proper}
 initial segment of $\cM$ is explicitly 
$(\alphaprime,\nprime)$-big  above $\beta$. But then, by definition,
$\cM$ is explicitly $(\alpha_0,0)$-big above $\beta$.) Now let
$\cN\unlhd\cM$ be the $\unlhd$-least initial segment of $\cM$ such
that $\cN$ is explicitly $(\alphaprime,\nprime)$-big above $\beta$.
Notice that $\R\intersect\cN[G]=\R\intersect\cM[G]$. So it suffices
to show that $\R\intersect\cN[G]$ is closed under $\OD^{\JalphaR}_{n+1}$.
We will in fact show that $\R\intersect\cN[G]$ is closed under 
$\OD^{\JofR{\alphaprime}}_{\nprime+1}$.
If $(\alphaprime,\nprime) = (\omega_1^{\omega_1}, 0)$, then we may
now apply Lemma \ref{TechnicalVersion} to $\cN$, and we are done.
 So suppose that 
$\alphaprime<\omega_1^{\omega_1}$. Then we have to worry about the
hypothesis in Lemma \ref{TechnicalVersion} about $\max(\code(\alphaprime))$.

The case $\cof(\alphaprime)\leq\omega$ and $\nprime=0$ is slightly awkward,
so let us handle the other cases first. In all other cases, there is
at least one total-on-$\cN$ extender from the $\cN$-sequence, with critical
point above $\beta$. Let $\mu$ be the least such critical point.
By Lemma  \ref{BignessIsBelowKappa},
$\max(\code(\alphaprime))<\mu$. Let $\betaprime$ be any ordinal
with $\max(\code(\alphaprime))<\betaprime<\mu$. By Lemma
\ref{BigUpToKappa} on page \pageref{BigUpToKappa},
 and the discussion preceding that lemma, we have that
$\cN$ is explicitly $(\alphaprime,\nprime)$-big above $\betaprime$.
So we may apply Lemma
\ref{TechnicalVersion} to $\cN$ and $\betaprime$ to conclude that
$\R\intersect\cN[G]$ is closed under
$\OD^{\JofR{\alphaprime}}_{\nprime+1}$.

Finally, let us consider the case $\cof(\alphaprime)\leq\omega$ and
$\nprime=0$. By proposition \ref{ODone} we have that
$\OD^{\JofR{\alphaprime}}_{1}$ = $A_{(\alphaprime,0)}$ =
$$\setof{x\in\R}{(\exists\gamma<\alphaprime)\, x\in\OD^{\JofR{\gamma}}}.$$
We have that for all $\gamma<\alphaprime$, and all $m$, 
$\cN$ is $(\gamma,m)$-big above $\beta$.
By the previous cases of the theorem we are proving, for all 
$\gamma<\alphaprime$, and all 
$m\not=0$, $\R\intersect\cN[G]$ is closed under 
$\OD^{\JofR{\gamma}}_{m+1}$. Thus $\R\intersect\cN[G]$ is closed under
$\OD^{\JofR{\alphaprime}}_{1}$.
\end{proof}

Now we turn to the proof of Lemma \ref{TechnicalVersion}.
We will prove the lemma by induction on $(\alpha,n)$.
There will be several different cases for the inductive step. Some
of these cases are direct analogues of cases in the proof by Steel and
Woodin in \cite{MaSt}. Other cases are unique to our setting.
Two of the unique cases are: $\cof(\alpha)\leq\omega$ with $n=1$, and $n=2$.
In these two cases, our primary tool will be the following lemma.
This lemma is the motivation for the definition of explicitly
$(\alpha,1)$-big, in the case $\cof(\alpha)\leq\omega$. The lemma
says the if $\cM$ is explicitly $(\alpha,1)$-big, then $\cM$ can tell
whether or not a given $\Sigma_2(\JalphaR)$ fact is true.

\begin{lemma}
\label{ASublemma}
Suppose $2\leq\alpha<\omega_1^{\omega_1}$, and that 
$\cof(\alpha)\leq\omega$. Assume $\Det(\JalphaR)$.
Let $\cM$ be a countable, realizable premouse. Let
$\P$ be a poset in $\cM$ and let $\kappa$ be such that 
$\P\subseteq J^{\cM}_{\kappa}$.
Suppose $\cM\models\kappa^+$ exists, and let $\beta=(\kappa^+)^{\cM}$.
Suppose that $\cM$ is explicitly $(\alpha,1)$-big above $\beta$
and no proper initial segment of $\cM$ is  explicitly $(\alpha,1)$-big
above $\beta$.
Let $\formulaphi$ be a $\Sigma_2$ formula in the language of $\JalphaR$,
with one free variable. Then there is a $\Pi_1$ formula $\psi$ in the
language of premice, and a parameter $a\in\cM$, such that
whenever $G$ is $\cM$-generic over $\P$ and $x\in\R\intersect\cM[G]$
we have that:
\begin{gather*}
\JalphaR\models\formulaphi[x]\\
\text{iff}\\
\cM[G]\models\psi[x,a,G].
\end{gather*}
\end{lemma}

Lemma \ref{TechnicalVersion} and Lemma \ref{ASublemma} will be proved
by simultaneous induction. For ease of readability  we will separate
these into two proofs. We begin with the proof of Lemma \ref{ASublemma}.

\begin{proof}[Proof of Lemma \ref{ASublemma}]
We assume that Lemma \ref{TechnicalVersion} holds for all
$(\alphaprime,n)$
with $1\leq\alphaprime<\alpha$.

Since no proper initial segment of $\cM$ is explicitly $(\alpha,1)$-big
above $\beta$, there is an ordinal $\rho\in\cM$ such $\beta<\rho$ and
$(\cM,\rho)$ witnesses that $\cM$ is explicitly $(\alpha,1)$-big
above $\beta$. Recall that this implies that $\rho$ is a limit
cardinal of $\cM$, and that $\cM$ is the first admissible structure
over $\cJ^{\cM}_{\rho}$. Also,
by Lemma \ref{LadderIsInMouse}, there is a function $f$ 
such that:
\begin{itemize}
\item[(a)] $f$ is an $\alpha$-ladder above $\beta$ relative to $\cM$, and
\item[(b)] $f\in\cM$, and
\item[(c)] $\rho=\sup(f)$.
\end{itemize}
Recall what it means that $f$ is an $\alpha$-ladder above 
$\beta$ relative to $\cM$.
There are two cases. First suppose that $\alpha=\alphaprime+1$ is
a successor ordinal. (See Figure \ref{fig:bigpremouse} on page 
\pageref{fig:bigpremouse}.) Then:
\begin{itemize}
\item[(i)] $\dom(f)=\omega$, and
\item[(ii)] for all $n$,  $f(n)=(\delta_n, \gamma_n)$, 
with $\beta<\delta_n<\gamma_n<\delta_{n+1}<\ORD^{\cM}$, and
\item[(iii)]
$\delta_n$ is a cardinal of $\cM$, and 
$J^{\cM}_{\gamma_n}\models$ ``$\delta_n$ is a Woodin cardinal,''
and $\gamma_n$ is least
such that $J^{\cM}_{\gamma_n}$ is explicitly $(\alphaprime,2n)$-big 
above $\delta_n$, and
\item[(iv)] $\max(\code(\alpha))<\delta_0$.
\end{itemize}
Next suppose that $\alpha$ is a limit ordinal of cofinality $\omega$.
Then:
\begin{itemize}
\item[(i)] $\dom(f)$ is a countable limit ordinal 
$\mu\leq\max(\code(\alpha))$, and 
\item[(ii)] for all $\xi\in\dom(f)$ we have that:
$f(\xi)=(\delta_{\xi}, \gamma_{\xi}, s_{\xi})$ 
with $\beta<\delta_{\xi}<\gamma_{\xi}<\ORD^{\cM}$, and
$\delta_{\xi}>\sup\setof{\gamma_{\eta}}{\eta<\xi}$, and
$s_{\xi}\in\CodeSet$ is a code for an ordinal less than $\omega_1^{\omega}$,
 and, letting $\alpha_{\xi}=|s_{\xi}|$,
we have that $\alpha_{\xi}<\alpha$ and
\item[(iii)]
$\delta_{\xi}$ is a cardinal of $\cM$, and
$J^{\cM}_{\gamma_{\xi}}\models$ ``$\delta_{\xi}$ is a Woodin cardinal,''
and $\gamma_{\xi}$ is least such that $J^{\cM}_{\gamma_{\xi}}$ is 
explicitly $(\alpha_{\xi},2)$-big 
above $\delta_{\xi}$, and
\item[(iv)] for all $\xi$, $\max(s_{\xi})\leq\max(\code(\alpha))<\delta_0$, and
\item[(v)] $\setof{\alpha_{\xi}}{\xi\in\dom(f)}$ is cofinal in $\alpha$.
\end{itemize}
Also recall that $\sup(f)=\sup\setof{\gamma_{\xi}}{\xi\in\dom(f)}=
\sup\setof{\delta_{\xi}}{\xi\in\dom(f)}$.

For notational simplicity we make the following definitions. Let
$\mu =\dom(f)$. So $\mu=\omega$ if $\alpha$ is a successor ordinal,
and if $\alpha$ is a limit ordinal of cofinality $\omega$, then
$\mu\leq\max(\code(\alpha))$ is some countable limit ordinal.
If $\alpha$ is a limit ordinal, then for $\xi<\mu$
the codes $s_{\xi}$ were defined
above, and we set $\alpha_{\xi}=|s_{\xi}|$. If $\alpha$ is a successor
ordinal, then for $\xi<\mu$ let $\alpha_{\xi} = \alpha-1$ and let
$s_{\xi}=\code(\alpha-1)$. Also, for $\xi<\mu$, let
$$
n_{\xi} =
\begin{cases}
2\xi & \text{if $\alpha$ is a successor ordinal},\\
2& \text{if $\alpha$ is a limit ordinal}.
\end{cases}
$$
Notice then that for all $\xi<\mu$,
 $J^{\cM}_{\gamma_{\xi}}\models\delta_{\xi}^+$
exists, and $J^{\cM}_{\gamma_{\xi}}$ is explicitly
$(\alpha_{\xi},n_{\xi})$-big
above $(\delta_{\xi}^+)^{J^{\cM}_{\gamma_{\xi}}}$.  This is important
because we want to apply Lemma \ref{TechnicalVersion} for
$(\alpha_{\xi},n_{\xi})$ to the premouse $J^{\cM}_{\gamma_{\xi}}$.

Now, fix a $\Sigma_2$ formula $\formulaphi$, as in the statement of the 
Lemma.
Recall that there is a $\Sigma_1(\JalphaR)$ partial function from $\R$ onto
$\JalphaR$. As a consequence of this,
$$\Sigma_2(\JalphaR)=\exists^{\R}[\Sigma_1(\JalphaR)
\AND \Pi_1(\JalphaR)].$$
So there are two $\Sigma_1$ formulas
$\formulaphi_1$ and $\formulaphi_2$ so that for all $x\in\R$
\begin{gather*}
\JalphaR\models\formulaphi[x]\\
\text{iff}\\
(\exists y\in\R)\; \JalphaR\models(\formulaphi_1\AND\neg\formulaphi_2)[x,y].
\end{gather*}
Fix such formulas $\formulaphi_1$ and $\formulaphi_2$.

We now fix two sequences of formulas
$\sequence{\theta^1_{\xi}}{\xi<\mu}$ and
$\sequence{\theta^2_{\xi}}{\xi<\mu}$
such that $(\forall x,y\in\R)$
\begin{gather*}
\JalphaR\models\formulaphi_1[x,y]\\
\text{iff}\\
(\exists\xi<\mu)\quad\JofR{\alpha_{\xi}}\models\theta^1_{\xi}[x,y]
\end{gather*}
and
\begin{gather*}
\JalphaR\models\formulaphi_2[x,y]\\
\text{iff}\\
(\exists\xi<\mu)\quad\JofR{\alpha_{\xi}}\models\theta^2_{\xi}[x,y].
\end{gather*}
If $\alpha$ is a limit ordinal, then we simply set
$\theta^1_{\xi}=\formulaphi_1$  and $\theta^2_{\xi}=\formulaphi_2$
for all $\xi<\alpha$. If $\alpha$ is
a successor ordinal, then
by Lemma \ref{SigmaOneIsUnion}, we may take
$\sequence{\theta^1_{\xi}}{\xi\in\omega}$ and
$\sequence{\theta^2_{\xi}}{\xi\in\omega}$ to be recursive.
Notice that in either case we have that $\sequence{\theta^1_{\xi}}{\xi<\mu}$
and $\sequence{\theta^2_{\xi}}{\xi<\mu}$ are in $J^{\cM}_{\rho}$.
Also, we may assume that for $\xi<\mu$,
$\theta^1_{\xi}$ and $\theta^2_{\xi}$ are both (logically equivalent to)
$\Sigma_{n_{\xi}}$ formulas.

Suppose for a moment that $G$ is $\cM$-generic over $\P$, and that
for each $\xi\in\dom(f)$,
$H_{\xi}$ is $\cJ^{\cM}_{\gamma_{\xi}}[G]$-generic over 
$\Col(\omega,\delta_{\xi})$, where
$\Col(\nu_1,\nu_2)$ is the partial order for collapsing
$\nu_2$ to $\nu_1$.
In this situation we will want to show later 
that $\cJ^{\cM}_{\gamma_{\xi}}[G][H_{\xi}]$ has enough 
``correctness'' that it can tell 
when $\theta^1_{\xi}$ and $\theta^2_{\xi}$ are true in the model
$\JofR{\alpha_{\xi}}$. To do this, we must consider the following technical point.
For $\xi<\mu$, let $\bar{\R}_\xi=\R\intersect 
J^{\cM}_{\gamma_{\xi}}[G][H_{\xi}]$. Notice that
$\omega_1^{L(\bar{\R}_{\xi})}=(\delta_{\xi}^+)^{J^{\cM}_{\gamma_{\xi}}}$.
Recall that $s_{\xi}=\code(\alpha_{\xi})$ and that $\max(s_{\xi})
<\delta_0\leq\delta_{\xi}$. Thus $\max(s_{\xi})<
\omega_1^{L(\bar{\R}_{\xi})}$. Thus it makes sense to consider the term
$|s_{\xi}|^{L(\bar{\R}_{\xi})}$.  This will be important below.
(cf. the discussion preceding 
Lemma \ref{MaxCodeIsGoodParameter} on page \pageref{MaxCodeIsGoodParameter}.)

For any $\cM$, $\P$ as in the statement of the lemma, and any ladder $f$
as above, and any $G$ which is $\cM$ generic over $\P$, let
let $S(\cM,f,\P,G)$ be the
collection of all $(x,y)\in\R^2$ such that:
\begin{quote}
There exists (in $V$)
some sequence $\lsequence{H_{\xi}}{\xi\in\dom(f)}\;$ s.t.
$(\forall\xi\in\dom(f))\Bigl[H_{\xi}$ is 
$J^{\cM}_{\gamma_{\xi}}[G]$-generic
over $\Col(\omega,\delta_{\xi})$ and  
$(x,y)\in J^{\cM}_{\gamma_{\xi}}[G][H_{\xi}]\;\Bigr]$ and,
letting $\bar{\R}_{\xi}=\R\intersect J^{\cM}_{\gamma_{\xi}}[G][H_{\xi}]$,
and $\bar{\alpha}_{\xi}=|s_{\xi}|^{L(\bar{\R}_{\xi})}$ we have that
$(\exists\xi\in\dom(f))
\bigl[J_{\bar{\alpha}_{\xi}}(\bar{\R}_{\xi})\models
\theta^1_{\xi}[x,y]\bigr]$ and
$(\forall\xi\in\dom(f))
\bigl[J_{\bar{\alpha}_{\xi}}(\bar{\R}_{\xi})\models
\neg\theta^2_{\xi}[x,y]\bigr]$.
\end{quote}
Then $S(\cM,f,\P,G)$ is $\Sigma^1_1$ uniformly in any real coding
$(\cJ^{\cM}_{\rho},f,\P,G)$, where $\rho=\sup(f)$.
That is, there is a recursive tree $T$ on $\omega\times\omega\times\omega$
such that for any $\cM,\P$ as in the statement of the lemma,
and any ladder $f$ as above, and any $G$ which is $\cM$ generic over $\P$, 
and any real $z$ coding
$(\cJ^{\cM}_{\rho},f,\P,G)$, where $\rho=\sup(f)$,
we have that
$$S(\cM,f,\P,G)=[T_z]=\setof{(x,y)}{(\forall i\in\omega)\;
(x\restr i, y\restr i, z\restr i)\in T}.$$
Fix such a recursive tree $T$.

Now we return to our fixed $\cM$, $\P$, $f$, and $\rho=\sup(f)$. 
(We have not yet fixed a $G$.) 
We will use the following notation: If $\tau$ is a term in a forcing
language, and $K$ is a generic filter, then $\tau_{K}$ denotes the
$K$-interpretation of $\tau$.
Now, let $\dot{z}$ be a $\P$-name in $\cM$
such that whenever $G$ is $\cM$-generic over
$\P$, $\dot{z}_{G}$ is a $\Col(w,\rho)$-name in $\cM[G]$, and
whenever $H$ is $\cM[G]$-generic over $\Col(\omega,\rho)$, $(\dot{z}_G)_H$
is a real coding $(\cJ^{\cM}_{\rho},f,\P,G)$.
Let $a=\angles{\rho,\dot{z}}$.  We claim that this $a$, in conjunction
with an appropriate $\psi$, satisfies the conclusion of the lemma.
Now we define the appropriate formula $\psi$.

For any reals $x,z$ let $T_{x,z}$ be the tree on $\omega$ given
by:
$$s\in T_{x,z} \Iff
(x\restr\length{s},s,z\restr\length{s})\in T.$$

\begin{claim}
There is a $\Pi_1$ formula $\psi$ such that 
whenever
$G$ is $\cM$ generic over $\P$,
and $x\in\R\intersect\cM[G]$, we have that:
\begin{gather*}
\cM[G]\models\psi[x,a,G]\\
\text{iff}\\
\text{whenever $H$ is $\cM[G]$-generic over $\Col(\omega,\rho)$, }
T_{x,(\dot{z}_G)_H}\text{ is ill-founded.}
\end{gather*}
\end{claim}
\begin{subproof}[Proof of Claim]
We can take $\psi$ to say:
$$\emptyset\forcein{\Col(\omega,\rho)}{\cM[G]}\neg[\exists
\text{ a rank function for }T_{\check{x},\dot{z}_G}].$$
More precisely, we can take $\psi$ to say: 
``For all forcing conditions $q$, and
all names $\sigma_1$ and $\sigma_2$, if $\sigma_2 = \dot{z}_G$, then
$\neg[q\Vdash \sigma_1$ is a rank function for
$T_{\check{x},\sigma_2}]$.''  We leave it to the reader to see that
this can be expressed over $\cM[G]$ as a $\Pi_1$ statement
about $x$, $a$, and $G$. 

Now, if $T_{x,(\dot{z}_G)_H}$ is ill-founded
whenever $H$ is $\cM[G]$-generic over $\Col(\omega,\rho)$, then
clearly $\cM[G]\models\psi[x,a,G]$. The converse is true because
for all such $H$, $\cM[G][H]$ is a model of $\KP$.
\end{subproof}

Fix $\psi$ as in the proof of the claim.  We claim that this $\psi$
satisfies the conclusion of the lemma. 
Let $G$ be $\cM$-generic over $\P$,
and let $x\in\R\intersect\cM[G]$. We will show that
\begin{gather*}
\JalphaR\models\formulaphi[x]\\
\text{iff}\\
\cM[G]\models\psi[x,a,G].
\end{gather*}

Suppose first that $\cM[G]\models\psi[x,a,G]$. Let $H$ be $\cM[G]$-generic
over $\Col(\omega,\rho)$. Let $z=(\dot{z}_G)_H$. Then
$\exists y\in\R$ such that $(x,y,z)\in[T]$. Fix such a $y$.
Then $(x,y)\in S(\cM,f,\P,G)$.
So let $\lsequence{H_{\xi}}{\xi\in\dom{f}}$
be such that:
\begin{quote}
$(\forall\xi\in\dom(f))\Bigl[H_{\xi}$ is 
$J^{\cM}_{\gamma_{\xi}}[G]$-generic
over $\Col(\omega,\delta_{\xi})$ and  
$(x,y)\in J^{\cM}_{\gamma_{\xi}}[G][H_{\xi}]\;\Bigr]$ and,
letting $\bar{\R}_{\xi}=\R\intersect J^{\cM}_{\gamma_{\xi}}[G][H_{\xi}]$,
and $\bar{\alpha}_{\xi}=|s_{\xi}|^{L(\bar{\R}_{\xi})}$ we have that
$(\exists\xi\in\dom(f))
\bigl[J_{\bar{\alpha}_{\xi}}(\bar{\R}_{\xi})\models
\theta^1_{\xi}[x,y]\bigr]$ and
$(\forall\xi\in\dom(f))
\bigl[J_{\bar{\alpha}_{\xi}}(\bar{\R}_{\xi})\models
\neg\theta^2_{\xi}[x,y]\bigr]$.
\end{quote}
For a moment, fix $\xi\in\dom{f}$.
Let $\bar{\R}_{\xi}=\R\intersect J^{\cM}_{\gamma_{\xi}}[G][H_{\xi}]$.
Let $\bar{\alpha}_{\xi}=|s_{\xi}|^{L(\bar{\R}_{\xi})}$.
We will show that $J_{\bar{\alpha}_{\xi}}(\bar{\R}_{\xi})$ is
$\Sigma_{n_{\xi}}$ elementarily embeddable into
$J_{\alpha_{\xi}}(\R)$.
There is a $K$ which is $J^{\cM}_{\gamma_{\xi}}$-generic over
$\Col(\omega,\delta_{\xi})$ such that
$\bar{\R}_{\xi}=\R\intersect J^{\cM}_{\gamma_{\xi}}[K]$.
As $J^{\cM}_{\gamma_{\xi}}$ is $(\alpha_{\xi},n_{\xi})$-big
above $(\delta_{\xi}^+)^{J^{\cM}_{\gamma_{\xi}}}$,
Lemma \ref{TechnicalVersion} for
$(\alpha_{\xi},n_{\xi})$ implies that $\bar{\R}_{\xi}$ is closed
under $\OD^{\JofR{\alpha_{\xi}}}_{n_{\xi}+1}$. In particular
$\bar{\R}_{\xi}$ is closed  downward under
$\Delta_{(\alpha_{\xi},n_{\xi})}$. 
As $J^{\cM}_{\gamma_{\xi}}$ has a measurable cardinal above
$\delta_{\xi}$, $\bar{\R}_{\xi}=\R\intersect L(\bar{\R}_{\xi})$.
Since $n_{\xi}$ is even, we can apply  Lemma \ref{MaxCodeIsGoodParameter}
on page \pageref{MaxCodeIsGoodParameter}
to conclude that $J_{\bar{\alpha}_{\xi}}(\bar{\R}_{\xi})$ is in fact
$\Sigma_{n_{\xi}+1}$ elementarily embeddable into
$J_{\alpha_{\xi}}(\R)$. (This is why we organized the definition of
an $\alpha$-ladder so that $n_{\xi}$ is even.) 

Now notice that we have:
\begin{gather*}
(\exists\xi\in\dom{f})
\bigl[J_{\bar{\alpha}_{\xi}}(\bar{\R}_{\xi})\models\theta^1_{\xi}[x,y]\bigr],\\
\text{and}\\
(\forall\xi\in\dom{f})
\bigl[J_{\bar{\alpha}_{\xi}}(\bar{\R}_{\xi})\models\neg\theta^2_{\xi}[x,y]\bigr].
\end{gather*}
Thus we have:
\begin{gather*}
(\exists\xi\in\dom{f})
\bigl[\JofR{\alpha_{\xi}}\models\theta^1_{\xi}[x,y]\bigr],\\
\text{and}\\
(\forall\xi\in\dom{f})
\bigl[\JofR{\alpha_{\xi}}\models\neg\theta^2_{\xi}[x,y]\bigr].
\end{gather*}
And so we have:
$$ \JalphaR\models(\formulaphi_1 \AND \neg \formulaphi_2)[x,y].$$
And so we have:
$$ \JalphaR\models\formulaphi[x].$$

Conversely, assume that $\JalphaR\models\formulaphi[x].$ We  will show that
$\cM[G]\models\psi[x,a,G]$.

Since $\JalphaR\models\formulaphi[x]$, $\exists y\in\R$
such that $\JalphaR\models(\formulaphi_1 \AND \neg \formulaphi_2)[x,y].$ Fix
such a $y$.
We are going to iterate $\cM$, yielding an iteration map
$i:\cM\map\bar{\cM}$, with $\crit(i)>\beta$,
so that for all $\xi\in\dom(f)$, $\angles{x,y}$
is $i(J^{\cM}_{\gamma_{\xi}})[G]$-generic over a partial order contained
in $i(J^{\cM}_{\delta_{\xi}})$.

We apply Lemma \ref{MakeRealGenericTwo}
iteratively, for each $\xi\in\dom{f}$.  This gives us, for each
$\xi\in\dom{f}$: a premouse $\cM_{\xi}$, and
a partial order $\Q_{\xi}\in\cM_{\xi}$,
and an iteration tree $\cT_{\xi}$ on $\cM_{\xi}$
of countable length $\theta_{\xi}+1$; and a system of commuting maps
$i_{\xi,\eta}:\cM_{\xi}\map\cM_{\eta}$ for $\xi<\eta\in\dom{f}$;
such that:
\begin{itemize}
\item[(a)] $\cM_0=\cM$, and
\item[(b)] for each $\xi\in\dom(f)$,
$\cM_{\xi+1}$ is the last model of $\cT_{\xi}$, and
\item[(c)] for each $\xi\in\dom(f)$,
$D^{\cT_{\xi}}=\emptyset$
and $i_{\xi,\xi+1}=i^{\cT_{\xi}}_{0,\theta_{\xi}}$, and
for $\xiprime<\xi$, $i_{\xiprime,\xi+1}=
i_{\xiprime,\xi}\circ i_{\xi,\xi+1}$, and
\item[(d)] for limit $\eta\in\dom(f)$, $\cM_{\eta}$ is the direct limit
of the $\cM_{\xi}$ under the $i_{\xi,\xiprime}$ for
$\xi<\xiprime<\eta$, and $i_{\xi,\eta}$ is the direct limit map, and
\item[(e)] for each $\xi\in\dom(f)$,
$\cM_{\xi}$ is realizable, and
\item[(f)] for each $\xi\in\dom(f)$,
$\length(E^{\cT_{\xi}}_{\nu})<
i^{\cT_{\xi}}_{0,\nu}\circ i_{0,\xi}(\delta_{\xi})$ for all
$\nu<\theta_{\xi}$, and
\item[(g)] for each $\xi\in\dom(f)$,
$\crit(E^{\cT_{\xi}}_{\nu}) >
\sup\lsetof{i_{0,\xi}(\gamma_{\xiprime})}%
{\xiprime<\xi}$ for all $\nu<\theta_{\xi}$
(in particular $\crit(i_{0,\xi+1})>\beta$ so that
$G$ is $\cM_{\xi+1}$-generic over $\P$), and
\item[(h)] for each $\xi\in\dom{f}$,
$\Q_{\xi}\subseteq i_{0,\xi}(J^{\cM}_{\delta_{\xi}})$
and
$\Q_{\xi}\in i_{0,\xi}(J^{\cM}_{\gamma_{\xi}})$
and
$\angles{x,y}$ is
$i_{0,\xi+1}(J^{\cM}_{\gamma_{\xi}})[G]$-generic over
$i_{\xi,\xi+1}(\Q_{\xi})$.
\end{itemize}

The partial order $\Q_{\xi}$ and the
tree $\cT_{\xi}$ comes from applying Lemma \ref{MakeRealGenericTwo}
with
\begin{itemize}
\item[(i)] $\cM=\cM_{\xi}$
\item[(ii)] $\kappa=\sup\lsetof{i_{0,\xi}(\gamma_{\xiprime})}%
{\xiprime<\xi}$, (In the case $\xi=0$ take $\kappa=\max(\code(\alpha)$.)
\item[(iii)] $\delta=i_{0,\xi}(\delta_{\xi})$,
\item[(iv)] $\gamma=i_{0,\xi}(\gamma_{\xi})$,
\item[(v)] $w=\angles{x,y}$,
\item[(vi)] $\P=\P$, and $G=G$.
\end{itemize}
The reason that we can have (ii) above is that
by the definition of an $\alpha$-ladder
$\delta_{\xi}>\sup\setof{\gamma_{\xiprime}}{\xiprime<\xi}$,
and so
$i_{0,\xi}(\delta_{\xi})>\sup\lsetof{i_{0,\xi}(\gamma_{\xiprime})}%
{\xiprime<\xi}$.

The above remarks clearly define by recursion a sequence
$$\bigsequence{\cM_{\xi},\cT_{\xi},\Q_{\xi}}{\xi\in\dom{f}}$$
and a system of maps
$$\lsetof{i_{\xi,\eta}}{\xi<\eta<\dom(f)}$$
such that properties (a) through (h) above hold.
The reason that we are able to maintain that (e) holds at limit
$\xi$ is that for all but finitely many $\eta_1 < \eta_2<\xi$, the 
realization
maps for $\cM_{\eta_1}$ and $\cM_{\eta_2}$ commute with each other.
(cf. the proof in \cite{St2} that if $\cM$ is a countable realizable mouse,
then player II wins the weak iteration game on $\cM$ of length 
$\omega_1$.)

Let $\bar{\cM}$ be the direct limit of the $\cM_{\xi}$ under the
$i_{\xi,\eta}$, and let $i_{\xi}:\cM_{\xi}\map\bar{\cM}$
be the direct limit map. Let $i=i_0$.
As in the previous paragraph, $\bar{\cM}$ is realizable. By (g)
we have that $\crit(i)>\beta$ so that $G$ is $\bar{\cM}$-generic
over $\P$.
And by (g) and (h) we have that for all $\xi\in\dom(f)$, $\angles{x,y}$
is $i(J^{\cM}_{\gamma_{\xi}})[G]$-generic over
$i_{\xi}(\Q_{\xi})$.

To complete the proof of Lemma \ref{ASublemma}, we must show
that $\cM[G]\models\psi[x,a,G]$. Let $\dot{x}$ and $\dot{G}$ be $\P$-names
in $J^{\cM}_{\beta}$ such that $\dot{x}_G=x$ and $\dot{G}_G=G$.
Since $\psi$ is a $\Pi_1$ formula,
the relation $R(p,\sigma_1,\sigma_2,\sigma_3)\Iff
p\Vdash\psi(\sigma_1,\sigma_2,\sigma_3)$ is $\Pi_1$ definable over $\cM$.
Then since $\crit(i)>\beta$, we have
\begin{gather*}
\cM[G]\models\psi[x,a,G]\\
\text{iff}\\
(\exists p\in G)\cM\models\text{``}
p\Vdash\psi(\dot{x},\check{a},\dot{G})\text{''}\\
\text{iff}\\
(\exists p\in G)\bar{\cM}\models\text{``}
p\Vdash\psi(\dot{x},i(\check{a}),\dot{G})\text{''}\\
\text{iff}\\
\bar{\cM}[G]\models\psi[x,i(a),G]
\end{gather*}
So it suffices to show that
$\bar{\cM}[G]\models\psi[x,i(a),G]$.
Now, $i(a)=\angles{i(\rho),i(\dot{z})}$.
Let $H$ be $\bar{\cM}[G]$-generic
over $\Col(\omega,i(\rho))$. Let $z=(i(\dot{z})_G)_H$. It suffices to
show that $\cM[G][H]\models$ ``There is no rank function for
$T_{x,z}$.'' For this it suffices to show that
$(x,y,z)\in[T]$. 

Since $\max(\code(\alpha))<\crit(i)$, $i(f)$ is an $\alpha$-ladder
above $\beta$ relative to $\bar{\cM}$.
Thus $S(\bar{\cM},i(f),\P,G)$ is defined.
It is not difficult to see that
$z$ codes $(J^{\bar{\cM}}_{i(\rho)},i(f),\P,G)$.
Thus $[T_z]=S(\bar{\cM},i(f),\P,G)$.
So to complete the proof it suffices to show
that $(x,y)\in S(\bar{\cM},i(f),\P,G)$.

Since $\dom(f)\leq\max(\code(\alpha))<\crit(i)$, $\dom(i(f))=\dom(f)$.
For $\xi\in\dom(f)$, let $(\bar{\delta}_{\xi},\bar{\gamma}_{\xi})=
i(\delta_{\xi},\gamma_{\xi})$. For each $\xi\in\dom{f}$,
$\angles{x,y}$ is  $J^{\bar{\cM}}_{\bar{\gamma}_{\xi}}[G]$-generic over
$i_{\xi}(\Q_{\xi})$, and $i_{\xi}(\Q_{\xi})\subseteq
J^{\bar{\cM}}_{\bar{\delta}_{\xi}}$. Let $H_{\xi}$ be
$J^{\bar{\cM}}_{\bar{\gamma}_{\xi}}[G]$-generic over
$\Col(\omega,\bar{\delta}_{\xi})$ with $\angles{x,y}\in
J^{\bar{\cM}}_{\bar{\gamma}_{\xi}}[G][H_{\xi}]$.

For a moment, fix $\xi\in\dom{f}$.
Let $\bar{\R}_{\xi}=\R\intersect
J^{\bar{\cM}}_{\bar{\gamma}_{\xi}}[G][H_{\xi}]$,
and let $\bar{\alpha}_{\xi}=|s_{\xi}|^{L(\bar{\R}_{\xi})}$.
Exactly as we did before, me may show that
$J_{\bar{\alpha}_{\xi}}(\bar{\R}_{\xi})$ is
$\Sigma_{n_{\xi}}$ elementarily embeddable into
$J_{\alpha_{\xi}}(\R)$.

Now we have:
$$ \JalphaR\models(\formulaphi_1 \AND \neg \formulaphi_2)[x,y].$$
So we have:
\begin{gather*}
(\exists\xi\in\dom{f})
\bigl[\JofR{\alpha_{\xi}}\models\theta^1_{\xi}[x,y]\bigr],\\
\text{and}\\
(\forall\xi\in\dom{f})
\bigl[\JofR{\alpha_{\xi}}\models\neg\theta^2_{\xi}[x,y]\bigr].
\end{gather*}
And so we have:
\begin{gather*}
(\exists\xi\in\dom{f})
\bigl[J_{\bar{\alpha}_{\xi}}(\bar{\R}_{\xi})\models\theta^1_{\xi}[x,y]\bigr],\\
\text{and}\\
(\forall\xi\in\dom{f})
\bigl[J_{\bar{\alpha}_{\xi}}(\bar{\R}_{\xi})\models\neg\theta^2_{\xi}[x,y]\bigr].
\end{gather*}
Thus the sequence $\lsequence{H_{\xi}}{\xi\in\dom{f}}$ witnesses that
$(x,y)\in S(\bar{\cM},i(f),\P,G)$.
\end{proof}

\skipbig

We now turn to the proof of Lemma \ref{TechnicalVersion}. The
statement of the lemma is on page \pageref{TechnicalVersion}.

\begin{proof}[Proof of Lemma \ref{TechnicalVersion}]
The lemma is proved by induction on $(\alpha,n)$. For $\alpha=1$
the lemma is essentially proved in Section 4 of \cite{St2}.
So we assume that $\alpha\geq2$. We also assume that  Lemma
\ref{ASublemma} on page \pageref{ASublemma}
holds for $\alpha$. 

Without loss of generality we may assume that there is no
proper initial segment $\cN\lhd\cM$ such that $\beta\in\cN$ and $\cN$ is
explicitly $(\alpha,n)$-big above $\beta$.
(If there were such an $\cN\lhd\cM$ then, since $\beta\geq(\kappa^+)^{\cM}$,
$\R\intersect\cN[G]=\R\intersect\cM[G]$, and so we may replace
$\cM$ with $\cN$.)

We split into 6 cases depending on $(\alpha,n)$.

\underline{CASE 1}\quad  $\cof(\alpha)\leq\omega$ and $n=0$.\qquad
By Proposition \ref{ODone} we have that
$\OD^{\JalphaR}_1 = \Aa{0}$ =
$$\setof{x\in\R}%
{(\exists\alphaprime<\alpha)\, x\in\OD^{\JofR{\alphaprime}}}.$$
We are assuming that
$\cM$ is explicitly $(\alpha,0)$-big above $\beta$.
This implies that for cofinally many $\alphaprime<\alpha$, and all $n$, 
$\cM$ is explicitly $(\alphaprime,n)$-big above $\beta$.
By Lemma \ref{CodeFacts} on page \pageref{CodeFacts}, for all
sufficiently large $\alphaprime<\alpha$, $\max(\code(\alphaprime))
\leq\max(\code(\alpha))$, and so $\max(\code(\alphaprime))<\beta$.
Thus we may apply our induction hypothesis and conclude that
for cofinally many $\alphaprime<\alpha$, and all $n$,
$\R\intersect\cM[G]$ is closed under $\OD^{\JofR{\alphaprime}}_{n+1}$.
Thus $\R\intersect\cM[G]$ is closed under $\OD^{\JalphaR}_1$.

\underline{CASE 2}\quad  $\cof(\alpha)>\omega$ and $n=0$.\qquad
We are assuming that $\cM$ is explicitly $(\alpha,0)$-big above $\beta$,
and no proper initial segment of $\cM$ is. This implies that
$\cM$ is active. Let $\kappa$ be the critical point of the last
extender on the $\cM$ sequence. Then $\kappa>\beta$.
Fix any $\alphaprime<\alpha$ and any $n\in\omega$.
As in the proof of Lemma  \ref{CriterionFourPrime}, 
 there is a countable premouse $\cM^{\prime}$
such that
\begin{itemize}
\item[(a)] $\cM^{\prime}$ is a linear iterate of $\cM$ via its last
extender, and
\item[(b)] letting $i:\cM\map\cM^{\prime}$ be the iteration map
we have that $i(\kappa)>\max(\code(\alphaprime))$, and
\item[(c)] $\cM^{\prime}$ is realizable, and
\item[(c)]$\cM^{\prime}$ is explicitly $(\alphaprime,n)$-big above 
$\beta$.
\end{itemize}
Let $\betaprime$ be any ordinal with 
$\max(\code(\alphaprime))<\betaprime<i(\kappa)$.
By Lemma \ref{BigUpToKappa}, 
$\cM^{\prime}$ is explicitly $(\alphaprime,n)$-big above
$\betaprime$. Thus we may apply the Lemma we are proving,
for $(\alphaprime,n)$, to $\cM^{\prime}$ and $\betaprime$,
to conclude that $\R\intersect\cM^{\prime}[G]$ is closed
under $\OD^{\JofR{\alphaprime}}_{n+1}$.

Notice that for any such $\cM^{\prime}$ as above,
$\R\intersect\cM^{\prime}[G]=
\R\intersect\cM[G]$. 
So we have that for all $\alphaprime<\alpha$ and
all $n$, $\R\intersect\cM[G]$ is closed under 
$\OD^{\JofR{\alphaprime}}_{n+1}$.
As in the previous case, this implies that
$\R\intersect\cM[G]$ is closed under $\OD^{\JalphaR}_1$.

\underline{CASE 3}\quad  $\cof(\alpha)\leq\omega$ and
$n=1$.\qquad  We are assuming that
$\cM$ is explicitly $(\alpha,1)$-big above $\beta$ and that no proper 
initial segment of $\cM$ is explicitly $(\alpha,1)$ big above $\beta$.
This means that there is an ordinal $\rho$ in $\cM$ such that $(\cM,\rho)$
witnesses that $\cM$ is explicitly $(\alpha,1)$-big above $\beta$.
Recall that this means that $\rho$ is a limit cardinal of $\cM$ and:
\begin{itemize}
\item[(i)] $\cM$ is the first admissible structure over $\cJ^{\cM}_{\rho}$, 
and
\item[(ii)] $\code(\alpha)\in J^{\cM}_{\rho}$, and
\item[(iii)] for cofinally many $\alphaprime<\alpha$, and all $n\in\omega$,
and cofinally many $\delta<\rho$ we have that
\begin{itemize}
\item[(a)] $\delta$ is a cardinal of $\cM$, and
\item[(b)] $\delta$ is locally $(\alphaprime,n)$ Woodin 
            in $\cJ^{\cM}_{\rho}$.
\end{itemize}
\end{itemize}

Let $x_0\in\R\intersect\cM[G]$, and 
let $y_0\in \OD^{\JalphaR}_2(x_0)$.  We will show that $y_0\in\cM[G]$.
Fix $\xi<\omega_1$ such that for all $w\in\WO$ with $|w|=\xi$,
$\singleton{y_0}\in\Sigma_{2}(\JalphaR,\singleton{x_0,w})$.

We begin with a claim that says that we may ignore the ordinal
parameter $\xi$.

\begin{claim}[Claim 1]
Without loss of generality we may assume that
$\singleton{y_0}\in\Sigma_{2}(\JalphaR,\singleton{x_0})$.
\end{claim}
\begin{subproof}[Proof of Claim 1]
Recall that we are assuming that $\max(\code(\alpha))<\beta$.
Notice that $\cM$ has a measurable cardinal $\mu$ such $\beta<\mu$,
and such that $\cM$ is explicitly
$(\alpha,1)$-big above $\mu$.
Let $j:\cM\map\bar{\cM}$ be the iteration map resulting from iterating
$\cM$, via $\mu$, $\xi+1$ times. Then $\xi<j(\mu)$.
By Lemma
\ref{IteratingPreservesBigness}, $\bar{\cM}$ is explicitly
$(\alpha,1)$-big above $j(u)$. Notice that $\R\intersect\cM[G]=
\R\intersect\bar{\cM}[G]$. So to see that $y_0\in\cM[G]$, it
suffices to see that $y_0\in\bar{\cM}[G]$.

In order to see that $y_0\in\bar{\cM}[G]$
it suffices to see that $y_0\in\bar{\cM}[G][H]$, whenever $H$ is 
$\bar{\cM}[G]$ generic over $\Col(\omega,j(\mu))$.  
Fix an arbitrary such $H$. Let $J=G\times H$. 
So it suffices to see that $y_0\in\bar{\cM}[J]$.  But
$\bar{\cM}[J]=\bar{\cM}[K]$, for some $K$ which is $\bar{\cM}$-generic over
$\Col(\omega,j(\mu))$.  So in order to show that $y_0\in\cM[G]$ it
suffices to show that $y_0\in\bar{\cM}[K]$ whenever $K$ is 
$\bar{\cM}$-generic
over $\Col(\omega,j(\mu))$ and $x_0\in\bar{\cM}[K]$. 
Since $\bar{\cM}\models (j(\mu))^+$ exists, and 
$\bar{\cM}$ is explicitly $(\alpha,1)$-big above $(j(\mu)^+)^{\bar{\cM}}$,
it suffices to replace $\cM$ with $\bar{\cM}$, and $\P$ with
$\Col(\omega,j(\mu))$, and $\kappa$ with $j(\mu)$, and
$\beta$ with $(j(\mu)^+)^{\bar{\cM}}$. In other words, without loss
of generality we may assume that $\P=\Col(\omega,\kappa)$ and that
$\xi<\kappa$. So we do assume this now.

Then there is a real $w$ in $\WO$ such that $w$ codes $\xi$, and
such that $w$ in $\cM[G]$. We are assuming that
$\singleton{y_0}\in\Sigma_{2}(\JalphaR,\singleton{x_0,w})$.
So we may replace $x_0$ by $\angles{x_0,w}$ to yield Claim 1.
\end{subproof}

So we do now assume that 
$\singleton{y_0}\in\Sigma_{2}(\JalphaR,\singleton{x_0})$. 
This means that as a
set of integers $y_0\in\Delta_2(\JalphaR,\singleton{x_0})$. By lemma
\ref{ASublemma}, $y_0$ is a $\BfDeltaOne(\cM[G])$ set of integers.
As $\cM[G]$ is admissible, $y_0\in\cM[G]$.

\underline{CASE 4}\quad  $\cof(\alpha)\leq\omega$ and $n=2$.\qquad  
We are assuming that
$\cM$ is explicitly $(\alpha,2)$-big above $\beta$ and that no proper 
initial
segment of $\cM$ is explicitly $(\alpha,2)$ big above $\beta$.
This means that there is an ordinal $\delta\in\cM$ with $\beta<\delta$
such that $\cM\models\text{``}\delta$ is a Woodin cardinal'',
and $\cM$ is explicitly $(\alpha,1)$ big above $\delta$.

Let $x_0\in\R\intersect\cM[G]$, and 
let $y_0\in \OD^{\JalphaR}_3(x_0)$.  We will show that $y_0\in\cM[G]$.
Fix $\xi<\omega_1$ such that for all $w\in\WO$ with $|w|=\xi$,
$\singleton{y_0}\in\Sigma_{3}(\JalphaR,\singleton{x_0,w})$.

Exactly as in Case 3,  we may now ignore the ordinal
parameter $\xi$ and assume that 
 $\singleton{y_0}\in\Sigma_{3}(\JalphaR,\singleton{x_0})$.
Recall that there is a $\Sigma_1(\JalphaR)$
partial function from $\R$ onto $\JalphaR$. This implies that
$\Sigma_3(\JalphaR)=\exists^{\R}\Pi_2(\JalphaR)$.
Fix a $\Pi_2$ formula $\formulaphi$ such that
$$\singleton{y_0}=
\setof{y}{(\exists z\in\R)\JalphaR\models\formulaphi[x_0,y,z]}.$$
Fix $z_0\in\R$ such that $\JalphaR\models\formulaphi[x_0,y_0,z_0]$.
We now apply Lemma \ref{MakeRealGeneric}. This gives us
a partial order $\Q\subseteq J^{\cM}_{\delta}$, with
$\Q\in \cM$, and a $0$-maximal
iteration tree $\cT$ on $\cM$ of countable length $\theta+1$ such that
\begin{itemize}
\item[(a)] $\cM^{\cT}_{\theta}$ is realizable, and
\item[(b)] $D^{\cT}=\emptyset$ so that $i^{\cT}_{0,\theta}$ is defined,
and
\item[(c)] $\length{E^{\cT}_{\xi}}<i^{\cT}_{0,\xi}(\delta)$ for all
$\xi<\theta$, and
\item[(d)] $\crit{E^{\cT}_{\xi}} > \beta$ for all $\xi<\theta$
(so $G$ is $\cM^{\cT}_{\theta}$-generic over $\P$), and
\item[(e)] $\angles{y_0,z_0}$ is $\cM^{\cT}_{\theta}[G]$-generic over
$i^{\cT}_{0,\theta}(\Q)$.
\end{itemize}

As in the proof of Claim 1, we may apply Lemma
\ref{IteratingPreservesBigness} to conclude that
$\cM^{\cT}_{\theta}$ is explicitly $(\alpha,1)$-big above
$i^{\cT}_{0,\theta}(\delta)$. Let $\bar{\cM}=\cM^{\cT}_{\theta}$,
and $\bar{\delta}=i^{\cT}_{0,\theta}(\delta)$.
Since $\R\intersect\bar{\cM}[G]=\R\intersect\cM[G]$, it suffices to
show that $y_0\in\bar{\cM}[G]$.

$\bar{\cM}$ is explicitly $(\alpha,1)$ big above $\bar{\delta}$, 
and no proper
initial segment of $\bar{\cM}$ is explicitly $(\alpha,1)$-big 
above $\bar{\delta}$.
This implies that $\bar{\cM}$ is explicitly $(\alpha,1)$ big above
$(\bar{\delta}^+)^{\bar{\cM}}$.
By Lemma \ref{ASublemma} applied to the formula $\neg\formulaphi$,
there is a $\Sigma_1$ formula
$\psi$,and there is
a parameter $a\in\bar{\cM}$ such that whenever $K$ is
$\bar{\cM}$-generic over $\P\times\Col(\omega,\bar{\delta})$,
and $x,y,z\in\R\intersect\bar{\cM}[K]$, we have
\begin{gather*}
\JalphaR\models\formulaphi[x,y,z]\\
\text{iff}\\
\bar{\cM}[K]\models\psi[x,y,z,a,K].
\end{gather*}
Fix such a formula $\psi$ and such a parameter $a$.

Since $\bar{\cM}[G]$ is admissible, to show that $y_0\in\bar{\cM}[G]$ it
is sufficient to show that $y_0$ is a $\BfDeltaOne(\bar{\cM}[G])$ set
of integers. That is the content of the following claim. Let $\dot{H}$
be the canonical $\Col(\omega,\bar{\delta})$-name in $\bar{\cM}[G]$
such that whenever $H$ is $\bar{\cM}[G]$-generic over
$\Col(\omega,\bar{\delta})$, $\dot{H}_H=H$.

\begin{claim}
Let $n,m\in\omega$. Then the following are equivalent:
\begin{itemize}
\item[(a)] $y_0(n)=m$
\item[(b)] $\bar{\cM}[G]\models\exists q\in\Col(\omega,\bar{\delta})$
s.t. $q\Vdash\text{``}(\exists y,z\in\R)
\big[\psi(\check{x_0},y,z,\check{a},\check{G}\times\dot{H})\;
\AND\\ y(n)=m\big]$.''
\item[(c)] $\bar{\cM}[G]\models\forall q\in\Col(\omega,\bar{\delta})$,
$q\Vdash\text{``}(\forall y,z\in\R)
\big[\psi(\check{x_0},y,z,\check{a},\check{G}\times\dot{H})\;
\implies\\ y(n)=m\big]$.''
\end{itemize}
\end{claim}
\begin{subproof}[Proof of Claim]

\noindent
(a) $\Implies$ (b).\quad Suppose $y_0(n)=m$. Let $H$ be
$\bar{\cM}[G]$-generic
over $\Col(\omega,\bar{\delta})$ and such that
$\angles{y_0,z_0}\in\bar{\cM}[G][H]$. Then $G\times H$ is
$\bar{\cM}$-generic over $\P\times\Col(\omega,\bar{\delta})$
and $\JalphaR\models\formulaphi[x_0,y_0,z_0]$, so
$\bar{\cM}[G][H]\models\psi[x_0,y_0,z_0,a,G\times H]$. From this
(b) easily follows.

\noindent
(b) $\Implies$ (a).\quad Fix $q$ witnessing that (b) holds.
Let $H$ be
$\bar{\cM}[G]$-generic
over $\Col(\omega,\bar{\delta})$ and such that
$q\in H$. Then fix $y,z\in\R$ such that $y(n)=m$ and
$\bar{\cM}[G][H]\models\psi[x_0,y,z,a,G\times H]$. Then
$\JalphaR\models\formulaphi[x_0,y,z]$, so $y=y_0$.

\noindent
(a) $\Implies$ (c).\quad Suppose $y_0(n)=m$. Let $H$ be
$\bar{\cM}[G]$-generic
over $\Col(\omega,\bar{\delta})$. Let $y,z\in\R\intersect\bar{\cM}[G][H]$,
and suppose that $\bar{\cM}[G][H]\models\psi[x_0,y,z,a,G\times H]$. Then
$\JalphaR\models\formulaphi[x_0,y,z]$, so $y=y_0$. So $y(n)=m$. Since
$H,y,z$ were arbitrary, (c) follows.

\noindent
(c) $\Implies$ (a).\quad Let $H$ be $\bar{\cM}[G]$-generic
over $\Col(\omega,\bar{\delta})$ and such that
$\angles{y_0,z_0}\in\bar{\cM}[G][H]$. Since
$\JalphaR\models\formulaphi[x_0,y_0,z_0]$, we have that
$\bar{\cM}[G][H]\models\psi[x_0,y_0,z_0,a,G\times H]$.
By (c), $y_0(n)=m$.
\end{subproof}
This completes CASE 4.

\underline{CASE 5}\quad  $n\geq1$ is odd, 
and if $\cof(\alpha)\leq\omega$, then
$n\geq 3$.
\qquad  We are assuming that
$\cM$ is explicitly $(\alpha,n)$-big above $\beta$ and that no proper 
initial segment of $\cM$ is explicitly $(\alpha,n)$-big above $\beta$.
If $\cof(\alpha)>\omega$ this means that $n\geq1$ and
 there are $n$ ordinals
$\delta_1,\dots,\delta_n\in\cM$ with
$\beta<\delta_1<\cdots<\delta_n$
such that for $1\leq i\leq n$,
$\cM\models\text{``}\delta_i$ is a Woodin cardinal'',
and $\cM$ is explicitly $(\alpha,0)$-big above $\delta_n$.
If $\cof(\alpha)\leq\omega$ this means that $n\geq3$ and
 there are $n-1$ ordinals
$\delta_1,\dots,\delta_{n-1}\in\cM$ with
$\beta<\delta_1<\cdots<\delta_{n-1}$
such that for $1\leq i\leq n-1$,
$\cM\models\text{``}\delta_i$ is a Woodin cardinal'',
and $\cM$ is explicitly $(\alpha,1)$-big above $\delta_{n-1}$.

Let $x_0\in\R\intersect\cM[G]$, and 
let $y_0\in \OD^{\JalphaR}_{n+1}(x_0)$.  We will show that $y_0\in\cM[G]$.
Fix $\xi<\omega_1$ such that for all $w\in\WO$ with $|w|=\xi$,
$\singleton{y_0}\in\Sigma_{n+1}(\JalphaR,\singleton{x_0,w})$.

Exactly as in Case 3,  we may now ignore the ordinal
parameter $\xi$ and assume that 
 $\singleton{y_0}\in\Sigma_{n+1}(\JalphaR,\singleton{x_0})$.
Recall that there is a $\Sigma_1(\JalphaR)$
partial function from $\R$ onto $\JalphaR$. This implies that
$\Sigma_{n+1}(\JalphaR)=\exists^{\R}\Pi_{n}(\JalphaR)$.
Fix a $\Pi_{n}$ formula $\formulaphi$ such that
$$\singleton{y_0}=
\setof{y}{(\exists z\in\R)\JalphaR\models\formulaphi[x_0,y,z]}.$$
Fix $z_0\in\R$ such that $\JalphaR\models\formulaphi[x_0,y_0,z_0]$.
We now apply Lemma \ref{MakeRealGeneric}. This gives us
a partial order $\Q\subseteq J^{\cM}_{\delta_1}$, with
$\Q\in \cM$, and a $0$-maximal
iteration tree $\cT$ on $\cM$ of countable length $\theta+1$ such that
\begin{itemize}
\item[(a)] $\cM^{\cT}_{\theta}$ is realizable, and
\item[(b)] $D^{\cT}=\emptyset$ so that $i^{\cT}_{0,\theta}$ is defined,
and
\item[(c)] $\length{E^{\cT}_{\xi}}<i^{\cT}_{0,\xi}(\delta_1)$ for all
$\xi<\theta$, and
\item[(d)] $\crit{E^{\cT}_{\xi}} > \beta$ for all $\xi<\theta$
(so $G$ is $\cM^{\cT}_{\theta}$-generic over $\P$), and
\item[(e)] $\angles{y_0,z_0}$ is $\cM^{\cT}_{\theta}[G]$-generic over
$i^{\cT}_{0,\theta}(\Q)$.
\end{itemize}

Let $\bar{\cM}=\cM^{\cT}_{\theta}$,
and let
$\bar{\delta}_i=i^{\cT}_{0,\theta}(\delta_i)$.
Since $\R\intersect\bar{\cM}[G]=\R\intersect\cM[G]$, it suffices to
show that $y_0\in\bar{\cM}[G]$.

Notice that $\bar{\cM}$ is explicitly $(\alpha,n-1)$-big above
$(\bar{\delta}_1^+)^{\bar{\cM}}$. 
(If $\cof(\alpha)\leq\omega$ then, as in the proof of 
Claim 1, we may apply Lemma
\ref{IteratingPreservesBigness} to conclude that
$\bar{\cM}$ is explicitly $(\alpha,1)$-big above
$\bar{\delta}_{n-1}$. If $\cof(\alpha)>\omega$ then it is easy to see
that $\bar{\cM}$ is explicitly $(\alpha,0)$-big above
$\bar{\delta}_{n}$.)  

Let $H$ be $\bar{\cM}[G]$ generic
over $\Col(\omega,\delta_1)$ and such that
$\angles{y_0,z_0}\in\bar{\cM}[G][H]$.
Let $\bar{\R}=\R\intersect\bar{\cM}[G][H]$.
By the lemma we are proving,
for $(\alpha,n-1)$, $\bar{\R}$ is closed under
$\OD^{\JalphaR}_{n}$. In particular $\bar{\R}$
is closed downward under $\Da{n-1}$. Also, clearly
$\bar{\R}=\R\intersect L(\bar{\R})$. 
By assumption, $\max(\code(\alpha))<\beta<\bar{\delta}_1<
(\omega_1)^{L(\bar{\R}}$.
Since $n-1$ is even, we can apply  Lemma \ref{MaxCodeIsGoodParameter}
to conclude that there is some  ordinal $\bar{\alpha}$ in $\bar{\cM}$
 such that
$J_{\bar{\alpha}}(\bar{\R})$ is 
$\Sigma_{n}$ elementarily embeddable into
$J_{\alpha}(\R)$.

Since $\formulaphi$ is a $\Pi_n$ formula and $y_0,z_0\in\bar{\R}$,
it follows that in $\bar{\cM}[G][H]$ the following is true:
$$\singleton{y_0}=
\setof{y}{(\exists z\in\R)\JofR{\bar{\alpha}}\models\formulaphi[x_0,y,z]}.$$
As the collapse forcing is homogeneous, this implies that
$y_0\in\bar{\cM}[G]$.

This completes CASE 5.

\underline{CASE 6}\quad  $n\geq2$ is even, and if $\cof(\alpha)\leq\omega$
then $n\geq 4$. \qquad  
We are assuming that
$\cM$ is explicitly $(\alpha,n)$-big above $\beta$ and that no proper 
initial segment of $\cM$ is explicitly $(\alpha,n)$-big above $\beta$.
This means that there are two ordinals, $\delta_1,\delta_2$, such that
$\beta<\delta_1<\delta_2$ and both $\delta_1$ and $\delta_2$ are Woodin
in $\cM$ and $\cM$ is explicitly $(\alpha,n-2)$-big above $\delta_2$.

Let $x_0\in\R\intersect\cM[G]$, and 
let $y_0\in \OD^{\JalphaR}_{n+1}(x_0)$.  We will show that $y_0\in\cM[G]$.
Fix $\xi<\omega_1$ such that for all $w\in\WO$ with $|w|=\xi$,
$\singleton{y_0}\in\Sigma_{n+1}(\JalphaR,\singleton{x_0,w})$.

Exactly as in Case 3,  we may now ignore the ordinal
parameter $\xi$ and assume that 
$\singleton{y_0}\in\Sigma_{n+1}(\JalphaR,\singleton{x_0})$.
The proof that we may ignore the ordinal parameter
$\xi$ also yields the following:

\begin{claim}
Without loss of generality we may assume that 
$\max(\code(\alpha))$ is countable in $\cM[G]$.
\end{claim}
\begin{subproof}[Proof of Claim]
Recall that $\beta\geq(\kappa^+)^{\cM}$, where 
$\P\subseteq J^{\cM}_{\kappa}$.
By assumption $\max(\code(\alpha))<\beta$.  
So, exactly as in the proof of Claim 1 in CASE 3, 
without loss of generality
we may assume that  $\max(\code(\alpha))<\kappa$ and
$\P=\Col(\omega,\kappa)$.
\end{subproof}

So we do now assume that $\max(\code(\alpha))$ is countable in $\cM[G]$.
We will now show that $\cM[G]$ is $\Sigma_{n+1}(\JalphaR)$ correct.
More precisely,
let $\bar{\R}=\R\intersect\cM[G]$. Let $s=\code(\alpha)$ and let
$\bar{\alpha}=|s|^{\cM[G]}$.
We will show that there is  a $\Sigma_{n+1}$ elementary embedding
$$j:J_{\bar{\alpha}}(\bar{\R})\map\JalphaR.$$
Since, $\singleton{y_0}\in\Sigma_{n+1}(\JalphaR,\singleton{x_0})$,
this will give us that $y_0\in\bar{\R}$.

\begin{remark}
The original $\cM[G]$ with which we started the proof may {\em not}
be $\Sigma_{n+1}(\JalphaR)$-correct. Because of our application of
the above Claim, we are now
working with a (possibly) modified $\cM[G]$ in which 
$\max(\code(\alpha))$ is countable. This is what will give us that
$\bar{\R}$ is correct.
\end{remark}

As $\cM$ is explicitly $(\alpha,n)$-big above $\beta$, in particular $\cM$ 
is explicitly $(\alpha,n-2)$-big above $\beta$.
By the lemma we are proving, for $(\alpha,n-2)$,
$\bar{\R}$ is closed under $\OD^{\JalphaR}_{n-1}$. In particular,
$\bar{\R}$ is closed down under $\Da{n-2}$. Also, 
$\bar{\R}=\R\intersect L(\bar{\R})$. Also, by assumption,
$\max(\code(\alpha))<\omega_1^{L(\bar{\R})}$.
By Lemma \ref{MaxCodeIsGoodParameter}, since $n-2$ is even,
there is a $\Sigma_{n-1}$ elementary
map
\mbox{$j:(J_{\bar{\alpha}}(\bar{\R});\;\in,\;\bar{\R})\map
(\JalphaR;\;\in,\;\R)$}.
We will complete our proof by showing that $j$ is $\Sigma_{n+1}$ elementary.

Recall that there is a $\Sigma_1(\JalphaR)$
partial function from $\R$ onto $\JalphaR$. 
So to show that $j$ is
$\Sigma_{n+1}(\JalphaR)$-elementary, it suffices to show that for
every $\Pi_{n}$ formula $\formulaphi$ and every $x\in\bar{\R}$, if
$\JalphaR\models(\exists y\in\R)\formulaphi[x,y]$, then
$J_{\bar{\alpha}}(\bar{\R})\models(\exists y\in\bar{\R})\formulaphi[x,y]$.
So fix a $\Pi_{n}$ formula $\formulaphi$, and fix an $x\in\bar{\R}$.
Now fix a $y\in\R$ and suppose that
$\JalphaR\models\formulaphi[x,y]$. we will show that
$J_{\bar{\alpha}}(\bar{\R})\models(\exists y\in\bar{\R})\formulaphi[x,y]$.

We now apply Lemma \ref{MakeRealGeneric}. This gives us
a partial order $\Q\subseteq J^{\cM}_{\delta_1}$, with
$\Q\in \cM$, and a $0$-maximal
iteration tree $\cT$ on $\cM$ of countable length $\theta+1$ such that
\begin{itemize}
\item[(a)] $\cM^{\cT}_{\theta}$ is realizable, and
\item[(b)] $D^{\cT}=\emptyset$ so that $i^{\cT}_{0,\theta}$ is defined,
and
\item[(c)] $\length{E^{\cT}_{\xi}}<i^{\cT}_{0,\xi}(\delta_1)$ for all
$\xi<\theta$, and
\item[(d)] $\crit{E^{\cT}_{\xi}} > \beta$ for all $\xi<\theta$
(so $G$ is $\cM^{\cT}_{\theta}$-generic over $\P$), and
\item[(e)] $y$ is $\cM^{\cT}_{\theta}[G]$-generic over
$i^{\cT}_{0,\theta}(\Q)$.
\end{itemize}
Let $\bar{\cM}=\cM^{\cT}_{\theta}$,
and for $1\leq i\leq 2$, let
$\bar{\delta}_i=i^{\cT}_{0,\theta}(\delta_i)$.
Notice that $\R\intersect\bar{\cM}[G]=\bar{\R}$.
As in the previous case, we have that
$\bar{\cM}$ is explicitly $(\alpha,n-2)$-big above
$\bar{\delta}_2$. As no proper initial segment of
$\bar{\cM}$ is explicitly $(\alpha,n-2)$-big above
$\bar{\delta}_2$, we have that
$\bar{\cM}$ is explicitly $(\alpha,n-2)$-big above
$(\bar{\delta}_2^+)^{\bar{\cM}}$.

Since $y$ is $\bar{\cM}[G]$ generic over $i^{\cT}_{0,\theta}(\Q)$, and
$i^{\cT}_{0,\theta}(\Q)\subset J^{\bar{\cM}}_{\bar{\delta}_1}$,
we can absorb $y$ into a generic object for the stationary tower
forcing over $\bar{\cM}[G]$ up to $\bar{\delta}_2$. (For a
description of Woodin's Stationary Tower Forcing, see
\cite{Ma2}.)
Let $H$ be this
generic object. Let $\cN$ be the generic ultrapower of
$\bar{\cM}[G]$ by $H$,
and let
$\pi:\bar{\cM}[G]\map\cN$ be the generic elementary embedding.
It is one of the fundamental properties of the Stationary Tower
forcing that $\R\intersect\cN=\R\intersect\bar{\cM}[G][H]$.
Let $\R^*=\R\intersect\bar{\cM}[G][H]$.

Recall   that $s=\code(\alpha)$, that
$\max(s)<\omega_1^{L(\bar{\R})}$, and  that 
$\bar{\alpha}=|s|^{L(\bar{\R})}$.
Let $\alpha^*=|s|^{L(\R^*)}$. Then 
$\pi(J_{\bar{\alpha}}(\bar{\R}))=J_{\alpha^*}(\R^*)$.

Now $\bar{\cM}[G][H]$ is a generic extension of $\bar{\cM}$ via
a partial order $\Z\subset J^{\bar{\cM}}_{\bar{\delta}_2}$.
Since $\bar{\cM}$ is explicitly $(\alpha,n-2)$ big above
$(\bar{\delta}_2^+)^{\bar{\cM}}$, we have by the lemma we are proving,
for $(\alpha,n-2)$, that $\R^*$ is closed under
$\OD^{\JalphaR}_{n-1}$. In particular, $\R^*$ is closed downward
under $\Da{n-2}$.
Also, 
$\R^*=\R\intersect L(\R^*)$. 
By Lemma \ref{MaxCodeIsGoodParameter}, there is a $\Sigma_{n-1}$ elementary
map
\mbox{$i:(J_{\alpha^*}(\R^*);\;\in,\;\R^*)\map
(\JalphaR;\;\in,\;\R)$}.

As $\formulaphi$ is a $\Pi_{n}$
formula, and $x,y\in\R^*$, 
$J_{\alpha^*}(\R^*)\models\formulaphi[x,y]$.
As $\pi$ is fully elementary 
and $x\in\bar{\R}$, we have that
$J_{\bar{\alpha}}(\bar{\R})\models(\exists y\in\bar{\R})\formulaphi[x,y]$.
This completes our proof.
\end{proof}

%


\skipbig

\section{The Ranked Iteration Game}

\label{section:iterability}

Suppose that $\cM$ is a countable, iterable premouse, and
$\cM$ is $(\alpha,n)$-big, but every proper initial segment
of $\cM$ is $(\alpha,n)$-small. In the previous section we
proved that $\Aan\subseteq\R\intersect\cM$. Our next major goal
is to prove that $\R\intersect\cM\subseteq\Aan$. This will give
us that $\Aan=\R\intersect\cM$, and so $\Aan$ is a mouse set.

Showing that $\R\intersect\cM\subseteq\Aan$ essentially amounts
to showing that every proper initial segment of $\cM$ is sufficiently
simply definable. Suppose $\cN\lhd\cM$ is a proper initial segment of
$\cM$. Suppose further that $\cN$ projects to $\omega$.
Why should $\cN$ be simply definable? The reason is that
$\cN$ is the unique premouse $\cN^{\prime}$
such that $\cN^{\prime}$ projects to $\omega$, $\cN^{\prime}$ is sound,
$\cN^{\prime}$ has
the same ordinals as $\cN$, and $\cN^{\prime}$ is fully iterable.
(This fact is demonstrated by performing a \emph{comparison} of
$\cN^{\prime}$ with $\cN$.) This would give us a simple definition
of $\cN$, except for the fact that being fully iterable is not
simply definable. 

The strategy for showing that $\cN$ is simply definable is the following:
Find a condition other than ``full iterability'' which is simply definable,
and yet is sufficient to allow us to perform a comparison. Such a condition
is called a definable ``iterability condition.'' That is the topic of
this section.

We follow \cite{MaSt} in using ``iteration games'' as the basis
for our iterability conditions.
 We work in the theory ZF + DC.

We will need a coding function $x\mapsto x^*$ from $\R$ onto $\HC$ = the
collection
of hereditarily countable sets.  For the sake of concreteness we
will describe a particular such coding function.  Let $\cT$ be the
collection of wellfounded trees on $\omega$. First we define
a function $h:\cT\surjection\HC$ by induction on the rank of
a tree $T\in\cT$. Let
$$T_{\angles{n}}=\setof{s\in\omega^{<\omega}}{\angles{n}\frown s\in T},$$
and set
$$h(T)=\setof{h(T_{\angles{n}})}{\angles{n}\in T}.$$
Now let
$\pi:\R\surjection\Powerset(\omega^{<\omega})$ be recursive,
and set $\WF = \setof{x}{\pi(x)\in\cT}$.
Finally, for $x\notin \WF$, set $x^*=\emptyset$, and
for $x\in \WF$, set $x^*=h(\pi(x))$. Then if $n\geq1$ and $\theta$
is any $\Sigma_n$ formula then the set
$$\setof{x\in\R}{\HC\models\theta[x^*]}$$
is $\Sigma^1_{n+1}$.

We will need several concepts from \cite{MaSt}. The following definition
is taken from that paper.

\begin{definition}
Let $\cM$ be a premouse, and $\delta<\ORD^{\cM}$.
\begin{itemize}
\item[(a)] $k(\cM,\delta)$ is the unique $k<\omega$ such that $\cM$ is
$k$-sound, $k+1$ solid, and $\rho_{k+1}(\cM)\leq\delta<\rho_k(\cM)$,
if such a $k$ exists. $k(\cM,\delta)$ is undefined otherwise.
\item[(b)] $\cM$ is a \emph{$\delta$-mouse} iff $k(\cM,\delta)$ is
defined, and letting $k=k(\cM,\delta)$, we have that
$\cM=\cH^{\cM}_{k+1}(\delta\union p_{k+1}(\cM))$.
\item[(c)] An iteration tree $\cT$ on $\cM$ is \emph{above} $\delta$
iff $\crit(E)\geq\delta$ for all extenders $E$ used on $\cT$.
\end{itemize}
\end{definition}

See \cite{MiSt} for the definitions of the various fine structural terms
used above.

\begin{definition}
If $x\in\R$ then we say that $x$ is a {\em mouse code} iff
$x^*=(\cM,\,\delta)$ where $\cM$ is a countable premouse
and $\delta\in\ORD$
and  either
\begin{itemize}
\item[(1)] $\delta=\ORD^{\cM}$, or
\item[(2)] $\delta\in\ORD^{\cM}$ and $\cM$ is a $\delta$-mouse.
\end{itemize}
\end{definition}

We will need a few more concepts from \cite{MaSt}. The following 
definition is  also taken from that paper.

\begin{definition}
Let $\cT$ be an iteration tree of limit length.
\begin{itemize}
\item[(a)] $\delta(\cT)$ is the supremum of the lengths of the extenders
used in $\cT$.
\item[(b)] Let $b$ be a cofinal wellfounded branch of $\cT$. Then
$Q(b,\cT)=\cJ^{\cM^{\cT}_b}_{\alpha}$, where $\alpha$ is the largest $\beta$
such that $\beta=\delta(\cT)$, or $\delta(\cT)<\beta\leq\ORD^{\cM^{\cT}_b}$
and $\cJ^{\cM^{\cT}_b}_{\beta}\models\delta(\cT)$ is a Woodin cardinal.
\end{itemize}
\end{definition}

In \cite{MaSt} it is shown that if $\cM$ is a $\delta$-mouse and
$\cT$ is on $\omega$-maximal tree on $\cM$ above $\delta$, and
$\delta(\cT)<\ORD^{Q(b,\cT)}$, then $Q(b,\cT)$ is a $\delta(\cT)$-mouse.

We will also need from \cite{MaSt} the notion of a \emph{putative iteration
tree}. This is a tree $\cT$
of premice which has all of the properties of an
iteration tree, except that if $\cT$ has successor length, then the last
model of $\cT$ is allowed to be illfounded.

\skipmed

Let $\alpha\geq1$.  Let $x$ be a mouse code.
We now describe $G^{\alpha}(x)$,
the {\em rank-$\alpha$ iteration game}
on $x$. $G^{\alpha}(x)$ is a clopen game on $\R\times\alpha$, so
a run of $G^{\alpha}(x)$ is completed after some finite number of rounds.

At the start of round $n$ with $1\leq n<\omega$, we have an ordinal
$\alpha_{n-1}$ and a real $x_{n-1}$ such that $x_{n-1}$ is a mouse
code.  (We get started by setting $\alpha_0=\alpha$ and
$x_0=x$.) Round $n$ begins with player $\ONE$ playing an ordinal
$\alpha_n<\alpha_{n-1}$. If $\alpha_{n-1}=0$ so that player $\ONE$ can
not play $\alpha_n$, then player $\ONE$ loses.  Otherwise let
$(\cM,\,\delta) = x_{n-1}^*$.  Then player $\ONE$ must play a real $y_n$
such that $y_n^*=\cT$ where $\cT$ is an $\omega$-maximal, putative
iteration tree on $\cM$ above $\delta$.  (We allow player $\ONE$ to play
the trivial tree with only one model, so that player $\ONE$ can
always meet this obligation.) If $\cT$ has a last model
which is wellfounded then player $\ONE$ loses. Otherwise, player
$\TWO$ must play a mouse code $x_n$ such that there is some
maximal wellfounded branch $b$ of $\cT$ so that
$x_n^*=\bigl\langle Q(b,\cT\restriction\sup(b)),\,
\delta\bigl(\cT\restriction\sup(b)\bigr)\,\bigr\rangle$.
If there is no
maximal wellfounded branch $b$ of $\cT$ then player $\TWO$ loses.
This completes round $n$.

Since there is no infinite decreasing sequence of ordinals, one
of the two players will win after a finite number of rounds.  So
this completes the description of $G^{\alpha}(x)$.
\begin{definition}
Let $\cM$ be a countable $\delta$-mouse. Let $\alpha\geq1$.
Then $\cM$ is {\em rank-$\alpha$ iterable above $\delta$} iff
for some (equivalently, for all) $x\in\R$ such that
$x^*=(\cM,\,\delta)$, there is a winning quasi-strategy for player
$\TWO$ in $G^{\alpha}(x)$.  (From now on we will abbreviate this
to ``$\TWO$ wins $G^{\alpha}(x)$''.)
\end{definition}

Our immediate goal is to explore the definability of the set
$$\setof{x\in\R}{\TWO\text{ wins }G^{\alpha}(x)}.$$
\begin{definition}
\label{def:lambda}
For $\alpha\geq2$ set
$$\lambda(\alpha)=
\begin{cases}
\omega\alpha& \text{if $\alpha\geq\omega$}\\
\omega(\alpha - 1)& \text{if $\alpha<\omega$}.
\end{cases}
$$
\end{definition}

In Corollary \ref{PiOneGame} below we show that if $\alpha\geq2$
and $\lambda=\lambda(\alpha)$, then the set
$$\setof{x\in\R}%
{x\text{ is a mouse code }\AND\TWO\text{ wins }G^{\lambda}(x)}$$
is $\Pi_1(\JalphaR)$ uniformly in $\alpha$.

If $x$ is a mouse code and
$p=\angles{y_1,x_1,\dots y_n,x_n}$ is a sequence of reals
of even length then we will say that $p$ is a {\em legal play
on $x$ with $\ONE$ to move} iff there are ordinals
$\alpha,\alpha_1\dots\alpha_n$ such that letting
$\pprime = \angles{\alpha_1,y_1,x_1,\dots\alpha_n,y_n,x_n}$ we have
that $\pprime$ is a legal play in $G^{\alpha}(x)$.
Notice that this is equivalent
to saying that if we let
$\pprime=\angles{n-1,y_1,x_1,\dots 0,y_n,x_n}$ then
$\pprime$ is a legal play in $G^n(x)$. Notice that the
relation ``$p$ is a legal play on $x$ with $\ONE$ to move''
is a simply definable relation of $p$ and $x$. (The relation
is $\Sigma^1_2$.)
The definition of a
{\em legal play on $x$ with $\TWO$ to move} is similar.

The following lemma gives an alternate characterization of
the property ``$\TWO$ wins $G^{\alpha}(x)$''. The lemma says that we could
also have
defined this relation on $x$ and $\alpha$ by induction on $\alpha$.
\begin{lemma}
\label{InductiveGame}
Let $x$ be a mouse code and $\alpha\geq1$.
\begin{itemize}
\item[(1)] If $\alpha=\beta+1$ then $\TWO$ wins $G^{\alpha}(x)$
iff $(\forall y_1\in\R)(\exists x_1\in\R)\big[\angles{y_1}$ is a legal
play on $x$ with $\TWO$ to move $\Implies (\angles{y_1,x_1}$ is a legal
play on $x$ with $\ONE$ to move, and $\TWO$ wins $G^{\beta}(x_1)\;)\,\big]$.
\item[(2)] If $\alpha$ is a limit ordinal then $\TWO$ wins
$G^{\alpha}(x)$ iff $(\forall \beta<\alpha)$ $\TWO$ wins
$G^{\beta}(x)$.
\end{itemize}
\end{lemma}
\begin{proof}
One direction of both (1) and (2) is trivial.  For the other direction,
the cautious reader might ask whether we need to use the Axiom of
Choice in order to pick winning quasi-strategies to piece together.
But since $G^{\beta}(x)$ is a clopen game, if player $\TWO$ wins
$G^{\beta}(x)$, then he has a canonical winning quasi-strategy. So the
Axiom of Choice is in fact not needed.
\end{proof}

Recall the definition of $\lambda(\alpha)$ in Definition
\ref{def:lambda} above. Notice that $\lambda(\alpha)\leq\ORD^{\JalphaR}$.
For $\gamma\in\ORD$, let $F(\gamma)=\{x\in\R\;:\; x$
is a mouse code and player $\TWO$ wins $G^{\gamma}(x) \; \}$.

\begin{lemma}
Let $\alpha\geq2$. Let $\lambda=\lambda(\alpha)$. 
\begin{itemize}
\item[(1)] For all $\gamma<\lambda$, $F\restr\gamma \in \JalphaR$.
\item[(2)] $F\restr\lambda$ is $\Delta_1(\JalphaR)$.
\end{itemize}
\end{lemma}
\begin{proof}
Let $P$ be the binary relation defined by 
$P(f,\gamma)$ iff $f=F\restr\gamma$.
Notice that $P$ is a $\Sigma_0$ relation of $f$ and $\gamma$, 
with $\R$ as a parameter.
This follows from the previous lemma.

First we will show that $(1) \implies (2)$. Assuming $(1)$, we
have that for $\gamma<\lambda$, $S = F(\gamma)$ iff there is
a function $f\in\JalphaR$ such that $P(f,\gamma+1)$ and
$S=f(\gamma)$. This shows that $F\restr\lambda$ is $\Sigma_1(\JalphaR)$.
Since $S=F(\gamma)$ iff $(\forall T)$ $[T=F(\gamma)\rightarrow S=T]$,
we have that $F\restr\lambda$ is also $\Pi_1(\JalphaR)$.

Next we will prove (1) by induction on $\alpha$.
As the basis of our induction, we start
with $\alpha=2$. Then $\lambda(\alpha)=\omega$. For $n<\omega$ it
is easy to see that $F\restr n$ is analytically definable, and so
$F\restr n\in J_2(\R)$. If $\alpha$ is a limit ordinal, then
$\lambda(\alpha)=\sup\setof{\lambda(\beta)}{\beta<\alpha}$, and so
the inductive step of the proof of (1) is trivial for limit $\alpha$.
Suppose $\alpha=\beta+1$ is a successor ordinal.  
Then $\lambda(\alpha)=\lambda(\beta)+\omega$.
By (2) for $\beta$, $F\restr\lambda(\beta)\in \JalphaR$. We show by
induction on $n$ that $F\restr(\lambda(\beta)+n)\in\JalphaR$.
Let $f=F\restr(\lambda(\beta)+n)$, and $g=F\restr(\lambda(\beta)+n+1)$.
Suppose $f\in\JalphaR$. Let $S=g(\lambda(\beta)+n)$. Then $S\subset \R$
and, by the previous
lemma, $S$ is $\Sigma_0$ definable with $f$ and $\R$ as parameters.
Thus $S\in\JalphaR$. As $g=f\union\singleton{\angles{\lambda(\beta)+n,S}}$,
we have that $g\in\JalphaR$.
\end{proof}

\begin{corollary}
\label{PiOneGame}
Let $\alpha\geq2$. Let $\lambda=\lambda(\alpha)$. Let $n\in\omega$.
Then the set
$$\setof{x}%
{x\text{ is a mouse code }\AND\TWO\text{ wins }G^{\lambda+n}(x)}$$
is $\Pa{2n}$, uniformly in $\alpha$.
\end{corollary}
\begin{proof}
For $n=0$,  this follows from part (2) of the previous Lemma
and part (2) of Lemma \ref{InductiveGame}.
For $n>0$ the result follows  by induction on $n$, using 
part (1) of Lemma \ref{InductiveGame}.
\end{proof}

We have defined rank-$\alpha$ iterability for all ordinals $\alpha$.
  Before
going on we should point out that the definition only has significance
for $\alpha\leq\kappa^{\R}$.  To explain this we need to define the
{\em length} $\omega$ (weak) iteration game on $x$, $WG(x)$.
We turn to this now. Instead of explicitly giving the definition of
$WG(x)$, we only describe how $WG(x)$ differs from the
\emph{ranked} iteration game described above:
The definition of $WG(x)$ is similar to the definition of
the ranked iteration game, except that player $\ONE$ does
not play any ordinals.  So a run of $WG(x)$ can last for $\omega$
rounds, in which case we say that $\TWO$ wins.  So
$WG(x)$ is a game on $\R$ with closed payoff for player $\TWO$.
This completes our description of $WG(x)$. We will also refer
to $WG(x)$ as ``the length $\omega$ weak iteration game on $x$.''

\begin{note}
Our definition of the weak iteration game differs slightly from
definitions given in other papers. See for example \cite{MaSt}
and \cite{St2}. The reason for the name ``weak'' is that there is something
else called the strong iteration game.
See \cite{IT}. We will have no need for the strong iteration game in this 
paper.
\end{note}

Let $\kappa^{\R}$ be the least ordinal $\kappa$ such that
$\JofR{\kappa}$ is admissible.
The following is a well-known fact about admissible sets and closed
games. The result
has nothing to do with the particular details of the iteration
games we have defined.

\begin{proposition}
Let $x$ be a mouse code.  Then the following are equivalent:
\begin{itemize}
\item[(1)] Player $\TWO$ wins $WG(x)$.
\item[(2)] For all ordinals $\alpha$ player $\TWO$ wins $G^{\alpha}(x)$.
\item[(3)] Player $\TWO$ wins $G^{\kappa^{\R}}(x)$.
\end{itemize}
\end{proposition}
\begin{proof}
Clearly $(1)\Implies(2)\Implies(3)$. We will show that $(3)\Implies(1)$.
For this it suffices to show that
for all mouse codes $x$ such that player $\TWO$ wins $G^{\kappa^{\R}}(x)$,
and for all legal first moves $y_0$ by player $\ONE$ in $WG(x)$,
there is an $x_0$, a legal response by player $\TWO$ in $WG(x)$,
such that player $\TWO$ wins the game $G^{\kappa^{\R}}(x_0)$. Suppose
this were false. Let $x$ and $y_0$ be a counterexample. For each
$x_0$ which is a legal response by player $\TWO$ in $WG(x)$, let $f(x_0)$ be
the least ordinal $\alpha$ such that player $\TWO$ does not win
$G^{\alpha}(x_0)$. Since player $\TWO$ wins $G^{\kappa^{\R}}(x)$, the
range of $f$ is cofinal in $\kappa^{\R}$.  But $f$ is $\BfSigmaOne$
definable over $\JofR{\kappa^{\R}}$. This contradicts the fact
that $\JofR{\kappa^{\R}}$ is admissible.
\end{proof}
\begin{definition}
Let $\cM$ be a countable $\delta$-mouse.
Then $\cM$ is {\em length-$\omega$ iterable above $\delta$} iff
for some (equivalently, for all) $x\in\R$ such that
$x^*=(\cM,\,\delta)$, there is a winning quasi-strategy for player
$\TWO$ in $WG(x)$.
\end{definition}

\skipmed

In the case $n=0$, Corollary \ref{PiOneGame} above tells us that
the set
$$\setof{x}%
{x\text{ is a mouse code }\AND\TWO\text{ wins }G^{\lambda(\alpha)}(x)}$$
is $\Pa{0}$ definable. This result is satisfactory for our applications
in the sequel. However, if $n\geq1$, then 
Corollary \ref{PiOneGame} tells us that
the set
$$\setof{x}%
{x\text{ is a mouse code }\AND\TWO\text{ wins }G^{\lambda(\alpha)+n}(x)}$$
is $\Pa{2n}$ definable. This is \emph{not} satisfactory for
our applications. Our solution is to replace $G^{\lambda(\alpha)+n}(x)$
with a different game in the case $n>1$.
Let $x$ be a mouse code. Suppose $2\leq\alpha<\kappa^{\R}$,
and let $\lambda=\lambda(\alpha)$.  Let $n\in\omega$.
As in \cite{MaSt} we now define a game $G_{(\alpha,n)}(x)$
which is an approximation of the the game $G^{\lambda+n}(x)$.
The reason we do this is that under the assumption of
$\Det(\JofR{\alpha+1})$, the set
$$\setof{x}%
{x\text{ is a mouse code }\AND\TWO\text{ wins }G_{(\alpha,n)}(x)}$$
will be $\Pa{n}$, uniformly in $\alpha$. This beats the conclusion
of Lemma \ref{PiOneGame} by a factor of 2.

As we mentioned above, 
for $n=0$ the game $G^{\lambda}(x)$ is perfectly satisfactory, and
so we set $G_{(\alpha,0)}(x)\defeq G^{\lambda}(x)$.  
For $n>0$, $G_{(\alpha,n)}(x)$ will not be exactly the same
as $G^{\lambda+n}(x)$.  The
definition of $G_{(\alpha,n)}(x)$ for $n>0$ depends on whether
$n$ is even or odd. The idea is to use the Bounded Quantification
Lemma, that is Lemma \ref{BoundedQuantification} on page
\pageref{BoundedQuantification}, in order to reduce the number of
quantifiers needed to define the game.

\begin{note}
The clever idea of using the Bounded Quantification Lemma to simplify the
definition of the iteration game is \emph{not} due to the author.
The idea appears in \cite{MaSt}.
\end{note}

First let $n\geq2$ be even. Then $G_{(\alpha,n)}(x)$ is played in
$n$ rounds. At the start of round $k$, where $1\leq k\leq n$,
we have a mouse code $x_{k-1}$.
(We get started by setting $x_0=x$.) The rules for round $k$ depend
on whether $k$ is even or odd. First suppose that $1\leq k\leq n-1$
and $k$ is odd. Round $k$ is played as follows:
Let $(\cM,\,\delta) = x_{k-1}^*$. Player $\ONE$ must play two reals,
$y_k$ and $z_k$,
such that $y_k^*=\cT$ where $\cT$ is an $\omega$-maximal, putative
iteration tree on $\cM$ above $\delta$, and $z_k$ is arbitrary.
(We allow player $\ONE$ to play
the trivial tree with only one model, so that player $\ONE$ can
always meet this obligation.) If $\cT$ has a last model
which is wellfounded then player $\ONE$ loses. Otherwise, player
$\TWO$ must play a real $x_k$ such that $x_k$ is a mouse code and
there is some maximal wellfounded
branch $b$ of $\cT$ so that
$x_k^*=\bigl\langle Q(b,\cT\restriction\sup(b)),\,
\delta\bigl(\cT\restriction\sup(b)\bigr)\,\bigr\rangle$.
If there is no maximal wellfounded branch of $\cT$
then player $\TWO$ loses.  This completes
round $k$ in the case that $k$ is odd.

Now suppose that $2\leq k\leq n$
and $k$ is even. Round $k$ is played as follows:
Let $(\cN,\,\gamma) = x_{k-1}^*$. Player $\ONE$ must play a real,
$y_k$
such that:
\begin{itemize}
\item[(i)] $y_k^*=\cU$ where $\cU$ is an $\omega$-maximal, putative
iteration tree on $\cN$ above $\gamma$.  (We allow player $\ONE$ to play
the trivial tree with only one model, so that player $\ONE$ can
always meet this obligation.)
And,
\item[(ii)] $y_k\in\Da{n-k+1}(\angles{x_{k-1},z_{k-1}})$.
\end{itemize}
If there is no such real $y_k$, or
if $\cU$ has a last model
which is wellfounded, then player $\ONE$ loses. Otherwise, player
$\TWO$ must play a real $x_k$ such that $x_k$ is a mouse code and
there is some maximal wellfounded
branch $c$ of $\cU$ so that
$x_k^*=\bigl\langle Q(c,\cU\restriction\sup(c)),\,
\delta\bigl(\cU\restriction\sup(c)\bigr)\,\bigr\rangle$.
If there is no maximal wellfounded branch of $\cU$
then player $\TWO$ loses.  This completes
round $k$ in the case that $k$ is even.

After $n$ rounds of play, if neither player has lost then the play
will have produced a real $x_n$.
Then player $\TWO$ wins this run of  $G_{(\alpha,n)}(x)$
iff  player $\TWO$ has a winning quasi-strategy in the
rank-$\lambda$ iteration game on $x_n$, $G^{\lambda}(x_n)$, where
$\lambda=\lambda(\alpha)$. This completes the description of
$G_{(\alpha,n)}(x)$ in the case that $n\geq2$ is even.

Notice that for even $n\geq 2$, the game $G_{(\alpha,n)}(x)$ demands more
of Player $\ONE$ than does the game $G^{\lambda+n}(x)$. Thus we have
the following.
\begin{proposition}
Let $n\geq2$ be even. If $\TWO$ wins
$G^{\lambda+n}(x)$ then $\TWO$ wins $G_{(\alpha,n)}(x)$.
\end{proposition}

Furthermore we have the following.
\begin{lemma}
\label{EvenDefinableGame}
Suppose $2\leq\alpha<\kappa^{\R}$.
Assume $\Det(\JofR{\alpha+1})$.
Let $n\geq2$ be even. Then the set
$$\setof{x}%
{x\text{ is a mouse code }\AND\TWO\text{ wins }G_{(\alpha,n)}(x)}$$
is $\Pa{n}$, uniformly in $\alpha$.
\end{lemma}
\begin{proof}
There is a relation $R\in\Pi_1(\JalphaR)$ such that if $x$ is a
mouse code, then $\TWO$ wins $G_{(\alpha,n)}(x)$ iff
\begin{multline*}
(\forall y_1,z_1)(\exists x_1)(\forall y_2\in\Da{n-1}(\angles{x_1,z_1})
(\exists x_2)\cdots\\
(\forall y_{n-1},z_{n-1})(\exists x_{n-1})
(\forall y_n\in\Da{1}(\angles{x_{n-1},z_{n-1}})
(\exists x_n)R(x,x_1,y_1,\dots,x_n,y_n).
\end{multline*}
Now $R\in\Pa{0}$, and for odd $k\geq 1$, $\Pa{0}\subseteq\Sa{k}$.
So $R\in\Sa{k}$ for all odd $k$. Using
the Bounded Quantification Lemma, it is easy to see that our claim is
true.
\end{proof}

For even $n\geq 2$, one should think of
the game $G_{(\alpha,n)}(x)$ defined above as an approximation
of the game $G^{\lambda+n}(x)$, where $\lambda=\lambda(\alpha)$.
Notice that the game $G^{\lambda}(x)$ is employed in the definition
of $G_{(\alpha,n)}(x)$. The relation
``$\TWO$ wins $G^{\lambda}(x)$'' is a $\Pa{0}$ relation of $x$, and
since $n$ is even $\Pa{0}\subseteq\Pan$.

Now suppose that $n\geq 1$ is odd. Then, recall that we do not
necessarily have that $\Pa{0}\subseteq\Pan$. This means that
we can not employ the game $G^{\lambda}(x)$ in our definition
of the game $G_{(\alpha,n)}(x)$. But yet we still want to have
that the game $G_{(\alpha,n)}(x)$ is an ``approximation'' of
the game $G^{\lambda+n}(x)$. So we must use a slightly different
strategy. We continue with our definition  of $G_{(\alpha,n)}(x)$ below.

Let $n\geq1$ be odd. Then $G_{(\alpha,n)}(x)$ is played in
$n+1$ rounds. At the start of round $k$, where $1\leq k\leq n+1$,
we have a mouse code $x_{k-1}$.
(We get started by setting $x_0=x$.) Round $k$ is played as follows:
Let $(\cM,\,\delta) = x_{k-1}^*$. Player $\ONE$ must play a real $y_k$
such that $y_k^*=\cT$ where $\cT$ is an $\omega$-maximal, putative
iteration tree on $\cM$ above $\delta$.  (We allow player $\ONE$ to play
the trivial tree with only one model, so that player $\ONE$ can
always meet this obligation.) If $\cT$ has a last model
which is wellfounded then player $\ONE$ loses. Otherwise, player
$\TWO$ must play a real $x_k$ such that:
\begin{itemize}
\item[(i)] $x_k$ is a mouse code and there is some maximal wellfounded
branch $b$ of $\cT$ so that
$x_k^*=\bigl\langle Q(b,\cT\restriction\sup(b)),\,
\delta\bigl(\cT\restriction\sup(b)\bigr)\,\bigr\rangle$.
And,
\item[(ii)] $x_k\in\Dan(y_k)$.
\end{itemize}
If there is no such $x_k$ then player $\TWO$ loses.  This completes
round $k$.

After $n+1$ rounds of play, if neither player has lost then the play
will have produced a real $x_{n+1}$.
Then player $\TWO$ wins this run of  $G_{(\alpha,n)}(x)$
iff $(\exists \sigma\in\JalphaR)$ such that $\sigma$ is a winning
quasi-strategy for player $\TWO$ in the {\em length} $\omega$
(weak) iteration game on $x_{n+1}$, $WG(x_{n+1})$. This completes the
description of $G_{(\alpha,n)}(x)$ in the case that $n$ is odd.

Notice that for odd $n\geq 1$, the game $G_{(\alpha,n)}(x)$ demands more
of Player $\TWO$ than does the game $G^{\lambda+n}(x)$. In fact, it demands
more of Player $\TWO$ than does the game $WG(x)$.
Thus we have the following.
\begin{proposition}
Let $n\geq1$ be odd. If $\TWO$ wins $G_{(\alpha,n)}(x)$ then
$\TWO$ wins $WG(x)$, and so in particular $\TWO$ wins
$G^{\lambda+n}(x)$.
\end{proposition}

Furthermore we have the following.
\begin{lemma}
\label{OddDefinableGame}
Suppose $2\leq\alpha<\kappa^{\R}$.
Assume $\Det(\JofR{\alpha+1})$.
Let $n\geq1$ be odd. Then the set
$$\setof{x}%
{x\text{ is a mouse code }\AND\TWO\text{ wins }G_{(\alpha,n)}(x)}$$
is $\Pa{n}$, uniformly in $\alpha$.
\end{lemma}
\begin{proof}
There is a relation $R\in\Sigma_1(\JalphaR)$ such that if $x$ is a
mouse code, then $\TWO$ wins $G_{(\alpha,n)}(x)$ iff
$$(\forall y_1)(\exists x_1\in\Dan(y_1))\cdots(\forall y_{n+1})
(\exists x_{n+1}\in\Dan(y_{n+1})) R(x,x_1,y_1,\dots,x_{n+1},y_{n+1}).$$
Since $n$ is odd, $R\in\Pan$. Our claim then follows from the
Bounded Quantification Lemma.
\end{proof}

\begin{remark}
The reader may wonder about our formulation of the game
$G_{(\alpha,n)}(x)$ in the case that $n$ is odd.  
In particular it may seem surprising that
the game is played in $n+1$ rounds instead of $n$.
The explanation is that we want to be able to show that player 
$\TWO$ can win
$G_{(\alpha,n)}(x)$ whenever $x$ is a mouse code and, writing
$x^*=(\cM,\delta)$, we have that
 $\cM$ is $(\alpha,n+1)$-small above $\delta$,
and $\cM$ is fully iterable. This fact will be needed in the proof
of Theorem \ref{ComparisonTheorem}. In order to be able to show that
player $\TWO$ can win $G_{(\alpha,n)}(x)$ in this situation, we
need to have the game played for $n+1$ rounds. Notice that we still
get Lemma \ref{OddDefinableGame} above even though $G_{(\alpha,n)}(x)$
is played in $n+1$ rounds.
\end{remark}

\begin{definition}
Let $\cM$ be a countable $\delta$-mouse. Let $\alpha\geq2$.
Let $n\in\omega$.
Then $\cM$ is {\em $\Pan$ iterable above $\delta$} iff
for some (equivalently, for all) $x\in\R$ such that
$x^*=(\cM,\,\delta)$, there is a winning quasi-strategy for player
$\TWO$ in $G_{(\alpha,n)}(x)$.
\end{definition}

For completeness we would like to define $\Pan$ iterability in
the case $\alpha=1$. We will say that $\cM$ is $\Pa[1]{n}$ iterable
above $\delta$ iff $\cM$ is $\Pi_{n-1}$ iterable above $\delta$ in the
sense of \cite{MaSt}.

%


\skipbig

\section{Petite Premice}

\label{section:petitemice}

Recall that our goal in this paper is to prove that $\Aan$ is a mouse set.
In the previous section we defined the notion ``$\cM$ is
$\Pan$-iterable.''
At this point, we would have liked to prove a certain
comparison theorem:

\begin{conjecture}
Suppose that $\cM$ and $\cN$ are countable premice. Suppose that
$\cN$ is $\Pan$-iterable and $(\alpha,n+1)$-small, and
$\cM$ is fully iterable and $(\alpha,n+1)$-small. Then $\cM$ and $\cN$
can be compared.
\end{conjecture}

Suppose that $\cM$ is a fully iterable, countable premouse, and that
$\cM$ is $(\alpha,n+1)$-big, but every proper initial segment of $\cM$
is $(\alpha,n+1)$-small.
If true, the above conjecture would imply that 
 $\R\intersect\cM\subseteq\Aa{n+1}$.
The results of Section \ref{section:correctness} imply that
$\Aa{n+1}\subseteq\R\intersect\cM$. So we would have that
$\Aa{n+1}=\R\intersect\cM$, and so $\Aa{n+1}$ is a mouse set.

Unfortunately we are not able to prove this conjecture---at least, not
for every $(\alpha,n)$. So we are
forced to retreat. First, we will define the notion ``$\cM$ is 
$(\alpha,n)$-petite.'' Then, in the next section, we will prove
a version of the 
conjecture above, with ``petite'' in place of ``small''.

Our goal in this section then is to define the notion 
``$\cM$ is $(\alpha,n)$-petite.'' This will only be a new notion
for
pairs $(\alpha,n)$ such that $\cof(\alpha)\leq\omega$ and $n\geq1$.
For other pairs,  we will define ``$(\alpha,n)$-petite'' to be
synonymous with ``$(\alpha,n)$-small''. But for pairs $(\alpha,n)$ such
that  $\cof(\alpha)\leq\omega$ and $n\geq1$, the condition of being 
$(\alpha,n)$-petite
will be a little more restrictive than the condition of being 
$(\alpha,n)$-small.
In other words, if $\cM$ is $(\alpha,n)$-petite then $\cM$ is 
$(\alpha,n)$-small, but the converse does not  hold. We believe
that the notion of $(\alpha,n)$-small is the ``natural'' notion,
and that our notion of $(\alpha,n)$-petite is an ad-hoc notion with which
we are unfortunately forced to deal.

Recall that a mouse $\cM$ is $(\alpha,n)$-small iff $\cM$ is not
$(\alpha,n)$-big. To parallel this situation, we will first define the
notion ``weakly $(\alpha,n)$-big'', and then we will say that
$\cM$ is $(\alpha,n)$-petite iff $\cM$ is not weakly
$(\alpha,n)$-big. Furthermore, recall that a mouse $\cM$ is
$(\alpha,n)$-big iff $\cM$ is explicitly $(\alphaprime,\nprime)$-big
for some $(\alphaprime,\nprime)$ with 
$(\alpha,n)\lexleq(\alphaprime,\nprime)$. To parallel this situation,
we will first define the notion ``weakly, explicitly $(\alpha,n)$-big'',
and then we will say that $\cM$ is weakly
$(\alpha,n)$-big iff $\cM$ is weakly, explicitly $(\alphaprime,\nprime)$-big
for some $(\alphaprime,\nprime)$ with 
$(\alpha,n)\lexleq(\alphaprime,\nprime)$. The name 
``weakly, explicitly $(\alpha,n)$-big'' is, of course, very awkward.
We choose to tolerate this awkwardness, among other reasons, because we
believe that the notion itself is ad-hoc and should be dispensable.

If $\cof(\alpha)>\omega$, or if $\alpha=1$,
or if $\cof(\alpha)\leq\omega$ and $n=0$,
then $\cM$ is \emph{weakly, explicitly $(\alpha,n)$-big above $\beta$}
iff $\cM$ is (strongly)
explicitly $(\alpha,n)$-big above $\beta$.  So we only
need give our definition in the case $\cof(\alpha)\leq\omega$ and $n\geq 1$.
We start with the case $n=1$. Compare the following
definition with Definition \ref{cofomeganisone} on page
\pageref{cofomeganisone}. The main difference between the two definitions
is that in the definition below we do not demand that the $\delta$'s
remain cardinals up to the next admissible.

\begin{definition}
\label{def:petitezero}
Let $\cM$ be a premouse and $\beta\in\ORD^{\cM}$. Let 
$2\leq\alpha<\omega_1^{\omega_1}$, and
suppose that $\alpha$ is a successor ordinal, or  a limit 
ordinal of cofinality $\omega$. Then $\cM$ is 
\emph{weakly, explicitly $(\alpha,1)$-big above $\beta$} iff there is an
initial segment $\cN\unlhd\cM$ with $\beta\in\cN$,
such that for cofinally many
$\alphaprime<\alpha$, and all $n\in\omega$, and cofinally many
$\delta<\ORD\intersect\cN$, we have that:
\begin{itemize}
\item[(a)] $\delta$ is a a cardinal of $\cN$, and
\item[(b)] $\delta$ is locally $(\alphaprime,n)$-Woodin in $\cN$.
\end{itemize}
\end{definition}

In Figure \ref{fig:bigpremouse} on page \pageref{fig:bigpremouse},
and in the remark immediately following that figure, we explored the
definition of (strongly) explicitly $(\alpha,1)$-big in the case that
$\alpha$ is a successor ordinal. Similarly,
in Figure \ref{fig:weaklybigpremouse}
and Remark \ref{rem:weaklytwoonebig} below, we explore the definition
of weakly, explicitly $(\alpha,1)$-big in the case that $\alpha$ is
a successor ordinal.

\begin{figure}[hbt]
\begin{center}
\input{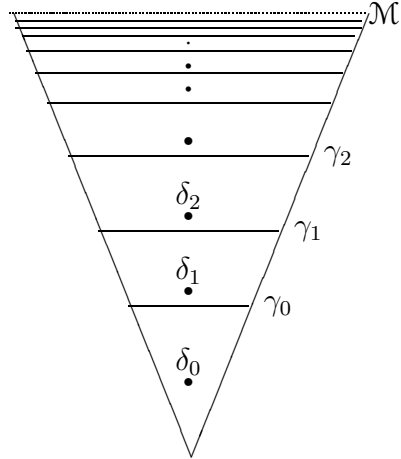}
\caption{$\alpha\geq2$ is a successor and $\cM$ is weakly, explicitly 
$(\alpha,1)$-big}
\label{fig:weaklybigpremouse}
\end{center}
\end{figure}

\begin{remark}
\label{rem:weaklytwoonebig}
Suppose $\alpha=\alphaprime+1$ with $\alphaprime\geq1$.
Suppose that $\cM$ is a countable premouse and 
$\cM$ is weakly, explicitly $(\alpha,1)$-big (above $0$)
 but no proper initial
segment of $\cM$ is. We will examine what $\cM$ looks like.
(See Figure \ref{fig:weaklybigpremouse}.)
There is a sequence of ordinals 
$\delta_0<\gamma_0<\delta_1<\gamma_1<\delta_2<\dots$ such that
 each
$\delta_n$ is a cardinal of $\cM$ 
and a Woodin cardinal of $\cJ^{\cM}_{\gamma_n}$, and
$\cJ^{\cM}_{\gamma_n}$ is explicitly $(\alphaprime,n)$-big above $\delta_n$.
Since no proper initial segment of $\cM$ is weakly, explicitly
$(\alpha,1)$-big,
we have that $\sup\singleton{\delta_n}=\sup\singleton{\gamma_n}=\ORD^{\cM}$.

Let $\rho=\ORD^{\cM}$, and let $\cM^*$ be the first admissible structure
over $\cM$. It is tempting to think that $(\cM^*,\rho)$ witnesses that
$\cM^*$ is (strongly) explicitly $(\alpha,1)$-big.
Note that in general \emph{this is not true.} The reason is that
in general $\rho$ will not be a cardinal in $\cM^*$.
\end{remark}

We now complete the inductive definition of weakly, explicitly 
$(\alpha,n)$-big.

\begin{definition}
Let $\cM$ be a premouse and $\beta\in\ORD^{\cM}$. Let 
$2\leq\alpha<\omega_1^{\omega_1}$, and
suppose that $\alpha$ is a successor ordinal, or  a limit 
ordinal of cofinality $\omega$. Let $n\geq 2$. Then $\cM$ is
\emph{weakly, explicitly $(\alpha,n)$-big above $\beta$} iff there is
an initial segment $\cN\unlhd\cM$, and an ordinal $\delta\in\cN$ with
$\beta<\delta$, such that $\delta$ is a Woodin cardinal of $\cN$, and
$\cN$ is weakly, explicitly $(\alpha,n-1)$-big above $\delta$.
\end{definition}

It is easy to see that if $\cM$ is  explicitly $(\alpha,n)$-big
then $\cM$ is weakly, explicitly $(\alpha,n)$-big. Conversely we 
have the following.

\begin{lemma}
\label{WeakPlusOneGivesStrong}
Let $\cM$ be a countable, realizable premouse, and $\beta\in\ORD^{\cM}$.
 Let  $2\leq\alpha<\omega_1^{\omega_1}$, and
suppose that $\alpha$ is a successor ordinal, or  a limit 
ordinal of cofinality $\omega$. Let $n\geq 0$. Suppose that $\cM$ is
weakly, explicitly $(\alpha,n+1)$-big above $\beta$. Furthermore, suppose
that $\max(\code(\alpha))<\beta$.
Then $\cM$
is (strongly) explicitly $(\alpha,n)$-big above $\beta$.
\end{lemma}
\begin{proof}
First consider the case $n=0$.
Suppose that $\cM$ is weakly, explicitly $(\alpha,1)$-big
above $\beta$. Then clearly $\cM$ is explicitly
$(\alpha,0)$-big above $\beta$. This does not use our hypothesis
that $\max(\code(\alpha))<\beta$.

Next consider the case $n=1$. Suppose $\cM$ is weakly, explicitly 
$(\alpha,2)$-big above $\beta$, and we will show that $\cM$ is
explicitly $(\alpha,1)$-big above $\beta$. Let $\cN\unlhd\cM$,
and let $\delta$ be an ordinal of $\cN$ such that $\delta$ is Woodin
in $\cN$ and $\cN$ is weakly, explicitly $(\alpha,1)$-big above $\delta$.
By the previous paragraph, $\cN$ is explicitly $(\alpha,0)$-big above
$\delta$. By our hypothesis, $\max(\code(\alpha))<\delta$.
By Lemma \ref{FixedDeltaIsTooBig} on page 
\pageref{FixedDeltaIsTooBig}, $\cJ^{\cN}_{\delta}$ is explicitly 
$(\alpha,1)$-big above $\beta$.

The result now follows trivially for $n\geq 2$.
\end{proof}

Finally, we get to the definition of $(\alpha,n)$-petite.

\begin{definition}
Let $\cM$ be a premouse and $\beta\in\ORD^{\cM}$. Let 
$2\leq\alpha\leq\omega_1^{\omega_1}$, and let $n\in\omega$.
Then $\cM$ is \emph{weakly $(\alpha,n)$-big above $\beta$}
iff $\cM$ is weakly, explicitly $(\alphaprime,\nprime)$-big above
$\beta$ for some $(\alphaprime,\nprime)$ with 
$(\alpha,n)\lexleq(\alphaprime,\nprime)$. Finally, $\cM$ is
\emph{$(\alpha,n)$-petite above $\beta$} iff $\cM$ is not
weakly $(\alpha,n)$-big above $\beta$.
\end{definition}

In the remainder of this section we will prove a few lemmas concerning
our new notion of $(\alpha,n)$-petite. We mentioned above that this is
only a new notion for certain pairs $(\alpha,n)$. Let us begin by
verifying that this is so:

\begin{lemma}
Suppose that $\cof(\alpha)>\omega$, or $n=0$. Then $(\alpha,n)$-petite
is synonymous with $(\alpha,n)$-small.
\end{lemma}
\begin{proof}
Let $(\alpha,n)$ be as in the statement of the lemma.
If $\cM$ is not $(\alpha,n)$-small above $\beta$, then it follows easily
that $\cM$ is not $(\alpha,n)$-petite above above $\beta$. Conversely,
suppose that $\cM$ is not $(\alpha,n)$-petite above above $\beta$.
Let $(\alphaprime,\nprime)$ be the lexicographically least pair such
that $(\alpha,n)\lexleq(\alphaprime,\nprime)$ and $\cM$ is weakly, 
explicitly $(\alphaprime,\nprime)$-big above $\beta$. If 
$(\alphaprime,\nprime)=(\alpha,n)$ then, because of our conditions
on $(\alpha,n)$, we have that by definition
$\cM$ is \emph{strongly} explicitly $(\alpha,n)$-big above $\beta$,
and so $\cM$ is not $(\alpha,n)$-small above $\beta$.
So suppose that $(\alpha,n)\lexless(\alphaprime,\nprime)$. Then
the proof of Lemma \ref{biggoesdownbyone} shows that $\alphaprime>\alpha$
and $\nprime=0$.
So again, by definition, $\cM$ is \emph{strongly} explicitly 
$(\alphaprime,\nprime)$-big above $\beta$. So $\cM$ is not 
$(\alpha,n)$-small above $\beta$.
\end{proof}

The following lemma is fairly obvious, but the proof is a bit convoluted.

\begin{lemma}
\label{ZeroPetiteIsPetiter}
Let $\cM$ be a countable, realizable premouse, and
$\beta\in\ORD^{\cM}$. Suppose $2\leq\alpha\leq\omega_1^{\omega_1}$. Then
$\cM$ is $(\alpha,0)$-petite above $\beta$ iff $\cM$ is 
$(\alphaprime,n)$-petite above $\beta$ for some pair
$(\alphaprime,n)\lexless(\alpha,0)$.
\end{lemma}
\begin{proof}
If $\cM$ is $(\alphaprime,n)$-petite above $\beta$ for some pair
$(\alphaprime,n)\lexless(\alpha,0)$, then of course
$\cM$ is $(\alpha,0)$-petite above $\beta$. Conversely,
suppose that for all $(\alphaprime,n)\lexless(\alpha,0)$, 
$\cM$ is not $(\alphaprime,n)$-petite above $\beta$. We will show that
$\cM$ is not $(\alpha,0)$-petite above $\beta$. Suppose towards a
contradiction that $\cM$ is $(\alpha,0)$-petite above $\beta$.

Since $\cM$ is
not $(\alphaprime,n)$-petite above $\beta$ for all 
$(\alphaprime,n)\lexless(\alpha,0)$, we have that either:\\
(i) $\cM$ is weakly, explicitly $(\alphaprime,n)$-big above
$\beta$ for cofinally many $(\alphaprime,n)\lexless(\alpha,0)$,\\
or\\
(ii) $\cM$ is weakly, explicitly $(\alphaprime,n)$-big above
$\beta$ for some $(\alphaprime,n)$ with
$(\alpha,0)\lexleq(\alphaprime,n)$.\\
In case (ii), we have that $\cM$ is not $(\alpha,0)$-petite above $\beta$.
So suppose that case (i) holds. First notice that case (i) can only hold
if $\cof(\alpha)\leq\omega$. This is because $\cM$ can only be weakly,
explicitly $(\alphaprime,n)$-big above $\beta$ for countably many
pairs $(\alphaprime,n)$. Thus we may assume that $\cof(\alpha)\leq\omega$.
We will show that $\cM$ is weakly, explicitly $(\alpha,0)$-big
above $\beta$. Recall that this is synonymous with strongly, explicitly
$(\alpha,0)$-big above $\beta$. By definition, we must prove the
following:
 
\begin{claim}
For cofinally many $\alphaprime<\alpha$, and 
\emph{all} $n$, there is a \emph{proper} initial segment of $\cM$  
which is is \emph{strongly}, explicitly $(\alphaprime,n)$-big
above $\beta$.
\end{claim}
\begin{subproof}[Proof of Claim]
First suppose that $\alpha$ is a limit ordinal of cofinality $\omega$.
From our assumption that (i) above holds, and using
 the proof of Lemma \ref{biggoesdownbyone}, we have that for cofinally
many $\alphaprime<\alpha$, $\cM$ is weakly explicitly $(\alphaprime,0)$-big
above $\beta$.
Let $\alpha_0<\alpha$ and let $n_0\in\omega$. We will find an 
$\alphaprime$ with $\alpha_0\leq\alphaprime<\alpha$, and 
an $\cN\properseg\cM$ such that
$\cN$ is strongly, explicitly  $(\alphaprime,n_0)$-big above $\beta$.

Let $\alpha_1$ be least such that $\alpha_0<\alpha_1<\alpha$ 
and $\cM$ is weakly, explicitly  $(\alpha_1,0)$-big above $\beta$. 
By definition, $\cM$ is \emph{strongly} explicitly
$(\alpha_1,0)$-big above $\beta$.
If $\alpha_1=\alpha_0+1$, then  by definition
there is a proper initial segment
of $\cM$ which is  explicitly $(\alpha_0,n_0)$-big above $\beta$.
So we are done if we set $\alphaprime=\alpha_0$.
So suppose that
$\alpha_1>\alpha_0+1$. Then it can not be true that $\alpha_1$ is a 
successor ordinal or a limit ordinal of cofinality $\omega$, so it must
be the case that $\alpha_1$ is a limit ordinal of
cofinality $\omega_1$. Let $\cN_1\initseg\cM$ be
 such that $\cN_1$ witnesses that $\cM$ is  explicitly 
$(\alpha_1,0)$-big above $\beta$. By definition, $\cN_1$ is active.
Let $\kappa_1$ be the critical point of the last extender from the
$\cN_1$-sequence. In review we have:
\begin{itemize}
\item[(a)]$\alpha_0<\alpha_1<\alpha$.
\item[(b)]$\cof(\alpha_1)=\omega_1$.
\item[(c)]$\cN_1\initseg\cM$ witnesses that $\cM$ is explicitly
$(\alpha_1,0)$-big above $\beta$.
\item[(d)] $\kappa_1$ is the critical point of the last extender from
the $\cN$-sequence.
\end{itemize}
By replacing $\alpha_1$ and $\cN_1$ if necessary, we may also assume
that $\kappa_1$ is least possible so that (a) through (d) above holds.

 Now let 
 $\alpha_2$ be least such that $\alpha_1<\alpha_2<\alpha$ 
and $\cM$ is weakly, explicitly  $(\alpha_2,0)$-big above $\beta$.
By definition, $\cM$ is \emph{strongly} explicitly  
$(\alpha_2,0)$-big above $\beta$.
If $\alpha_2=\alpha_1+1$, then  by definition
there is a proper initial segment
of $\cM$ which is  explicitly $(\alpha_1,n_0)$-big above $\beta$.
So we are done if we set $\alphaprime=\alpha_1$.
So suppose that
$\alpha_2>\alpha_1+1$. Then it can not be true that $\alpha_2$ is a 
successor ordinal or a limit ordinal of cofinality $\omega$, so it must
be the case that $\alpha_2$ is a limit ordinal of
cofinality $\omega_1$. Let $\cN_2\initseg\cM$ be
 such that $\cN_2$ witnesses that $\cM$ is  explicitly
$(\alpha_2,0)$-big above $\beta$. By definition, $\cN_2$ is active.
Let $\kappa_2$ be the critical point of the last extender from the
$\cN_2$-sequence. By our minimality hypothesis on $\kappa_1$, we
have that $\kappa_1\leq\kappa_2$. By definition, 
$\max(\code(\alpha_1))<\kappa_1$. Thus $\max(\code(\alpha_1))<\kappa_2$.
Then, since $\alpha_1<\alpha_2$ and $\max(\code(\alpha_1))<\kappa_2$,
there is a proper initial segment $\cN\properseg J^{\cM}_{\kappa_2}$ such
that $\cN$ is explicitly $(\alpha_1,n_0)$-big above $\beta$. Thus we
have accomplished our task, with $\alphaprime=\alpha_1$.

Now we must prove our claim in the case that $\alpha$ is a successor
ordinal. Suppose $\alpha=\alpha_0+1$. From our assumption that (i) above 
holds, and using the proof of Lemma \ref{biggoesdownbyone}, we have that 
$\cM$ is weakly, explicitly $(\alpha_0,n)$-big above $\beta$ for
all $n\in\omega$. Fix $n_0\in\omega$. We must show that there is
a proper initial segment $\cN\properseg\cM$ such that $\cN$ is
\emph{strongly}, explicitly $(\alpha_0,n_0)$-big above $\beta$.
Let $\cN_0\initseg\cM$ be the $\initseg$-least initial segment of
$\cM$ which is weakly, explicitly $(\alpha_0,0)$-big above $\beta$.
By Proposition \ref{OmegaWoodinsIsBig}, if $\cN_0$ had 
$\omega$ Woodin cardinals above $\beta$, 
then $\cN_0$ would be $(\omega_1^{\omega_1},0)$-big above $\beta$.
Since we are assuming that $\cM$ is $(\alpha,0)$-petite above $\beta$,
$\cN_0$ does not have $\omega$ Woodin cardinals above $\beta$.
Let $m$ be such that $0\leq m <\omega$ and $\cN_0$ has 
$m$ Woodin cardinals, but not $m+1$ Woodin cardinals above $\beta$.
Let $n=m+n_0+10$, and let $\cN_1\initseg\cM$ be the $\initseg$-least 
initial
segment of $\cM$ which is weakly, explicitly $(\alpha_0,n)$-big above
$\beta$. Since $\cN_1$ is also weakly, explicitly $(\alpha_0,0)$-big
above $\beta$ we have that $\cN_0\initseg\cN_1$. Since $\cN_0$ does
not have $n$ Woodin cardinals, we have that $\cN_0$ is a proper initial
segment of $\cN_1$. Since $n>m+n_0+1$, $\cN_1$ is weakly explicitly
$(\alpha_0,n_0+1)$-big above $\ORD^{\cN_0}$.
By  Lemma \ref{bigcontainscode}, $\max(\code(\alpha_0))\leq\ORD^{\cN_0}$.
Then, by the proof of Lemma \ref{WeakPlusOneGivesStrong}, 
there is a proper initial segment
of $\cN_1$ which is \emph{strongly} explicitly $(\alpha_0,n_0)$-big
above $\ORD^{\cN_0}$, and so above $\beta$.
\end{subproof}

This completes the proof of the lemma.
\end{proof}

The following lemma is also fairly obvious. The proof just requires
untangling definitions. We will use this lemma in the next section.

\begin{lemma}
\label{PetiteAboveWoodin}
Let $\cM$ be a countable, realizable premouse. Let $\beta<\delta$
be ordinals in $\cM$, and suppose that $\delta$ is a Woodin 
cardinal of $\cM$. Let $\alpha\geq 2$.
\begin{itemize}
\item[(a)] Let $n\geq1$ and suppose that $\cM$ is $(\alpha,n+1)$-petite
above $\beta$. Then $\cM$ is $(\alpha,n)$-petite above $\delta$.
\item[(b)] Suppose that $\cM$ is $(\alpha,1)$-petite above
$\beta$. Then $\cM$ is $(\alphaprime,m)$-petite above $\delta$ for
some $(\alphaprime,m)$ with $1\leq\alphaprime<\alpha$.
\end{itemize}
\end{lemma}
\begin{proof}

(a) Suppose $\cM$ is not $(\alpha,n)$ petite above $\delta$.
Since $\cM$ is $(\alpha,n+1)$-petite above $\delta$, we
must have that $\cM$ is weakly, explicitly $(\alpha,n)$-big
above $\delta$. As $\cM\models\delta$ is Woodin,
 $\cM$ is weakly,  explicitly $(\alpha,n+1)$-big above $\beta$. 
This is a contradiction.

(b) Suppose that for all $\alphaprime<\alpha$ and all
$m\in\omega$, $\cM$ is not $(\alphaprime,m)$-petite above
$\delta$. By the previous lemma, $\cM$ is not $(\alpha,0)$-petite
above $\delta$. Since $\cM$ is $(\alpha,1)$-petite above $\delta$,
we must have that $\cM$ is weakly, explicitly $(\alpha,0)$-big above
$\delta$.  Recall that  weakly, explicitly $(\alpha,0)$-big
is synonymous with (strongly) explicitly $(\alpha,0)$-big.

First consider the case that $\cof(\alpha)>\omega$. Since
$\cM$ is  explicitly $(\alpha,0)$-big
above $\delta$ and
$\cM\models\delta$ is Woodin, $\cM$ is explicitly 
$(\alpha,1)$-big above $\beta$. This is a contradiction.

Next consider the case that $\cof(\alpha)\leq\omega$.

\begin{claim}
$\max(\code(\alpha))<\delta$.
\end{claim}
\begin{subproof}[Proof of Claim]
Let $k$ be the least measurable cardinal of $\cM$ which is greater
than $\beta$. So we have that $\beta<\kappa<\delta$. Since
$\cM$ is explicitly $(\alpha,0)$-big above $\beta$ but
$\cM$ is not $(\alpha^*,0)$-big above $\beta$ for any $\alpha^*>\alpha$,
Lemma \ref{BignessIsBelowKappa} on page 
\pageref{BignessIsBelowKappa} implies that $\max(\code(\alpha))<\kappa$.
\end{subproof}
We may now apply Lemma \ref{FixedDeltaIsTooBig}
on page \pageref{FixedDeltaIsTooBig} to conclude that
$\cJ^{\cM}_{\delta}$ is explicitly
$(\alpha,1)$-big above $\beta$. This is a contradiction.
\end{proof}

In Section \ref{section:correctness}, 
it was essential to our proof that an iterate of a big
premouse is big. In fact this requirement was one of the central motivating
factors behind our definition of $(\alpha,n)$-big. Similarly, in
the next section it will be essential to our proof that an iterate of
a petite premouse is petite. That is the content of the next lemma.

\begin{lemma}
\label{PreservationOfSmallness}
Let $\cM$ be a premouse, and $\beta\in\ORD^{\cM}$. Let $\cT$ be an
iteration tree on $\cM$, with $\cT$ above $\beta$. Suppose that $\cM$
is $(\alpha,n)$-petite above $\beta$. Then $\forall\xi<\length(\cT)$,
$\cM^{\cT}_{\xi}$ is $(\alpha,n)$-petite above $\beta$.
\end{lemma}
\begin{proof}
By induction on $\xi$.
For $\xi<\length(\cT)$ we will write $\cM_{\xi}$ instead
of $\cM^{\cT}_{\xi}$, and $E_{\xi}$ instead of $E^{\cT}_{\xi}$.

First suppose $\xi=\eta+1$ is a successor ordinal. Let
$\gamma=T\pred(\xi)$. Let $k=\deg(\xi)$. Then there is an
initial segment $\cM^*_{\xi}\unlhd\cM_{\gamma}$ such that
$\cM_{\xi}=\Ult_k(\cM^*_{\xi}, E_{\eta})$.
Let $j:\cM^*_{\xi}\map\cM_{\xi}$ be the canonical $k$-embedding. Let
$\kappa=\crit(j)=\crit(E_{\eta})$. 

By induction, $\cM_{\gamma}$ is $(\alpha,n)$-petite above $\beta$,
and so in particular, $\cM^*_{\xi}$ is $(\alpha,n)$-petite
above $\beta$.
Let $(\alpha_0,n_0)$ be the  lexicographically greatest pair such that
$\cM^*_{\xi}$ is weakly, explicitly $(\alpha_0,n_0)$-big above $\beta$.
As in the proof of Theorem \ref{BigMiceAreClosed} on 
page \pageref{BigMiceAreClosed}, we can see that there is such a pair
$(\alpha_0,n_0)$. We have that $(\alpha_0,n_0)\lexless(\alpha,n)$,
and that $\cM^*_{\xi}$ is $(\alpha_0,n_0+1)$-petite above $\beta$.
We will show that $\cM_{\xi}$ is $(\alpha_0,n_0+1)$-petite above $\beta$.

\begin{claim}[Claim 1]
There is an ordinal $\mu\leq\kappa$ such that
$\mu$ is the critical point of a total-on-$\cM^*_{\xi}$ extender
from the $\cM^*_{\xi}$ sequence.
\end{claim}
\begin{subproof}[Proof of Claim 1]
  If $\gamma=\eta$, then $E_{\eta}$
is on the $\cM^*_{\xi}$ sequence, and so we may take $\mu=\kappa$.
So suppose that $\gamma<\eta$. Then we have that 
$\kappa<\length(E_{\gamma})$, and $\length(E_{\gamma})$
is a cardinal of $\cM_{\eta}$, and $\cM^*_{\xi}$ and $\cM_{\eta}$
agree below $\length(E_{\gamma})$.
If $\crit(E_{\gamma})\leq\kappa<\length(E_{\gamma})$, then
it follows that $E_{\gamma}$ is total-on-$\cM^*_{\xi}$ and so
we may take $\mu=\crit(E^{\cT}_{\gamma})$. 
So suppose
$\kappa<\crit{E_{\gamma}}$.
Then we have that
$(\kappa^{++})^{\cM_{\eta}} < \crit(E_{\gamma})$. 
By the initial segment condition on $E_{\eta}$, there is an
initial segment of $E_{\eta}$ on the $\cM_{\eta}$ sequence
below $\crit(E_{\gamma})$,
and so on the $M^*_{\xi}$ sequence. So we may take $\mu=\kappa$.
\end{subproof}

\begin{claim}[Claim 2]
$\max(\code(\alpha_0))<\kappa$.
\end{claim}
\begin{subproof}[Proof of Claim 2]
Let $\mu<\kappa$ be as in Claim 1.
By the proof of Lemma \ref{BignessIsBelowKappa} on page 
\pageref{BignessIsBelowKappa}, since $\cM^*_{\xi}$ is not
explicitly $(\alpha^*,0)$-big for any $\alpha^*>\alpha_0$,
$\max(\code(\alpha_0))<\mu$.
\end{subproof}

Now, suppose towards a contradiction that $\cM_{\xi}$ were not
$(\alpha_0,n_0+1)$-petite above $\beta$. Then there is some
$(\alphaprime,\nprime)$ with $(\alpha_0,n_0)\lexless(\alphaprime,\nprime)$,
such that $\cM_{\xi}$ is weakly, explicitly $(\alphaprime,\nprime)$-big
above $\beta$. Let $\sprime=\code(\alphaprime)$. Let $s_0=\code(\alpha_0)$.
Then $\cM_{\xi}$ satisfies the statement:
\begin{quote}
There exists an $\sprime$ such that 
$(|s_0|,n_0)\lexless(|\sprime|,\nprime)$, and I am weakly, explicitly
$(|\sprime|,n_0)$-big above $\beta$.
\end{quote}
Since $\beta<\crit(j)$, and $\max(s_0)<\crit(j)$, and $j$ is either 
$\Sigma_2$-elementary, or $\Sigma_1$-elementary and cofinal, it follows
that $\cM^*_{\xi}$ satisfies the same statement.
Since this is a contradiction, we have that $\cM_{\xi}$ is
$(\alpha_0,n_0+1)$-petite above $\beta$.
(There is one subtle,
technical point here. If we were using the notion of ``$(\alpha,n)$-small''
instead of the notion of ``$(\alpha,n)$-petite'' in this lemma, then at
this point, in order to take care of one particular case,
 we would need to show that  if $\cM_{\xi}$ were admissible, then
$\cM^*_{\xi}$ would be admissible too. The contrapositive of this would
require us to show that \emph{non-admissibility} is preserved by
ultrapowers.  Recall that earlier in this paper we were required to show
that admissibility was preserved by ultrapowers. We do in fact have a proof
that non-admissibility is preserved by ultrapowers. However, since we
are using the notion of ``$(\alpha,n)$-petite'' instead of the notion
of ``$(\alpha,n)$-small'' in this lemma, we will not need this proof here.)

Finally, we must consider the case that $\xi$ is a limit ordinal.
Let $\gamma<_{T}\xi$ be such that $i^{\cT}_{\gamma,\xi}$ is defined.
Let $j=i^{\cT}_{\gamma,\xi}$ and let $\kappa=\crit(j)$.
We may assume that $\gamma$ is greatest so that $\gamma<_{T}\xi$
and $\kappa=\crit(i^{\cT}_{\gamma,\xi})$.
As in Claim 1 we have that there is an ordinal $\mu\leq\kappa$ such that
$\mu$ is the critical point of a total-on-$\cM_{\gamma}$ extender
from the $\cM_{\gamma}$ sequence. Let $(\alpha_0,n_0)$ be the  
lexicographically greatest pair such that
$\cM_{\gamma}$ is weakly, explicitly $(\alpha_0,n_0)$-big above $\beta$.
As in Claim 2 we have that $\max(\code(\alpha_0))<\kappa$. The rest
of the proof is the same as the successor case.
\end{proof}

%


\skipbig

\section{A Comparison Theorem}

\label{section:comparison}

In this section we prove the following comparison theorem:
 If $\cM$ and $\cN$ are $\omega$-mice,
and $\cN$ is $\Pan$-iterable, and $\cM$ is fully iterable and
$(\alpha,n+1)$-petite, then $\cN\unlhd\cM$ or $\cM\unlhd\cN$.
 (See part (3) Theorem \ref{ComparisonTheorem}  below.)
Our theorem is an extension of Lemma 2.2 from \cite{MaSt}.
We employ the same architecture here as in the proof of that lemma.
Consequently, we will be borrowing many definitions and ideas
from the paper \cite{MaSt}. We will attempt to make this section
self-contained by giving definitions (at least informally) for all
of the terms we use. Having a copy of the paper \cite{MaSt} will be
useful to the reader, but should not be absolutely necessary.
It may also prove useful to the reader to have some familiarity
with the paper \cite{St2}.
We work in the theory ZFC.

We begin with a short discussion about our method of coding hereditarily
countable sets by reals. Recall that if $x\in\R$ then $x^*$ is the
hereditarily countable set coded by $x$. Suppose that $\cT$ is a countable
iteration tree on a countable premouse, and $x\in\R$ and $x^*=\cT$.
Suppose that $b$ is a cofinal, wellfounded branch of $\cT$, and that
$b$ is the unique cofinal branch of $\cT$ with some given property.
Thus $\cM^{\cT}_{b}$ is definable from $\cT$. Consider the following
question. Is there a real $y$ such that $y^*=\cM^{\cT}_b$ and
$y$ is definable from $x$? The answer is yes. This is because, given
$\cT$ and $\cM^{\cT}_b$, the coding
of $\cT$ by $x$ induces, in a simple way, a coding of $\cM^{\cT}_b$.
Below we make this idea precise.

\begin{definition}
Let $P\subset\HC\times\HC$ be a binary relation on $\HC$. Let 
$f:P\times\R\map\R$.
Then we say that {\em $f$ is a code transformation for $P$} iff
for all $\angles{a,b,x}\in\dom(f)$, if $x^*=a$ then $(f(a,b,x))^*=b$.
\end{definition}
\begin{definition}
\label{DefinableCodeTransformation}
Let $P\subset\HC\times\HC$ be a binary relation on $\HC$. Suppose that 
$P$ is definable
over $\HC$.  Then we say that {\em $P$ admits a definable code 
transformation} iff
there is a code transformation $f$ for $P$ such that $f$ is definable 
over $\HC$.
\end{definition}

For example, suppose that $P(a,b)\Iff b\subseteq a$. Then $P$ admits a 
simple definable
code transformation, namely,  given a wellfounded tree coding $a$, give the 
subtree which
codes $b$.

More to the point here is the following proposition.

\begin{proposition}
\label{DefCodTransProp}
Let $P$ be the binary relation on $\HC$ defined by
 $P(\cT,\angles{b,\cM})\Iff$ $\cT$ is an iteration
tree and $b$ is a cofinal wellfounded branch of $\cT$ and
$\cM=\cM^{\cT}_b$. Then  $P$ admits
a definable code transformation. 
\end{proposition}
\begin{proof}
Suppose that $P(\cT,\angles{b,\cM})$ holds, and 
$x\in\R$ and $x^*=\cT$. We must 
define a real $y$ such that $y^*=\angles{b,\cM}$. 
Our definition is allowed to employ $x$, $\cT$, $\cM$, and $b$ as 
parameters. We will define two reals, $y_1$ and $y_2$, such that
$y_1^*=b$ and $y_2^*=\cM$. From $y_1$ and $y_2$ it is easy to define $y$.

To begin with, there is some ambiguity
about exactly which hereditarily countable set $\cT$ is. For concreteness,
let us define an iteration tree $\cT$ to be an ordered $8$-tuple.
That is $\cT=\angles{\alpha,T,\vec{E},\vec{M},\vec{M^*},\vec{i},D,\deg}$ 
where 
$\alpha=\length(\cT)$, $T$ is the binary relation on $\alpha$ which is
the tree ordering associated with $\cT$, $D\subseteq\alpha$ is the
drop set associated with $\cT$, $\deg:\alpha\map\omega+1$ 
is the degree function associated
with $\cT$, $\vec{E}$ is the sequence of extenders
used on $\cT$, $\vec{M}$ is the sequence of models on  $\cT$,
$\vec{M^*}$ is the sequence of models which are dropped to, and $\vec{i}$
is the sequence of embeddings associated with $\cT$. (This definition of
an iteration tree differs slightly from the one in \cite{MiSt}. In that
paper an iteration tree is defined as an ordered $5$-tuple. The extra three
components of $\cT$ which we have here, $\alpha$, $\vec{\cM}$, and 
$\vec{i}$, are definable from the other 5 components. Thus
they are redundant. We have included them here to simplify our discussion
of definable code transformations.)

Now, since $b\subseteq\alpha$ and $x$ gives us a coding of $\alpha$, it is
easy to see how to define $y_1$. So we concentrate on $y_2$.
To define a
real $y_2$ such that $y_2^*=\cM$, it suffices to define an enumeration
$f:\omega\bijection\cM$. Since the transitive collapse map is
definable over $\HC$, it suffices to define an enumeration of
 the pre-collapsed version of the direct limit model $\cM^{\cT}_b$. 
Thus it suffices to define an enumeration of the set 
$\Union{\gamma\in b}\singleton{\gamma}\times\cM^{\cT}_{\gamma}$.
It is easy to see that we can define such an enumeration over $\HC$,
given $x$, $\cT$,
and $b$ as parameters.
\end{proof}

We will need to borrow several notions from the paper \cite{MaSt}.
The first such notion is that of
a \emph{coiteration} of two premice. If $\cM$ and $\cN$ are premice,
then a coiteration of $\cM$ and $\cN$ is the sequence of iteration trees
that results from comparing $\cM$ with $\cN$ for some number of steps. 
The precise definition below
is taken word-for-word from \cite{MaSt}.

\begin{definition}
Let $\cM$ and $\cN$ be premice. A \emph{coiteration} of $\cM$ and $\cN$ is
a sequence $\angles{(\cT_{\alpha},\cU_{\alpha}) \mid \alpha<\theta}$ such
that 
\begin{itemize}
\item[(1)] $\cT_{\alpha}$ and $\cU_{\alpha}$ are $\omega$-maximal iteration
trees on $\cM$ and $\cN$ respectively.
\item[(2)] $\alpha<\beta$ $\implies$ $\cT_{\beta}$ extends $\cT_{\alpha}$
and $\cU_{\beta}$ extends $\cU_{\alpha}$.
\item[(3)] $\alpha$ limit $\implies$ $\Bigl(\cT_{\alpha}=
\Union{\beta<\alpha}\cT_{\beta} \text{ and } \cU_{\alpha}=
\Union{\beta<\alpha}\cU_{\beta} \Bigr)$.
\item[(4)] $\length(\cT_{\alpha})$ a limit $\implies$
$\length(\cT_{\alpha+1})=\length(\cT_{\alpha})+1$ and
$\length(\cU_{\alpha})$ a limit $\implies$
$\length(\cU_{\alpha+1})=\length(\cU_{\alpha})+1$.
\item[(5)] $\Bigl(\length(\cT_{\alpha}) \text{ a limit and }
\length(\cU_{\alpha}) \text{ a successor } \Bigr) \implies
\length(\cU_{\alpha+1}) = \length(\cU_{\alpha})$, and \\
$\Bigl(\length(\cU_{\alpha}) \text{ a limit and }
\length(\cT_{\alpha}) \text{ a successor } \Bigr) \implies
\length(\cT_{\alpha+1}) = \length(\cT_{\alpha})$.
\item[(6)] If $\length(\cT_{\alpha})$ and $\length(\cU_{\alpha})$ are
both successor ordinals, and $\alpha+1<\theta$, then $\cT_{\alpha+1}$
and $\cU_{\alpha+1}$ are determined by ``iterating the least disagreement''
between the last models in $\cT_{\alpha}$ and $\cU_{\alpha}$, and by the
rules for $\omega$-maximal iteration trees. (See \S7 of \cite{MiSt}.)
\end{itemize}
\end{definition}

A coiteration of $\cM$ and $\cN$ is determined by the cofinal wellfounded
branches of $\cT_{\alpha}$ or $\cU_{\alpha}$ used to produce
$\cT_{\alpha+1}$ or $\cU_{\alpha+1}$ in the case $\length(\cT_{\alpha})$
or $\length(\cU_{\alpha})$ is a limit. This means that
$\Union{\alpha<\theta}\cT_{\alpha}$ and $\Union{\alpha<\theta}\cU_{\alpha}$
determine $\angles{(\cT_{\alpha},\cU_{\alpha}) \mid \alpha<\theta}$,
so we can identify the two, and speak of an appropriate pair $(\cT,\cU)$
of iteration trees on $\cM$ and $\cN$ respectively as a coiteration.

\begin{definition}
Let $(\cT,\cU)$ be a coiteration of $\cM$ and $\cN$. We say that
$(\cT,\cU)$ is \emph{terminal} iff either
\begin{itemize}
\item[(1)] $\length(\cT)$ is a limit and $\cT$ has no cofinal wellfounded
branch, or $\length(\cU)$ is a limit and $\cU$ has no cofinal wellfounded
branch, or
\item[(2)] $\cT$ and $\cU$ have last models $\cP$ and $\cQ$,
$\cP\notinitseg\cQ$ and $\cQ\notinitseg\cP$, and one of the 
ultrapowers 
determined by iterating the least disagreement between $\cP$ and $\cQ$,
according to the rules for $\omega$-maximal trees, is illfounded, or
\item[(3)] $\cT$ and $\cU$ have last models $\cP$ and $\cQ$, and
$\cP\unlhd\cQ$ or $\cQ\unlhd\cP$.
\end{itemize}
\end{definition}
We shall call a coiteration \emph{successful} just in case it is terminal
and case (3) holds above.

\skipmed

The main theorems in this section and the next
section include a large cardinal hypothesis---we
will assume that there exists infinitely many Woodin cardinals.
We will use this large cardinal hypothesis in two ways. Firstly, the
hypothesis will give us the determinacy of sufficiently many games
in $\LofR$. This determinacy is needed in order to be able to apply
some of the results from earlier sections of this paper. 
The other way in which
the large cardinal hypothesis will be useful to us is that it will
give us a certain amount of \emph{generic absoluteness}. We need this
generic absoluteness in order to show that our coiterations terminate
in countably many steps. (See Lemma \ref{CoiterationsAreCountable} below.)

Both applications of our large cardinal hypothesis---determinacy and
generic absoluteness---are proved using the machinery of Woodin's
\emph{stationary tower forcing.}
  See Chapter 9 of \cite{Ma2} for a thorough
treatment of this machinery. Below we state two facts from the theory
of stationary tower forcing. Then we will use these facts to
derive the specific results we shall need.

\begin{proposition}[Woodin]
\label{firstprop}
Let $\angles{\delta_n \mid n\in\omega}$ be a strictly increasing sequence
of Woodin cardinals. Let $\lambda = \sup_{n} \kappa_n$. Let 
$\gamma > \kappa$ be any ordinal. Then there is a generic extension
of the universe, $V[G]$, such that in $V[G]$ there is an elementary
embedding $j:V\map N$ such that
\begin{itemize}
\item[(i)] $\gamma$ is in the wellfounded part of $\cN$,
\item[(ii)] $\crit(j)=\omega_1$, and $j(\omega_1)=\lambda$, and
\item[(iii)] letting $\R^*$ be the reals of $N$, we have that $\R^*$
is also the set of reals in a symmetric collapse of $V$ up to $\lambda$.
\end{itemize}
\end{proposition}

See section 9.5 of \cite{Ma2} for a proof of this proposition. 

\begin{proposition}[Woodin]
\label{secondprop}
Let $\angles{\delta_n \mid n\in\omega}$ be a strictly increasing sequence
of Woodin cardinals. Let $\lambda = \sup_{n} \kappa_n$.
Let $\R^*$ be the reals of a symmetric collapse of $V$ up to $\lambda$.
Then $L(\R^*)\models\AD$.
\end{proposition}

See section 9.6 of \cite{Ma2} for a proof of this proposition.

From the existence of infinitely many Woodin cardinals, one can not
prove $\AD^{\LofR}$. However, the above two propositions do yield that
\emph{some} of the games in $\LofR$ are determined. 
For our purposes in this
paper, it will be more than sufficient to know that every game
in $\JofR{\kappa^{\R}}$ is determined. We give that argument next.

\begin{corollary}
\label{determinacyholds}
Let $\angles{\delta_n \mid n\in\omega}$ be a strictly increasing sequence
of Woodin cardinals. 
Let $k$ be the least ordinal such that $\JofR{\kappa}$ is admissible.
Then every game in $\JofR{\kappa}$ is determined.
\end{corollary}
\begin{proof}
Let $\lambda = \sup_{n} \kappa_n$.
Let $\R^{\prime}$ be the set of reals in some symmetric collapse of $V$
up to $\lambda$. Then there is some ordinal $\alpha$ such that
$J_{\alpha}(\R^{\prime})$ is admissible. Let $\alpha_0$ be the least
such ordinal. Since the collapse forcing is homogeneous, $\alpha_0$
does not depend on our particular choice of $\R^{\prime}$.

Let $\gamma>\alpha_0$ be any ordinal.
Let $j:V\map N$ be as in Proposition \ref{firstprop}.
 Let $\R^*$ be the reals of $N$.
By  Proposition \ref{firstprop}, $\R^*$ is also the set of reals in a symmetric
collapse of $V$ up to $\lambda$. Thus $\alpha_0$ is the least
$\alpha$ such that $J_{\alpha}(\R^*)$ is admissible.
 Now $\alpha_0$ is in the wellfounded part of $N$.
So $N\models$``$\alpha_0$ is the least
$\alpha$ such that $\JalphaR$ is admissible.'' Thus $j(\kappa)=\alpha_0$.

By Proposition \ref{secondprop} $L(\R^*)\models\AD$, and so in
particular $J_{\alpha_0}(\R^*)\models\AD$. By the elementarity of
$j$, $\JofR{\kappa}\models\AD$. So every game in $\JofR{\kappa}$ is 
determined.
\end{proof}

The next lemma will tell us that our coiterations terminate in countably 
many steps.
\begin{lemma}
\label{CoiterationsAreCountable}
Assume that there exists infinitely many Woodin cardinals. 
Let $\kappa$ be the least ordinal such that $\JofR{\kappa}$
is admissible.
Let $\cM$ and $\cN$ be countable premice. Let
$\lsequence{(\cT_{\alpha},\cU_{\alpha})}{\alpha<\theta}$ be a coiteration
of $\cM$ and $\cN$ which is in $\JofR{\kappa}$.
Then $\theta<\omega_1$.
\end{lemma}
\begin{proof}
Suppose otherwise. Let $\cT$ on $\cM$ and $\cU$ on $\cN$ be the
iteration trees of length $\omega_1$ associated with our coiteration.
For a contradiction it is enough to show that both $\cT$ and $\cU$ have
cofinal branches.  (For then we can get the usual contradiction involving
stationary subsets of $\omega_1$. See the proof of the Comparison
Theorem in \cite{MiSt}.)

As in the proof of Corollary \ref{determinacyholds} above,
there is a generic elementary embedding 
\mbox{$j:\JofR{\kappa}\map J_{\alpha_0}(\R^*)$} with $\crit(j)=\omega_1$.
As $\kappa$ begins a $\Sigma_1$-gap,
$\lsequence{(\cT_{\alpha},\cU_{\alpha})}{\alpha<\theta}$ is 
$\Sigma_1$ definable in $\JofR{\kappa}$ from some real $x$.

Now $j(\lsequence{(\cT_{\alpha},\cU_{\alpha})}{\alpha<\omega_1})$ is
a coiteration of $\cM$ and $\cN$ in $J_{\alpha_0}(\R^*)$
which properly
extends $\lsequence{(\cT_{\alpha},\cU_{\alpha})}{\alpha<\omega_1}$.
Furthermore
$j(\lsequence{(\cT_{\alpha},\cU_{\alpha})}{\alpha<\omega_1})$ is 
$\Sigma_1$ definable in $J_{\alpha_0}(\R^*)$ from the real $x$.
 Thus there are cofinal
branches $b$ of $\cT$ and $c$ of $\cU$ such that both $b$ and $c$
are  definable in $L(\R^*)$ from $x$ and $\alpha_0$. As the collapse
forcing is homogeneous and $x\in V$, $b$ and $c$ are in $V$.
\end{proof}

\skipmed

We will need to borrow a few more ideas from \cite{MaSt}. First a
definition.

\begin{definition}
Let $\cM$ be a premouse and $\delta\in\ORD^{\cM}$. We say that
$\delta$ is a \emph{cutpoint} of $\cM$ iff for no extender $E$ on the 
$\cM$-sequence do we have that $\crit(E)<\delta\leq\length(E)$.
\end{definition}

The following definition comes from \cite{St2}.

\begin{definition}
Let $\cM$ be  premouse. Then $\cM$ is \emph{tame} iff for all
initial segments $\cN\unlhd\cM$, and all ordinals
$\delta\in\ORD\intersect\cN$, if $\delta$ is a Woodin cardinal of $\cN$,
then for no extender $E$ on the $\cM$ sequence do we have that
$\crit(E)\leq\delta\leq\length(E)$. 
\end{definition}
In particular, if $\cM$ is tame and $\delta$ is a Woodin cardinal of $\cM$,
then $\delta$ is a cutpoint of $\cM$.

Recall that in Section \ref{section:prereq} we gave the definition
of a \emph{meek} premouse. An easy reflection argument shows that
a meek premouse is tame. In particular, if $\cM$ is
premouse  with the property that whenever $\cN\initseg\cM$, 
$\cN$ does not have $\omega$ Woodin cardinals, then $\cM$ is tame.
So all of the premice we are working with in this paper are tame.
The reason that we brought up the notion of tameness is that we
want to highlight the following fact:

\begin{proposition}
Suppose $\cM$ is a tame premouse and $\cT$ is an iteration tree
on $\cM$ of limit length and $b$ is a cofinal wellfounded branch of
$\cT$. Let $\delta=\delta(\cT)$ and $\cQ=Q(b,\cT)$. Suppose
that $\delta\in\cQ$. Then $\delta$ is a cutpoint of $\cQ$.
\end{proposition}

We continue with some more ideas from \cite{MaSt}
The following is Lemma 1.10 from that paper.  The lemma says that
if you can successfully carry out a comparison of two $\omega$-mice
$\cM$ and $\cN$, then in fact the comparison must have ended
immediately as we already have that $\cM\unlhd\cN$ or $\cN\unlhd\cM$.
\begin{lemma}
\label{AlreadyCompared}
Let $\cM$ and $\cN$ be $\delta$-mice, where $\delta$ is a cutpoint of
$\cM$ and $\cN$.  Suppose that $\cJ^{\cM}_{\delta}=\cJ^{\cN}_{\delta}$.
Let $(\cT,\cU)$ be a coiteration of $\cM$ and $\cN$, with $\cP$ the last 
model
of $\cT$, and $\cQ$ the last model of $\cU$.  Suppose $\cP\unlhd\cQ$.
Then $\cM\unlhd\cN$.
\end{lemma}

Let us say that two premice $\cM$ and $\cN$ are {\em comparable} iff
either $\cM\unlhd\cN$ or $\cN\unlhd\cM$. If $\cM$ and $\cN$ are not
comparable, then we will say that they are incomparable.
We will also need the following, which is Lemma 1.12 from
\cite{MaSt}.

\begin{lemma}
\label{IncomparableBranches}
Let $\cM$ be a $\delta$-mouse, and $\cT$ and $\omega$-maximal
iteration tree on $\cM$ above $\delta$, with $\cT$ of limit
length. Let $b$ and $c$ be distinct cofinal wellfounded
branches of $\cT$. Then $\cM^{\cT}_b$ and $\cM^{\cT}_c$
are incomparable.
\end{lemma}

The following lemma and its proof are implicit in \cite{MaSt}.

\begin{lemma}
\label{IncomparableQStructures}
Let $\cM$ be a $\delta$-mouse, and $\cT$ an $\omega$-maximal
iteration tree on $\cM$ above $\delta$, with $\cT$ of limit length.
Suppose that $b$ and $c$ are distinct cofinal wellfounded
branches of $\cT$. Then:
\begin{itemize}
\item[(a)] $Q(b,\cT)$ and $Q(c,\cT)$ are incomparable.
\item[(b)] Let $\deltaprime=\delta(\cT)$. Let $\cP$ be any premouse
such that $\cJ^{\cP}_{\deltaprime}
=\cJ^{\cM^{\cT}_b}_{\deltaprime}=\cJ^{\cM^{\cT}_c}_{\deltaprime}$.
Let $\cR=\cJ^{\cP}_{\alpha}$ where $\alpha$ is largest so that
$\deltaprime=\alpha$ or $\deltaprime$ is Woodin in $\cJ^{\cP}_{\alpha}$.
Suppose that $\cR$ is comparable with  both $Q(b,\cT)$ and $Q(c,\cT)$.
Then $\cP$ is a proper initial segment of both
$Q(b,\cT)$ and $Q(c,\cT)$.
\end{itemize}
\end{lemma}
\begin{proof}
(a) Suppose for example that $Q(b,\cT)\unlhd Q(c,\cT)$.
If $Q(b,\cT)\lhd Q(c,\cT)$, then since $\delta(\cT)$ is Woodin
in $Q(c,\cT)$, $Q(b,\cT)=\cM^{\cT}_b$ and so
$\cM^{\cT}_b\lhd\cM^{\cT}_c$. This contradicts Lemma 
\ref{IncomparableBranches} above.
So we must have $Q(b,\cT)=Q(c,\cT)$ and $Q(b,\cT)\not=\cM^{\cT}_b$ and
$Q(c,\cT)\not=\cM^{\cT}_c$.  Then $\delta(\cT)$ is Woodin in
$\cM^{\cT}_b$ with respect to functions definable over $Q(b,\cT)$ 
(since such
functions are in $\cM^{\cT}_b\intersect\cM^{\cT}_c$.)  But
this contradicts the definition of $Q(b,\cT)$.

(b) Since $Q(b,\cT)$ and $Q(c,\cT)$ are incomparable, we must
have that $\cR$ is a proper initial segment of both
$Q(b,\cT)$ and $Q(c,\cT)$. But then $\cR$ cannot define a function
witnessing the failure of the Woodiness of $\deltaprime$. So $\cR=\cP$.
\end{proof}

Finally, we get to the main content of this section.
The notion of a $k$-realizable premouse, and a 
$k$-realizable branch of an iteration tree are mentioned in the 
statement of the next theorem. See \cite{St2} for these definitions.

\skipsmall

\begin{theorem}
\label{ComparisonTheorem}
Assume that there exists $\omega$ Woodin cardinals. 
Let $(\alpha,n)$ be such that
$2\leq\alpha<\omega_1^{\omega_1}$, and $n\geq 0$.
\begin{itemize}
\item[(1)] Let $\cM$ be a countable $\delta$-mouse which is
$(\alpha,n+1)$-petite  above $\delta$, where $\delta$ is
a cutpoint of $\cM$. Suppose $\cM$ is $k(\cM,\delta)$-realizable via
$\pi$. Let $\cT$ be an $\omega$-maximal iteration tree on $\cM$ of countable
limit length, with $\cT$ above $\delta$. Suppose $\cM^{\cT}_{\xi}$ is
$(\pi,\cT)$-realizable for all $\xi<\length(\cT)$. Let $b$ be the
unique cofinal $(\pi,\cT)$-realizable branch of $\cT$. If $n\geq1$ then
$b$ is the unique cofinal branch $c$ of $\cT$ such that $\cM^{\cT}_c$
is wellfounded and $Q(c,\cT)$ is $\Pa{n-1}$-iterable above $\delta(\cT)$.
If $n=0$ then 
$b$ is the unique cofinal branch $c$ of
$\cT$ such that $M^{\cT}_c$ is wellfounded and $Q(c,\cT)$ is
$\Pa{0}$ iterable above $\delta(\cT)$.

\skipsmall

\item[(2)] Let $\cM$ be a countable $\delta$-mouse which is
$(\alpha,n+1)$-petite above $\delta$, where $\delta$
is a cutpoint of $\cM$. Suppose $\cM$ is $k(\cM,\delta)$ realizable.
Then $\cM$ is $\Pan$ iterable above $\delta$.

\skipsmall

\item[(3)] Suppose $\cM$ and $\cN$ are countable $\delta$-mice, where
$\delta$ is a cutpoint of $\cM$ and $\cN$.  Suppose $\cJ^{\cM}_{\delta}=
\cJ^{\cN}_{\delta}$.  Suppose $\cM$ is $k(\cM,\delta)$ realizable and
$(\alpha,n+1)$-petite  above $\delta$, and that $\cN$ is
$\Pan$ iterable above $\delta$. Then $\cM\unlhd\cN$ or $\cN\unlhd\cM$.
\end{itemize}
\end{theorem}

\begin{note}
We conjecture that the above theorem is also true if the word
``petite'' is replaced everywhere in the statement of the theorem
with the word ``small.'' Unfortunately, we are unable to prove
this.  During the following
proof, we will point out exactly where the proof of the stronger
version of the theorem---the version with ``small'' instead of
``petite''---breaks down. A proof of the stronger version of the
theorem would be a very nice result. In particular, it would imply
that $\Aan$ is a mouse set in the case $\cof(\alpha)\leq\omega$ and
$n>0$. Currently, we are not even able to show that $\Aa[2]{1}$
is a mouse set.
\end{note}

\begin{proof}[Proof of the theorem]

The theorem is proved by induction on $(\alpha,n)$.
We begin with (1) for $(\alpha,n)$. Let
$\cM,\delta,\pi$ and $\cT$ be as in (1). Let $b$ be a cofinal
$(\pi,\cT)$ realizable branch of $\cT$; it is shown in \cite{St2}
that there is a unique such branch.
First we will show that $b$ has the iterability property stated in
(1).  Then we will show that $b$ is the unique branch with this
property.

First suppose $n=0$.
Since $Q(b,\cT)$ is realizable, $Q(b,\cT)$ is length-$\omega$ iterable
above $\delta(\cT)$. (See Corollary 1.9 from \cite{St2}.) This trivially
implies that $Q(b,\cT)$ is $\Pa{0}$ iterable above $\delta(\cT)$.

Now suppose that $n\geq1$.
By Lemma \ref{PreservationOfSmallness}, $\cM^{\cT}_b$ is 
$(\alpha,n+1)$-petite above $\delta$. 
If $\delta(\cT)\notin Q(b,\cT)$
then $Q(b,\cT)$ is trivially iterable above $\delta(\cT)$. So suppose that
$\delta(\cT)\in Q(b,\cT)$. Then $Q(b,\cT)\models\delta(\cT)$ is Woodin.
By part (a) of Lemma \ref{PetiteAboveWoodin}
on page \pageref{PetiteAboveWoodin},
$Q(b,\cT)$ is $(\alpha,n)$ petite above $\delta(\cT)$. 
So by part (2) of the theorem we are proving,
for $(\alpha,n-1)$, $Q(b,\cT)$ is $\Pa{n-1}$ iterable
above $\delta(\cT)$. 

Now we show that $b$ is unique. Suppose towards a contradiction that there 
is
a second cofinal branch $c$ of $\cT$ such that
$\cM^{\cT}_c$
is wellfounded, and in the case $n\geq1$, $Q(c,\cT)$ is $\Pa{n-1}$-iterable
above $\delta(\cT)$, and in the case $n=0$,
$Q(c,\cT)$ is $\Pa{0}$-iterable above $\delta(\cT)$.

First suppose $n\geq 1$.
As we pointed out above,
$Q(b,\cT)$ is $(\alpha,n)$ petite above $\delta(\cT)$.
By part (3) of the theorem we are proving, for
$(\alpha,n-1)$,  we have that
$Q(b,\cT)\unlhd Q(c,\cT)$ or $Q(c,\cT)\unlhd Q(b,\cT)$. This contradicts
part (a) of Lemma \ref{IncomparableQStructures}.

Now suppose $n=0$ and there is a second cofinal branch $c$ of $\cT$ such
that $\cM^{\cT}_c$ is wellfounded and $Q(c,\cT)$ is 
$\Pa{0}$-iterable above $\delta(\cT)$. Again we will derive a contradiction
by showing that $Q(b,\cT)\unlhd Q(c,\cT)$ or $Q(c,\cT)\unlhd Q(b,\cT)$.
If $\delta(\cT)\notin Q(b,\cT)$
then $Q(b,\cT)\unlhd Q(c,\cT)$. So suppose that
$\delta(\cT)\in Q(b,\cT)$. Then $Q(b,\cT)\models\delta(\cT)$ is Woodin.
By part (b) of Lemma \ref{PetiteAboveWoodin}
on page \pageref{PetiteAboveWoodin}, there is some $\alphaprime<\alpha$
and some $m\in\omega$ such that 
$Q(b,\cT)$ is $(\alphaprime,m+1)$-petite above
$\delta(\cT)$. Fix such an $\alphaprime$ and an $m$ with $m$ even.
Since $m$ is even, $Q(c,\cT)$ is $\Pa[\alphaprime]{m}$-iterable
above $\delta(\cT)$. By part (3) of the theorem we are proving, for
$(\alphaprime,m)$, (or by \cite{MaSt} in the case $\alphaprime=1$)
we have that
$Q(b,\cT)\unlhd Q(c,\cT)$ or $Q(c,\cT)\unlhd Q(b,\cT)$. Again,
this contradicts
part (a) of Lemma \ref{IncomparableQStructures}.
This completes the proof of (1) for $(\alpha,n)$.

\skipmed

Next we prove (2) for $(\alpha,n)$. Let $\cM,\delta$ be as in (2). Since
$\cM$ is realizable, $\cM$ is length-$\omega$ iterable above $\delta$.
(See Corollary 1.9 from \cite{St2}.)
If $n$ is even this implies that $\cM$ is $\Pan$-iterable above $\delta$.
So suppose $n\geq1$ is odd. Let $x_0\in\R$ be a mouse code for
$(\cM,\delta)$. We will now describe a winning strategy for player
$\TWO$ in $G_{(\alpha,n)}(x_0)$.

Let $1\leq k\leq n+1$. At round $k$ of the game we have a real
$x_{k-1}$ such that $x_{k-1}$ is a mouse code for some
$(\cM^{\prime},\deltaprime)$. (If $k=1$ then
$(\cM^{\prime},\deltaprime) = (\cM,\delta)$.)
We maintain by induction on $k$ the
following hypothesis, which we shall refer to using the symbol $\dagger$:
\begin{align}
&\cM^{\prime}\text{ is }(\alpha,n+2-k)\text{-petite %
 above }\deltaprime,\tag*{$(\dagger)_{k-1}$}\\
&\text{and }\deltaprime\text{ is a cutpoint of }\cM^{\prime}
\text{ and }\cM^{\prime}\text{ is }k(\cM^{\prime},\deltaprime)
\text{ realizable.}\notag
\end{align}

Notice that with $k=1$ we do have that $(\dagger)_0$ holds.

At the start of round $k$, player
$\ONE$ will play a real $y_k$ so that $y_k^*$ is an $\omega$-maximal
putative iteration tree $\cT$ on $\cM^{\prime}$ above $\deltaprime$.
We will define a real $x_k$ which player $\TWO$ should play in response.
If $\cT$ has a last wellfounded model then player $\ONE$ loses. So we
assume otherwise. Let $\pi$ be such that $\cM^{\prime}$ is
$k(\cM^{\prime},\deltaprime)$-realizable via $\pi$.
Then $\cT$ has a maximal, wellfounded, $(\pi,\cT)$-realizable branch.
(This is proved in \cite{St2}.)
Let $b$ be such a branch of minimal sup.
Let $\bar{\cT}=\cT\restriction\sup(b)$. Then $b$ is the
unique cofinal $(\pi,\bar{\cT})$-realizable branch of $\bar{\cT}$
and $\cM^{\bar{\cT}}_{\xi}$ is $(\pi,\bar{\cT})$-realizable
for all $\xi<\length(\bar{\cT})$. So $\cM^{\prime},\deltaprime,
\bar{\cT}, b$ satisfy the hypotheses of part (1) of this theorem
for $(\alpha,n)$.
So $b$ is the unique cofinal branch $c$ of $\bar{\cT}$
such that $\cM^{\bar{\cT}}_c$
is wellfounded and $Q(c,\bar{\cT})$ is $\Pa{n-1}$-iterable above
$\delta(\bar{\cT})$. (Notice that
$n\geq 1$ since we are assuming that $n$ is odd.)

Let $\cQ=Q(b,\bar{\cT})$ and $\gamma=\delta(\bar{\cT})$.
Since $\cM^{\prime}$ is $(\alpha,n+2-k)$-petite above $\deltaprime$,
by Lemma \ref{PreservationOfSmallness},
$\cM^{\bar{\cT}}_b$ is $(\alpha,n+2-k)$-petite above $\deltaprime$.
Then by Lemma \ref{PetiteAboveWoodin},
$\cQ$ is $(\alpha,n+2-(k+1))$-petite  above $\gamma$.
We will have
player $\TWO$ play some $x_k$ such that $x_k$ is a mouse code for
$(\cQ,\gamma)$. Notice that if we can do this then we will have
maintained our inductive hypotheses $(\dagger)_k$
about $x_k$. The main problem now
is to show that there is such an $x_k$ with $x_k\in\Dan(y_k)$.

Let $\bar{y}$ be a real such that $\bar{y}^* =\bar{\cT}$ and
$\bar{y}$ is recursive in $y_k$.
Let $P$ be the binary relation on $\HC$ given by
$P(\cT^{\prime},\angles{\cQ^{\prime},\delta^{\prime}})\Iff$ 
$\cT^{\prime}$ is
an iteration tree of limit length, and 
$\delta^{\prime}=\delta(\cT^{\prime})$,
and there exists exactly one cofinal wellfounded branch $\bprime$ 
of $\cT^{\prime}$ such that $\cQ^{\prime}=Q(b^{\prime},\cT^{\prime})$.
Then $P$ admits a definable code transformation.
(See the proof of Proposition \ref{DefCodTransProp}.) 
Let $f$ be a definable code
transformation for $P$. 
Notice that $P(\bar{\cT},\angles{\cQ,\gamma})$ holds.
Let $x_k=f(\bar{\cT},\angles{\cQ,\gamma},\bar{y})$. 
Then $x_k^*=(\cQ,\gamma)$. Furthermore, $x_k$ is the unique real
$x$ such that
\begin{quote}
$(\exists z,w,r\in\R)$ s.t. $z^*$ is a cofinal wellfounded branch of
$\bar{y}^*$, and $w^* = Q(z^*,\bar{y}^*)$ and $r^*=\delta(\bar{y}^*)$
 and $x^*=(w^*,r^*)$
and player $\TWO$ wins $G_{(\alpha,n-1)}(x)$ and
$x=f(\bar{y}^*,x^*,\bar{y})$.
\end{quote}
Since the relation ``player $\TWO$ wins $G_{(\alpha,n-1)}(x)$'' is
$\Pa{n-1}$, we have that
$x_k\in\Dan(\bar{y})$.

Clearly,
if player $\TWO$ plays this $x_k$ he will not lose $G_{(\alpha,n)}(x_0)$
at round $k$.

Playing in this way for $n+1$ rounds we eventually produce a real
$x_{n+1}$.  We must now show that $\exists \sigma\in\JalphaR$
such that $\sigma$ is a winning quasi-strategy for player $\TWO$ in
the  length-$\omega$ iteration game $WG(x_{n+1})$. By our induction
hypothesis $(\dagger)_{n}$,
 $x_n$ is a mouse code for some
$(\cM^{\prime},\deltaprime)$ such that $\cM^{\prime}$ is $(\alpha,1)$-petite
above $\deltaprime$. 
Thus
 $x_{n+1}$ is a mouse code for some
$(\cM^{\prime\prime},\delta^{\prime\prime})$ such that
$\cM^{\prime\prime}$  is $(\alpha,1)$-petite above $\deltaprime$, and
so by Lemma \ref{PetiteAboveWoodin},
for some $\beta<\alpha$ and some $m\in\omega$,
$\cM^{\prime\prime}$  is $(\beta,m+1)$ petite above $\delta^{\prime\prime}$.
Let $\sigma$ be the quasi-strategy for player $\TWO$ in the game
$WG(x_{n+1})$ which dictates that whenever player $\ONE$ plays a tree,
player $\TWO$ should play the unique maximal realizable branch of minimal
sup. This is a winning quasi-strategy for player $\TWO$, and by
(1) for $(\beta,m)$, $\sigma$ is definable over $\JbetaR$ and so
$\sigma\in\JalphaR$.

\skipmed

Finally we prove (3) for $(\alpha,n)$. There are three cases: $n=0$,
$n\geq 1$ odd, and $n\geq 2$ even.

\underline{CASE 1}\quad  $n=0$.\qquad 
Suppose $\cM$ and $\cN$ are countable $\delta$-mice, where
$\delta$ is a cutpoint of $\cM$ and $\cN$.  Suppose $\cJ^{\cM}_{\delta}=
\cJ^{\cN}_{\delta}$.  Suppose $\cM$ is $k(\cM,\delta)$ realizable and
$(\alpha,1)$-petite  above $\delta$, and that $\cN$ is
$\Pa{0}$ iterable above $\delta$. Recall that this this means that
$\cN$ is rank-$\lambda(\alpha)$ iterable above $\delta$ where
$$\lambda(\alpha)=
\begin{cases}
\omega\alpha& \text{if $\alpha\geq\omega$,}\\
\omega(\alpha - 1)& \text{if $\alpha<\omega$}.
\end{cases}
$$
Suppose that
$\cM\notinitseg\cN$ and $\cN\notinitseg\cM$. 
We will derive a contradiction.

\begin{note}
It is interesting to note that the proof we are giving works in both
the case $\cof(\alpha)\leq\omega$ and the case $\cof(\alpha)>\omega$.
\end{note}

We define a coiteration $\bigsequence{(\cT_{\gamma},\cU_{\gamma})}%
{\gamma<\theta}$ of $\cM$ and $\cN$ by induction on $\gamma$. For this it
suffices to define $\cT_{\gamma+1}$ in the case $\length(\cT_{\gamma})$
is a limit ordinal, and similarly with $\cU_{\gamma+1}$. Fix a
$k(\cM,\delta)$ realization $\pi$ of $\cM$. If
$\length(\cT_{\gamma})$ is a limit ordinal then let $b$ be the unique
cofinal $(\pi,\cT)$ realizable branch of $\cT_{\gamma}$. It is shown
in \cite{St2} that
there always is such a $b$. Let
$\cT_{\gamma+1}$ be the iteration tree extending $\cT_{\gamma}$
with length $\length(\cT_{\gamma})+1$ and such that
$b=[0,\length(\cT_{\gamma})]_{T_{\gamma+1}}$.

If $\length(\cU_{\gamma})$ is a limit ordinal then let $c$ be the unique
cofinal wellfounded branch of $\cU_{\gamma}$ such that $Q(c,\cU_{\gamma})$
is rank-$\lambda(\alpha)$ iterable above $\delta(\cU_{\gamma})$, if there
is a unique such branch. If $c$ is defined then let
$\cU_{\gamma+1}$ be the iteration tree extending $\cU_{\gamma}$
with length $\length(\cU_{\gamma})+1$ and such that
$c=[0,\length(\cU_{\gamma})]_{U_{\gamma+1}}$. If $c$ is not defined then
stop the construction.

Suppose this construction produces a coiteration
$\bigsequence{(\cT_{\gamma},\cU_{\gamma})}{\gamma<\omega_1}$ of length
$\omega_1$. By (1) for $(\alpha,0)$ the branch $b$ of $\cT_{\gamma}$
chosen when $\length(\cT_{\gamma})$ is a limit is
$\Pi_1(\JalphaR,\singleton{\cT_{\gamma}})$ uniformly in $\gamma$.
Clearly the branch $c$ of $\cU_{\gamma}$
chosen when $\length(\cU_{\gamma})$ is a limit is
$\Pi_1(\JalphaR,\singleton{\cU_{\gamma}})$ uniformly in $\gamma$.
It follows that $\bigsequence{(\cT_{\gamma},\cU_{\gamma})}{\gamma<\omega_1}
\in\JofR{\alpha+\omega}$. 
This contradicts Lemma \ref{CoiterationsAreCountable}. Therefore
the construction must produce a coiteration
\mbox{$\bigsequence{(\cT_{\gamma},\cU_{\gamma})}{\gamma\leq\theta}$} with
$\theta<\omega_1$, such that either the construction stops at stage
$\theta$ for the reason described above, or
$\bigsequence{(\cT_{\gamma},\cU_{\gamma})}{\gamma\leq\theta}$ is terminal.
If the construction produces a terminal coiteration, and if this
coiteration is successful, then by Lemma \ref{AlreadyCompared},
$\cM\unlhd\cN$ or $\cN\unlhd\cM$. This is a contradiction. So
we assume that either the construction produces a terminal,
non-successful coiteration, or the construction stops for the reason
described above. This leaves us with three possibilities:
\begin{itemize}
\item[(a)] $\cU_{\theta}$ has successor length and the coiteration is
terminal because the next ultrapower
which is prescribed by the definition of a coiteration
in order to define $\cU_{\theta+1}$ is illfounded. Or,
\item[(b)] $\cU_{\theta}$ has limit length but has no cofinal wellfounded
branch $c$ such that $Q(c,\cU_{\theta})$ is rank-$\lambda(\alpha)$
iterable above $\delta(\cU_{\theta})$. Or,
\item[(c)]$\cU_{\theta}$ has limit length and has two cofinal wellfounded
branches $c_1$ and $c_2$ such that both $Q(c_1,\cU_{\theta})$ and
$Q(c_2,\cU_{\theta})$ are rank-$\lambda(\alpha)$ iterable above
$\delta(\cU_{\theta})$.
\end{itemize}
We will show that each of these possibilities leads to a contradiction.
We begin with the easiest, which is the third.

Suppose towards
a contradiction that possibility (c) holds. Let
$\cQ_1=Q(c_1,\cU_{\theta})$ and $\cQ_2=Q(c_2,\cU_{\theta})$.
Let $\cR$ be the corresponding  realizable $Q$-structure from
$\cT_{\theta}$. (Below we explain what we mean by this. After that we
will feel free to use the phrase ``the corresponding realizable
$Q$-structure'' without further explanation.) That is,
let $\cP$ be the last model from $\cT_{\theta}$, if $\cT_{\theta}$ has
successor length, otherwise let $\cP$ be the direct limit model of the
unique cofinal realizable branch of $\cT_{\theta}$. Let
$\gamma=\delta(\cU_{\theta})$. Notice that $\gamma=\delta(\cT_{\theta})$
if $\cT_{\theta}$ has limit length. In any case we know that
$\cJ^{\cQ_1}_{\gamma}=\cJ^{\cQ_2}_{\gamma}=\cJ^{\cP}_{\gamma}$. Finally
let $\cR=\cJ^{\cP}_{\xi}$ where $\xi$ is largest such that $\xi=\gamma$
or $\cJ^{\cP}_{\xi}\models\gamma$ is Woodin. This is what we mean by
the corresponding realizable $Q$-structure from $\cT_{\theta}$. Notice
that if $\cT_{\theta}$ has limit length and if $b$ is the unique cofinal
realizable branch of $\cT_{\theta}$  then $\cR=Q(b,\cT_{\theta})$. Now,
since $\cM$ is $(\alpha,1)$-petite above $\delta$, by Lemma
\ref{PreservationOfSmallness},
$\cP$ is $(\alpha,1)$-petite above $\delta$, and
so, by Lemma \ref{PetiteAboveWoodin},
there is an $\alphaprime$ with $1\leq\alphaprime<\alpha$, and an
$m\in\omega$, such that $\cR$ is $(\alphaprime,m)$-petite above $\gamma$.
Fix such a pair $(\alphaprime,m)$ with $m$ even. Since $m$ is even and
$\cQ_1$ and $\cQ_2$ are rank-$\lambda(\alpha)$ iterable above $\gamma$,
we have that
$\cQ_1$ and $\cQ_2$ are $\Pa[\alphaprime]{m}$ iterable above $\gamma$.
So by part (3) of this theorem for $(\alphaprime,m)$ we have that
$\cQ_1\unlhd\cR$ or $\cR\unlhd\cQ_1$ and  also that
$\cQ_2\unlhd\cR$ or $\cR\unlhd\cQ_2$. By part (b) of Lemma
\ref{IncomparableQStructures}, we have that $\cP\lhd\cQ_1$.
But then by Lemma \ref{AlreadyCompared},
$\cM\unlhd\cN$. This is a contradiction.

Now we assume towards a contradiction
that either possibility (a) or possibility (b) has occurred
with our coiteration
$\bigsequence{(\cT_{\gamma},\cU_{\gamma})}{\gamma\leq\theta}$.
Let $x_0$
be a mouse code for $(\cN,\delta)$. Since $\cN$ is rank $\lambda(\alpha)$
iterable above $\delta$, player $\TWO$ wins the rank-$\lambda(\alpha)$
iteration game $G^{\lambda(\alpha)}(x_0)$. Fix a winning strategy $\sigma$
for player $\TWO$ in this game.
If possibility
(a) holds then let $\cU$ be the putative iteration tree extending
$\cU_{\theta}$ whose last model is illfounded and which results from
$\cU_{\theta}$ by taking the ultrapower which is mentioned in possibility
(a). If possibility (b) holds then let $\cU=\cU_{\theta}$. In either case
let $y_1$ be a real such that $y_1^* = \cU$.  Notice that
for each $\beta$ with $1\leq\beta<\alpha$, and each $m\in\omega$,
$(\lambda(\beta)+m,y_1)$ is a legal first move for player $\ONE$ in
$G^{\lambda(\alpha)}(x_0)$, and that if player $\ONE$ plays this
move he does not immediately lose. For each such pair $(\beta,m)$,
let $x_1^{(\beta,m)}$ be
$\sigma$'s response when player $\ONE$ plays $(\lambda(\beta)+m,y_1)$
as his first move in the game $G^{\lambda(\alpha)}(x_0)$. Since $\sigma$
is a winning strategy, for each pair $(\beta,m)$
there is some maximal wellfounded branch $c$ of $\cU$ so that
$(x_1^{(\beta,m)})^*=\bigl\langle Q(c,\cU\restriction\sup(c)),\,
\delta\bigl(\cU\restriction\sup(c)\bigr)\,\bigr\rangle$. 
Let $c_{(\beta,m)}$ be
such a maximal wellfounded branch. Let
$\cQ_{(\beta,m)}$, $\delta_{(\beta,m)}$
be such that $(x_1^{(\beta,m)})^*=
\bigl\langle\cQ_{(\beta,m)},\, \delta_{(\beta,m)}\,\bigr\rangle$.  
Since $\sigma$ is a
winning strategy, $\cQ_{(\beta,m)}$ is rank-$(\lambda(\beta)+m)$
iterable above $\delta_{(\beta,m)}$.
In summary, for each $(\beta,m)$ such that
$(1,0)\lexleq(\beta,m)\lexless(\alpha,0)$, we have the following:
\begin{itemize}
\item[(i)] $c_{(\beta,m)}$ is a maximal wellfounded branch of $\cU$.
\item[(ii)]$\cQ_{(\beta,m)}=
Q(c_{(\beta,m)},\cU\restriction\sup(c_{(\beta,m)}))$
\item[(iii)]$\delta_{(\beta,m)}=
\delta(\cU\restriction\sup(c_{(\beta,m)})).$
\item[(iv)]$\cQ_{(\beta,m)}$ is rank-$(\lambda(\beta)+m)$ iterable above
$\delta_{(\beta,m)}$.
\end{itemize}
By discarding some of the $c_{(\beta,m)}$ if necessary, and then re-indexing,
we may maintain
(i) through (iv) above and in addition we may assume
\begin{itemize}
\item[(v)] If $(\beta,m)<_{\text{lex}}(\betaprime,\mprime)$ then
$\sup(c_{(\beta,m)})\leq\sup(c_{(\betaprime,\mprime)})$.
\end{itemize}
It follows from this that if $(\beta,m)<_{\text{lex}}(\betaprime,\mprime)$
then $\delta_{(\beta,m)}\leq\delta_{(\betaprime,\mprime)}$.

Now, for each $(\beta,m)$ as above,
let $\gamma_{(\beta,m)}\leq\theta$ be that stage of the coiteration
construction
such that $c_{(\beta,m)}$ is
cofinal in $\cU_{\gamma_{(\beta,m)}}$. By property (v) above,
if $(\beta,m)<_{\text{lex}}(\betaprime,\mprime)$ then
$\gamma_{(\beta,m)}\leq\gamma_{(\betaprime,\mprime)}$.

Fix $(\beta,m)$ for a moment.
Notice that
$$\delta_{(\beta,m)}= \delta(\cU_{\gamma_{(\beta,m)}})$$
and
$$\cQ_{(\beta,m)}= Q(c_{(\beta,m)},\cU_{\gamma_{(\beta,m)}}).$$
Now let $\cR_{(\beta,m)}$ be
the corresponding realizable $Q$-structure from $\cT_{\gamma_{(\beta,m)}}$.
(We explained what we meant by this phrase in the proof above that
possibility (c) cannot occur.)

\begin{claim}[Claim 1]
Suppose $1\leq\beta<\alpha$ and $m\in\omega$ and $m$ is even. Then
$\cR_{(\beta,m)}$ is not $(\beta,m+1)$-petite above $\delta_{(\beta,m)}$.
\end{claim}
\begin{subproof}[Proof of Claim 1]
Suppose towards a contradiction that $\cR_{(\beta,m)}$ is
$(\beta,m+1)$-petite above $\delta_{(\beta,m)}$. 
Now, since $\cQ_{(\beta,m)}$ is rank-$(\lambda(\beta)+m)$ iterable above
$\delta_{(\beta,m)}$ and $m$ is even, $\cQ_{(\beta,m)}$ is
$\Pa[\beta]{m}$ iterable above $\delta_{(\beta,m)}$.
By part (3) of this theorem for $(\beta,m)$,
$\cR_{(\beta,m)}\unlhd\cQ_{(\beta,m)}$ or vice-versa.
Now, $c_{(\beta,m)}$ is a cofinal branch of $\cU_{\gamma_{(\beta,m)}}$,
and $\cQ_{(\beta,m)}= Q(c_{(\beta,m)},\cU_{\gamma_{(\beta,m)}})$.
To derive our contradiction we will show that there is a second cofinal
branch $c$ of $\cU_{\gamma_{(\beta,m)}}$ such that
$Q(c,\cU_{\gamma_{(\beta,m)}})$ is also rank-$(\lambda(\beta)+m)$
iterable above $\delta_{(\beta,m)}$. To see that such a $c$ exists consider
two cases. First suppose that $\gamma_{(\beta,m)}<\theta$. Then
$c_{(\beta,m)}$ is a maximal, non-cofinal branch of $\cU_{\theta}$.
In this case let $c$ be the unique cofinal branch of
$\cU_{\gamma_{(\beta,m)}}$ which is a non-maximal branch of $\cU$. That
is, let $c$ be the branch which was chosen by our construction. Then
by the definition of our construction, $Q(c,\cU_{\gamma_{(\beta,m)}})$
is rank-$\lambda(\alpha)$ iterable above $\delta_{(\beta,m)}$.
Next consider the case that $\gamma_{(\beta,m)}=\theta$.  Then
$c_{(\beta,m)}$ is a cofinal branch of $\cU_{\theta}$ and so
it must be possibility (b) which occurred. That is,
$\cU_{\theta}$ has limit length but has no cofinal wellfounded
branch $c$ such that $Q(c,\cU_{\theta})$ is rank-$\lambda(\alpha)$
iterable above $\delta(\cU_{\theta})$.
Now we also have that
$\gamma_{(\betaprime,\mprime)}=\theta$ for all
$(\betaprime,\mprime)$ with
$(\beta,m)<_{\text{lex}}(\betaprime,\mprime)<_{\text{lex}}(\alpha,0)$.
So in this case we can let
$c=c_{(\betaprime,\mprime)}$ where $(\betaprime,\mprime)$ is
lexicographically
least such that $(\beta,m)<_{\text{lex}}(\betaprime,\mprime)$ and
$c_{(\betaprime,\mprime)}\not= c_{(\beta,m)}$.  There is such a
$(\betaprime,\mprime)$ because otherwise $\cQ_{(\beta,m)}$ is
rank-$\lambda(\alpha)$ iterable above
$\delta_{(\beta,m)} = \delta(\cU)$, and this violates the fact that
possibility (b) has occurred.
So in either case we have seen that there is a second
cofinal branch $c$ of $\cU_{\gamma_{(\beta,m)}}$ such that
$Q(c,\cU_{\gamma_{(\beta,m)}})$ is rank-$(\lambda(\beta)+m)$ iterable
above $\delta_{(\beta,m)}$.
So exactly
as with $\cQ_{(\beta,m)}$ above, we see that
$\cR_{(\beta,m)}\unlhd Q(c,\cU_{\gamma_{(\beta,m)}})$, or vice-versa.
Now, as in the proof above that possibility (c) cannot occur, we
see that part (b) of Lemma
\ref{IncomparableQStructures} implies that at stage $\gamma_{(\beta,m)}$
we attain a successful coiteration of $\cM$ and $\cN$, and so
by Lemma \ref{AlreadyCompared}, $\cM\unlhd\cN$.
This is a contradiction.
\end{subproof}

If $\cT_{\theta}$ has successor length, let $\cP$ be the last model of
$\cT_{\theta}$.  Otherwise, let $\cP$ be the direct limit model of the
unique cofinal realizable branch of $\cT_{\theta}$. We will derive our
contradiction by showing that $\cP$ is not $(\alpha,1)$-petite above
$\delta$.  As $\cM$ is $(\alpha,1)$-petite above $\delta$, this will
contradict Lemma \ref{PreservationOfSmallness}.

\begin{claim}[Claim 2]
Suppose $1\leq\beta<\alpha$ and $m\in\omega$. Then
$\cR_{(\beta,m)}\unlhd \cP$.
\end{claim}
\begin{subproof}[Proof of Claim 2]
Let $\cW$
be the last model model from $\cT_{\gamma_{(\beta,m)}}$, if
$\cT_{\gamma_{(\beta,m)}}$ has successor length, otherwise, let
$\cW$ be the direct limit model of the unique cofinal realizable
branch of $\cT_{\gamma_{(\beta,m)}}$. By definition of
$\cR_{(\beta,m)}$, $\cR_{(\beta,m)}\unlhd\cW$.
If $\gamma_{(\beta,m)}=\theta$ then $\cW=\cP$ and so
$\cR_{(\beta,m)}\unlhd \cP$. So suppose that $\gamma_{(\beta,m)}<\theta$.
Then $\cW$ is a model on the tree $\cT_{\theta}$.

Now, let $c$ be the
unique cofinal branch of $\cU_{\gamma_{(\beta,m)}}$ which is a
non-maximal branch of $\cU_{\theta}$. That
is, let $c$ be the branch which was chosen by our construction. Then
by the definition of our construction, $Q(c,\cU_{\gamma_{(\beta,m)}})$
is rank-$\lambda(\alpha)$ iterable above $\delta_{(\beta,m)}$.
Let $\cS=Q(c,\cU_{\gamma_{(\beta,m)}})$. Since
$\cR_{(\beta,m)}$ is $(\alpha,1)$-petite above $\delta$,
by Lemma \ref{PetiteAboveWoodin}, there is some
$\betaprime<\alpha$ and some $\mprime\in\omega$ such that
$\cR_{(\beta,m)}$ is $(\betaprime,\mprime+1)$-petite above
$\delta_{(\beta,m)}$. Fix such a $(\betaprime,\mprime)$ with $\mprime$
even. Since $\cS$ is rank-$\lambda(\alpha)$ iterable above
$\delta_{(\beta,m)}$ and $\mprime$ is even, $\cS$ is
$\Pa[\betaprime]{\mprime}$ iterable above
$\delta_{(\beta,m)}$.
Then by part (3) of this theorem for $(\betaprime,\mprime)$,
$\cR_{(\beta,m)}\unlhd\cS$ or $\cS\unlhd\cR_{(\beta,m)}$. 
Suppose for a moment that $\cR_{(\beta,m)}\lhd\cS$.
If $\cR_{(\beta,m)}\lhd\cW$,
 then $\delta_{(\beta,m)}$ is Woodin
in $\cW$ with respect to functions definable over $\cR_{(\beta,m)}$.
This contradicts the definition of $\cR_{(\beta,m)}$, and
so we must have that $\cR_{(\beta,m)}=\cW$.
 So $\cW\lhd\cS$, and so by Lemma \ref{AlreadyCompared},
$\cM\unlhd\cN$. This is a contradiction.
So $\cR_{(\beta,m)}\notproperseg\cS$. Similarly,
$\cS\notproperseg\cR_{(\beta,m)}$. Thus $\cR_{(\beta,m)}=\cS$.

Now let $\xi$ be such that $\cW=\cM^{\cT_{\theta}}_{\xi}$. If
$\cW$ is the last model of $\cT_{\theta}$, then $\cW=\cP$ and
so $\cR_{(\beta,m)}\unlhd \cP$. So suppose that $\cW$ is not the
last model of $\cT_{\theta}$. Then the extender $E^{\cT_{\theta}}_{\xi}$
is defined.
By the definition of a coiteration,
$\length(E^{\cT_{\theta}}_{\xi})$ is the least ordinal $\eta$ such
that $\cJ^{\cW}_{\eta}\not=\cJ^{\cN_c}_{\eta}$.
Since $\cR_{(\beta,m)}=\cS$, $\length(E^{\cT_{\theta}}_{\xi}) >
\ORD\intersect\cR_{(\beta,m)}$. By the properties of an iteration tree,
$\cW$ and $\cP$ agree below the
length of this extender.  Thus $\cR_{(\beta,m)}\lhd \cP$.
\end{subproof}

One consequence of  Claim 2 is that for all
$(\beta,m),(\betaprime,\mprime)$, we have that
$\cR_{(\beta,m)}\unlhd\cR_{(\betaprime,\mprime)}$ or
$\cR_{(\betaprime,\mprime)}\unlhd\cR_{(\beta,m)}$.
If $\gamma_{(\beta,m)}=\gamma_{(\betaprime,\mprime)}$, then
by definition $\delta_{(\beta,m)}=\delta_{(\betaprime,\mprime)}$ and
$\cR_{(\beta,m)}=\cR_{(\betaprime,\mprime)}$. We also have the following.

\begin{claim}[Claim 3]
Suppose $(\beta,m)$ and $(\betaprime,\mprime)$ are such
that $\gamma_{(\beta,m)}<\gamma_{(\betaprime,\mprime)}<\theta$.
Then
$\cR_{(\beta,m)}\lhd\cR_{(\betaprime,\mprime)}$ and
and $\ORD\intersect\cR_{(\beta,m)}<\delta_{(\betaprime,\mprime)}$
and
$\delta_{(\beta,m)}$ is a cardinal in $\cR_{(\betaprime,\mprime)}$.
\end{claim}
\begin{subproof}[Proof of Claim 3]
Fix $(\beta,m),(\betaprime,\mprime)$ as in the statement of the
claim.
Let $c$ be the
unique cofinal branch of $\cU_{\gamma_{(\beta,m)}}$ which is a
non-maximal branch of $\cU_{\theta}$. That
is, let $c$ be the branch which was chosen by our construction
in order to define $\cU_{\gamma_{(\beta,m)}+1}$.
Let $\cS=Q(c,\cU_{\gamma_{(\beta,m)}})$.
Let $\cH=\cN_c$. Then $\cH$ is the last model of
$\cU_{\gamma_{(\beta,m)}+1}$. Also $\cH$ is on the tree $\cU_{\theta}$.
Let $\xi$ be such that $\cH=\cN^{\cU_{\theta}}_{\xi}$. $\cH$ is not
the last model of $\cU_{\theta}$ and so the extender
$E^{\cU_{\theta}}_{\xi}$ is defined. As in the proof of Claim 2 we
see that $\cR_{(\beta,m)}=\cS$, and so
$\length(E^{\cU_{\theta}}_{\xi}) > \ORD\intersect\cS$.
Let $\lambda=\length(E^{\cU_{\theta}}_{\xi})$.

Now  let $\cprime$ be the
unique cofinal branch of $\cU_{\gamma_{(\betaprime,\mprime)}}$ which is a
non-maximal branch of $\cU_{\theta}$. That
is, let $\cprime$ be the branch which was chosen by our construction
in order to define $\cU_{\gamma_{(\betaprime,\mprime)}+1}$.
Let $\cS^{\prime}=Q(\cprime,\cU_{\gamma_{(\betaprime,\mprime)}})$.
Let $\cH^{\prime}=\cN_{\cprime}$. Then $\cH^{\prime}$ is the last model of
$\cU_{\gamma_{(\betaprime,\mprime)}+1}$.
Also $\cH^{\prime}$ is on the tree $\cU_{\theta}$.
As in the proof of Claim 2 we see that $\cR_{(\betaprime,\mprime)}=
\cS^{\prime}$.

Let $\xiprime$ be such that $\cH^{\prime}=\cN^{\cU_{\theta}}_{\xiprime}$.
Clearly $\xi<\xiprime$.
So by the properties of an iteration tree we have that
$J^{\cH}_{\lambda}=J^{\cH^{\prime}}_{\lambda}$ and $\lambda$ is a
cardinal in $\cH^{\prime}$. Also $\lambda<\delta_{(\betaprime,\mprime)}$
and so $J^{\cH}_{\lambda}\lhd\cS^{\prime}$. So we have
$\cR_{(\beta,m)}=\cS\lhd J^{\cH}_{\lambda}\lhd\cS^{\prime}=
\cR_{(\betaprime,\mprime)}$. Also, $\lambda$ is a cardinal in
$\cS^{\prime}=\cR_{(\betaprime,\mprime)}$ and $\delta_{(\beta,m)}$
is a cardinal in $\cH$, so $\delta_{(\beta,m)}$ is a cardinal in
$J^{\cH}_{\lambda}=J^{\cS^{\prime}}_{\lambda}$ so
$\delta_{(\beta,m)}$ is a cardinal in $\cR_{(\betaprime,\mprime)}$.
\end{subproof}

To obtain our contradiction,
we will show that $\cP$ is not
 $(\alpha,1)$-petite above $\delta$. Our proof of this
fact has two cases.

Suppose first that there is some $(\beta_0,m_0)$ such that
for all $(\beta,m)$ with
$(\beta_0,m_0)\lexleq(\beta,m)\lexless(\alpha,0)$, we have that
$\gamma_{(\beta,m)}=\gamma_{(\beta_0,m_0)}$. 
That is, suppose that all but boundedly many of our branches
are cofinal in the same initial segment of $\cU_{\theta}$.
(Notice that this first case \emph{must} be true if
$\cof(\alpha)>\omega$. This is because $\theta$ is a countable
ordinal.)
Then for all sufficiently large $(\beta,m)$, 
$\delta_{(\beta_0,m_0)}=\delta_{(\beta,m)}$, and
$\cR_{(\beta_0,m_0)}=\cR_{(\beta,m)}$.
Write $\cR=\cR_{(\beta_0,m_0)}$ and
 $\deltaprime=\delta_{(\beta_0,m_0)}$. Then $\deltaprime$ is Woodin in 
$\cR$ and, by Claim 1,
for all $(\beta,m)$ with
$(\beta_0,m_0)\lexleq(\beta,m)\lexless(\alpha,0)$, $\cR$
is not $(\beta,m)$-petite above $\deltaprime$. 
By part (b) of Lemma \ref{PetiteAboveWoodin}, $\cR$
is not $(\alpha,1)$-petite above $\delta$. 
By Claim 2, $\cR\initseg\cP$. Thus $\cP$ is not
 $(\alpha,1)$-petite above $\delta$.
This is a contradiction.

\begin{note}
Notice that in the situation described in the previous paragraph,
if $\cof(\alpha)\leq\omega$ then
we have something stronger than the fact that $\cR$
is not $(\alpha,1)$-petite above $\delta$.
We actually have that $\cR$ is not $(\alpha,1)$-\emph{small}
above $\delta$. This is relevant because we conjecture that
the theorem we are currently proving is also true if the
word ``petite'' is replaced by the word ``small'' in the statement
of the theorem. Up to this point in the proof, we have not made any
use of the distinction between ``small'' and ``petite.'' We shall
do so in the next paragraph.
\end{note}

Next suppose that there is no such pair $(\beta_0,m_0)$ as above.
As we pointed out above, this case can only occur if 
$\cof(\alpha)\leq\omega$.
Then we can choose an $\omega$-sequence
of pairs
$(\beta_0,m_0)\lexless(\beta_1,m_1)\lexless(\beta_2,m_2)\lexless\cdots$
such that $\sup\setof{(\beta_i,m_i)}{i\in\omega}=(\alpha,0)$, and
such that $\gamma_{(\beta_i,m_i)}<\gamma_{(\beta_{i+1},m_{i+1})}$.
Let $\eta=\sup_i\delta_{(\beta_i,m_i)}$.
It follows from Claims 2 and 3 that for each $i$, $\delta_{(\beta_i,m_i)}$
is a cardinal in $J^{\cP}_{\eta}$, and $\cR_{(\beta_i,m_i)}$ is
a proper initial segment of $J^{\cP}_{\eta}$. Let us write
$\cR_i=\cR_{(\beta_i,m_i)}$ and
$\delta_i=\delta_{(\beta_i,m_i)}$.
If there is some $i$ such that $\cR_i$
is not $(\alpha,1)$-petite above $\delta_i$, then we have our contradiction.
So suppose otherwise. Then, since $\cR_i$ is not $(\beta_i,m_i)$-petite
above $\delta_i$, there must be some $(\betaprime_i,\mprime_i)$
with $(\beta_i,m_i)\lexleq(\betaprime_i,\mprime_i)\lexless(\alpha,0)$
such that $\cR_i$ is weakly, explicitly $(\betaprime_i,\mprime_i)$-big
above $\delta_i$. Let $\kappa$ be the least measurable cardinal
of $J^{\cP}_{\eta}$ greater than $\delta$.
Then $\delta<\kappa<\delta_0$. (In fact $\delta_0$ is a limit of measurable
cardinals in $J^{\cP}_{\eta}$ since $\delta_0$ is a cardinal in
$J^{\cP}_{\eta}$ and $\delta_0$ is Woodin in $\cR_0$.)
By the proof of Lemma \ref{BignessIsBelowKappa},
we may assume that $\max(\code(\betaprime_i))<\kappa$ for all $i$.
Then, by Lemma \ref{WeakPlusOneGivesStrong}, we have that
$\cR_{i+1}$ is \emph{strongly} explicitly $(\betaprime_i,\mprime_i)$-big 
above
$\delta_{i+1}$. Thus, $\delta_{i+1}$ is locally
$(\betaprime_i,\mprime_i)$-Woodin in $J^{\cP}_{\eta}$.
Since the $\delta_i$ are cofinal in the ordinals of
$J^{\cP}_{\eta}$ and the $(\betaprime_i,\mprime_i)$ are cofinal in
$(\alpha,0)$, it follows that $J^{\cP}_{\eta}$ is weakly, explicitly
$(\alpha,1)$-big above $\delta$. Thus $\cP$ is not $(\alpha,1)$-petite
above $\delta$.
As this is a contradiction, this completes the proof
of (3) for $(\alpha,n)$ in the case that $n=0$.

\begin{note}
Let us consider for a moment how we might prove the current
theorem if it were stated with the word ``small'' instead of
the word ``petite.'' The previous paragraph is the only place in
the proof where we  use the distinction between ``small'' and 
``petite''.  In the previous paragraph we showed
that $J^{\cP}_{\eta}$ is not $(\alpha,1)$-petite above $\delta$, and
so $\cP$ is not $(\alpha,1)$-petite above $\delta$. 
This is a contradiction since
we assumed that $\cM$ is $(\alpha,1)$-petite above $\delta$.
We did not show that $\cP$ is not $(\alpha,1)$-small above $\delta$,
and so if we were only assuming that $\cM$ was $(\alpha,1)$-small
above $\delta$, we would not yet have a contradiction. One possible
strategy for proving the stronger version of the current theorem---the
version with ``small'' instead of ``petite''---is to now
get a contradiction
by showing that our coiteration is successful.
For a moment, suppose that we are assuming only that
$\cM$ is $(\alpha,1)$-small above $\delta$, rather than making
the stronger assumption that $\cM$ is $(\alpha,1)$-petite above $\delta$.
Using the notation from the previous paragraph, let
$\gamma=\sup\setof{\gamma_{(\beta_i,m_i)}}{i\in\omega}$.
Then $\gamma\leq\theta$, and $\length(\cU_{\gamma})$ is a limit ordinal.
Notice that 
$\sup\setof{\delta_i}{i\in\omega}=\ORD\intersect J^{\cP}_{\eta}=\eta
=\delta(\cU_{\gamma})$. Suppose that $\cU_{\gamma}$ has a cofinal,
wellfounded branch, $c$. (This is true if $\gamma<\theta$.)
Let $\cH=\cN_c^{\cU_{\gamma}}$.
Notice that each $\delta_i$ is a cardinal in $\cH$. Notice also that
$\cJ^{\cP}_{\eta}=\cJ^{\cH}_{\eta}$. 
If $\cH$ is at least as tall as the first
admissible structure over $J^{\cP}_{\eta}$, then since $\cP$ is
$(\alpha,1)$-small above $\delta$, we must have that
$\cP\lhd\cH$, and so our coiteration is successful. If $\cH$ is not as
tall as the first admissible structure over $J^{\cP}_{\eta}$, 
then $\ORD\intersect\cH$ is less than the least $\xi>\eta$ such that
$\xi$ is the index of an extender on the $\cP$-sequence. So again we
have $\cH\unlhd\cP$ or $\cP\unlhd\cH$. So our coiteration is successful,
and this is a contradiction.

Notice that in the previous paragraph, we don't need that $\cH$ is
fully wellfounded, but only that $\eta$ is in the wellfounded part 
of $\cH$. For then if $\cH$ is not wellfounded, then the wellfounded
part of $\cH$ is at least as tall as the first admissible over 
$J^{\cP}_{\eta}$, and so $\cP\lhd\cH$.

So our inability to prove the stronger version of the current theorem---the
version with ``small'' instead of ``petite''---boils down to our inability
to prove that if $\gamma=\theta$, then $\cU_{\theta}$ has a cofinal
branch $c$, that this branch $c$ has only boundedly many drops so
that the direct limit model $\cN^{\cU_{\theta}}_c$ is defined, and
that $\eta$ is in the wellfounded part of $\cN^{\cU_{\theta}}_c$.
Actually, it is not too hard to see that $\cU_{\theta}$ does have a cofinal
branch $c$. If $c$ has only boundedly many drops so that the direct limit
model $\cN^{\cU_{\theta}}_c$ is defined, then it is not too hard to
see that $\eta$ is in the wellfounded part of $\cN^{\cU_{\theta}}_c$.
So the problem boils down to showing that $c$ has only boundedly many
drops. This is what we are unable to do!
\end{note}

\skipmed

\underline{CASE 2}\quad  $n\geq1$ is odd.\qquad
Suppose $\cM$ and $\cN$ are countable $\delta$-mice, where
$\delta$ is a cutpoint of $\cM$ and $\cN$.  Suppose $\cJ^{\cM}_{\delta}=
\cJ^{\cN}_{\delta}$.  Suppose $\cM$ is $k(\cM,\delta)$ realizable and
$(\alpha,n+1)$-petite  above $\delta$, and that $\cN$ is
$\Pa{n}$ iterable above $\delta$. Suppose that
$\cM\notinitseg\cN$ and $\cN\notinitseg\cM$. We will derive a contradiction.

Since  $\cN$ is $\Pa{n}$ iterable above $\delta$ and $n$ is odd,
we actually have that
$\cN$ is length-$\omega$ iterable above $\delta$.

Again we define a coiteration $\bigsequence{(\cT_{\gamma},\cU_{\gamma})}%
{\gamma<\theta}$ of $\cM$ and $\cN$ by induction on $\gamma$. For this it
suffices to define $\cT_{\gamma+1}$ in the case $\length(\cT_{\gamma})$
is a limit ordinal, and similarly with $\cU_{\gamma+1}$. Fix a
$k(\cM,\delta)$ realization $\pi$ of $\cM$. If
$\length(\cT_{\gamma})$ is a limit ordinal then let $b$ be the unique
cofinal $(\pi,\cT)$ realizable branch of $\cT_{\gamma}$. It is shown
in \cite{St2} that
there always is such a $b$.  Let
$\cT_{\gamma+1}$ be the iteration tree extending $\cT_{\gamma}$
with length $\length(\cT_{\gamma})+1$ and such that
$b=[0,\length(\cT_{\gamma})]_{T_{\gamma+1}}$.

If $\length(\cU_{\gamma})$ is a limit ordinal then let $c$ be the unique
cofinal wellfounded branch of $\cU_{\gamma}$ such that $Q(c,\cU_{\gamma})$
is $\Pa{n-1}$ iterable
above $\delta(\cU_{\gamma})$, if there
is a unique such branch. If $c$ is defined then let
$\cU_{\gamma+1}$ be the iteration tree extending $\cU_{\gamma}$
with length $\length(\cU_{\gamma})+1$ and such that
$c=[0,\length(\cU_{\gamma})]_{U_{\gamma+1}}$. If $c$ is not defined then
stop the construction.  Note that there always is at least one cofinal 
wellfounded
branch $\cprime$  of $\cU_{\gamma}$ such that $Q(\cprime,\cU_{\gamma})$
is 
$\Pa{n-1}$ iterable
above $\delta(\cU_{\gamma})$. 
(Proof: We are assuming that $\cN$ is length-$\omega$ iterable above 
$\delta$.  As $n-1$ is even, length-$\omega$ iterability implies
$\Pa{n-1}$ iterability. This implies that for all limit 
$\gammaprime<\gamma$, we actually chose the unique length-$\omega$ iterable
branch of $\cU_{\gammaprime}$. This implies that there is a cofinal
branch $\cprime$ of $\cU_{\gamma}$ such that
$Q(\cprime,\cU_{\gamma})$ is length-$\omega$ iterable
above $\delta(\cU_{\gamma})$.)
So the only way
that the construction can stop is if there are two such branches.

The construction cannot produce a coiteration
$\bigsequence{(\cT_{\gamma},\cU_{\gamma})}{\gamma<\omega_1}$ of length
$\omega_1$.  This is because, by part (a) of this theorem for $(\alpha,n)$,
such a coiteration would be definable over $\JalphaR$,
 contradicting  Lemma \ref{CoiterationsAreCountable}.
Therefore
the construction must produce a coiteration
\mbox{$\bigsequence{(\cT_{\gamma},\cU_{\gamma})}{\gamma\leq\theta}$} with
$\theta<\omega_1$, such that either the construction stops at stage
$\theta$ for the reason described above, or
$\bigsequence{(\cT_{\gamma},\cU_{\gamma})}{\gamma\leq\theta}$ is terminal.
If the construction produces a terminal coiteration, and if this
coiteration is successful, then by Lemma \ref{AlreadyCompared},
$\cM\unlhd\cN$ or $\cN\unlhd\cM$. This is a contradiction.  So we assume
that the construction does not produce a successful coiteration. But
our construction makes it impossible to produce a non-successful,
terminal coiteration. So
the only possibility left is that the construction stopped for the reason
described above. That is: $\cU_{\theta}$ has limit length and has two
cofinal wellfounded branches $c_1$ and $c_2$ such that both $Q(c_1,\cU_{\theta})$ and
$Q(c_2,\cU_{\theta})$ are $\Pa{n-1}$
iterable above  $\delta(\cU_{\theta})$.
Let  $\cQ_1=Q(c_1,\cU_{\theta})$ and $\cQ_2=Q(c_2,\cU_{\theta})$.
Let $\cR$ be the corresponding  realizable $Q$-structure from
$\cT_{\theta}$.  Let $\gamma=\delta(\cU_{\theta})$.
Since $\cM$ is $(\alpha,n)$-petite  above $\delta$, $\cR$ is 
$(\alpha,n-1)$-petite above $\gamma$.
So by part (3) of this theorem for $(\alpha,n-1)$ we have that
$\cQ_1\unlhd\cR$ or $\cR\unlhd\cQ_1$ and  also that
$\cQ_2\unlhd\cR$ or $\cR\unlhd\cQ_2$.  As in the proof that possibility (c)
cannot occur, in the proof of Case 1 above, this is a contradiction. 
This completes the proof
of (3) for $(\alpha,n)$ in the case that $n\geq1$ is odd.

\skipmed

\underline{CASE 3}\quad  $n\geq2$ is even.\qquad
Suppose $\cM$ and $\cN$ are countable $\delta$-mice, where
$\delta$ is a cutpoint of $\cM$ and $\cN$.  Suppose $\cJ^{\cM}_{\delta}=
\cJ^{\cN}_{\delta}$.  Suppose $\cM$ is $k(\cM,\delta)$ realizable and
$(\alpha,n+1)$-petite  above $\delta$, and that $\cN$ is
$\Pa{n}$ iterable above $\delta$. Suppose that
$\cM\notinitseg\cN$ and $\cN\notinitseg\cM$. We will derive a contradiction.

Fix a mouse code $x_0$ for $(\cN,\delta)$.
Again we define a coiteration $\bigsequence{(\cT_{\gamma},\cU_{\gamma})}%
{\gamma<\theta}$ of $\cM$ and $\cN$ by induction on $\gamma$. For this it
suffices to define $\cT_{\gamma+1}$ in the case $\length(\cT_{\gamma})$
is a limit ordinal, and similarly with $\cU_{\gamma+1}$. Fix a
$k(\cM,\delta)$ realization $\pi$ of $\cM$. If
$\length(\cT_{\gamma})$ is a limit ordinal then let $b$ be the unique
cofinal $(\pi,\cT)$ realizable branch of $\cT_{\gamma}$. We know that
there always is such a $b$. Let
$\cT_{\gamma+1}$ be the iteration tree extending $\cT_{\gamma}$
with length $\length(\cT_{\gamma})+1$ and such that
$b=[0,\length(\cT_{\gamma})]_{T_{\gamma+1}}$.

Now suppose that $\length(\cU_{\gamma})$ is a limit ordinal.
Let $\cR_{\gamma}$ be the corresponding realizable $Q$-structure from
$\cT_{\gamma}$.   As $\cN$ is
$\Pa{n}$ iterable above $\delta$, player $\TWO$ wins the game
$G_{(\alpha,n)}(x_0)$.  Let $A_{\gamma}$ be the set of all reals $x_1$ such that:
\begin{quote}
$(\exists y_1,z_1\in\R)$ s.t. $y_1^*=\cU_{\gamma}$ and $z_1^*=\cR_{\gamma}$
and
$\angles{y_1,z_1,x_1}$ is a winning position for player $\TWO$ in
$G_{(\alpha,n)}(x_0)$.
\end{quote}
Now let  $c$ be the unique
cofinal wellfounded branch of $\cU_{\gamma}$ such that  $(\exists x_1\in
 A_{\gamma})$ with $x_1^*=
\bigl(\,Q(c,\cU_{\gamma}),\delta(\cU_{\gamma})\,\bigr)$,  if there
is a unique such branch. If $c$ is defined then let
$\cU_{\gamma+1}$ be the iteration tree extending $\cU_{\gamma}$
with length $\length(\cU_{\gamma})+1$ and such that
$c=[0,\length(\cU_{\gamma})]_{U_{\gamma+1}}$. If $c$ is not defined then
stop the construction.

The construction cannot produce a coiteration
$\bigsequence{(\cT_{\gamma},\cU_{\gamma})}{\gamma<\omega_1}$ of length
$\omega_1$.  This is because, by part (a) of this theorem for $(\alpha,n)$,
such a coiteration would be definable over $\JalphaR$
 contradicting  Lemma \ref{CoiterationsAreCountable}.
Therefore
the construction must produce a coiteration
\mbox{$\bigsequence{(\cT_{\gamma},\cU_{\gamma})}{\gamma\leq\theta}$} with
$\theta<\omega_1$, such that either the construction stops at stage
$\theta$ for the reason described above, or
$\bigsequence{(\cT_{\gamma},\cU_{\gamma})}{\gamma\leq\theta}$ is terminal.
If the construction produces a terminal coiteration, and if this
coiteration is successful, then by Lemma \ref{AlreadyCompared},
$\cM\unlhd\cN$ or $\cN\unlhd\cM$. This is a contradiction.
So
we assume that either the construction produces a terminal,
non-successful coiteration, or the construction stops for the reason
described above. This leaves us with three possibilities:
\begin{itemize}
\item[(a)] $\cU_{\theta}$ has successor length and the coiteration is
terminal because the next ultrapower
which is prescribed by the definition of a coiteration
in order to define $\cU_{\theta+1}$ is illfounded. Or,
\item[(b)] $\cU_{\theta}$ has limit length but has no cofinal wellfounded
branch $c$ such that $(\exists x_1\in A_{\theta})$ with $x_1^*=
\bigl(\,Q(c,\cU_{\theta}),\delta(\cU_{\theta})\,\bigr)$. Or,
\item[(c)]$\cU_{\theta}$ has limit length and has two cofinal wellfounded
branches $c_1$ and $c_2$ such that both $c_1$ and $c_2$ satisfy the 
conditions on $c$ in
possibility (b) above.
\end{itemize}
We will show that each of these possibilities leads to a contradiction.

First suppose that either possibility (a) or possibility (b) holds.
If possibility
(a) holds then let $\cU$ be the putative iteration tree extending
$\cU_{\theta}$ whose last model is illfounded and which results from
$\cU_{\theta}$ by taking the ultrapower which is mentioned in possibility
(a). If possibility (b) holds then let $\cU=\cU_{\theta}$. In either case
let $y_1$ be a real such that $y_1^* = \cU$.

Now, if $\cT_{\theta}$ has
limit length, then let $b$ be the unique cofinal realizable branch of
$\cT_{\theta}$ and let $X=Q(b,\cT_{\theta})$. Otherwise let $X=\emptyset$.
In either case let  $z_1$ be any real such that $z_1^*=$ the transitive 
closure of
$\singleton{\cT_{\theta},X}$.

Notice that
$\angles{y_1,z_1}$ is a legal first move for player $\ONE$ in
$G_{(\alpha,n)}(x_0)$, and that if player $\ONE$ plays this
move he does not immediately lose. As player $\TWO$ wins
$G_{(\alpha,n)}(x_0)$ there is a real $x_1$ such that
$\angles{y_1,z_1,x_1}$ is a winning position for player $\TWO$. Fix such
an $x_1$. Then there is some maximal wellfounded branch $c$ of $\cU$ so that
$x_1^*=
\bigl\langle Q(c,\cU\restriction\sup(c)),\,
 \delta\bigl(\cU\restriction\sup(c)\bigr)\,\bigr\rangle$.

First we will show that $c$ is not a cofinal branch of $\cU$.
Suppose towards a contradiction  that $c$ is a cofinal branch of $\cU$. Then
$\cU$ has limit length and so $\cU=\cU_{\theta}$. Let $\cR_{\theta}$ be the
corresponding realizable $Q$-structure from $\cT_{\theta}$. Then
$\cR_{\theta}$ is in the transitive closure of $\singleton{\cT_{\theta},X}$.
 So there is a real $z$
such that $z^*=\cR_{\theta}$ and $z$ is recursive
in $z_1$.  Examining the rules $G_{(\alpha,n)}(x_0)$ for even $n$ we wee 
that
$\angles{y_1,z,x_1}$ is also a winning position for player $\TWO$.   This
means that $x_1\in A_{\theta}$.
But this violates the hypothesis that either possibility (a) or possibility
 (b) occurred. So $c$ is not
a cofinal branch of $\cU$.

So let $\gamma<\theta$ be such that $c$ is cofinal in $\cU_{\gamma}$.
Then $\cU_{\gamma}$ has limit length.   Now
$\cR_{\gamma}$ is in the transitive closure of $\singleton{\cT_{\theta},X}$.
 So there is a real $z$
such that $z^*=\cR_{\gamma}$ and $z$ is recursive
in $z_1$. Fix such a  $z$. Let $y$ be any real such that $y^*=\cU_{\gamma}$.
  By examining
the rules of $G_{(\alpha,n)}(x_0)$ for even $n$ we wee that 
$\angles{y,z,x_1}$
is also a winning position for player $\TWO$. Thus $x_1\in A_{\gamma}$.
Now let $\cprime$ be the unique cofinal branch of $\cU_{\gamma}$ which is
a non-maximal branch of $\cU_{\theta}$.  Then $c$ and $\cprime$ witness that
we should have stopped the coiteration construction at stage $\gamma$ 
because
possibility (c) occurred. This is a contradiction.  So we have shown that 
neither
possibility (a) nor possibility (b) can occur.

Finally, assume towards a contradiction that possibility (c) has occurred.
Let $c_1$ and $c_2$ be two cofinal wellfounded
branches of $\cU_{\theta}$ such that 
$(\exists x_1\in A_{\theta})$ with $x_1^*=
\bigl(\,Q(c_1,\cU_{\theta}),\delta(\cU_{\theta})\,\bigr)$,
and similarly with $c_2$.
Let $\delta^*=\delta(\cU_{\theta})$.
Let $\cQ=Q(c_1,\cU_{\theta})$, and let $\cQ^{\prime}=Q(c_2,\cU_{\theta})$.
Let $\cR$ be the corresponding realizable $Q$-structure from $\cT_{\theta}$.
To obtain a contradiction, we will show that $\cR$ is comparable with both
$\cQ$ and $\cQ^{\prime}$.  Below we will show that $\cR$ is comparable
with $\cQ$.  The proof that $\cR$ is comparable with $\cQ^{\prime}$ is
identical.

So suppose towards a contradiction that $\cR\notinitseg\cQ$ and
$\cQ\notinitseg\cR$.  Now we define a coiteration
$\bigsequence{(\bar{\cT}_{\gamma},\bar{\cU}_{\gamma})}%
{\gamma<\bar{\theta}}$ of $\cR$ and $\cQ$ by induction on $\gamma$.
 For this it
suffices to define $\bar{\cT}_{\gamma+1}$ in the case
$\length(\bar{\cT}_{\gamma})$
is a limit ordinal, and similarly with $\bar{\cU}_{\gamma+1}$. Fix a
$k(\cR,\delta)$ realization $\bar{\pi}$ of $\cR$. If
$\length(\bar{\cT}_{\gamma})$ is a limit ordinal then let $b$ be the unique
cofinal $(\bar{\pi},\bar{\cT})$ realizable branch of $\bar{\cT}_{\gamma}$. 
We know that
there always is such a $b$. Let
$\bar{\cT}_{\gamma+1}$ be the iteration tree extending $\bar{\cT}_{\gamma}$
with length $\length(\bar{\cT}_{\gamma})+1$ and such that
$b=[0,\length(\bar{\cT}_{\gamma})]_{\bar{T}_{\gamma+1}}$.

If $\length(\bar{\cU}_{\gamma})$ is a limit ordinal then let $c$ be the 
unique
cofinal wellfounded branch of $\bar{\cU}_{\gamma}$ such that
$Q(c,\bar{\cU}_{\gamma})$
is $\Pa{n-2}$ iterable above $\delta(\bar{\cU}_{\gamma})$, if there
is a unique such branch. If $c$ is defined then let
$\bar{\cU}_{\gamma+1}$ be the iteration tree extending $\bar{\cU}_{\gamma}$
with length $\length(\bar{\cU}_{\gamma})+1$ and such that
$c=[0,\length(\bar{\cU}_{\gamma})]_{\bar{U}_{\gamma+1}}$.
If $c$ is not defined then
stop the construction.

The construction cannot produce a coiteration
$\bigsequence{(\bar{\cT}_{\gamma},\bar{\cU}_{\gamma})}{\gamma<\omega_1}$
of length
$\omega_1$.  This is because, by part (a) of this theorem for $(\alpha,n-1)$,
such a coiteration would be definable over $\JalphaR$,
 contradicting  Lemma \ref{CoiterationsAreCountable}.
Therefore
the construction must produce a coiteration
\mbox{$\bigsequence{(\bar{\cT}_{\gamma},\bar{\cU}_{\gamma})}%
{\gamma\leq\bar{\theta}}$} with
$\bar{\theta}<\omega_1$, such that either the construction stops at stage
$\bar{\theta}$ for the reason described above, or
$\bigsequence{(\bar{\cT}_{\gamma},\bar{\cU}_{\gamma})}%
{\gamma\leq\bar{\theta}}$ is terminal.
If the construction produces a terminal coiteration, and if this
coiteration is successful, then by Lemma \ref{AlreadyCompared},
$\cR\unlhd\cQ$ or $\cQ\unlhd\cR$. This is a contradiction.
So
we assume that either the construction produces a terminal,
non-successful coiteration, or the construction stops for the reason
described above. This leaves us with three possibilities:
\begin{itemize}
\item[(i)] $\bar{\cU}_{\bar{\theta}}$ has successor length and the 
coiteration is
terminal because the next ultrapower
which is prescribed by the definition of a coiteration
in order to define $\bar{\cU}_{\bar{\theta}+1}$ is illfounded. Or,
\item[(ii)] $\bar{\cU}_{\bar{\theta}}$ has limit length but has no cofinal 
wellfounded
branch $c$ such that  $Q(c,\bar{\cU}_{\bar{\theta}})$
is $\Pa{n-2}$ iterable above $\delta(\bar{\cU}_{\bar{\theta}})$.  Or,
\item[(iii)]$\bar{\cU}_{\bar{\theta}}$ has limit length and has two cofinal wellfounded
branches $c$ and $\cprime$ such that both  $Q(c,\bar{\cU}_{\bar{\theta}})$
and $Q(\cprime,\bar{\cU}_{\bar{\theta}})$
are $\Pa{n-2}$ iterable above $\delta(\bar{\cU}_{\bar{\theta}})$
\end{itemize}
We will show that each of these possibilities leads to a contradiction.
We begin with the easiest, which is the third.

Suppose towards
a contradiction that possibility (iii) holds.
Let $\bar{\delta}=\delta(\bar{\cU}_{\bar{\theta}})$.  Let
$\cQ_1=Q(c,\bar{\cU}_{\bar{\theta}})$ and
$\cQ_2=Q(\cprime,\bar{\cU}_{\bar{\theta}})$.
Let $\bar{\cR}$ be the corresponding  realizable $Q$-structure from
$\bar{\cT}_{\bar{\theta}}$.
Since $\cM$ is $(\alpha,n+1)$-petite  above $\delta$,
$\cR$ is $(\alpha,n-1)$-petite and  above $\delta^*$,
and so
$\bar{\cR}$ is $(\alpha,n-2)$-petite and  above $\bar{\delta}$.
Also
$\cQ_1$ and $\cQ_2$ are $\Pa{n-2}$ iterable above $\bar{\delta}$.
So by part (3) of this theorem for $(\alpha,n-2)$ we have that
$\cQ_1\unlhd\bar{\cR}$ or $\bar{\cR}\unlhd\cQ_1$ and  also that
$\cQ_2\unlhd\bar{\cR}$ or $\bar{\cR}\unlhd\cQ_2$. This is a contradiction.
(See the proof that possibility (c) cannot occur in the proof of
Case 1 above.)

Now we assume towards a contradiction
that either possibility (i) or possibility (ii) has occurred.
Recall that by our current case hypotheses we have that
$(\exists x_1\in A_{\theta})$ with $x_1^*=(\cQ,\delta^*)$. 
Fix such an $x_1$.
By the definition of $A_{\theta}$, $(\exists y_1,z_1\in\R)$ s.t.
$y_1^*=\cU_{\theta}$ and $z_1^*=\cR$, and
$\angles{y_1,z_1,x_1}$ is a winning position for player $\TWO$ in
$G_{(\alpha,n)}(x_0)$.  (Recall that $x_0$ is our fixed mouse code for
$(\cN,\delta)$.) Fix such reals $y_1,z_1$.

If possibility
(i) holds then let $\cU$ be the putative iteration tree extending
$\bar{\cU}_{\bar{\theta}}$ whose last model is illfounded and which results
 from
$\bar{\cU}_{\bar{\theta}}$ by taking the ultrapower which is mentioned in
possibility
(i). If possibility (ii) holds then let $\cU=\bar{\cU}_{\bar{\theta}}$.

\begin{claim}
There is a real $y_2$ such that  $y_2^*=\cU$ and
$y_2\in\Da{n-1}(x_1,z_1)$.
\end{claim}
\begin{subproof}[Proof of Claim]
We will prove the claim in the case that possibility (ii) has occurred,
that is assuming that $\cU=\bar{\cU}_{\bar{\theta}}$.  The proof
in the case that possibility (i) has occurred is similar.

Let $\rho=\lh(\cU)$. First we will show that there is a real
$w\in\WO$ such that $|w|=\rho$ and $w\in\Da{n-1}(x_1,z_1)$.
For $y,w\in\R$ put $S(y,w)$ iff $y^*$ is an iteration tree on $\cQ_1$,
and $w\in\WO$ and $|w|$ is equal to the length of the tree, and the
tree results from a comparison of $\cQ_1$ with $\cR$ in which all branches
chosen are $\Pa{n-2}$ iterable. More precisely, $S(y,w)$ iff
\begin{quote}
$w\in\WO$ and there is some real $v$ such that
$v^*$ is some coiteration
$\lsequence{(\cT^{\prime}_{\gamma},\cU^{\prime}_{\gamma})}%
{\gamma\leq\mu}$ of some two premice $\cR^{\prime}$, $\cQ^{\prime}$
and $(x_1^*)_0=\cQ^{\prime}$ and
$z_1^*=\cR^{\prime}$ and $y^*=\cU^{\prime}_{\mu}$ and
$|w|=\lh(y^*)$
and $(\forall a,b,c,d,e\in\R)$ if there is some
$\gamma<\mu$ s.t. (either $a^*=\cT^{\prime}_{\gamma}$, or
$a^*=\cU^{\prime}_{\gamma}$)  and
$a^*$ has limit length and $b^*$ is the unique cofinal branch of $a^*$
which was chosen by the construction, and $c^*=Q(b^*,a^*)$ and
$d^*=\delta(a^*)$ and $e^*=(c^*,d^*)$ then player $\TWO$ wins the game
$G_{(\alpha,n-2)}(e)$.
\end{quote}

Then $S\in\Pa{n-2}(x_1,z_1)$.  Now set $R(w)\Iff(\exists y)S(y,w)$.  Then
$R\in\Sa{n-1}(x_1,z_1)$ and $\rho$  is the largest element of the set
$\dotsetof{|w|}{R(w)}$.  Since $n-1$ is odd, we may apply
Kechris's Boundedness Theorem for Subsets of $\omega_1$,
Lemma \ref{BoundedSetofOrdinals} on page
\pageref{BoundedSetofOrdinals}, to conclude that  there is a $w\in\WO$
such that $|w|=\rho$ and $w\in\Da{n-1}(x_1,z_1)$. Fix such a $w$.

Now, let $P$ be the binary relation on $\HC$ given by
$P(\angles{\cM^{\prime},\rho^{\prime}},\cT^{\prime})\Iff$ $\cM^{\prime}$ is
 a premouse
and $\cT^{\prime}$ is an iteration tree on $\cM^{\prime}$ of length 
$\rho^{\prime}$.
Then $P$ admits a definable code transformation,  $f$.
[Proof:  Suppose $P(\angles{\cM^{\prime},\rho^{\prime}},\cT^{\prime})$ 
holds.
Suppose $(\xprime)^*=\cM^{\prime}$ and $\wprime\in\WO$ and
$|\wprime|=\rhoprime$.
As in the proof of Proposition \ref{DefCodTransProp}, 
 it suffices to see that we can define (over $\HC$, using
$\cM^{\prime}$, $\rho^{\prime}$, $\cT^{\prime}$, $\xprime$ and $\wprime$
as parameters)
a bijection from $\omega$ onto $\cT^{\prime}$. Since $\xprime$ gives
us a bijection from $\omega$ on $\cM^{\prime}$, and $\wprime$ gives
us a bijection from $\omega$ onto $\rhoprime$, it is easy to see
how to define a bijection from $\omega$ onto $\cT^{\prime}$.]

Now $P(\angles{\cQ_1,\rho},\cU)$ holds.
Let $y_2=f(\angles{\cQ_1,\rho},\cU,\angles{x_1,w})$.  Then $y_2^*=\cU$ and
$y_2$ is the unique real $y$ such that:
\begin{quote}
$S(y,w)$ and $y=f(\angles{x_1^*,|w|},y^*,\angles{x_1,w})$.
\end{quote}

As $w\in\Da{n-1}(x_1,z_1)$, we have that $y_2\in\Da{n-1}(x_1,z_1)$.
\end{subproof}

Fix a real $y_2$ as in the claim.  Notice that $\angles{y_1,z_1,x_1,y_2}$ is
 a
legal position for player $\ONE$ in $G_{(\alpha,n)}(x_0)$ which is not an
immediate loss for player $\ONE$.  Since $\angles{y_1,z_1,x_1}$ is a winning
position for player $\TWO$, there is a real $x_2$ such that
$\angles{y_1,z_1,x_1,y_2,x_2}$ is a wining position for player $\TWO$. By
examining the rules of $G_{(\alpha,n)}(x_0)$ we see that this means that
there is a maximal branch  $c$ of $\cU$ such that
$x_2^*=\bigl\langle Q(c,\cU\restriction\sup(c)),\,
\delta\bigl(\cU\restriction\sup(c)\bigr)\,\bigr\rangle$
and also such that
player $\TWO$ wins the game $G_{(\alpha,n-2)}(x_2)$.  The branch $c$ can
not be cofinal in $\cU$ because this would violate our hypothesis that
either possibility (i) or possibility (ii) occurred.  So let 
$\gamma<\bar{\theta}$
be such that $c$ is cofinal in $\bar{\cU}_{\gamma}$. Then at stage $\gamma$
of our coiteration we should have stopped the construction because 
possibility
(iii) occurred.  This contradiction completes the proof of
Theorem  \ref{ComparisonTheorem}.
\end{proof}

\skipbig

We conclude this section with the one and only application of the
previous theorem that we shall need.

\begin{corollary}
\label{AlphaPetiteImpliesAlpha}
Assume that there exists $\omega$ Woodin cardinals.
 Suppose $2\leq\alpha<\omega_1^{\omega_1}$, and $n\in\omega$.
Let $\cM$ be a countable, realizable,  sound premouse. Suppose that
 $\cM$ is $(\alpha,n+1)$-petite. 
 Then $\R\intersect\cM\subseteq \Aa{n+1}$.
\end{corollary}
\begin{proof}
Recall that $(\alpha,n+1)$-petite means
$(\alpha,n+1)$-petite above $0$, or equivalently above $\omega$.
Let $x\in\R\intersect\cM$. Let $\xi$ be the rank of $x$ in the order of 
construction
of $\cM$. Let $\cN\initseg\cM$ be the $\initseg$-least initial segment
of $\cM$ containing $x$. Then $\cN$ projects to $\omega$, and since
$\cN$ is sound, $\cN$ is an $\omega$-mouse. By part (3) of
Theorem  \ref{ComparisonTheorem}, if $\cN^{\prime}$ is any $\omega$-mouse,
and $\cN^{\prime}$ is $\Pan$ iterable, then $\cN\initseg\cN^{\prime}$
or $\cN^{\prime}\initseg\cN$. (Notice that the requirements that
$\cN$ and $\cN^{\prime}$ agree through $\omega$, and that $\omega$ be
a cutpoint of $\cN$ and $\cN^{\prime}$ are met trivially.) Thus,
$x$ is the unique real such that there is some
$\Pan$ iterable $\omega$-mouse $\cN^{\prime}$ such that $x$ is the 
$\xi^{\text{th}}$ real in the order of construction of
of $\cN^{\prime}$.
Let $w\in\WO$ with $|w|=\xi$. Then $x$ is the unique real $\xprime$
such that:
\begin{quote}
$(\exists y,z\in\R)$ s.t. $z$ is a mouse code  and player $\TWO$ wins the 
game
$G_{(\alpha,n)}(z)$ and $(z^*)_1=\omega$  and $y^*=(z^*)_0$ and 
$\xprime\in y^*$
and the rank of $\xprime$ in the order of construction of $y^*$ is $|w|$.
\end{quote}
Thus $x\in\Da{n+1}(w)$. Since $w$ was arbitrary, $x\in\Da{n+1}(\xi)$.
So $x\in\Aa{n+1}$.
\end{proof}

%


\skipbig

\section{Conclusions}

\label{section:concl}

In this section we put together the results of Section
\ref{section:correctness} and Section \ref{section:comparison}
to obtain our main results.

\begin{theorem}
\label{MainTheoremOne}
Assume that there exists $\omega$ Woodin cardinals. Suppose that
$(2,0)\lexleq(\alpha,n)\lexleq(\omega_1^{\omega_1},0)$,
and either $\cof(\alpha)=\omega_1$, or $n=0$.
 Then $\Aan$ is a mouse set.
\end{theorem}
\begin{proof}
Let $\cM$ be a countable, realizable premouse such that $\cM$ is
$(\alpha,n)$-big, but every proper initial segment of $\cM$ is
$(\alpha,n)$-small. We will show that $\Aan=\R\intersect\cM$.

First we will show that there exists a premouse $\cM$ as in the
previous paragraph.  In fact we will show that there exists a
``canonical'' such $\cM$.
It is shown in \cite{St2} that our large cardinal hypothesis implies
that there is a proper class premouse $L[\vec{E}]$ such that
$L[\vec{E}]\models$``$\exists$ $\omega$ Woodin cardinals.'' 
Let $\cM^*$ be the $\initseg$-least  initial segment of $L[\vec{E}]$ such 
that $\cM^*$ has $\omega$ Woodin cardinals which are cofinal in the 
ordinals of $\cM^*$. It is easy to see that $\cM^*$ projects to
$\omega$. (This uses a simple condensation argument, similar to the proof
of Corollary \ref{PhiMinimalProjects}.)
Thus $\cM^*$ is an $\omega$-mouse. It follows 
 that $\cM^*$ is
the unique $\cM^{\prime}$ such that $\cM^{\prime}$ is a realizable,
$\omega$-mouse, $\cM^{\prime}\models$``There are $\omega$ Woodin
cardinals cofinal in the ordinals,'' and no proper initial segment of
$\cM^{\prime}$ models this sentence. (Use Lemma \ref{AlreadyCompared}
on page \pageref{AlreadyCompared} to see this.)
Thus we may think of $\cM^*$ 
as the canonical inner model for the theory ``There are $\omega$ Woodin
cardinals cofinal in the ordinals.'' 
By Proposition \ref{OmegaWoodinsIsBig} on page
\pageref{OmegaWoodinsIsBig}, $\cM^*$ is
$(\omega_1^{\omega_1},0)$-big. So in particular, $\cM^*$ is 
$(\alpha,n)$-big. Let $\cM$ be the $\initseg$-least initial segment
of $\cM^*$ which is $(\alpha,n)$-big. We may think of $\cM$ as the
``canonical'' $(\alpha,n)$-big mouse.

 Our large cardinal hypothesis
implies that every game in $\JofR{\kappa}$ is determined, where $\kappa$
is the least $\R$-admissible ordinal. 
(See  Corollary \ref{determinacyholds}.)
In particular we have that
$\Det(\BfSigmaOne(\JofR{\omega_1^{\omega_1}}))$ holds. So we can apply
our main correctness result,
Corollary \ref{BigMouseContainsAan} on page \pageref{BigMouseContainsAan},
to conclude that 
$\Aan\subseteq\R\intersect\cM$.

Conversely, let $x\in\R\intersect\cM$ and we will show that 
$x\in\Aan$. Let $\cN$ be the $\initseg$-least initial segment of $\cM$
such that $x\in\cN$. Then $\cN$ is a proper initial segment of $\cM$,
so by hypothesis $\cN$ is $(\alpha,n)$-small.

First suppose that $n\geq 1$. Then by hypothesis $\cof(\alpha)>\omega$.
By definition $(\alpha,n)$-petite is synonymous with $(\alpha,n)$-small, so we have that
$\cN$ is $(\alpha,n)$-petite. By the corollary to our main
comparison theorem,
Corollary \ref{AlphaPetiteImpliesAlpha}
on page \pageref{AlphaPetiteImpliesAlpha}, $\R\intersect\cN\subseteq\Aan$. 
Thus $x\in\Aan$.

Now suppose that $n=0$. Recall that
$$\Aa{0}=\setof{z\in\R}{(\exists\beta<\alpha)\; z\in\OD^{\JbetaR}}.$$
So we must show that $x\in\OD^{\JbetaR}$ for some $\beta<\alpha$.
We have that $\cN$ is $(\alpha,0)$-small. By definition this is
synonymous with $(\alpha,0)$-petite, so $\cN$ is $(\alpha,0)$-petite.
It follows that
$\cN$ is $(\beta,m)$-petite
for some $(\beta,m)\lexless(\alpha,0)$. 
(See Lemma \ref{ZeroPetiteIsPetiter} on page \pageref{ZeroPetiteIsPetiter}.)
Fix such a $(\beta,m)$ with
$m\geq1$. As in the previous paragraph we have that
$\R\intersect\cN\subseteq\Aa[\beta]{m}$. 
Thus $x\in\Aa[\beta]{m}$. In particular, $x\in\OD^{\JbetaR}$.
\end{proof}

\begin{remark}
Let $\cM^*$ be as in the proof above. We believe that
$\R\intersect\cM^*=\Aa[\kappa]{0}$, where $\kappa=\kappa^{\R}$ is the
least $\R$-admissible ordinal. It would take us too far afield to prove
this here, so we will save this for another paper.
$\Aa[\kappa]{0}$ is an important set in descriptive set theory. It is the
set of reals which are hyperprojective in a countable ordinal. It is
also the largest countable inductive set of reals.
\end{remark}

\begin{remark}
One specific consequence of Theorem \ref{MainTheoremOne} is that for
$\alpha<\omega_1^{\omega_1}$, $\OD^{\JalphaR}$ is a mouse set.
This is because $\OD^{\JalphaR}=\Aa[\alpha+1]{0}$.
\end{remark}

Recall that on page \pageref{aquestion} we discussed the relationship
between the sets $\Aan$ and $\OD^{\JalphaR}_{n+1}$. As promised
we can now answer Questions 1 and 2 from that discussion affirmatively,
at least for certain pairs $(\alpha,n)$. We begin with Question 1.

\begin{theorem}
Assume that there exists $\omega$ Woodin cardinals. Suppose that
$2<\alpha<\omega_1^{\omega_1}$,  $\cof(\alpha)=\omega_1$, and 
$n\geq 1$. Then $\Aan=\OD^{\JalphaR}_{n+1}$.
\end{theorem}
\begin{proof}
The proof of Theorem \ref{MainTheoremOne} actually showed this,
since our main correctness result,
Corollary \ref{BigMouseContainsAan}, not only gives us that
$\Aan\subseteq\R\intersect\cM$, but actually that
$\OD^{\JalphaR}_{n+1}\subseteq\R\intersect\cM$.
\end{proof}

The next theorem answers Question 2 from our discussion on page
\pageref{aquestion}, at least for ordinals $\alpha<\omega_1^{\omega_1}$.

\begin{theorem}
Assume that there exists $\omega$ Woodin cardinals. Suppose that
$2\leq\alpha<\omega_1^{\omega_1}$. Then $\Union{n}\Aan=\OD^{\JalphaR}$.
\end{theorem}
\begin{proof}
Fix $\alpha$. Notice that $\OD^{\JalphaR}=\Aa[\alpha+1]{0}$.
Let $\cM$ be a countable, realize premouse such that
$\cM$ is $(\alpha+1,0)$-big, but every proper initial segment of
$\cM$ is $(\alpha+1,0)$-small.
The proof of Theorem \ref{MainTheoremOne} actually showed that
$\Aa[\alpha+1]{0}\subseteq\R\intersect\cM$, and
$\R\intersect\cM\subseteq\Union{n}\Aan$.
\end{proof}

\bigskip 

By showing that $\Aan$ is a mouse set, we have demonstrated a connection
between inner model theory and descriptive set theory. Connections between
two different areas of mathematics are useful, because they allow us to 
use the tools from one area to tackle problems from another area. In the
rest of this section we will explore in a little more detail the connection
between mice and reals which are ordinal definable in $\LofR$.

If $\cM$ is a mouse and $\Aan=\R\intersect\cM$, then certain aspects of
the structure of $\cM$ can give us information about $\Aan$. For example,
we can learn something about $\Aan$ by looking at the ordinal 
$\omega_1^{\cM}$. We turn to this idea now.

Suppose $\alpha\leq\omega_1^{\omega_1}$ and $\cof(\alpha)>\omega$.
Under the hypothesis that there are $\omega$ Woodin cardinals, we
have enough determinacy to conclude that the set
$$\Aa{0}=\Union{\beta<\alpha}\OD^{\JbetaR}$$
is countable. This implies in particular that $\Aa{0}=\Aa[\beta]{0}$
for some $\beta<\alpha$. Let $\beta$ be the least such ordinal. Then
we will say that $\beta$ is the \emph{$\OD$ cutoff ordinal} below $\alpha$.

\begin{proposition}
\label{MainTheoremTwo}
Assume that there exists $\omega$ Woodin cardinals. Suppose that
$\alpha\leq\omega_1^{\omega_1}$ and $\cof(\alpha)=\omega_1$.
Let $\cM$ be a countable, realizable premouse such that $\cM$ is
$(\alpha,0)$-big, but every proper initial segment of $\cM$ is
$(\alpha,0)$-small. Then $\cM\models$``$\omega_1$ exists.'' Furthermore,
letting $\gamma=\omega_1^{\cM}$ and letting 
$$\beta=\sup\setof{\alphaprime<\alpha}{\max(\code(\alphaprime))<\gamma},$$
we have that $\beta$ is the $\OD$ cutoff ordinal
below $\alpha$.
\end{proposition}
\begin{proof}
Let $\cM$ be as in the statement of the proposition. It is shown in the
proof of Theorem \ref{MainTheoremOne} above
that $\Aa{0}=\R\intersect\cM$.
If $\cM$ is \emph{explicitly} $(\alpha,0)$-big then, since no proper
initial segment of $\cM$ is explicitly $(\alpha,0)$-big, we have that
$\cM\models$``$\omega_1$ exists.'' Note also that in this case $\cM$ is
explicitly $(\alpha,0)$-big above $\omega_1^{\cM}$.
Now suppose that $\cM$ is not explicitly
$(\alpha,0)$-big, and let $\alphaprime>\alpha$ be least such that
$\cM$ is explicitly $(\alphaprime,0)$-big. Then $\alphaprime$ cannot
be a successor ordinal or a limit ordinal of cofinality $\omega$, so
$\alphaprime$ must be a limit ordinal of cofinality $\omega_1$. So,
as no proper initial segment of $\cM$ is $(\alphaprime,0)$-big,
$\cM\models$``$\omega_1$ exists.'' Note also that in this case $\cM$ is
explicitly $(\alphaprime,0)$-big above $\omega_1^{\cM}$. In any case
we have that $\cM\models$``$\omega_1$ exists.''

Now let $\gamma=\omega_1^{\cM}$ and let 
$$\beta=\sup\setof{\alphaprime<\alpha}{\max(\code(\alphaprime))<\gamma}.$$
It is easy to see that $\max(\code(\beta))=\gamma$, and that $\beta$
is a limit ordinal of cofinality $\omega$. 
Let $\delta$ be the $\OD$ cutoff ordinal below $\alpha$. We will show
that $\beta=\delta$.

First we will show that $\delta\leq\beta$. Let $\alphaprime<\delta$, and
we will show that $\alphaprime<\beta$. Since $\alphaprime<\delta$, there
is an $x\in\R$ such that $x\in\Aa{0}$ but $x\notin\Aa[\alphaprime]{0}$.
Since $\Aa{0}=\R\intersect\cM$, there is a proper initial segment 
$\cN\properseg\cM$ such that $x\in\R\intersect\cN$ and $\ORD^{\cN}<\gamma$.
Fix such an $\cN$. We claim that $\cN$ is $(\beta,0)$-petite. For suppose
$\cN$ were not $(\beta,0)$-petite. Since $\cN$ is $(\alpha,0)$-petite,
there must be some $(\alpha_0,n_0)$ with 
$(\beta,0)\lexleq(\alpha_0,n_0)\lexless(\alpha,0)$ such that $\cN$ is
explicitly $(\alpha_0,n_0)$-big. Since $\beta\leq\alpha_0<\alpha$,
we have that $\max(\code(\alpha_0))\geq\gamma$. But this is impossible
because according to Lemma \ref{bigcontainscode},
 $\max(\code(\alpha_0))\leq\ORD^{\cN}$. Thus $\cN$ is $(\beta,0)$-petite.
Then, as in the proof of Theorem \ref{MainTheoremOne} above, it follows 
that $x\in\Aa[\beta]{0}$. Since $x\notin\Aa[\alphaprime]{0}$ it follows
that $\alphaprime<\beta$.

Now we will show that $\beta\leq\delta$. Let $\alphaprime<\beta$ and
we will show that $\alphaprime<\delta$. Since $\alphaprime<\beta$, there
is an ordinal $\xi$ with $\alphaprime\leq\xi<\alpha$ such that
$\max(\code(\xi))<\gamma$. Fix such a $\xi$. As we pointed out in the
first paragraph of this proof, $\cM$ is explicitly $(\alpha_0,0)$-big
above $\gamma$, for some $\alpha_0\geq\alpha$. Since $\xi<\alpha$ and
$\max(\code(\xi))<\gamma$, we have by Lemma \ref{ProperSegmentsAreBig}
that for all $n\in\omega$, 
there is a proper initial segment $\cN\properseg\cM$ 
such that $\cN$ is explicitly $(\xi,n)$-big above $\gamma$. In particular
there is a proper initial segment $\cN\properseg\cM$ 
such that $\cN$ is explicitly $(\xi,1)$-big (above $0$.) Let $\cN$ be
the $\initseg$-least such initial segment. By our second corollary to
the Condensation Theorem, Corollary \ref{PhiMinimalIsCountable} on page 
\pageref{PhiMinimalIsCountable},
$\ORD^{\cN}<\gamma$. Thus there is a real $x$ such that $x\in\cM$ but
$x\notin\cN$. Fix such a real $x$. Since $x\in\cM$, $x\in\Aa{0}$.
Since $x\notin\cN$, $x\notin\Aa[\xi]{0}$. Thus $\xi<\delta$.
Since $\alphaprime\leq\xi$, we have shown that $\alphaprime<\delta$.
\end{proof}

Below we highlight two particular concrete examples of the previous 
proposition: $\alpha=\omega_1$, and $\alpha=\omega_1^{\omega_1}$.

\begin{corollary}
Assume that there exists $\omega$ Woodin cardinals. 
Let $\cM$ be a countable, realizable premouse such that $\cM$ is
$(\omega_1,0)$-big, but every proper initial segment of $\cM$ is
$(\omega_1,0)$-small. Let $\gamma=\omega_1^{\cM}$. Then
$\gamma$ is the $\OD$ cutoff ordinal below $\omega_1$.
\end{corollary}
\begin{proof} Let
$$\beta=\sup\setof{\alphaprime<\omega_1}{\max(\code(\alphaprime))<\gamma}.$$
If $\alphaprime<\omega_1$ then 
$\code(\alphaprime)=\singleton{(0,\alphaprime)}$, and so
$\max(\code(\alphaprime))=\alphaprime$. Thus $\beta=\gamma$.
\end{proof}

\begin{corollary}
Assume that there exists $\omega$ Woodin cardinals. 
Let $\cM$ be a countable, realizable premouse such that $\cM$ is
$(\omega_1^{\omega_1},0)$-big, but every proper initial segment of $\cM$ is
$(\omega_1^{\omega_1},0)$-small. Let $\gamma=\omega_1^{\cM}$. Then
$\omega_1^{\gamma}$ is the $\OD$ cutoff ordinal below $\omega_1^{\omega_1}$.
\end{corollary}
\begin{proof}
Let
$$\beta=\sup\setof{\alphaprime<\omega_1^{\omega_1}}%
{\max(\code(\alphaprime))<\gamma}.$$
$\code(\omega_1^{\gamma})=\singleton{(\gamma,1)}$. If 
$\max(\code(\alphaprime))<\gamma$ then it follows that 
$\alphaprime<\omega_1^{\gamma}$. Furthermore, there are cofinally many
$\alphaprime<\omega_1^{\gamma}$ such that $\max(\code(\alphaprime))<\gamma$.
Thus $\beta=\omega_1^{\gamma}$.
\end{proof}

We continue with our program of exploring the connections between the
structure of certain mice, and certain sets of ordinal definable
reals. In this paper we have defined two distinct notions: the notion
of $(\alpha,n)$-big, and the notion of \emph{explicitly} $(\alpha,n)$-big.
One might ask, what is the significance of a mouse $\cM$ which is
$(\alpha,n)$-big but not explicitly $(\alpha,n)$-big? The next proposition
sheds some light on this question. The proposition says that $\alpha$ must
be in the midst of a long gap in which no new reals are ordinal defined
in $\LofR$.

\begin{proposition}
Assume that there exists $\omega$ Woodin cardinals. 
Suppose $(2,0)\lexleq(\alpha,n)\lexless(\omega_1^{\omega_1},0)$.
Let $\cM$ be a countable, realizable premouse such that $\cM$ is
$(\alpha,n)$-big, but every proper initial segment of $\cM$ is
$(\alpha,n)$-small. Suppose that $\cM$ is not explicitly $(\alpha,n)$-big.
Then there is an ordinal $\alpha^*>\alpha$ such that 
$\cof(\alpha^*)=\omega_1$, and such that $\alpha$ is greater than
the $\OD$ cutoff ordinal below $\alpha^*$.
\end{proposition}
\begin{proof}
Let $(\alpha^*,n^*)$ be the lexicographically least pair such that
$(\alpha,n)\lexleq(\alpha^*,n^*)$,  and $\cM$ is explicitly 
$(\alpha^*,n^*)$-big. We are assuming that 
$(\alpha,n)\lexless(\alpha^*,n^*)$. It follows that $n^*=0$ and
$\alpha^*>\alpha$.
Then $\alpha^*$ must be a limit ordinal of cofinality
$\omega_1$. 
By definition, $\cM$ is active. Let $\kappa$ be the critical
point of the last extender on the $\cM$ sequence. Since $\cM$ is not
explicitly $(\alpha,n)$-big, we must have that 
$\max(\code(\alpha))\geq\kappa$. In particular 
$\max(\code(\alpha))>(\omega_1)^{\cM}$. 
Now $\cM$ is the $\initseg$-least initial segment of $\cM$
which is $(\alpha^*,0)$-big. So by Proposition \ref{MainTheoremTwo} 
above, $\alpha$ is greater than
the $\OD$ cutoff ordinal below $\alpha^*$.
\end{proof}

We will mention one last result about the connection between the
structure of mice and ordinal definable reals.
As we have pointed out,
assuming determinacy there are only countably many reals which are
ordinal definable in $\LofR$. So there will be many pairs
$(\alpha,n)$ so that $\Aan=\Aa{n+1}$. If $\Aan\not=\Aa{n+1}$, let
us say that $\Aa{n+1}$ is \emph{non-trivial.} Let $\cM_1$ be the
least mouse which is $(\alpha,n)$-big, and let $\cM_2$ be the least
mouse which is $(\alpha,n+1)$-big. Of course we can tell whether or not
$\Aa{n+1}$ is trivial by checking whether or not 
$\R\intersect\cM_2\subseteq\cM_1$. Our last proposition says that we
can also check whether or not $\Aa{n+1}$ is trivial by comparing
$\code(\alpha)$ with $\omega_1^{\cM_2}$.

\begin{proposition} 
Assume that there exists $\omega$ Woodin cardinals. Suppose that
$\alpha<\omega_1^{\omega_1}$,  $\cof(\alpha)=\omega_1$,
and $n\geq 0$.
Let $\cM$ be a countable, realizable premouse such that $\cM$ is
 $(\alpha,n+1)$-big, but every proper initial segment of $\cM$ is
$(\alpha,n+1)$-small. 
Then $\Aa{n+1}$ is non-trivial, if and only if
$\max(\code(\alpha))<\omega_1^{\cM}$.
\end{proposition}
\begin{proof}
First suppose that $\max(\code(\alpha))<\omega_1^{\cM}$. We claim that
$\cM$ is \emph{explicitly} $(\alpha,n+1)$-big. For otherwise,
as in the proof of the previous proposition, we must have that
$\cM$ is active and
$\max(\code(\alpha))$ is greater than the critical point of the last
extender on the $\cM$ sequence. But this contradicts our assumption that
$\max(\code(\alpha))<\omega_1^{\cM}$. Thus 
$\cM$ is  explicitly $(\alpha,n+1)$-big.
Since no proper initial segment of $\cM$ is  explicitly
$(\alpha,n+1)$-big, $\cM$ is 
explicitly $(\alpha,n+1)$-big above
$\omega_1^{\cM}$. Since $\max(\code(\alpha))<\omega_1^{\cM}$, 
 Lemma \ref{ProperSegmentsAreBig} tells us that there is a \emph{proper} 
initial segment of 
$\cM$ which is  explicitly $(\alpha,n)$-big above
$\omega_1^{\cM}$. Let $\cN$ be the $\initseg$-least initial segment of
$\cM$ which is  explicitly $(\alpha,n)$-big (above $0$.) 
Our second corollary to the  Condensation Theorem, Corollary 
\ref{PhiMinimalIsCountable}, implies
that  $\ORD^{\cN}<\omega_1^{\cM}$.
So there is a real $x$ such that $x\in\cM$ but $x\notin\cN$. 
Since $x\in\cM$, $x\in\Aa{n+1}$. Since $x\notin\cN$, $x\notin\Aan$. 
So $\Aa{n+1}$ is non-trivial.

Conversely, suppose that $\Aa{n+1}$ is non-trivial. Let $x\in\Aa{n+1}-\Aan$.
Let $\cN\properseg\cM$ be such that $x\in\cN$ and 
$\ORD^{\cN}<\omega_1^{\cM}$. Then $\cN$ is not $(\alpha,n)$-petite.
(If $\cN$ were $(\alpha,n)$-petite then we would have that $x\in\Aan$.)
Since $\cof(\alpha)=\omega_1$, $(\alpha,n)$-petite is synonymous with
$(\alpha,n)$-small. So $\cN$ is not $(\alpha,n)$-small. But $\cN$ is
$(\alpha,n+1)$-small. Thus $\cN$ must be \emph{explicitly}
$(\alpha,n)$-big. By Lemma \ref{bigcontainscode}, 
$\max(\code(\alpha))\leq\ORD^{\cN}$.
Thus $\max(\code(\alpha))<\omega_1^{\cM}$.
\end{proof}

%



\end{document}